\documentclass[12pt]{article} 
\usepackage[final]{graphicx}
\usepackage[dvips]{color}
\usepackage{graphicx}
\usepackage{amsmath}
\usepackage{amssymb}
\usepackage{ascmac}
\usepackage{amsthm}
\usepackage{graphicx}
\usepackage{bm}
\usepackage{url}
\def\Bbb{\bf} 

\newcommand\C{{ \Bbb C}}

\newcommand\Z{{\Bbb Z}}
\newcommand\Q{{\Bbb Q}}

\newcommand\G{{\Gamma}} 

\newcommand*\pFqskip{8mu}
\catcode`,\active
\newcommand*\pFq{\begingroup
        \catcode`\,\active
        \def ,{\mskip\pFqskip\relax}%
        \dopFq
}
\catcode`\,12
\def\dopFq#1#2#3#4#5{%
        {}_{#1}F_{#2}\biggl(\genfrac..{0pt}{}{#3}{#4};#5\biggr)%
        \endgroup
}
\newtheorem{lmm}{Lemma}[section]

\newtheorem{rmk}[lmm]{Remark}
\newtheorem{prp}[lmm]{Proposition}

\def\comment#1{ }

\makeatletter
    
    \@addtoreset{equation}{section}
  \makeatother
\allowdisplaybreaks
\begin{document}
\title{Special values of the hypergeometric series}
\author{Akihito Ebisu
}

\maketitle
\begin{abstract}
In this paper, we present a new method for finding identities for hypergeoemtric series,
such as the (Gauss) hypergeometric series, the generalized hypergeometric series
and the Appell-Lauricella hypergeometric series.
Furthermore, using this method, we get identities for the hypergeometric series $F(a,b;c;x)$;
we show that values of $F(a,b;c;x)$ at some points $x$  
can be expressed in terms of gamma functions, together with certain elementary functions.
We tabulate the values of $F(a,b;c;x)$ 
that can be obtained with this method. 
We find that this set includes almost all previously known values
and many previously unknown values. 


\end{abstract}
\tableofcontents
\section{Introduction}
The series
\begin{gather*}
F(a;-;x)
:=\sum^{\infty}_{n=0}\frac{(a,n)}{(1,n)}\,x^n,
\label{binomial_series}
\end{gather*}
where $(a,n):=\Gamma (a+n)/\Gamma (a)$, can be expressed as
\begin{gather}
F(a;-;x)=(1-x)^{-a}\quad (|x|<1).
\label{binomial_thm}
\end{gather}
This identity (\ref{binomial_thm}) is called the binomial theorem.
The binomial theorem (\ref{binomial_thm}) has many applications.
For example, because
\begin{gather*}
F\left(\frac{1}{2};-;\frac{1}{50}\right)=\frac{5\sqrt{2}}{7}
\label{sqrt2}
\end{gather*}
from (\ref{binomial_thm}), we can perform numerical calculation of $\sqrt{2}$ effectively 
(cf. 292 page in Chapter 4 in Part II of [Eu]):
\begin{gather*}
\sqrt{2}=
\frac{7}{5}\left\{1+\frac{1}{100}+\frac{1\cdot 3}{2!}\left(\frac{1}{100}\right)^2
+\frac{1\cdot 3\cdot 5}{3!}\left(\frac{1}{100}\right)^3+\cdots\right\}.
\label{sqrt21}
\end{gather*}
Using (\ref{binomial_thm}) like this, we are able to do numerical calculations of algebraic numbers.
Moreover, we obtain, for non-negative integer $k$,
\begin{gather*}
F(-k;-;-1)=\sum ^{k} _{n=0}\binom {k}{n}=2^k.
\label{combinatorial}
\end{gather*}
Thus, we can obtain combinatorial identities from (\ref{binomial_thm}).
As seen above, the binomial theorem is very useful.
Analogically, if hypergeometric series which are generalizations of $F(a;-;x)$
can be expressed in terms of well-known functions,
then it is expected that those hypergeometric identities have many applications. 
Indeed, transcendental numbers are expressed in term of hypergeometric series as we see later.
Further, it is known that most of combinatorial identities appearing in [Goul]
are special cases of those identities (cf. [An]).

That is why many methods for obtaining hypergeometric identities have been constructed.
In particular, in the last several decades, 
many methods that exploit progress in computer technology have been formulated:
Gosper's algorithm, the W-Z method, Zeilberger's algorithm, etc. (cf. [Ko] and [PWZ]).
These algorithms have been used to obtain and prove hypergeometric identities.
In this article, 
we present a new method for finding identities for hypergeoemtric series,
such as the (Gauss) hypergeometric series, the generalized hypergeometric series
and the Appell-Lauricella hypergeometric series.
Moreover, as an example of the application of this method, we give identities for the hypergeometric series;
we see that values of the hypergeometric series at some arguments  
can be expressed in terms of gamma functions, together with certain elementary functions.
Because our method systematically yields almost all known
identities of this type besides many new identities, for completeness, 
we tabulate the entire set.
In the appendix of this paper,
we also make use of our method to have some hypergeometric identities for
generalized hypergeometric series and Appell-Lauricella hypergeometric series.
In particular, hypergeometric identities for latter series have never been obtained.
The reason comes from the fact that
it is difficult to apply known methods directly 
to Appell-Lauricella hypergeometric series
because these are multiple series.
Therefore, our method will become a powerful tool 
for investigating hypergeometric identities
for many kinds of hypergeometric series, 
especially for Appell-Lauricella hypergeometric series. 

Now, we present a new method for finding identities for hypergeoemtric series.
To begin with, we give a proof of the binomial theorem (\ref{binomial_thm}), which is the model of our method. 
Put $\partial :=\frac{d}{dx}$. Recalling that $x\partial x^n=nx^n$, we see
\begin{align*}
\left(x\partial +a\right)F(a;-;x)&=\sum ^{\infty}_{n=0}\frac{(a,n)}{(1,n)}(x\partial +a)x^n
=\sum ^{\infty}_{n=0}\frac{(a,n)}{(1,n)}(n +a)x^n\\
&=\sum ^{\infty}_{n=0}\frac{a(a+1,n)}{(1,n)}x^n=aF(a+1;-;x).
\end{align*}
We also get
\begin{align*}
\partial F(a;-;x)&=\sum ^{\infty}_{n=1}\frac{(a,n)}{(1,n)}nx^{n-1}
=\sum ^{\infty}_{n=0}\frac{(a,n+1)}{(1,n)}x^n\\
&=\sum ^{\infty}_{n=0}\frac{a(a+1,n)}{(1,n)}x^n=aF(a+1;-;x).
\end{align*}
The above two formulae lead to
\begin{gather}
F(a;-;x)=(1-x)F(a+1;-;x).
\label{proof_binom}
\end{gather}
Thus, we find that $F(a;-;x)$ and $F(a+1;-;x)$ are linearly related over rational functions.
Furthermore, we can regard (\ref{proof_binom}) as a first order difference equation 
with rational functions coefficients.
Now, we solve this difference equation (\ref{proof_binom}). 
Because
\begin{gather*}
F(a+n-1;-;x)=(1-x) F(a+n;-;x)
\label{proof_binom1}
\end{gather*}
holds for any positive integer $n$ from (\ref{proof_binom}), we see
\begin{align}
\begin{split}
F(a;-;x)&=(1-x)F(a+1;-;x)=(1-x)^2F(a+2;-;x)=\cdots\\
&=(1-x)^nF(a+n;-;x).
\end{split}
\label{proof_binom2}
\end{align}
Substituting $a=-n$ in (\ref{proof_binom2}), we get the following identity:
\begin{gather}
F(-n;-;x)=(1-x)^nF(0;-;x)=(1-x)^n.
\label{proof_binom3}
\end{gather}
Since the identity (\ref{proof_binom3}) is valid for any non-negative integer $n$,
we find that (\ref{proof_binom3}) holds for any complex number $n$
by virtue of Carlson's theorem (cf. 5.3 in [Ba] or Lemma 3.1 in this article).
Thus, we get the binomial theorem.
Our method for finding hypergeometric identities is a generalization of this procedure.
It is known that several hypergeometric series with the same parameters up to additive integers
are linearly related over rational functions.
Suppose that this linear relation degenerates into a linear relation between two hypergeometric series
under a condition,
and further, this degenerate relation can be regarded as the first order difference equation 
with rational functions coefficients. 
Then, solving this difference equation as we get the binomial theorem,
we will obtain a hypergeometric identity. 
In this paper, applying this method to the hypergeometric series
\begin{gather*}
F(a,b;c;x)
:=\sum^{\infty}_{n=0}\frac{(a,n)(b,n)}{(c,n)(1,n)}\,x^n,
\label{hypergeometric}
\end{gather*}
we actually obtain almost all previously known hypergeometric identities
and, of course, many new hypergeometric identities.
As stated previously, 
we will also make use of this method to 
yield some hypergeometric identities for generalized hypergeometric series
and Appell-Lauricella hypergeometric series in the appendix of this paper.  

It is known that 
for a given triple of integers $(k,l,m)\in\Z^3,$
there exists a unique pair of rational functions 
$(Q(x), R(x)) \in (\Q (a,b,c,x))^2$ 
satisfying
\begin{gather}
F(a+k,b+l;c+m;x)=Q(x)F(a+1,b+1;c+1;x)+R(x)F(a,b;c;x),
\label{three_term}
\end{gather}
where $\Q (a,b,c,x)$ is the field generated over $\Q$ by $a$, $b$, $c$ and $x$ (cf. \S 1 in [Eb1]).
This relation is called the {\bf three term relation} of the hypergeometric series. 
From this, we obtain
\begin{align*}
\begin{split}
&F(a+nk,b+nl;c+nm;x)\\
&=Q^{(n)}(x)F(a+(n-1)k+1,b+(n-1)l+1;c+(n-1)m+1;x)\\
&+R^{(n)}(x)F(a+(n-1)k,b+(n-1)l,c+(n-1)m;x),
\end{split}
\end{align*}
where 
\begin{align*}
Q^{(n)}(x)&:=Q(x)|_{(a,b,c)\rightarrow (a+(n-1)k,b+(n-1)l,c+(n-1)m)}, \\
R^{(n)}(x)&:=R(x)|_{(a,b,c)\rightarrow (a+(n-1)k,b+(n-1)l,c+(n-1)m)}.
\end{align*}
Let $(a,b,c,x)$ be a solution of the system
\begin{gather}
Q^{(n)}(x)=0\quad {\text {for $n=1,2,3,\ldots$}}.
\label{condition}
\end{gather}
Then, for this $(a,b,c,x)$, we have the first order difference equation
\begin{align}
F(a,b;c;x)=\frac{1}{R^{(1)}(x)}\times F(a+k,b+l;c+m;x),
\label{degenerated0}
\end{align}
or we have
\begin{align}
F(a,b;c;x)=\frac{1}{R^{(1)}(x)R^{(2)}(x)\cdots R^{(n)}(x)}\times F(a+nk,b+nl;c+nm;x).
\label{degenerated1}
\end{align}
Solving (\ref{degenerated0}) or (\ref{degenerated1}), 
we obtain hypergeometric identities;
for example, if $k \neq 0$, then
by substituting $a=-nk-k'$, where $k'$ is an integer satisfying
$0\leq k'< |k|$, into (\ref{degenerated1}),
we find the hypergeometric identity
\begin{align}
\begin{split}
&F(-nk-k',b;c;x)\\
&=\left. \left(\frac{1}{R^{(1)}(x)R^{(2)}(x)\cdots R^{(n)}(x)}
\right) \right|
_{a\rightarrow -nk-k'}\times F(-k',b+nl;c+nm;x).
\end{split}
\label{-nk-k'}
\end{align}
In this way, we can find values of $F(a,b;c;x)$.
In this paper, we call values of the hypergeometric series obtained by our method
{\bf special values} of the hypergeometric series. 
We summarize how to get special values of the hypergeometric series:
\begin{itemize}
\item[] Step 1: Give $(k,l,m)\in \Z ^3$. 
\item[] Step 2: Obtain $(a,b,c,x)$ satisfying the system (\ref{condition}).
\item[] Step 3: Solve (\ref{degenerated0}) or (\ref{degenerated1}). 
\end{itemize} 
We remark that the numerator of $Q^{(n)}(x)$ is a polynomial in $n$ over $\Z [a,b,c,x]$
for a given $(k,l,m)\in \Z ^3$.
Therefore, in order to obtain  $(a,b,c,x)$ satisfying the system (\ref{condition}),
we only have to seek the $(a,b,c,x)$ that eliminate all of the coefficients of this polynomial.
Namely, we solve polynomial systems in $a,b,c,x$ to implement Step 2.
The Gr\"obner basis theory is useful for this.
Thus, we see that (\ref{degenerated1}) and (\ref{-nk-k'}) are valid for any integer $n$
from the constitution method of $(a,b,c,x)$ satisfying (\ref{condition}).

As an example, we give a special value for $(k,l,m)=(0,0,1)$.
Because the three term relation for $(k,l,m)=(0,0,1)$ is
\begin{align*}
&F(a,b;c+1;x)\\
&=\frac{ab(1-x)}{(c-a)(c-b)}{}F(a+1,b+1;c+1;x)+\frac{c(c-a-b)}{(c-a)(c-b)}{}F(a,b;c;x),
\end{align*}
we have
\begin{gather*}
Q^{(n)}(x)=\frac{ab(1-x)}{(c+n-1-a)(c+n-1-b)}.
\end{gather*}
Hence, $(a,b,c,x)$ satisfying (\ref{condition}) are $(0,b,c,x), (a,0,c,x), (a,b,c,1)$.
However, since the special values in the cases that $(a,b,c,x)=(0,b,c,x), (a,0,c,x)$ 
are obvious from the definition of $F(a,b;c;x)$,
we only consider the special value in the case that $(a,b,c,x)=(a,b,c,1)$.
In this case, we have the first order difference equation
\begin{gather}
F(a,b;c;1)=\frac{(c-a)(c-b)}{c(c-a-b)}{}F(a,b;c+1;1),
\label{gauss2}
\end{gather}
or we have
\begin{gather}
F(a,b;c;1)=\frac{(c-a,n)(c-b,n)}{(c,n)(c-a-b,n)}{}F(a,b;c+n;1)
\label{gauss3}
\end{gather}
for any positive number $n$,
where we assume $\Re (c-a-b)>0$ which is the convergence condition of $F(a,b;c;1)$.
Noticing that
\begin{align*}
&\lim_{n \to +\infty}\frac{(c-a,n)(c-b,n)}{(c,n)(c-a-b,n)}
=\frac{\Gamma(c)\Gamma(c-a-b)}{\Gamma(c-a)\Gamma(c-b)},\\
&\lim_{n \to +\infty}{}F(a,b;c+n;1)=1,
\end{align*}
we obtain 
\begin{gather}
F(a,b;c;1)=\frac{\Gamma \left(c\right)\Gamma \left(c-a-b\right)}
{\Gamma \left(c-a\right)\Gamma \left(c-b\right)}
\label{gauss}
\end{gather}
for $\Re (c-a-b)>0$. 
This is called the Gauss summation formula (cf. [48] in 24 of [Ga]).
Moreover, as a particular case of (\ref{gauss}), we have
\begin{gather}
F(a,-n;c;1)=\frac{(c-a,n)}{(c,n)},
\label{chu}
\end{gather}
where $n \in \Z _{\geq 0}$.
Though (\ref{chu}) has been known since the 13th century,
today, this is called the Chu-Vandermonde equality (cf. Corollary 2.2.3 in [AAR]).
The special values for other lattice points $(k,l,m)$ are considered in Sections 3 and 4.
\begin{rmk}
In [48] in 24 of [Ga], Gauss obtains (\ref{gauss}) using (\ref{gauss2}) and (\ref{gauss3})
which are got by substituting $x=1$ into the relation
\begin{align*}
0&=c(c-1-(2c-a-b-1)x)F(a,b;c;x)+(c-a)(c-b)xF(a,b;c+1;x)\\
&-c(c-1)(1-x)F(a,b;c-1;x).
\end{align*}
Therefore, our method is applied by Gauss in the special case.
\end{rmk}

In the above example, we obtained the Gauss summation formula (\ref{gauss}).
This is the general expression of $F(a,b;c;x)$ at $x=1$.
Hence, we only have to consider values of $F(a,b;c;x)$ at other points.
So, in Step 2, we seek $(a,b,c,x)$ satisfying the following:
let $(a,b,c)$ be a triple such that the number of solutions of (\ref{condition}) is finite, 
and let $x$ $(\neq 0,1)$ be one of its solutions (cf. Remark 2.4).
We call such a quadruple $(a,b,c,x)$ an {\bf admissible quadruple}. 

Although we need to investigate $(k,l,m)\in \Z^3$ for obtaining the special values of $F(a,b;c;x)$
as seen in Step 1,
we actually only have to investigate $(k,l,m)$ contained in
\begin{gather}
\{(k,l,m)\in \Z ^3 \mid 0 \leq k+l-m \leq l-k \leq m\}.
\label{quotient}
\end{gather}
We close this section by giving an outline of this statement 
(see Section 2 for details).
The hypergeometric equation $E(a,b,c): L(a,b,c)y=0$, where
\begin{gather*}
L(a,b,c):=\partial ^2+\frac{c-(a+b+1)x}{x(1-x)}\partial -\frac{ab}{x(1-x)},
\label{hge}
\end{gather*}
admits 23 hypergeometric solutions in addition to $F(a,b;c;x)$.
For each of these solutions, there exists
a relation similar to (\ref{degenerated1}).
For example, for the solution $x^{-a}F(a,a+1-c;a+1-b;1/x)$, we have
\begin{align*}
\begin{split}
&x^{-a-k}F(a+k,a+1-c+k-m;a+1-b+k-l;1/x)\\
&=\frac{(-1)^{m-k-l+1}c(b,l)(a+1-b,k-l)}{b(c,m)(a+1-c,k-m)}Q(x)x^{-a-1}F(a+1,a+1-c;a+1-b;1/x)\\
&+\frac{(-1)^{m-k-l}(b,l)(a+1-b,k-l)}{(c,m)(a+1-c,k-m)}R(x)x^{-a}F(a,a+1-c;a+1-b;1/x),
\end{split}
\end{align*}
where $Q(x)$ and $R(x)$ appear in (\ref{three_term}).
Thus, because 
\begin{align}
\begin{split}
&F(a,a+1-c;a+1-b;1/x)\\
&=\frac{(-1)^{n(m-k-l)}(c,nm)(a+1-c,n(k-m))}{(b,nl)(a+1-b,n(k-l))}\frac{1}{R^{(1)}(x)R^{(2)}(x)\cdots R^{(n)}(x)}\\
&\times x^{-nk}F(a+nk,a+1-c+n(k-m);a+1-b+n(k-l);1/x)
\end{split}
\label{degenerated2}
\end{align}
for an admissible quadruple $(a,b,c,x)$,
we can get special values of $F(a,a+1-c;a+1-b;1/x)$
with this quadruple. 
The same can be done for the other 22 solutions.
Thus, for the lattice point $(k,l,m)$,
we are able to obtain special values 
of 24 hypergeometric series with the above quadruple.

The identity (\ref{degenerated2}) implies that the special values 
of the 24 hypergeometric series mentioned above 
for the lattice points $(k,k-m,k-l)$
coincide with those for the lattice point $(k,l,m)$
(see Subsection 2.3 for details).
In other words, $(k,k-m,k-l)$ is equivalent to $(k,l,m)$ 
with respect to the obtained special values.
In fact, this holds generally,
as all 24 lattice points represented by triples 
corresponding to the 24 hypergeomtric solutions
are equivalent with respect to these special values.
In addition, from the relation $F(a,b;c;x)=F(b,a;c;x)$, 
which implies that $(k,l,m)$ is equivalent to $(l,k,m)$,  
it can also be shown that 
the $48$ $(=24\cdot 2)$ lattice points are equivalent 
with respect to these special values. 
These 48 lattice points form the orbit of $(k,l,m)$
under the action of the group $G$ on $\Z ^3$,
where $G$ is the group generated by following mappings
\begin{align*}
&\sigma _1 : (k,l,m)\rightarrow (m-k,l,m),&
&\sigma _2 : (k,l,m)\rightarrow (k,l,k+l-m),\\
&\sigma _3 : (k,l,m)\rightarrow (l,k,m),&
&\sigma _4 : (k,l,m)\rightarrow (m-k,m-l,m),\\
&\sigma _5 : (k,l,m)\rightarrow (-k,-l,-m).
\end{align*}
We remark that $G=\langle \sigma _1 , \sigma _2\rangle
\ltimes\left( \langle \sigma _3 \rangle 
\times \langle \sigma _4 \rangle 
\times \langle \sigma _5 \rangle\right) 
= S_3 \ltimes \left( S_2 \times S_2 \times S_2\right)$,
where $S_n$ is the symmetric group of degree $n$.
From this, regarding these 48 lattice points as equivalent, 
we can take (\ref{quotient}) as a complete system of representatives of
the quotient $G{\backslash} \Z ^3$ of this action
(see Subsection 2.3 for details).
Thus, we only need to investigate the lattice points contained in (\ref{quotient})
to obtain the special values of the hypergeometric series.
In this paper, we tabulate the special values 
for $(k,l,m)$ satisfying $0 \leq k+l-m \leq l-k \leq m \leq 6$.

\section{Preliminaries}
As stated in the previous section, for a given $(k,l,m)$,
the 24 hypergeometric solutions of $E(a,b,c)$ 
with an admissible quadruple $(a,b,c,x)$
have two term relations as (\ref{degenerated1}) and (\ref{degenerated2}), say degenerate relations.
In this section, we list these degenerate relations explicitly.
This list will be used subsequently, when we evaluate 
special values of the corresponding hypergeometric series.
After presenting this list, we give a proof 
that
(\ref{quotient}) can be taken as a complete system of representatives of $G{\backslash} \Z ^3$.
\subsection{Contiguity operators}
In this subsection,
we introduce contiguity operators.

The contiguity operators of the hypergeometric series are defined as 
\begin{align*}
&H_1(a,b,c):=\vartheta +a,\, B_1(a,b,c):=\frac{1}{(1-a)(c-a)}(-x(1-x)\partial +(bx+a-c)),\\
&H_2(a,b,c):=\vartheta +b,\, B_2(a,b,c):=\frac{1}{(1-b)(c-b)}(-x(1-x)\partial +(ax+b-c)),\\
&H_3(a,b,c):=\frac{1}{(c-a)(c-b)}((1-x)\partial +c-a-b),\, B_3(a,b,c):=\vartheta +c-1,
\end{align*}
where $\vartheta := x\partial$. 
These operators satisfy the following relations:
\begin{align*}
H_1(a,b,c){}F(a,b;c;x)&=a\,{}F(a+1,b;c;x),\\
B_1(a,b,c){}F(a,b;c;x)&=1/(a-1)\cdot {}F(a-1,b;c;x),\\
H_2(a,b,c){}F(a,b;c;x)&=b\,{}F(a,b+1;c;x),\\
B_2(a,b,c){}F(a,b;c;x)&=1/(b-1)\cdot {}F(a,b-1;c;x),\\
H_3(a,b,c){}F(a,b;c;x)&=1/c\cdot{}F(a,b;c+1;x),\\
B_3(a,b,c){}F(a,b;c;x)&=(c-1){}F(a,b;c-1;x)
\end{align*}
(cf. Theorems 2.1.1 and 2.1.3 in [IKSY]).
Let $H(k,l,m)$ be a composition of these contiguity operators such that
\begin{gather}
H(k,l,m)F(a,b;c;x)=\frac{(a,k)(b,l)}{(c,m)}{}F(a+k,b+l;c+m;x).
\label{three_term_f}
\end{gather}
Also, $H(k,l,m)$ can be expressed as
\begin{gather} 
H(k,l,m)=p(\partial)L(a,b,c)+q(x)\partial +r(x),
\label{composition}
\end{gather}
where $p(\partial)$ is an element of the ring of the differential operators in $x$ over $\Q (a,b,c,x)$,
and $q(x), r(x)\in \Q (a,b,c,x)$
(cf. (2.4) in [Eb1]). 

\subsection{Degenerate relations}
In this subsection, we list degenerate relations.
\begin{lmm}
$E(a,b,c)$ admits the following 24 hypergeometric solutions (cf. 2.9 in [Erd]):
\begin{align*}
y_1(a,b,c,x)&:=F(a,b;c;x)\\
&=(1-x)^{c-a-b}F(c-a,c-b;c;x)\\
&=(1-x)^{-a}F(a,c-b,c;x/(x-1))\\
&=(1-x)^{-b}F(c-a,b;c;x/(x-1)),\\
y_2(a,b,c,x)&:=F(a,b;a+b+1-c;1-x)\\
&=x^{1-c}F(a+1-c,b+1-c;a+b+1-c;1-x)\\
&=x^{-a}F(a,a+1-c;a+b+1-c;1-x^{-1})\\
&=x^{-b}F(b+1-c,b;a+b+1-c;1-x^{-1}),\\
y_3(a,b,c,x)&:=x^{-a}F(a,a+1-c;a+1-b;1/x)\\
&=(-1)^{a}(-x)^{b-c}(1-x)^{c-a-b}F(1-b,c-b;a+1-b;1/x)\\
&=(-1)^{a}(1-x)^{-a}F(a,c-b;a+1-b;(1-x)^{-1})\\
&=(-1)^{a}(-x)^{1-c}(1-x)^{c-a-1}F(a+1-c,1-b;a+1-b;(1-x)^{-1}),\\
y_4(a,b,c,x)&:=x^{-b}F(b+1-c,b;b+1-a;1/x)\\
&=(-1)^{b}(-x)^{a-c}(1-x)^{c-a-b}F(1-a,c-a;b+1-a;1/x)\\
&=(-1)^{b}(1-x)^{-b}F(b,c-a;b+1-a;(1-x)^{-1})\\
&=(-1)^{b}(-x)^{1-c}(1-x)^{c-b-1}F(b+1-c,1-a;b+1-a;(1-x)^{-1}),\\
y_5(a,b,c,x)&:=x^{1-c}F(a+1-c,b+1-c;2-c;x)\\
&=x^{1-c}(1-x)^{c-a-b}F(1-a,1-b;2-c;x)\\
&=x^{1-c}(1-x)^{c-a-1}F(a+1-c,1-b;2-c;x/(x-1))\\
&=x^{1-c}(1-x)^{c-b-1}F(b+1-c,1-a,2-c;x/(x-1)),\\
y_6(a,b,c,x)&:=(1-x)^{c-a-b}F(c-a,c-b;c+1-a-b;1-x)\\
&=x^{1-c}(1-x)^{c-a-b}F(1-a,1-b;c+1-a-b;1-x)\\
&=x^{a-c}(1-x)^{c-a-b}F(c-a,1-a;c+1-a-b;1-x^{-1})\\
&=x^{b-c}(1-x)^{c-a-b}F(c-b,1-b;c+1-a-b;1-x^{-1}),
\end{align*}
In the above, we must take the appropriate 
branches of $(-1)^a$ and $(-1)^b$.
\end{lmm}
We now apply both sides of (\ref{composition}) 
to $y_i$ $(i=1,2,\cdots ,6)$.
First, noting the relation
\begin{gather*}
\partial F(a,b;c;x)=\frac{ab}{c}F(a+1,b+1;c+1;x),
\end{gather*}
we have 
\begin{align*}
y_1(a+k,b+l,c+m,x)&=\frac{ab(c,m)}{c(a,k)(b,l)}q(x)y_1(a+1,b+1,c+1,x)\\
&+\frac{(c,m)}{(a,k)(b,l)}r(x)y_1(a,b,c,x)
\end{align*}
from (\ref{three_term_f}) and (\ref{composition}).
This implies that
\begin{gather}
Q(x)=\frac{ab(c,m)}{c(a,k)(b,l)}q(x), \ R(x)=\frac{(c,m)}{(a,k)(b,l)}r(x).
\label{QR}
\end{gather}
Next, we apply both sides of (\ref{composition}) to $y_2$.
This yields
\begin{align*}
&\frac{(a,k)(b,l)(c-a-b,m-k-l)}{(c-a,m-k)(c-b,m-l)}y_2(a+k,b+l,c+m,x)\\
&=-\frac{ab}{a+b+1-c}q(x)y_2(a+1,b+1,c+1,x)+r(x)y_2(a,b,c,x)
\end{align*}
(cf. (3.2) and Lemma 2.2 in [Eb1]).
Combining this and (\ref{QR}), we obtain
\begin{align*}
&y_2(a+k,b+l,c+m,x)\\
&=\frac{c(c-a,m-k)(c-b,m-l)}{(c-a-b-1)(c-a-b,m-k-l)(c,m)}Q(x)y_2(a+1,b+1,c+1,x)\\
&+\frac{(c-a,m-k)(c-b,m-l)}{(c-a-b,m-k-l)(c,m)}R(x)y_2(a,b,c,x)
\end{align*}
Finally, 
applying both sides of (\ref{composition}) to $y_i$ $(i=3,4,5,6)$,
we have the following:
\begin{lmm}
We define $Q(x)$ and $R(x)$ as (\ref{three_term}). Then, we have
\begin{align*}
&y_1(a+k,b+l,c+m,x)=Q(x)y_1(a+1,b+1,c+1,x)+R(x)y_1(a,b,c,x),\\
&y_2(a+k,b+l,c+m,x)\\
&=\frac{c(c-a,m-k)(c-b,m-l)}{(c-a-b-1)(c,m)(c-a-b,m-k-l)}Q(x)y_2(a+1,b+1,c+1,x)\\
&+\frac{(c-a,m-k)(c-b,m-l)}{(c,m)(c-a-b,m-k-l)}R(x)y_2(a,b,c,x),\\
&y_3(a+k,b+l,c+m,x)\\
&=\frac{(-1)^{m+1-k-l}c(b,l)(a+1-b,k-l)}{b(c,m)(a+1-c,k-m)}Q(x)y_3(a+1,b+1,c+1,x)\\
&+\frac{(-1)^{m-k-l}(b,l)(a+1-b,k-l)}{(c,m)(a+1-c,k-m)}R(x)y_3(a,b,c,x),\\
&y_4(a+k,b+l,c+m,x)\\
&=\frac{(-1)^{m+1-k-l}c(a,k)(b+1-a,l-k)}{a(c,m)(b+1-c,l-m)}Q(x)y_4(a+1,b+1,c+1,x)\\
&+\frac{(-1)^{m-k-l}(a,k)(b+1-a,l-k)}{(c,m)(b+1-c,l-m)}R(x)y_4(a,b,c,x),\\
&y_5(a+k,b+l,c+m,x)\\
&=\frac{c(1-c)(a,k)(b,l)(2-c,-m)}{ab(c,m)(a+1-c,k-m)(b+1-c,l-m)}Q(x)y_5(a+1,b+1,c+1,x)\\
&+\frac{(a,k)(b,l)(2-c,-m)}{(c,m)(a+1-c,k-m)(b+1-c,l-m)}R(x)y_5(a,b,c,x),\\
&y_6(a+k,b+l,c+m,x)\\
&=\frac{c(a+b-c)(a,k)(b,l)}{ab(c,m)(a+b-c,k+l-m)}Q(x)y_6(a+1,b+1,c+1,x)\\
&+\frac{(a,k)(b,l)}{(c,m)(a+b-c,k+l-m)}R(x)y_6(a,b,c,x).
\end{align*}
\end{lmm}
Using Lemmas 2.1 and 2.2,
we are able to obtain 24 degenerate relations.
For example, substituting
$y_1(a,b,c,x)=(1-x)^{c-a-b}F(c-a,c-b;c;x)$ into 
the top formula in Lemma 2.2,
we obtain
\begin{multline}
(1-x)^{m-k-l}F(c-a+m-k,c-b+m-l;c+m;x)\\
=Q(x)(1-x)^{-1}F(c-a,c-b;c+1;x)+R(x)F(c-a,c-b;c;x).
\label{(m-k,m-l,m)}
\end{multline} 
Therefore, defining 
\begin{gather*}
S^{(n)}:=\frac{1}{R^{(1)}(x)R^{(2)}(x)\cdots R^{(n)}(x)},
\end{gather*}
where $R^{(n)}(x):=R(x)|_{(a,b,c)\mapsto (a+(n-1)k,b+(n-1)l,c+(n-1)m)}$,
we find that
\begin{align*}
&F(c-a,c-b;c;x)\\
&=S^{(n)}(1-x)^{n(m-k-l)}F(c-a+n(m-k),c-b+n(m-l);c+nm;x)
\end{align*}
for an admissible quadruple $(a,b,c,x)$.
The following is obtained similarly from Lemmas 2.1 and 2.2.
\begin{prp}
Fix $(k,l,m)\in \Z ^3$. 
For an admissible quadruple $(a,b,c,x)$,
we obtain the following 24 degenerate relations:
\begin{align*}
&{\rm (i)}\, F(a,b;c;x)=S^{(n)}F(a+nk,b+nl;c+nm;x),\\
&{\rm (ii)}\,F(c-a,c-b;c;x)\notag \\
&=S^{(n)}(1-x)^{n(m-k-l)}F(c-a+n(m-k),c-b+n(m-l);c+nm;x),\\
&{\rm (iii)}\,F(a,c-b;c;x/(x-1))\notag \\
&=S^{(n)}(1-x)^{-nk}F(a+nk,c-b+n(m-l);c+nm;x/(x-1)),\\
&{\rm (iv)}\,F(c-a,b;c;x/(x-1))\notag \\
&=S^{(n)}(1-x)^{-nl}F(c-a+n(m-k),b+nl;c+nm;x/(x-1)),\\
&{\rm (v)}\,F(a,b;a+b+1-c;1-x)\notag \\
&=\frac{(c,nm)(c-a-b,n(m-k-l))}{(c-a,n(m-k))(c-b,n(m-l))}S^{(n)}\notag \\
&\times F(a+nk,b+nl;a+b+1-c+n(k+l-m);1-x),\\
&{\rm (vi)}\,F(a+1-c,b+1-c;a+b+1-c;1-x)\notag \\
&=\frac{(c,nm)(c-a-b,n(m-k-l))}{(c-a,n(m-k))(c-b,n(m-l))}S^{(n)}x^{-nm}\notag \\
&\times F(a+1-c+n(k-m),b+1-c+n(l-m)\\
&\quad ;a+b+1-c+n(k+l-m);1-x),\\
&{\rm (vii)}\,F(a,a+1-c;a+b+1-c;1-x^{-1})\notag \\
&=\frac{(c,nm)(c-a-b,n(m-k-l))}{(c-a,n(m-k))(c-b,n(m-l))}S^{(n)}x^{-nk}\notag \\
&\times F(a+nk,a+1-c+n(k-m);a+b+1-c+n(k+l-m);1-x^{-1}),\\
&{\rm (viii)}\,F(b+1-c,b;a+b+1-c;1-x^{-1})\notag \\
&=\frac{(c,nm)(c-a-b,n(m-k-l))}{(c-a,n(m-k))(c-b,n(m-l))}S^{(n)}x^{-nl}\notag \\
&\times F(b+1-c+n(l-m),b+nl;a+b+1-c+n(k+l-m);1-x^{-1}),\\
&{\rm (ix)}\,F(a,a+1-c;a+1-b;1/x)\notag \\
&=\frac{(-1)^{n(m-k-l)}(c,nm)(a+1-c,n(k-m))}{(b,nl)(a+1-b,n(k-l))}S^{(n)}x^{-nk}\notag \\
&\times F(a+nk,a+1-c+n(k-m);a+1-b+n(k-l);1/x),\\
&{\rm (x)}\,F(1-b,c-b;a+1-b;1/x)\notag \\
&=\frac{(-1)^{n(m-l)}(c,nm)(a+1-c,n(k-m))}{(b,nl)(a+1-b,n(k-l))}S^{(n)}(-x)^{n(l-m)}(1-x)^{n(m-k-l)}\notag \\
&\times F(1-b-nl,c-b+n(m-l);a+1-b+n(k-l);1/x),\\
&{\rm (xi)}\,F(a,c-b;a+1-b;(1-x)^{-1})\notag \\
&=\frac{(-1)^{n(m-l)}(c,nm)(a+1-c,n(k-m))}{(b,nl)(a+1-b,n(k-l))}S^{(n)}(1-x)^{-nk}\notag \\
&\times F(a+nk,c-b+n(m-l);a+1-b+n(k-l);(1-x)^{-1}),\\
&{\rm (xii)}\,F(a+1-c,1-b;a+1-b;(1-x)^{-1})\notag \\
&=\frac{(-1)^{n(m-l)}(c,nm)(a+1-c,n(k-m))}{(b,nl)(a+1-b,n(k-l))}S^{(n)}(-x)^{-nm}(1-x)^{n(m-k)}\notag \\
&\times F(a+1-c+n(k-m),1-b-nl;a+1-b+n(k-l);(1-x)^{-1}),\\
&{\rm (xiii)}\,F(b+1-c,b;b+1-a;1/x)\notag \\
&=\frac{(-1)^{n(m-k-l)}(c,nm)(b+1-c,n(l-m))}{(a,nk)(b+1-a,n(l-k))}S^{(n)}x^{-nl}\notag \\
&\times F(b+1-c+n(l-m),b+nl;b+1-a+n(l-k);1/x),\\
&{\rm (xiv)}\,F(1-a,c-a;b+1-a;1/x)\notag \\
&=\frac{(-1)^{n(m-k)}(c,nm)(b+1-c,n(l-m))}{(a,nk)(b+1-a,n(l-k))}S^{(n)}(-x)^{n(k-m)}(1-x)^{n(m-k-l)}\notag \\
&\times F(1-a-nk,c-a+n(m-k);b+1-a+n(l-k);1/x),\\
&{\rm (xv)}\,F(b,c-a;b+1-a;(1-x)^{-1})\notag \\
&=\frac{(-1)^{n(m-k)}(c,nm)(b+1-c,n(l-m))}{(a,nk)(b+1-a,n(l-k))}S^{(n)}(1-x)^{-nl}\notag \\
&\times F(b+nl,c-a+n(m-k);b+1-a+n(l-k);(1-x)^{-1}),\\
&{\rm (xvi)}\,F(b+1-c,1-a;b+1-a;(1-x)^{-1})\notag \\
&=\frac{(-1)^{n(m-k)}(c,nm)(b+1-c,n(l-m))}{(a,nk)(b+1-a,n(l-k))}S^{(n)}(-x)^{-nm}(1-x)^{n(m-l)}\notag \\
&\times F(b+1-c+n(l-m),1-a-nk;b+1-a+n(l-k);(1-x)^{-1}),\\
&{\rm (xvii)}\,F(a+1-c,b+1-c;2-c;x)\notag \\
&=\frac{(c,nm)(a+1-c,n(k-m))(b+1-c,n(l-m))}{(a,nk)(b,nl)(2-c,-nm)}S^{(n)}x^{-nm}\notag \\
&\times F(a+1-c+n(k-m),b+1-c+n(l-m);2-c-nm;x),\\
&{\rm (xviii)}\,F(1-a,1-b;2-c;x)\notag \\
&=\frac{(c,nm)(a+1-c,n(k-m))(b+1-c,n(l-m))}{(a,nk)(b,nl)(2-c,-nm)}S^{(n)}x^{-nm}(1-x)^{n(m-k-l)}\notag \\
&\times F(1-a-nk,1-b-nl;2-c-nm;x),\\
&{\rm (xix)}\,F(a+1-c,1-b;2-c;x/(x-1))\notag \\
&=\frac{(c,nm)(a+1-c,n(k-m))(b+1-c,n(l-m))}{(a,nk)(b,nl)(2-c,-nm)}S^{(n)}x^{-nm}(1-x)^{n(m-k)}\notag \\
&\times F(a+1-c+n(k-m),1-b-nl;2-c-nm;x/(x-1)),\\
&{\rm (xx)}\,F(b+1-c,1-a;2-c;x/(x-1))\notag \\
&=\frac{(c,nm)(a+1-c,n(k-m))(b+1-c,n(l-m))}{(a,nk)(b,nl)(2-c,-nm)}S^{(n)}x^{-nm}(1-x)^{n(m-l)}\notag \\
&\times F(b+1-c+n(l-m),1-a-nk;2-c-nm;x/(x-1)),\\
&{\rm (xxi)}\,F(c-a,c-b;c+1-a-b;1-x)\notag \\
&=\frac{(c,nm)(a+b-c,n(k+l-m))}{(a,nk)(b,nl)}S^{(n)}(1-x)^{n(m-k-l)}\notag \\
&\times F(c-a+n(m-k),c-b+n(m-l);c+1-a-b+n(m-k-l);1-x),\\
&{\rm (xxii)}\,F(1-a,1-b;c+1-a-b;1-x)\notag \\
&=\frac{(c,nm)(a+b-c,n(k+l-m))}{(a,nk)(b,nl)}S^{(n)}x^{-nm}(1-x)^{n(m-k-l)}\notag \\
&\times F(1-a-nk,1-b-nl;c+1-a-b+n(m-k-l);1-x),\\
&{\rm (xxiii)}\,F(c-a,1-a;c+1-a-b;1-x^{-1})\notag \\
&=\frac{(c,nm)(a+b-c,n(k+l-m))}{(a,nk)(b,nl)}S^{(n)}x^{n(k-m)}(1-x)^{n(m-k-l)}\notag \\
&\times F(c-a+n(m-k),1-a-nk;c+1-a-b+n(m-k-l);1-x^{-1}),\\
&{\rm (xxiv)}\,F(c-b,1-b;c+1-a-b;1-x^{-1})\notag \\
&=\frac{(c,nm)(a+b-c,n(k+l-m))}{(a,nk)(b,nl)}S^{(n)}x^{n(l-m)}(1-x)^{n(m-k-l)}\notag \\
&\times F(c-b+n(m-l),1-b-nl;c+1-a-b+n(m-k-l);1-x^{-1}).
\end{align*}
\end{prp}
\begin{rmk}
From the treatment presented in Section 1,
it is understood that $x$ in an admissible quadruple 
is not a free parameter.
The reason for this is the following.
If $x$ were a free parameter, 
then some of the 24 degenerate relations may be incorrect,
while the special values obtained 
from the remaining correct relations would be trivial.
(Here `trivial' means obvious from the definition of $F(a,b;c;x)$.)
For example, in the case $(k,l,m)=(0,1,1)$,
we have  
\begin{gather} 
F \left( a,b+1;c+1;x \right) ={\frac {a \left( 1-x \right) F \left( a+
1,b+1;c+1;x \right) }{a-c}}-{\frac {cF \left( a,b;c;x \right) }{a-c}}.
\label{(0,1,1)}
\end{gather}
This implies that 
\begin{gather*}
 Q^{(n)}(x)=\frac{a(1-x)}{a-c-n+1},
\end{gather*} 
and, therefore, the $(a,b,c,x)$ satisfying (\ref{condition}) are $(a,b,c,x)=(0,b,c,x)$ and $(a,b,c,1)$.
In the former case,
we find that the 
12 degenerate relations {\rm (xiii)-(xxiv)} are incorrect
because the denominator of the coefficient of $y_i(a+1,b+1,c+1,x)$ $(i=4,5,6)$ in Lemma 2.2
contains $a$ as a factor. 
In addition, the special values obtained from the 12 correct degenerate relations are trivial.
For this reason, we exclude cases in which $x$ is a free parameter from consideration. 
\end{rmk}
\subsection{A complete system of representatives of $G{\backslash}\Z ^3$}
In this subsection, we prove 
that we only need to investigate the lattice points contained in (\ref{quotient})
to obtain the special values of the hypergeometric series.

First, we show that
the 48 lattice points
that form 
the orbit of $(k,l,m)$ under the action of $G$ on $\Z ^3$ (cf. Section 1)  
are equivalent with respect to the obtained special values.
We do this by considering two particular examples,
$(m-k,m-l,m)$ and $(k,m-l,m)$,
and demonstrating that they are equivalent
to $(k,l,m)$ with respect to the special values.
The general result follows by analogy.

To begin with, we consider the case of $(m-k,m-l,m)$.
We start by replacing $(a,b,c,x)$ with $(c-a,c-b,c,x)$.
Then, the three term relation for $(m-k,m-l,m)$ is expressed as
\begin{align*}
&F(c-a+m-k,c-b+m-l;c+m;x)\notag \\
&=Q'(x)F(c-a+1;c-b+1;c+1;x)+R'(x)F(c-a,c-b;c;x),
\end{align*}
and that for $(0,0,1)$ is
\begin{align*}
F(c-a,c-b;c+1;x)
&=\frac{(c-a)(c-b)(1-x)}{ab}F(c-a+1;c-b+1;c+1;x)\\
&+\frac{c(a+b-c)}{ab} F(c-a,c-b;c;x).
\end{align*}
These two three term relations lead to
\begin{align}
&\frac{(c-a)(c-b)(1-x)}{ab}F(c-a+m-k,c-b+m-l;c+m;x)\notag \\
&=Q'(x)F(c-a,c-b,c+1;x)\notag \\
&+\left\{\frac{(c-a)(c-b)(1-x)}{ab}R'(x)-\frac{(c-a)(a+b-c)}{ab}Q'(x)\right\}F(c-a,c-b;c;x).
\label{(m-k,m-l,m)-2}
\end{align}
Equating (\ref{(m-k,m-l,m)}) with (\ref{(m-k,m-l,m)-2}), we have 
\begin{gather}
Q(x)=\frac{ab}{(c-a)(c-b)}(1-x)^{m-k-l}Q'(x).
\label{qq'}
\end{gather}
Hence, the admissible quadruples for $(k,l,m)$ coincide with
those for $(m-k,m-l,m)$.
This implies that 
the special values for $(k,l,m)$
coincide with those for $(m-k,m-l,m)$;
that is,
$(m-k,m-l,m)$ is equivalent to $(k,l,m)$ 
with respect to the obtained special values.

We next show that 
$(k,m-l,m)$ is equivalent to $(k,l,m)$ 
with respect to the obtained special values.
From the top formula in Lemma 2.2, we have
\begin{align}
&(1-x)^{-k}F(a+k,c-b+m-l;c+m;x/(x-1))\notag \\
&=Q(x)(1-x)^{-1}F(a+1,c-b;c+1;x/(x-1))
+R(x)F(a,c-b;c;x/(x-1)).
\label{k,m-l,m}
\end{align}
Also, making the replacement
$(a,b,c,x)\rightarrow (a,c-b,c,x/(x-1))$,
we find that
the three term relation for $(k,m-l,m)$ is expressible as
\begin{align}
&F(a+k,c-b+m-l;c+m;x/(x-1))\notag \\
&=Q''(x)F(a+1,c-b+1;c+1;x/(x-1))
+R''(x)F(a,c-b;c;x/(x-1)).
\label{k,m-l,m-2}
\end{align}
Then, using 
\begin{align*}
F(a+1,c-b,c+1;x/(x-1))&=\frac{(b-c)F(a+1,c-b+1;c+1;x/(x-1))}{b(1-x)}\\
&+\frac{cF(a,c-b;c;x/(x-1))}{b}
\end{align*}
and (\ref{k,m-l,m-2}), we get
\begin{align}
&\frac{b-c}{b(1-x)}F(a+k,c-b+m-l;c+m;x/(x-1))\notag \\
&=Q''(x)F(a+1,c-b;c+1;x/(x-1))\notag \\
&+\left\{\frac{(b-c)}{b(1-x)}R''(x)-\frac{c}{b}Q''(x)\right\}
F(a,c-b;c;x/(x-1)).
\label{k,k-m,m-2}
\end{align}
Finally, equating (\ref{k,m-l,m}) and (\ref{k,k-m,m-2}),
we obtain
\begin{gather}
Q(x)=\frac{b(1-x)^{2-k}}{b-c}Q''(x).
\label{k,k-m,m-3}
\end{gather}
Therefore, the admissible quadruples for $(k,l,m)$ coincide with
those for $(k,m-l,m)$.
This means that 
$(k,m-l,m)$ is equivalent to $(k,l,m)$ 
with respect to the obtained special values.

We can show that 
the other 46 lattice points 
are equivalent to $(k,l,m)$
with respect to the obtained special values analogously.
\begin{rmk}
The reason that we assumed $x\neq 0, 1$ 
in the definition of an admissible quadruple $(a,b,c,x)$ is that
if we do not make this assumption,
there are cases in which
the admissible quadruples for $(k,l,m)$ 
do not coincide with those for the other 47 points 
(cf. (\ref{qq'}), (\ref{k,k-m,m-3}) and Subsection 3.1).
\end{rmk}

Next, we determine 
a complete system of representatives of the quotient of the action of $G$ on $\Z ^3$.
Recall that 
$G=\langle \sigma _1 , \sigma _2\rangle
\ltimes \left(\langle \sigma _3 \rangle 
\times \langle \sigma _4 \rangle 
\times \langle \sigma _5 \rangle\right) 
=S_3 \ltimes \left( S_2 \times S_2 \times S_2 \right)$.
To begin with, we find
a complete system of representatives of the quotient of the action
of $\left(\langle \sigma _3 \rangle 
\times \langle \sigma _4 \rangle 
\times \langle \sigma _5 \rangle\right)$ on $\Z ^3$.
From $\sigma _5$, we see that $(-k,-l,-m)$ is identical to $(k,l,m)$.
Hence, we can assume $0\leq m$.
Further,
because $(m-k,m-l,m)$ is identical to $(k,l,m)$ by $\sigma _4$,
we can assume $0\leq k+l-m$.
In addition, 
because $(l,k,m)$ is identical to $(k,l,m)$ by $\sigma _3$,
we can assume $0\leq l-k$.
Thus, we can take $\{(k,l,m)\in \Z ^3 \mid 0\leq k+l-m, l-k, m\}$,
which we call $D$, as 
a complete system of representatives of the quotient of the action
of $\left(\langle \sigma _3 \rangle 
\times \langle \sigma _4 \rangle 
\times \langle \sigma _5 \rangle\right)$ on $\Z ^3$.
We now proceed to obtain 
a complete system of representatives of the quotient of the action of $G$ on $\Z ^3$.
If $(k,l,m) \in D$, then $D$ also contains 
$(m-k,l,m)$, $(k,l,k+l-m)$, $(l-m,l,k+l-m)$, $(l-m,l,l-k)$ and $(m-k,l,l-k)$.
Then, with
\begin{gather*}
\alpha =k+l-m,\ \beta =l-k,\ \gamma =m \\
{\bf e_1}=(1/2,1/2,0),\ {\bf e_2}=(-1/2,1/2,0), \ {\bf e_3}=(1/2,1/2,1),
\end{gather*}
we have the following: 
\begin{align*}
(k,l,m)&=\alpha {\bf e_1}+\beta {\bf e_2}+\gamma {\bf e_3},\
(m-k,l,m)=\alpha {\bf e_2}+\beta {\bf e_1}+\gamma {\bf e_3},\\
(k,l,k+l-m)&=\alpha {\bf e_3}+\beta {\bf e_2}+\gamma {\bf e_1},\
(l-m,l,k+l-m)=\alpha {\bf e_3}+\beta {\bf e_1}+\gamma {\bf e_2},\\
(l-m,l,l-k)&=\alpha {\bf e_1}+\beta {\bf e_3}+\gamma {\bf e_2},\
(m-k,l,l-k)=\alpha {\bf e_2}+\beta {\bf e_3}+\gamma {\bf e_1}.
\end{align*}
Together, these imply
that we are able to adopt (\ref{quotient})
as 
a complete system of the quotient of the action of $G$ on $\Z ^3$.
Thus, we only have to investigate the lattice points contained in (\ref{quotient})
to obtain the special values of the hypergeometric series.

\section{Derivation of special values}
In this section, we derive the special values for some lattice points as examples.

Although we implicitly assumed 
that the parameter $c$ of $F(a,b;c;x)$
is an element of $\C \setminus \{0,-1,-2,\cdots \}$,
for the rest of paper,
we expand the definition of $F(a,b;c;x)$.
From this point, even if the parameter $c$ is a non-positive integer,
we define $F(a,b;c;x)$ as follows 
if the parameter $a$ is a non-positive integer satisfying $c<a$: 
\begin{gather*}
F(a,b;c;x):=\sum _{n=0} ^{|a|}\frac{(a,n)(b,n)}{(c,n)(1,n)}x^n.
\end{gather*}
\subsection{Example 1: $(k,l,m)=(0,1,1)$}
We first consider the case $m=1$.
Then, the only $(k,l,m)$ satisfying (\ref{quotient}) is $(0,1,1)$.
For this point, the corresponding three term relation is (\ref{(0,1,1)}),
and this leads to
\begin{gather*}
 Q^{(n)}(x)=\frac{a(1-x)}{a-c-n+1}.
\end{gather*} 
Hence, $(a,b,c,x)=(a,b,c,1)$ satisfies (\ref{condition}) 
(and so does $(0,b,c,x)$,
but for the reason explained in Remark 2.4,
we do not consider this case).

In the case $(a,b,c,x)=(a,b,c,1)$, because 
\begin{gather*}
S^{(n)}=\frac{(c-a,n)}{(c,n)},
\end{gather*}
we have 
\begin{gather*}
F(a,b;c;1)=\frac{(c-a,n)}{(c,n)}F(a,b+n;c+n;1).
\label{chu2}
\end{gather*}
From this, substituting $b=-n$, where $n\in \Z _{\geq 0}$, we get the Chu-Vandermonde equality, (\ref{chu}).
Note that this equality holds even if $c$ is a non-positive integer with $c<b=-n$.

Because $(a,b,c,1)$ is not an admissible quadruple,
we should investigate the other 47 points.
For example, considering the case $(k,l,m)=(0,0,1)$, we get the Gauss summation formula, (\ref{gauss}).
Similarly, we only get the Chu-Vandermonde equality and the Gauss summation formula from the other 46 points.

\subsection{Example 2: $(k,l,m)=(1,2,2)$}
In the case $(k,l,m)=(1,2,2)$, because
\begin{align*}
&F \left( a+1,b+2,c+2,x \right)\\
&={\frac { \left( c+1 \right)  \left( -c
+xa \right) F \left( a+1,b+1,c+1,x \right) }{x \left( a-c \right) 
 \left( b+1 \right) }}+{\frac {c \left( c+1 \right) F \left( a,b,c,x
 \right) }{x \left( a-c \right)  \left( b+1 \right) }},
\end{align*}
the numerator of $Q^{(n)}(x)$ is
\begin{gather}
 \left( -4+2\,x \right) {n}^{2}+ \left(  \left( c+2\,a-3 \right) x-4\,
c+6 \right) n+ \left( 1-c+ca-a \right) x-{c}^{2}+3\,c-2.
\label{consti_method}
\end{gather}
Therefore, the admissible quadruple $(a,b,c,x)$ is $(a,b,2a,2)$, and 
\begin{gather*}
S^{(n)}={\frac { \left( -1 \right) ^{n} \left( 1/2\,b+1/2,n
 \right) }{ \left( a+1/2,n \right) }},
\end{gather*}
where we denote $\frac{1}{2}b$ by $1/2\,b$, 
and we also denote thus for the rest of this article.
Hence, {\rm (i)} in Proposition 2.3 leads  to
\begin{align}
{\mbox{$F$}(a,b;\,2\,a;\,2)}={\frac { \left( -1 \right) ^{n} \left( 1/2\,b+1/2,n
 \right) }{
 \left( a+1/2,n \right) }{\mbox{$F$}(a+n,b+2\,n;\,2\,a+2\,n;\,2)}}
\label{(1,2,2)}
\end{align}
Although (\ref{(1,2,2)}) is valid by virtue of analytic continuation,
this equality regarded as 
an infinite series expression does not make sense.  
For this reason, we carry out the degeneration of this 
into a finite series expression that does make sense.
This is done separately in the following cases:
$a=-n-1$, $b=-2\, n$ and $b=-2\, n-1$,
where $n\in \Z _{\geq 0}$.
If $a=-n-1$, then the right hand side of (\ref{(1,2,2)}) becomes
\begin{gather*}
 {\frac { \left( -1 \right) ^{n} \left( 1/2\,b+1/2,n
 \right)  F(-1,b+2\, n;-2;2) }{ \left( -1/2-n,n
 \right) }}
={\frac { \left( 2\,n+b+1 \right)  \left( 1/2\,b+1/2,n
 \right) }{ \left( 3/2,n \right) }}.
\end{gather*}
If $b=-2\, n$, then the right hand side of (\ref{(1,2,2)}) is
\begin{gather*}
\dfrac{(1/2,n)}{(a+1/2,n)}.
\end{gather*}
If $b=-2\, n-1$, then the right hand side of (\ref{(1,2,2)}) is equal to 
\begin{gather*}
\dfrac{\left(1,n\right)F\left(a+n, -1; 2\,a+2\,n; 2\right)}
{\left(a+1/2,n\right)}=0.
\end{gather*}
The special values obtained from {\rm (ii), (iii) and (iv)} are identical to the above.

Next, using {\rm (v)} in Proposition 2.3, we have
\begin{align}
{\mbox{$F$}(a,b;\,b+1-a;\,-1)}=
{\frac {{2}^{2\,n}\left( 1/2\,b+1/2,n \right) 
}{ \left( b+1-a,n \right) }{\mbox{$F$}(a+n,b+2\,n;\,b+1-a+n;\,-1)}}.
\label{(1,2,2)-v}
\end{align}
Substituting $a=-n$ into (\ref{(1,2,2)-v}),
this becomes
\begin{align}
{\mbox{$F$}(-n,b;\,b+1+n;\,-1)}
={\frac {{2}^{2\, n}\Gamma  \left( 1/2\,b+1/2+n \right) \Gamma  \left( n+b+
1 \right) }{\Gamma  \left( 1/2\,b+1/2 \right) \Gamma  \left( 2\,n+b+1
 \right) }}.
 \label{(1,2,2)-v1}
\end{align}
Note that the validity of (\ref{(1,2,2)-v1}) for any integer $n$
follows from
the constitution method of $(a,b,c,x)$ satisfying (\ref{condition})
(cf. (\ref{consti_method})).
The following lemma is known (cf. 5.3 in [Ba]):
\begin{lmm}
(Carlson's theorem)
We assume that $f(z)$ and $g(z)$ are regular and of the form $O(e^{k|z|})$, where $k<\pi$, for $\Re z \geq 0$,
and $f(z)=g(z)$ for $z=0,1,2\cdots$. 
Then, $f(z)=g(z)$ on $\{z\mid \,\Re z \geq 0\}$.
\end{lmm}
It is easily confirmed that both sides of (\ref{(1,2,2)-v1}) satisfy the assumption of the above lemma.
Hence, (\ref{(1,2,2)-v1}) holds for any complex number $n$ 
for which the left hand side is meaningful.
Resultingly, substituting $n=-a$ into (\ref{(1,2,2)-v1}), we find
\begin{align*}
{\mbox{$F$}(a,b;\,b+1-a;\,-1)}
={\frac {{2}^{-2a}\Gamma  \left( 1/2\,b+1/2-a \right) \Gamma  \left( b+
1-a \right) }{\Gamma  \left( 1/2\,b+1/2 \right) \Gamma  \left( b+1-2\,a
 \right) }}.
\end{align*}

Now, we consider the case $a$, $b+1-a\in \Z _{\leq 0}$ with $b+1-a<a$.
Specifically, 
we consider the case $a=-n, b=-2\,n-n_1-2$, where $n$, $n_1\in \Z _{\geq 0}$.
Then, (\ref{(1,2,2)-v}) leads to
\begin{align*}
F(a,b;b+1-a;-1)={\dfrac {{2}^{ 2\,n} \left( 1/2\, n_1+3/2,
n \right) }{ \left( n_1+2, n \right) }}.
\end{align*}

Next, we consider the case $b$, $b+1-a\in \Z _{\leq 0}$ with $b+1-a<b$.
Then, if $a=n_1+2$, $b=-2\,n$, where $n_1$, $n\in \Z _{\geq 0}$,
(\ref{(1,2,2)-v}) implies 
\begin{gather*}
\mbox{$F$}(a,b;\,b+1-a;\,-1)=\dfrac{(1/2,n)(n_1+2,n)}{(1/2\,n_1+1,n)(1/2\,n_1+3/2,n)}.
\end{gather*}
If $a=n_1+2$, $b=-2\,n-1$, where $n_1$, $n\in \Z _{\geq 0}$,
we find
\begin{gather*}
\mbox{$F$}(a,b;\,b+1-a;\,-1)=0
\end{gather*}
from (\ref{(1,2,2)-v}).

We are able to obtain special values from the other 19 degenerate relations similarly.
Those values are tabulated in the following section.

\subsection{Example 3: $(k,l,m)=(1,2,3)$}
For an admissible quadruple $(a,2a-1/3,3a,9)$, we have
\begin{gather*}
S^{(n)}={\frac { (-4)^n\left( a+1/2,n \right) }{
 \left( a+2/3,n \right) }}.
\end{gather*}
Now, we consider the special values obtained from $\rm (vii)$ in Proposition 2.3.
Because $\rm (vii)$ in this case is 
\begin{gather*}
{\mbox{$F$}(a,1-2\,a;\,2/3;\,8/9)}
=(-3)^n\,{\mbox{$F$}(a+n,-2\,a+1-2\,n;\,2/3;\,8/9)},
\label{(1,2,3)-7}
\end{gather*}
we find
\begin{align*}
{\mbox{$F$}(a,1-2\,a;\,2/3;\,8/9)}=
\begin{cases}
(-3)^n 
\quad {\text {if $a=-n$}}, \\
(-3)^{-n}
\quad {\text {if $a=n+1/2$}}, \\
(-3)^{-n-1}
\quad {\text {if $a=n+1$}}.
\end{cases}
\end{align*}
However, we can not directly apply Lemma 3.1 to the above identity,
because its left hand side does not satisfy the assumption of that lemma.
For this reason, we get the special value of ${\mbox{$F$}(a,1-2\,a;\,2/3;\,8/9)}$ using 
the following algebraic transformation of the hypergeometric series:
\begin{align}
F\left(a,1-2\,a;\,2/3;\,x\right)=\left( 1-x \right) ^{1/2\,a-1/4}
{F\left(1/2\,a-1/12,1/4-1/2\,a;\,2/3;
\,u(x)\right)},
\label{algebraic}
\end{align}
with
\begin{gather*}
u(x)={\frac {1}{64}}\,{\frac {x \left( -8+9\,x \right) ^{3}}{x-1}}
\end{gather*}
(see formula (130) in [Gour]).
We remark that this holds near $x=0$.
Now, we carry out an analytic continuation 
of each side of (\ref{algebraic}) 
along a curve starting at $x=0$ and ending at $x=8/9$,
as depicted in Figure 1.
\begin{figure}[htbp]
 \begin{minipage}{0.5\hsize}
  \begin{center}
   \includegraphics[width=70mm]{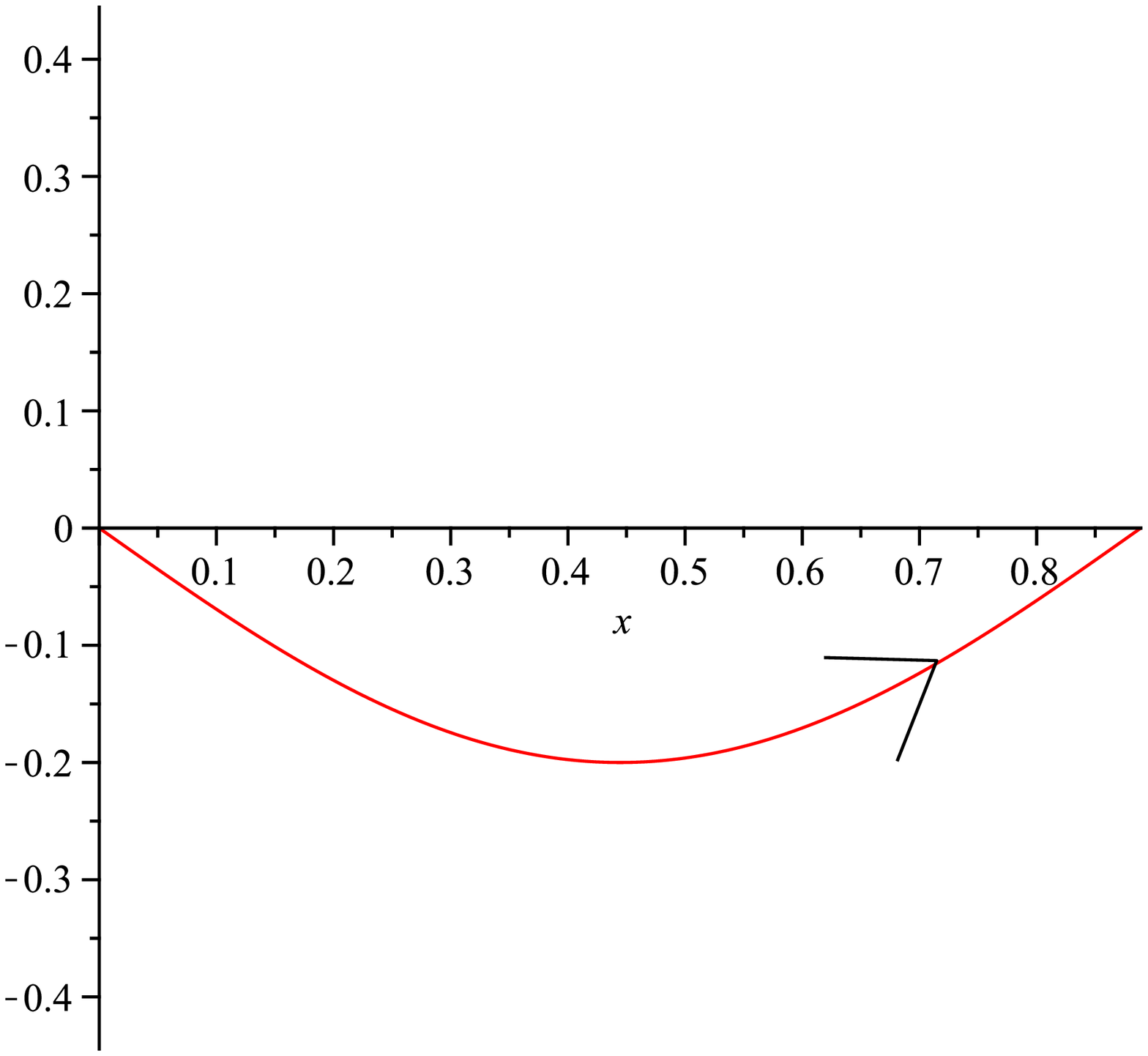}
  \end{center}
  \caption{$x$-plane}
  \label{fig:one}
 \end{minipage}
 \begin{minipage}{0.5\hsize}
  \begin{center}
   \includegraphics[width=70mm]{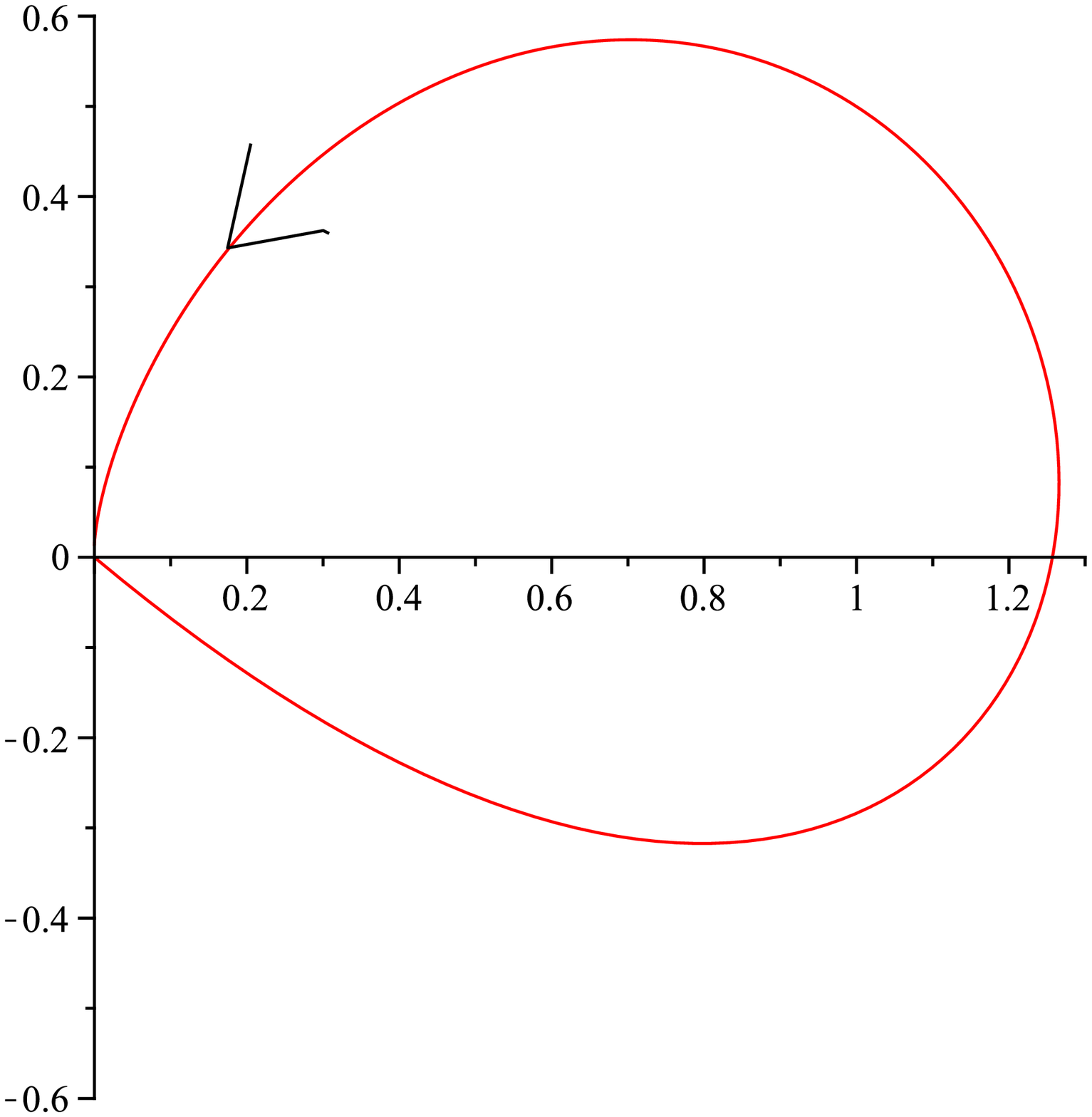}
  \end{center}
  \caption{$u$-plane}
  \label{fig:two}
 \end{minipage}
\end{figure}
Noting that $u(x)$ encircles $u=1$ once counterclockwise,
as shown in Figure 2,
we obtain  
\begin{gather*}
F\left(a,1-2\,a;\,2/3;\,8/9\right)= 2\cdot {3}^{-a}\sin \left(  \left( 5/6-a \right) \pi  \right)
\end{gather*}
using the following lemma
(cf. Theorem 4.7.2 in [IKSY]):
\begin{lmm}
Let $\gamma _1$ be a loop 
starting and ending at $x=a$, where $0<a<1$, 
and encircling $x=1$ once in the counterclockwise direction.
Then, analytic continuation of
$(y_1(a,b,c,x),y_5(a,b,c,x))$  
along $\gamma _1$ is given by
\begin{gather*}
(y_1(a,b,c,x),y_5(a,b,c,x))P^{-1}AP,
\end{gather*}
where
\begin{gather*}
P= \left[ \begin {array}{cc} {\frac {\Gamma  \left( c \right) \Gamma 
 \left( c-a-b \right) }{\Gamma  \left( c-a \right) \Gamma  \left( c-b
 \right) }}&{\frac {\Gamma  \left( 2-c \right)  \Gamma \left( c-a-b \right) }
{\Gamma  \left( 1-a \right) \Gamma  \left( 1-b \right) }}
\\ \noalign{\medskip}{\frac {\Gamma  \left( c \right) \Gamma  \left( a
+b-c \right) }{\Gamma  \left( a \right) \Gamma  \left( b \right) }}&{
\frac {\Gamma  \left( a+b-c \right) \Gamma  \left( 2-c \right) }{
\Gamma  \left( a+1-c \right) \Gamma  \left( b+1-c \right) }}
\end {array} \right],\
A=\left[ \begin {array}{cc} 1&0\\ \noalign{\medskip}0&{{\rm e}^{2\,i
\pi \, \left( c-a-b \right) }}\end {array} \right] .
\end{gather*}
\end{lmm}
\section{Tables of special values}
In this section, using Proposition 3.2.3,
we tabulate the special values for $(k,l,m)$ satisfying $0\leq k+l-m \leq l-k \leq m\leq 6$.
However, we exclude the Chu-Vandermonde  equality, (\ref{chu})
and the Gauss summation formula, (\ref{gauss}).

Let $n$, $n_1$, $n_2\in \Z _{\geq 0}$.
Then, recall that
\begin{gather}
F(a,b;c;x):=\sum _{n=0} ^{|a|}\frac{(a,n)(b,n)}{(c,n)(1,n)}x^n
\label{extend_hgf}
\end{gather}
when $a$ and $c$ are non-positive integers satisfying $c< a$.
Further, recall that we denote 
$\frac{1}{2}a$ by $1/2\,a$, and so on.

For a given $(k,l,m)$, 
we use the expression {\rm (**)}$\leq${\rm (*)} to indicate that
the values obtained from {\rm (**)} in Proposition 3.2.3 
coincide with or are contained in those obtained from {\rm (*)}.
\subsection{$m=1$}
\subsubsection{$(k,l,m)=(0,1,1)$}
In this case, there is no admissible quadruple
(cf. Subsection 3.1).

\subsection{$m=2$}

\subsubsection{$(k,l,m)=(1,2,2)$}
In the case that $(k,l,m)=(1,2,2)$, we obtain
\begin{gather}
 (a,b,c,x)=(a,b,2\, a,2),\
 S^{(n)}=\frac{(-1)^n(1/2\,b+1/2,n)}{(a+1/2,n)}.
\tag{1,2,2-1}
\end{gather}
\paragraph{(1,2,2-1)}  
\begin{flalign*}
&{\rm {(i)}}\ F(a,b;2\,a;2)=
\begin{cases}
{\dfrac { \left( b+1+2\,n \right) \left( 1/2\,b+1/2,n
 \right) }{ \left( 3/2,n \right) }}
 \quad {\text {if $a=-n-1$}},\\
 \dfrac{(1/2,n)}{(a+1/2,n)}\quad {\text {if $b=-2\, n$}},\\
0 \quad {\text {if $b=-2\,n-1$}}
\end{cases}&
\end{flalign*}
(The first case is identical to (4.11) in [Ge]).
We find that {\rm {(ii), (iii), (iv)} $\leq$ {\rm (i)}.
\begin{flalign*}
&{\rm {(v)}}\ F(a,b;b+1-a;-1)&\\
&=
\begin{cases}
{\dfrac {{2}^{-2\,a}\Gamma  \left( 1/2\,b+1/2-a \right) \Gamma  \left( b+1-a
 \right) }{\Gamma  \left( 1/2\,b+1/2 \right) \Gamma  \left( b+1-2\,
a \right) }},\\
\dfrac{(1/2,n_2)(n_1+2,n_2)}{(1/2\,n_1+1,n_2)(1/2\,n_1+3/2,n_2)}
\quad {\text {if $a=n_1+2$, $b=-2\,n_2$}},\\
0 
\quad {\text {if $a=n_1+2$, $b=-2\,n_2-1$}},\\
{\dfrac {{2}^{ 2\,n_1} \left( 1/2\, n_2+3/2,
n_1 \right) }{ \left( n_2+2, n_1 \right) }}\quad {\text{if $a=-n_1, b=-2\,n_1-n_2-2$}}
\end{cases}&
\end{flalign*}
(The first case is identical to 2.8(47) in [Erd]). We find that {\rm (vi)}$\leq${\rm (v)}.
\begin{flalign*}
&{\rm (vii)}\ {\mbox{$F$}(a,1-a;\,b+1-a;\,1/2)}&\\
&=\begin{cases}
{\dfrac {{2}^{-a}\Gamma  \left( 1/2\,b+1/2-a \right) \Gamma  \left( b+1-a
\right) }{\Gamma  \left( 1/2\,b+1/2 \right) \Gamma  \left( b+1-2\,
a \right) }},\\
{\dfrac {{2}^{ n_1} \left( 1/2\, n_2+3/2, 
n_1 \right) }{ \left(n_2+2, n_1 \right) }}
\quad {\text {if $a=-n_1$, $b=-2\,n_1-n_2-2$ }},\\
{\dfrac { \left(  n_1+ n_2+2, n_1 \right) }{ 2^{n_1}\left( 1/2\, n_2+1, n_1
 \right) }}\ 
{\text {if $a=n_1+1$, $b=-n_2-1$}}
\end{cases}&
\end{flalign*}
(The first case is identical to 2.8(51) in [Erd]).
\begin{flalign*}
&{\rm (viii)}\ {\mbox{$F$}(b,b+1-2\,a;\,b+1-a;\,1/2)}&\\
&=\begin{cases}
{\dfrac {{2}^{b-2\,a}\Gamma  \left( 1/2\,b+1/2-a \right) \Gamma 
 \left( b+1-a \right) }{\Gamma  \left( 1/2\,b+1/2 \right) \Gamma 
 \left( b+1-2\,a \right) }},\\
{\dfrac { \left( 1/2, n_2 \right) 
 \left( n_1+2,n_2 \right) }{{2}^{ 2\,n_2}
 \left( 1/2\,n_1+3/2,n_2 \right)  \left( 1/2
\, n_1+1,n_2 \right) }}
\quad {\text {if $a=n_1+2$, $b=-2\,n_2$ }},\\
0 \quad {\text {if $a=n_1+2$, $b=-2\,n_2-1$ }}
\end{cases}&
\end{flalign*}
(The first case is identical to 2.8(50) in [Erd]).
The special values obtained from {\rm (ix)-(xxiv)} 
are contained in the above. 

\subsubsection{$(k,l,m)=(0,2,2)$}
When $(k,l,m)=(0,2,2)$, we obtain 
\begin{gather}
\begin{cases}
(a,b,c,x)=(a,b,b+1-a,-1),\,\\
S^{(n)}={\dfrac {\left( 1/2\,b+1/2,n \right)
 \left( 1/2\,b+1-a,n \right) }{ \left( 1/2\,b+1/2-1/2
\,a,n \right) \left( 1/2\,b+1-1/2\,a,n \right) }}.
\tag{0,2,2-1}
\end{cases}
\end{gather}
\paragraph{(0,2,2-1)}
The special values obtained from (0,2,2-1) are contained in those obtained from (1,2,2-1).

\subsubsection{$(k,l,m)=(1,1,2)$}
In this case, there is no admissible quadruple.

\subsubsection{$(k,l,m)=(1,3,2)$}
In this case, we have
\begin{align}
(a,b,c,x)&=(a,3\,a-1,2\,a,1/2+1/2\,i\sqrt {3}), \ S^{(n)}=
{\frac { \left( -3/4\,i\sqrt {3} \right) ^{n} \left( a
+1/3,n \right) }{ \left( a+1/2,n \right) }}\tag {1,3,2-1},\\
(a,b,c,x)&=(a,3\,a-1,2\,a,1/2-1/2\,i\sqrt {3}), \ S^{(n)}=
{\frac { \left( 3/4\,i\sqrt {3} \right) ^{n} \left( a
+1/3,n \right) }{ \left( a+1/2,n \right) }}.\tag {1,3,2-2}
\end{align}
\paragraph{(1,3,2-1)}
\begin{flalign*}
&{\rm (i)}\, {\mbox{$F$}(a,3\,a-1;\,2\,a;\,1/2+1/2\,i\sqrt {3})}&\\
&=\begin{cases}
{\dfrac {{2}^{2\,a-2/3}{{\rm e}^{1/6\,i \left( 3\,a-1 \right) \pi }}
\Gamma  \left( 2/3 \right) \Gamma  \left( a+1/2 \right) }{{3}^{3/2\,a-
1/2} \Gamma  \left( 5/6 \right)\Gamma  \left( a+1/3 \right) }},\\
 {\dfrac { \left( -i \right) ^{n+1}{3}^{3/2\,n+1/2}
 \left( 5/3,n \right) }{{2}^{2\,n} \left( 3/2,n
 \right) }}
\quad {\text {if $a=-1-n$}}
\end{cases}&
\end{flalign*}
(The first case is identical to 2.8(55) in [Erd]
and  the second case is identical to Theorem 11 in [Ek]).
\begin{flalign*}
&{\rm (ii)}\, {\mbox{$F$}(a,1-a;\,2\,a;\,1/2+1/2\,i\sqrt {3})}&\\
&=\begin{cases}
{\dfrac {{2}^{2\,a-2/3}{{\rm e}^{1/6\,i \left( 1-a \right) \pi }}
\Gamma  \left( 2/3 \right) \Gamma  \left( a+1/2 \right) }{{3}^{3/2\,a-
1/2}\Gamma  \left( 5/6 \right)\Gamma  \left( a+1/3 \right)  }}
,\\
{\dfrac { -\left( -i \right) ^{n+1}{3}^{3/2\,n+1/2} \left( 1+i\sqrt {3} \right) ^
{2\,n-1} \left( 5/3,n \right) }{{2}^{4
\,n-1}\left( 3/2,n \right) }}
\quad {\text {if $a=-1-n$}}.
\end{cases}&
\end{flalign*}
\begin{flalign*}
&{\rm (iii)}\, {\mbox{$F$}(a,1-a;\,2\,a;\,1/2-1/2\,i\sqrt {3})}&\\
&=\begin{cases}
{\dfrac {{2}^{2\,a-2/3}{{\rm e}^{1/6\,i \left( a-1 \right) \pi }}
\Gamma  \left( 2/3 \right) \Gamma  \left( a+1/2 \right) }{{3}^{3/2\,a-
1/2}\Gamma  \left( 5/6 \right)\Gamma  \left( a+1/3 \right)  }},
\\
{\dfrac { -i ^{n+1}{3}^{3/2\,n+1/2} \left( 1-i\sqrt {3} \right) ^
{2\,n-1} \left( 5/3,n \right) }{{2}^{4
\,n-1}\left( 3/2,n \right) }}
\quad {\text {if $a=-1-n$}}.
\end{cases}&
\end{flalign*}
\begin{flalign*}
&{\rm (iv)}\,{\mbox{$F$}(a,3\,a-1;\,2\,a;\,1/2-1/2\,i\sqrt {3})}&\\
&=\begin{cases}
{\dfrac {{2}^{2\,a-2/3}{{\rm e}^{1/6\,i \left( 1-3\,a \right) \pi }}
\Gamma  \left( 2/3 \right) \Gamma  \left( a+1/2 \right) }{{3}^{3/2\,a-
1/2}\Gamma  \left( 5/6 \right) \Gamma  \left( a+1/3 \right) }}
,\\
{\dfrac {{i}^{n+1}{3}^{3/2\,n+1/2} \left( 5/3,n
 \right) }{{2}^{2\,n} \left( 3/2,n \right) }}
\quad {\text {if $a=-1-n$}}
\end{cases}&
\end{flalign*}
(The first case is identical to 2.8(56) in [Erd]
and the second case is identical to Theorem 11 in [Ek]).
The special values obtained from {\rm (v)-(xxiv)} are contained in the above.
\paragraph{(1,3,2-2)}
The special values obtained from (1,3,2-2) coincide with
those obtained from (1,3,2-1). 
\subsection{$m=3$}
\subsubsection{$(k,l,m)=(0,3,3)$}
In this case, we get
\begin{align}
&(a,b,c,x)=(1,b,b,-1/2+1/2\,i\sqrt {3}),\, S^{(n)}=1,\tag{0,3,3-1}\\
&(a,b,c,x)=(1,b,b,-1/2-1/2\,i\sqrt {3}),\, S^{(n)}=1,\tag{0,3,3-2}\\
&(a,b,c,x)=(0,b,b+1,-1/2+1/2\,i\sqrt {3}),\, S^{(n)}=1,\tag{0,3,3-3}\\
&(a,b,c,x)=(0,b,b+1,-1/2-1/2\,i\sqrt {3}),\, S^{(n)}=1.\tag{0,3,3-4}
\end{align}
\paragraph{(0,3,3-1)}
The special values obtained from (0,3,3-1) except trivial values
are special cases of 
\begin{gather}
F(1,b;2;x)=\frac{(1-x)^{1-b}-1}{(b-1)x}.
\label{F(1,b;2;x)}
\end{gather}

\paragraph{(0,3,3-2), (0,3,3-3), (0,3,3-4)}
The special values obtained from (0,3,3-2), (0,3,3-3) and (0,3,3-4) coincide with those obtained from (0,3,3-1). 
\subsubsection{$(k,l,m)=(1,2,3)$}
In this case, we have
\begin{align}
&(a,b,c,x)=(a,2\,a-1/3,3\,a,9),\,
S^{(n)}={\frac { (-4)^n\left( a+1/2,n \right) }{
 \left( a+2/3,n \right) }}\tag{1,2,3-1},\\
&(a,b,c,x)=(a,2\,a-2/3,3\,a-1,9),\,
S^{(n)}={\frac { (-4)^n\left( a+1/2,n \right) }{
 \left( a+2/3,n \right) }}\tag{1,2,3-2}.
\end{align}
\paragraph{(1,2,3-1)}
\begin{flalign*}
&{\rm (i)}\, {\mbox{$F$}(a,2\,a-1/3;\,3\,a;\,9)}
=\begin{cases}
{\dfrac { \left( -6 \right)({-4})^{n} \left( 3/2
,n \right) }{\left( 4/3,n \right) }}
\quad {\text {if $a=-1-n$}},\\
{\dfrac { \left( -4 \right) ^{n} \left( 1/3,n \right) }
{ \left( 1/6,n \right) }}
\quad {\text {if $a=1/6-n$}},\\
{\dfrac {10 \left( -4 \right) ^{n} \left( 11/6,n
 \right) }{ \left( 5/3,n \right) }}
\quad {\text {if $a=-4/3-n$}}
\end{cases}&
\end{flalign*}
(The first case is identical to Theorem 17 in [Ek]).
\begin{flalign*}
&{\rm (ii)}\, {\mbox{$F$}(2\,a,a+1/3;\,3\,a;\,9)}
=\begin{cases}
{\dfrac { (-6)\left( -4 \right) ^{n} \left( 3/2,n
 \right) }{ \left( 4/3,n \right) }}
\quad {\text {if $a=-1-n$}},\\
0
\quad {\text {if $a=-1/2-n$}},\\
{\dfrac { \left( -4 \right) ^{n} \left( 5/6,n \right) }
{\left( 2/3,n \right) }}
\quad {\text {if $a=-1/3-n$}}
\end{cases}&
\end{flalign*}
(The third case is identical to Theorem 18 in [Ek]).
\begin{flalign*}
&{\rm (iii)}\, {\mbox{$F$}(a,a+1/3;\,3\,a;\,9/8)}
=\begin{cases}
{\dfrac {3\left( 3/2,n \right) }{ {2}^{n+2} \left( 4/3,n \right) }}
\quad {\text {if $a=-1-n$}},\\
{\dfrac { \left( 5/6,n \right) }{{2}^{n} \left( 2/3,n \right) }}
\quad {\text {if $a=-1/3-n$}}.
\end{cases}&
\end{flalign*}
\begin{flalign*}
&{\rm (iv)}\, {\mbox{$F$}(2\,a,2\,a-1/3;\,3\,a;\,9/8)}
=\begin{cases}
{\dfrac {3}{2}}\,{\dfrac {\left( 3/2,n \right) }{
 \left( -16 \right) ^{n+1} \left( 4/3,n \right) }}
\quad {\text {if $a=-1-n$}},\\
0
\quad {\text {if $a=-1/2-n$}},\\
{\dfrac { \left( 1/3,n \right) }{ \left( -16 \right) ^{
n} \left( 1/6,n \right) }}
\quad {\text {if $a=1/6-n$}},\\
{\dfrac { -5\left( 11/6,n \right) }{
 \left( -16 \right) ^{n+2} \left( 5/3,n \right) }}
\quad {\text {if $a=-4/3-n$}}.\\
\end{cases}&
\end{flalign*}
\begin{flalign*}
&{\rm (v)}\, {\mbox{$F$}(a,2\,a-1/3;\,2/3;\,-8)}
=\begin{cases}
(-3)^{3\,n}
\quad {\text {if $a=-n$}},\\
(-3)^{3\,n}
\quad {\text {if $a=1/6-n$}},\\
\left( -3 \right) ^{3\,n+1}
\quad {\text {if $a=-1/3-n$}}
\end{cases}&
\end{flalign*}
(The second case is identical to (3.12) in [GS] and
the first case is identical to (5.23) in [GS]).
We find that {\rm (vi)}$\leq ${\rm (v)}.
\begin{flalign*}
&{\rm (vii)}\, F\left(a,1-2\,a;\,2/3;\,8/9\right)= 2\cdot {3}^{-a}\sin \left(  \left( 5/6-a \right) \pi  \right)&
\end{flalign*}
(The above is identical to (8/9.1) in [Gos] and (3.2) in [Ka]).
We derive (vii) using (\ref{algebraic}) and 
connection formulae for the hypergeometric series(see (25)--(44) in 2.9 in [Erd]).
We find that {\rm (viii)}$\leq ${\rm (vii)}.
\begin{flalign*}
&{\rm (ix)}\, {\mbox{$F$}(a,1-2\,a;\,4/3-a;\,1/9)}=
{\dfrac {{3}^{-a}\sqrt {\pi }\Gamma  \left( 4/3-a \right) }{{2}^{1/3-2
\,a}\Gamma  \left( 2/3 \right) \Gamma  \left( 7/6-a \right) }}&
\end{flalign*}
(The above is identical to (1/9.4) in [Gos] and (1.2) in [Ka]).
\begin{flalign*}
&{\rm (x)}\, {\mbox{$F$}(4/3-2\,a,a+1/3;\,4/3-a;\,1/9)}=
{\dfrac {{3}^{2/3-a}\sqrt {\pi }\Gamma  \left( 4/3-a \right) }{{2}^{4/3-2
\,a}\Gamma  \left( 2/3 \right) \Gamma  \left( 7/6-a \right) }}&
&
\end{flalign*}
(The above is identical to (1/9.5) in [Gos] and (1.3) in [Ka]).
\begin{flalign*}
&{\rm (xi)}\, {\mbox{$F$}(a,a+1/3;\,4/3-a;\,-1/8)}=
{\dfrac {{3}^{-3\,a}\sqrt {\pi }\Gamma  \left( 4/3-a \right) }{{2}^{1/3-5
\,a}\Gamma  \left( 2/3 \right) \Gamma  \left( 7/6-a \right) }}&.
\end{flalign*}
\begin{flalign*}
&{\rm (xii)}\, {\mbox{$F$}(1-2\,a,4/3-2\,a;\,4/3-a;\,-1/8)}=
{\dfrac {{3}^{3\,a-2}\sqrt {\pi }\Gamma  \left( 4/3-a \right) }{{2}^{4
\,a-8/3}\Gamma  \left( 2/3 \right) \Gamma  \left( 7/6-a \right) }}&
\end{flalign*}
(The above is a generalization of Theorem 32 in [Ek]).
We find {\rm (xiii)}$\leq${\rm (ix)}, {\rm (xiv)}$\leq${\rm (x)},
{\rm (xv)}$\leq${\rm (xii)}, {\rm (xvi)}$\leq${\rm (xi)},
{\rm (xvii)}$\leq${\rm (i)}, {\rm (xviii)}$\leq${\rm (ii)},
{\rm (xix)}$\leq${\rm (iv)}, {\rm (xx)}$\leq${\rm (iii)}.
\begin{flalign*}
&{\rm (xxi)}\, {\mbox{$F$}(2\,a,a+1/3;\,4/3;\,-8)}
=\begin{cases}
{\dfrac { \left( -3 \right) ^{3\,n} \left( 1/2,n \right) 
}{ \left( 7/6,n \right) }}
\quad {\text {if $a=-n$}},\\
0
\quad {\text {if $a=-1/2-n$}},\\
{\dfrac { \left( -3 \right) ^{3\,n} \left( 5/6,n \right) 
}{ \left( 3/2,n \right) }}
\quad {\text {if $a=-1/3-n$}}\\
\end{cases}&
\end{flalign*}
(The third case is identical to (3.7) in [GS]).
We find {\rm (xxii)}$\leq${\rm (xxi)}.
\begin{flalign*}
&{\rm (xxiii)}\, {\mbox{$F$}(2\,a,1-a;\,4/3;\,8/9)}
={\frac {{3}^{a-1}\sqrt {\pi }\Gamma  \left( 1/6 \right) }{2\Gamma 
 \left( 7/6-a \right) \Gamma  \left( a+1/2 \right) }}
&
\end{flalign*}
(The above is identical to (8/9.2) in [Gos], (5.24) in [DS] and (3.3) in [Ka]).
We derive (xxiii) by using formula (132) in [Gour]
and connection formulae for the hypergeometric series(see (25)--(44) in 2.9 in [Erd]).
We find {\rm (xxiv)}$\leq${\rm (xxiii)}.
\paragraph{(1,2,3-2)}
The special values obtained from (1,2,3-2) coincide with 
those obtained from (1,2,3-1).
\subsubsection{$(k,l,m)=(1,3,3)$}
In this case, we get
\begin{align}
&(a,b,c,x)=(a,3\,a-1/2,3\,a,-3),\, S^{(n)}=
{\frac {{2}^{2\,n}   \left( a+1/2,n 
 \right) ^{2}}{ \left( a+1/3,n \right) \left( a+2/3,n \right) }}
,\tag{1,3,3-1}\\
&(a,b,c,x)=(a,3\,a-3/2,3\,a-1,-3),\, S^{(n)}=
{\frac {{2}^{2\,n} \left( a-1/6,n \right)  \left( a+1/6,n \right) }{ \left( a-1/3,n \right)  \left( a+1/3,n \right) }}
,\tag{1,3,3-2}\\
&(a,b,c,x)=(a,3\,a+1,3\,a,3/2),\, S^{(n)}=
{\frac { \left( -3 \right) ^{3\,n}\left( a+1,n
 \right)  \left( 2\,a,2\,n \right) }{{2}^{3\,n} \left( 3\,a,3\,n \right)}}
,\tag{1,3,3-3}\\
 &(a,b,c,x)=(a,3\,a-3,3\,a-1,3/2),\, S^{(n)}=
 {\frac {  \left( 3\,a-2 \right) }{({-2})^{n}
 \left( 3\,a-2+3\,n \right) }}
\tag{1,3,3-4}.
\end{align}
\paragraph{(1,3,3-1)}
\begin{flalign*}
&{\rm (i)}\, {\mbox{$F$}(a,3\,a-1/2;\,3\,a;\,-3)}
=\begin{cases}
{\dfrac {9\cdot {2}^{2\,n-1}   \left( 3/2,n 
 \right) ^{2}}{\left( 5/3,n \right) 
 \left( 4/3,n \right) }}
\quad {\text {if $a=-1-n$}},\\
{\dfrac {{2}^{2\,n} \left(  1/3,n \right) 
^{2}}{ \left( 1/2,n \right) \left( 1/
6,n \right) }}
\quad {\text {if $a=1/6-n$}},\\
{\dfrac { {2}^{2\,n+1}   \left( 2/3,n 
 \right) ^{2}}{ \left( 5/6,n \right) 
 \left( 1/2,n \right) }}
\quad {\text {if $a=-1/6-n$}},\\
0
\quad {\text {if $a=-1/2-n$}}.\\
\end{cases}&
\end{flalign*}
\begin{flalign*}
&{\rm (ii)}\, {\mbox{$F$}(1/2,2\,a;\,3\,a;\,-3)}
=\begin{cases}
{\dfrac {9}{8}}\,{\dfrac {   \left( 3/2,n \right) 
  ^{2}}{ \left( 5/3,n \right)
 \left( 4/3,n \right) }}
\quad {\text {if $a=-1-n$}},\\
0
\quad {\text {if $a=-1/2-n$}}.\\
\end{cases}&
\end{flalign*}
\begin{flalign*}
&{\rm (iii)}\, {\mbox{$F$}(1/2,a;\,3\,a;\,3/4)}
=\begin{cases}
\dfrac{2}{\sqrt{3}}\,{\dfrac {\Gamma  \left( a+1/3 \right) \Gamma  \left( a+2/3
 \right) }{ \left( \Gamma  \left( a+1/2 \right)  \right) ^{2}
}},\\
{\dfrac {9}{8}}\,{\dfrac {   \left( 3/2,n 
 \right) ^{2}}{ \left( 5/3,n \right)
 \left( 4/3,n \right) }}
\quad {\text {if $a=-1-n$}}
\end{cases}&
\end{flalign*}
(The first case is identical to (3/4.1) in [Gos]).
\begin{flalign*}
&{\rm (iv)}\, {\mbox{$F$}(2\,a,3\,a-1/2;\,3\,a;\,3/4)}
=\begin{cases}
\dfrac{2^{4\,a}}{\sqrt{3}}\,{\dfrac {\Gamma  \left( a+1/3 \right) \Gamma 
 \left( a+2/3 \right) }{ \left( \Gamma  \left( a+1/2 \right)  \right) 
^{2}}},\\
{\dfrac { 9\left(  3/2,n
  \right) ^{2}}{{2}^{4\,n+7}\left( 5/3,n \right) 
 \left( 4/3,n \right) }}
\quad {\text {if $a=-1-n$}}.
\end{cases}&
\end{flalign*}
\begin{flalign*}
&{\rm (v)}\, {\mbox{$F$}(a,3\,a-1/2;\,a+1/2;\,4)}
=\begin{cases}
(-3)^{3\,n}
\quad {\text {if $a=-n$}},\\
(-3)^{3\,n}
\quad {\text {if $a=1/6-n$}},\\
-\left( -3 \right) ^{3\,n+1}
\quad {\text {if $a=-1/6-n$}}
\end{cases}&
\end{flalign*}
(The first case is identical to Theorem 12 in [Ek]).
\begin{flalign*}
&{\rm (vi)}\, {\mbox{$F$}(1/2,1-2\,a;\,a+1/2;\,4)}
=\begin{cases}
1
\quad {\text {if $a=1/2+n$}},\\
-1/3
\quad {\text {if $a=1+n$}}.
\end{cases}&
\end{flalign*}
\begin{flalign*}
&{\rm (vii)}\, {\mbox{$F$}(a,1-2\,a;\,a+1/2;\,4/3)}
=\begin{cases}
3^{2\,n}
\quad {\text {if $a=-n$}},\\
3^{-2\,n}
\quad {\text {if $a=1/2+n$}},\\
3^{-2\,n-2}
\quad {\text {if $a=1+n$}}.
\end{cases}&
\end{flalign*}
\begin{flalign*}
&{\rm (viii)}\, {\mbox{$F$}(1/2,3\,a-1/2;\,a+1/2;\,4/3)}
=\begin{cases}
1
\quad {\text {if $a=1/6-n$}},\\
-1
\quad {\text {if $a=-1/6-n$}}.\\
\end{cases}&
\end{flalign*}
\begin{flalign*}
&{\rm (ix)}\, {\mbox{$F$}(a,1-2\,a;\,3/2-2\,a;\,-1/3)}
={\frac {{2}^{3/2-2\,a}\Gamma  \left( 5/4-a \right) \Gamma  \left( 3/4
-a \right) }{{3}^{1-a}\Gamma  \left( 7/6-a \right) \Gamma  \left( 5/6
-a \right) }}.&
\end{flalign*}
\begin{flalign*}
&{\rm (x)}\, {\mbox{$F$}(1/2,3/2-3\,a;\,3/2-2\,a;\,-1/3)}
=\sqrt{\dfrac{2}{3}}\,\dfrac{\Gamma  \left( 5/4-a \right) \Gamma  \left( 3/4
-a \right) }{\Gamma  \left( 7/6-a \right) \Gamma  \left( 5/6
-a \right) }
&
\end{flalign*}
(The above is a generalization of (28.1) in [Ge]).
\begin{flalign*}
&{\rm (xi)}\, {\mbox{$F$}(1/2,a;\,3/2-2\,a;\,1/4)}
=\dfrac{2^{3/2}\,\Gamma  \left( 5/4-a \right) \Gamma  \left( 3/4
-a \right) }{3\,\Gamma  \left( 7/6-a \right) \Gamma  \left( 5/6
-a \right) }
&
\end{flalign*}
(The above is identical to (1/4.1) in [Gos] and (5.22) in [GS]).
\begin{flalign*}
&{\rm (xii)}\, {\mbox{$F$}(1-2\,a,3/2-3\,a;\,3/2-2\,a;\,1/4)}
=\dfrac{2^{7/2-6\,a}\Gamma  \left( 5/4-a \right) \Gamma  \left( 3/4
-a \right) }{3^{2-3\,a}\Gamma  \left( 7/6-a \right) \Gamma  \left( 5/6
-a \right) }
&
\end{flalign*}
(The above is a generalization of Theorem 30 in [Ek]).
\begin{flalign*}
&{\rm (xiii)}\, {\mbox{$F$}(1/2,3\,a-1/2;\,2\,a+1/2;\,-1/3)}
={\frac {{2}^{1-2\,a} \sqrt{\pi}
\Gamma  \left( 2\,a+1/2 \right) }{ \sqrt{3}\left( \Gamma  \left( a+1/2
 \right)  \right) ^{2} }}.
&
\end{flalign*}
\begin{flalign*}
&{\rm (xiv)}\, {\mbox{$F$}(2\,a,1-a;\,2\,a+1/2;\,-1/3)}
={\frac { \sqrt{\pi}
\Gamma  \left( 2\,a+1/2 \right) }{ 3^a\left( \Gamma  \left( a+1/2
 \right)  \right) ^{2} }}.
&
\end{flalign*}
\begin{flalign*}
&{\rm (xv)}\, {\mbox{$F$}(2\,a,3\,a-1/2;\,2\,a+1/2;\,1/4)}
={\frac {2^{4\,a} \sqrt{\pi}
\Gamma  \left( 2\,a+1/2 \right) }{ 3^{3\,a}\left( \Gamma  \left( a+1/2
 \right)  \right) ^{2} }}.
&
\end{flalign*}
\begin{flalign*}
&{\rm (xvi)}\, {\mbox{$F$}(1/2,1-a;\,2\,a+1/2;\,1/4)}
={\frac {2^{2-2\,a} \sqrt{\pi}
\Gamma  \left( 2\,a+1/2 \right) }{ 3\left( \Gamma  \left( a+1/2
 \right)  \right) ^{2} }}
&
\end{flalign*}
(The above is identical to (1/4.2) in [Gos]).
\begin{flalign*}
&{\rm (xvii)}\, {\mbox{$F$}(1/2,1-2\,a;\,2-3\,a;\,-3)}
=\begin{cases}
{\dfrac { \left( 1/3,n \right) 
 \left( 2/3,n \right) }{\left( 1/6,n \right)  \left( 5/6,n \right) }}
\quad {\text {if $a=1/2+n$}},\\
-{\dfrac {35}{64}}\,{\dfrac { \left( 11/6,n \right)  \left( 13/6,n \right) }{ \left( 5/3,n
 \right)  \left( 7/3,n \right) }}
\quad {\text {if $a=2+n$}}.
\end{cases}&
\end{flalign*}
\begin{flalign*}
&{\rm (xviii)}\,{\mbox{$F$}(1-a,3/2-3\,a;\,2-3\,a;\,-3)}&\\
&=\begin{cases}
{\dfrac {{2}^{2\,n}\left( 5/6,n \right)  \left( 7/
6,n \right) }{ \left( 2/3,n \right)  \left( 4/3,n \right) }}
\quad {\text {if $a=1+n$}},\\
{\dfrac {{2}^{2\,n} \left( 1/3,n \right) 
 \left( 2/3, n\right) }{ \left( 1/6,n \right)  \left( 5/6,n \right) }}
\quad {\text {if $a=1/2+n$}},\\
0
\quad {\text {if $a=5/6+n$}},\\
0
\quad {\text {if $a=7/6+n$}}.
\end{cases}&
\end{flalign*}
\begin{flalign*}
&{\rm (xix)}\,{\mbox{$F$}(1-2\,a,3/2-3\,a;\,2-3\,a;\,3/4)}&\\
&=\begin{cases}
{\dfrac {{2}^{2-4\,a}\Gamma  \left( 4/3-a \right) \Gamma 
 \left( 2/3-a \right) }{\sqrt {3}\,\Gamma  \left( 7/6-a \right) \Gamma  \left( 5/
6-a \right) }},\\
\dfrac{-35\left(11/6,n\right)\left(13/6,n\right)}{2^{4\,n+12}
\left(5/3,n\right)\left(7/3,n\right)}
\quad {\text {if $a=2+n$}}.\\
\end{cases}&
\end{flalign*}
\begin{flalign*}
&{\rm (xx)}\, {\mbox{$F$}(1/2,1-a;\,2-3\,a;\,3/4)}&\\
&=\begin{cases}
{\dfrac {{2}\,\Gamma  \left( 4/3-a \right) \Gamma 
 \left( 2/3-a \right) }{\sqrt {3}\,\Gamma  \left( 7/6-a \right) \Gamma  \left( 5/
6-a \right) }},\\
\dfrac{\left(5/6,n\right)\left(7/6,n\right)}{
\left(2/3,n\right)\left(4/3,n\right)}
\quad {\text {if $a=1+n$}}
\end{cases}&
\end{flalign*}
(The first case is identical to (3/4.2) in [Gos]).
\begin{flalign*}
&{\rm (xxi)}\, {\mbox{$F$}(1/2,2\,a;\,3/2-a;\,4)}
=\begin{cases}
{\dfrac { \left( 3/2,n \right)  \left( 
1/2,n \right) }{ \left( 5/6,n \right) 
 \left( 7/6,n \right) }}
\quad {\text {if $a=-n$}},\\
0
\quad {\text {if $a=-1/2-n$}}
\end{cases}&
\end{flalign*}
(The above are identical to (5.25) in [GS]).
\begin{flalign*}
&{\rm (xxii)}\, {\mbox{$F$}(1-a,3/2-3\,a;\,3/2-a;\,4)}&\\
&=\begin{cases}
{\dfrac { \left( -3 \right) ^{3\,n} \left( 5/6,n
 \right)  \left( 7/6,n \right) }{
 \left( 3/2,n \right)  \left( 1/2,n \right) }}
\quad {\text {if $a=1+n$}},\\
0
\quad {\text {if $a=5/6+n$}},\\
0
\quad {\text {if $a=7/6+n$}}
\end{cases}&
\end{flalign*}
(The first case is identical to Theorem 13 in [Ek]).
\begin{flalign*}
&{\rm (xxiii)}\, {\mbox{$F$}(2\,a,1-a;\,3/2-a;\,4/3)}
=\begin{cases}
{\dfrac { \left( 3/2,n \right)  \left( 
1/2,n \right) }{{3}^{2\,n} \left( 5/6,n \right) \left( 7/6,n \right) }}
\quad {\text {if $a=-n$}},\\
0
\quad {\text {if $a=-1/2-n$}},\\
{\dfrac {{3}^{2\,n} \left( 5/6,n \right)  \left( 7/6,n \right) }{ \left( 3/2,n
 \right)  \left( 1/2,n \right) }}
\quad {\text {if $a=1+n$}}.
\end{cases}&
\end{flalign*}
\begin{flalign*}
&{\rm (xxiv)}\, {\mbox{$F$}(1/2,3/2-3\,a;\,3/2-a;\,4/3)}
=\begin{cases}
0
\quad {\text {if $a=5/6+n$}},\\
0
\quad {\text {if $a=7/6+n$}}.
\end{cases}&
\end{flalign*}
\paragraph{(1,3,3-2)}
The special values obtained from (1,3,3-2) coincide with 
those obtained from (1,3,3-1).
\paragraph{(1,3,3-3)}
\begin{flalign*}
&{\rm (i)}\,
{\mbox{$F$}(a,3\,a+1;\,3\,a;\,3/2)}
=\begin{cases}
0
\quad {\text {if $a=-1-n$}},\\
{\dfrac { \left( -3 \right) ^{3\,n} \left( 1/3,n
 \right) \left( 5/3,2\,n \right) }{{2}^{3\,n} \left( 2,3\,n \right) }}
\quad {\text {if $a=-1/3-n$}},\\
{\dfrac { \left( -3 \right) ^{3\,n} \left( 2/3,n
 \right) \left( 7/3,2\,n \right) }{{2}^{3\,n+1} \left( 3,3\,n \right) }}
\quad {\text {if $a=-2/3-n$}}.
\end{cases}&
\end{flalign*}
The special values obtained from {\rm (ii)} and {\rm (iii)} are trivial.
\begin{flalign*}
&{\rm (iv)}\,
{\mbox{$F$}(2\,a,3\,a+1;\,3\,a;\,3)}
=\begin{cases}
0
\quad {\text {if $a=-1/2-1/2\,n$}},\\
{\dfrac {{3}^{3\,n} \left( 1/3,n \right)  \left( 5/3,2\,n \right) }{ \left( 2,3\,n \right) }}
\quad {\text {if $a=-1/3-n$}},\\
{\dfrac {-{3}^{3\,n} \left( 2/3,n \right)  \left( 7/3,2\,n \right) }{ \left( 3,3\,n \right) }}
\quad {\text {if $a=-2/3-n$}}.\\
\end{cases}&
\end{flalign*}
\begin{flalign*}
&{\rm (v)}\, {\mbox{$F$}(a,3\,a+1;\,a+2;\,-1/2)}
={2}^{3\,a}{3}^{-3\,a} \left( a+1 \right) &
\end{flalign*}
(The above is a special case of 
\begin{gather}
F\left(c-a, c-2, c;(c-2)/(a-1)\right)=(c-1)\left(\frac{a+1-c}{a-1}\right)^{a+1-c}.
\label{1.6eb2}
\end{gather}
This is given in (a/b.1) in [Gos] and (1.6) in [Eb2]).
\begin{flalign*}
&{\rm (vi)}\, {\mbox{$F$}(2,1-2\,a;\,a+2;\,-1/2)}
=2/3\,a+2/3&
\end{flalign*}
(The above is a special case of 
\begin{gather}
F\left(a, 2, c;(c-2)/(a-1)\right)=\frac{(a-1)(c-1)}{a+1-c}.
\label{1.5eb2}
\end{gather}
This is given in (1.5) in [Eb2]).
\begin{flalign*}
&{\rm (vii)}\,
{\mbox{$F$}(a,1-2\,a;\,a+2;\,1/3)}
={2}^{2\,a}{3}^{-2\,a} \left( a+1 \right)
&
\end{flalign*}
(The above is a special case of (\ref{1.6eb2})).
\begin{flalign*}
&{\rm (viii)}\, {\mbox{$F$}(2,3\,a+1;\,a+2;\,1/3)}
=3/2\,a+3/2
&
\end{flalign*}
(The above is a special case of (\ref{1.5eb2})).
\begin{flalign*}
&{\rm (ix)}\,{\mbox{$F$}(a,1-2\,a;\,-2\,a;\,2/3)}
=\begin{cases}
0,\\
{\dfrac {{2}^{2\,n} \left( 3/2,3\,n \right) }{{3}^{2\,n
} \left( 3/2,n \right)  \left( 1,2\,n \right) }}
\quad {\text {if $a=1/2+n$}},\\
{\dfrac {{2}^{2\,n+2} \left( 3,3\,n \right) }{{3}^{2\,n+1
} \left( 2,n \right)  \left( 2,2\,n \right) }}
\quad {\text {if $a=1+n$}}.\\
\end{cases}
&
\end{flalign*}
The special values obtained from {\rm (x)} and {\rm (xi)} are trivial.
\begin{flalign*}
&{\rm (xii)}\,{\mbox{$F$}(-3\,a,1-2\,a;\,-2\,a;\,-2)}
=\begin{cases}
0
\quad {\text {if $a=1/3+n$}},\\
0
\quad {\text {if $a=2/3+n$}},\\
{\dfrac {{2}^{2\,n}\left( 3/2,3\,n \right) }{ \left( 3/2,n \right)  \left( 1,2\,n \right) }}
\quad {\text {if $a=1/2+n$}},\\
{\dfrac {{2}^{2\,n+2}\left( 3,3\,n \right) }{ \left( 2,n \right)  \left( 2,2\,n \right) }}
\quad {\text {if $a=1+n$}}.\\
\end{cases}
&
\end{flalign*}

\begin{flalign*}
&{\rm (xiii)}\,
{\mbox{$F$}(2,3\,a+1;\,2\,a+2;\,2/3)}
=6\,a+3&
\end{flalign*}
(The above is a special case of (\ref{1.5eb2})).
\begin{flalign*}
&{\rm (xiv)}\,
{\mbox{$F$}(2\,a,1-a;\,2\,a+2;\,2/3)}
={3}^{-a} \left( 2\,a+1 \right) 
&
\end{flalign*}
(The above is a special case of (\ref{1.6eb2}})).
\begin{flalign*}
&{\rm (xv)}\,
{\mbox{$F$}(2\,a,3\,a+1;\,2\,a+2;\,-2)}=
\begin{cases}
 {3}^{3\,n}\left( 1-6\,n \right)
\quad {\text {if $a=-1/3-n$}},\\
{3}^{3\,n+1} \left( -1-6\,n \right)
\quad {\text {if $a=-2/3-n$}}
\end{cases}&
\end{flalign*}
(The above are special cases of (\ref{1.6eb2})).
\begin{flalign*}
&{\rm (xvi)}\,
{\mbox{$F$}(2,1-a;\,2\,a+2;\,-2)}
=2/3\,n+1
\quad {\text {if $a=1+n$}}
&
\end{flalign*}
(The above is a special case of (\ref{1.5eb2})).
\begin{flalign*}
&{\rm (xvii)}\,
{\mbox{$F$}(2,1-2\,a;\,2-3\,a;\,3/2)}
=\begin{cases}
6\,n+1
\quad {\text {if $a=1/2+n$}},\\
6\,n+10
\quad {\text {if $a=2+n$}}
\end{cases}&
\end{flalign*}
(The above are special cases of (\ref{1.5eb2})).
\begin{flalign*}
&{\rm (xviii)}\,
{\mbox{$F$}(-3\,a,1-a;\,2-3\,a;\,3/2)}
= \left( -2 \right) ^{-1-n} \left( -3\,n-2 \right) 
\quad {\text {if $a=1+n$}}
&
\end{flalign*}
(The above is a special case of (\ref{1.6eb2})).
\begin{flalign*}
&{\rm (xix)}\,
{\mbox{$F$}(-3\,a,1-2\,a;\,2-3\,a;\,3)}
=\begin{cases}
{2}^{2\,n}\left( 6\,n+1 \right) 
\quad {\text {if $a=1/2+n$}},\\
- {2}^{2\,n+3}\left( 6\,n+10 \right) 
\quad {\text {if $a=2+n$}}
\end{cases}&
\end{flalign*}
(The above are special cases of (\ref{1.6eb2})).
\begin{flalign*}
&{\rm (xx)}\,
{\mbox{$F$}(2,1-a;\,2-3\,a;\,3)}
=3/2\,n+1
\quad {\text {if $a=1+n$}}&
\end{flalign*}
The above is a special case of (\ref{1.5eb2})).
The special values obtained from {\rm (xxi)} are trivial.
\begin{flalign*}
&{\rm (xxii)}\,
{\mbox{$F$}(-3\,a,1-a;\,-a;\,-1/2)}
=\begin{cases}
0,\\
{\dfrac { \left( 3,3\,n \right) }{{2}^{n}\left( 2,n \right)  \left( 2,2\,n \right) }}
\quad {\text {if $a=1+n$}}.
\end{cases}&
\end{flalign*}
\begin{flalign*}
&{\rm (xxiii)}\,
{\mbox{$F$}(2\,a,1-a;\,-a;\,1/3)}
=\begin{cases}
0,\\
{\dfrac { \left( 3,3\,n \right) }{{3}^{n} \left( 2,n \right) \left( 2,2\,n \right) }}
\quad {\text {if $a=1+n$}}.
\end{cases}
&
\end{flalign*}
The special values obtained from {\rm (xxiv)} are trivial.
\paragraph{(1,3,3-4)}
The special values obtained from (1,3,3-4) coincide with 
those obtained from (1,3,3-3).

\subsubsection{$(k,l,m)=(1,4,3)$}
In this case, we have
\begin{flalign}
&(a,b,c,x)=(a,b,b+1-a,-1),
\tag{1,4,3-1}\\
&(a,b,c,x)=(a,4\,a-1/2,3\,a,-1),\,
S^{(n)}={\frac {{2}^{6\,n} \left( a+3/8,n \right)  \left( a+5/8,n \right) }{{3}^{3\,n}
 \left( a+1/3,n \right)  \left( a+2/3,n \right) }}\tag{1,4,3-2},&\\
&(a,b,c,x)=(a,4\,a-5/2,3\,a-1,-1),
S^{(n)}={\frac {{2}^{6\,n} \left( a-1/8,n \right)  \left( a+1/8,n \right) }{{3}^{3\,n}
 \left( a-1/3,n \right)  \left( a+1/3,n \right) }}\tag{1,4,3-3}.&
\end{flalign}
\paragraph{(1,4,3-1)}
The special values obtained from 
(1,4,3-1) are evaluated in paragraphs (1,2,2-1) and (0,2,2-1),
\paragraph{(1,4,3-2)}
\begin{flalign*}
&{\rm (i)}\, {\mbox{$F$}(a,4\,a-1/2;\,3\,a;\,-1)}&\\
&=\begin{cases}
{\dfrac {{2}^{3/4-6\,a}\sqrt {\pi }\Gamma  \left( 3/4 \right) \Gamma 
 \left( a+1/3 \right) \Gamma  \left( a+2/3 \right) }{{3}^{3/8-3\,a}
\Gamma \left(  11/24 \right)
\Gamma  \left( 19/24 \right) \Gamma  \left( a+3/8 \right) 
\Gamma  \left( a+5/8 \right) }},\\
{\dfrac {5\cdot {2}^{6\,n-1} \left( 13/8,n \right)  \left( 11/8,n \right) }{{3}^{3\,n} \left( 5/3
,n \right)  \left( 4/3,n \right) }}
\quad {\text {if $a=-1-n$}}
\end{cases}
&
\end{flalign*}
(The first case is a generalization of Theorem 1 in [Ek]).
\begin{flalign*}
&{\rm (ii)}\, {\mbox{$F$}(2\,a,1/2-a;\,3\,a;\,-1)}&\\
&=\begin{cases}
{\dfrac {{2}^{1/4-4\,a}\sqrt {\pi }\Gamma  \left( 3/4 \right) \Gamma 
 \left( a+1/3 \right) \Gamma  \left( a+2/3 \right) }{{3}^{3/8-3\,a}
\Gamma \left(  11/24 \right)
\Gamma  \left( 19/24 \right) \Gamma  \left( a+3/8 \right) 
\Gamma  \left( a+5/8 \right) }},\\
{\dfrac {5\cdot {2}^{4\,n-3} \left( 13/8,n \right)  \left( 11/8,n \right) }{{3}^{3\,n} \left( 5/3
,n \right)  \left( 4/3,n \right) }}
\quad {\text {if $a=-1-n$}}.
\end{cases}
&
\end{flalign*}
\begin{flalign*}
&{\rm (iii)}\,{\mbox{$F$}(a,1/2-a;\,3\,a;\,1/2)}\\
&=\begin{cases}
{\dfrac {{2}^{3/4-5\,a}\sqrt {\pi }\Gamma  \left( 3/4 \right) \Gamma 
 \left( a+1/3 \right) \Gamma  \left( a+2/3 \right) }{{3}^{3/8-3\,a}
\Gamma \left(  11/24 \right)
\Gamma  \left( 19/24 \right) \Gamma  \left( a+3/8 \right) 
\Gamma  \left( a+5/8 \right) }},\\
{\dfrac {5\cdot {2}^{5\,n-2} \left( 13/8,n \right)  \left( 11/8,n \right) }{{3}^{3\,n} \left( 5/3
,n \right)  \left( 4/3,n \right) }}
\quad {\text {if $a=-1-n$}}.
\end{cases}
&
\end{flalign*}
\begin{flalign*}
&{\rm (iv)}\,{\mbox{$F$}(2\,a,4\,a-1/2;\,3\,a;\,1/2)}\\
&=\begin{cases}
{\dfrac {{2}^{1/4-2\,a}\sqrt {\pi }\Gamma  \left( 3/4 \right) \Gamma 
 \left( a+1/3 \right) \Gamma  \left( a+2/3 \right) }{{3}^{3/8-3\,a}
\Gamma \left(  11/24 \right)
\Gamma  \left( 19/24 \right) \Gamma  \left( a+3/8 \right) 
\Gamma  \left( a+5/8 \right) }},\\
{\dfrac {5\cdot {2}^{2\,n-5} \left( 13/8,n \right)  \left( 11/8,n \right) }{{3}^{3\,n} \left( 5/3
,n \right)  \left( 4/3,n \right) }}
\quad {\text {if $a=-1-n$}}.
\end{cases}
&
\end{flalign*}
\begin{flalign*}
&{\rm (v)}\, {\mbox{$F$}(a,4\,a-1/2;\,2\,a+1/2;\,2)}&\\
&=\begin{cases}
{\dfrac { \left( -4 \right) ^{n} \left( 5/8,n \right) 
 \left( 3/8,n \right) }{ \left( 3/4,
\,n \right)\left( 1/4,\,n \right) }}
\quad {\text {if $a=-n$}},\\
{\dfrac { \left( -4\right) ^{n} \left( 1/2,n \right) 
 \left( 1/4,n \right) }{ \left( 5/8,
\,n \right)\left( 1/8,\,n \right) }}
\quad {\text {if $a=1/8-n$}},\\
{\dfrac { 2\left( -4 \right) ^{n} \left( 3/4,n \right) 
 \left( 1/2,n \right) }{ \left( 7/8,
\,n \right)\left( 3/8,\,n \right) }}
\quad {\text {if $a=-1/8-n$}},\\
0
\quad {\text {if $a=-3/8-n$}},\\
0
\quad {\text {if $a=-5/8-n$}}.
\end{cases}
&
\end{flalign*}
\begin{flalign*}
&{\rm (vi)}\, {\mbox{$F$}(1-2\,a,a+1/2;\,2\,a+1/2;\,2)}&\\
&=\begin{cases}
\dfrac{ \left( 3/4,n \right)\left( 5/4,n \right)}{2^{2\,n}  \left( 7/8,n \right)  \left( 9/8,n \right) }
\quad {\text {if $a=1/2+n$}},\\
\dfrac{- \left( 5/4,n \right)\left( 7/4,n \right)}{  5\cdot 2^{2\,n}\left( 11/8,n \right)  \left( 13/8,n \right) }
\quad {\text {if $a=1+n$}},\\
{\dfrac {  2  ^{2\,n} \left( 9/8,n
 \right)  \left( 7/8,n \right) }{
 \left( 5/4,n \right)\left( 3/4,n \right) }}
\quad {\text {if $a=-1/2-n$}}.
\end{cases}
&
\end{flalign*}
\begin{flalign*}
&{\rm (vii)}\,{\mbox{$F$}(a,1-2\,a;\,2\,a+1/2;\,2)}
=\begin{cases}
{\dfrac {  {2}^{2\,n}\left( 5/8,n \right)  \left( 
3/8,n \right)}{ \left( 3/4,n \right)\left( 1/4,n \right) }}
\quad {\text {if $a=-n$}},\\
\dfrac{ \left( 3/4,n \right) \left( 5/4,n \right)}{ 2^{2\,n} \left( 7/8,n \right)  \left( 9/8,n \right) }
\quad {\text {if $a=1/2+n$}},\\
\dfrac{ \left( 5/4,n \right)\left( 7/4,n \right)}{  5\cdot 2^{2\,n }\left( 11/8,n \right)  \left( 13/8,n \right) }
\quad {\text {if $a=1+n$}}.
\end{cases}
&
\end{flalign*}
\begin{flalign*}
&{\rm (viii)}\,{\mbox{$F$}(4\,a-1/2,a+1/2;\,2\,a+1/2;\,2)}&\\
&=\begin{cases}
{\dfrac { \left( -4 \right) ^{n}\left( 1/2,n \right) 
 \left( 1/4,n \right) }{\left( 5/8,n \right) \left( 1/8,n \right) }}
\quad {\text {if $a=1/8-n$}},\\
{\dfrac { \left( -4 \right) ^{n+1}\left( 3/4,n
 \right)  \left( 1/2,n \right) }{2
 \left( 7/8,n \right)\left( 3/8,n \right) }}
\quad {\text {if $a=-1/8-n$}},\\
0
\quad {\text {if $a=-3/8-n$}},\\
0
\quad {\text {if $a=-5/8-n$}},\\
{\dfrac { \left( -4 \right) ^{n}\left( 9/8,n \right) 
 \left( 7/8,n \right) }{ \left( 5/4,n \right)\left( 3/4,n \right) }}
\quad {\text {if $a=-1/2-n$}}.\\
\end{cases}
&
\end{flalign*}
We find that {\rm (ix)}$\leq${\rm (ii)},
{\rm (x)}$\leq${\rm (i)},
{\rm (xi)}$\leq${\rm (iii)},
{\rm (xii)}$\leq${\rm (iv)}.
\begin{flalign*}
&{\rm (xiii)}\,{\mbox{$F$}(4\,a-1/2,a+1/2;\,3\,a+1/2;\,-1)}&\\
&=\begin{cases}
{\dfrac { {3}^{3\,a-3/8}\sqrt {\pi }\Gamma 
 \left( 3/4 \right) \Gamma  \left( a+1/6 \right)\Gamma  \left( a+5/6 \right)  }{ 2^{6\,a-3/4}\Gamma\left( 7/24 \right)\Gamma\left( 23/24 \right) \Gamma  \left( a+3/8 \right) \Gamma  \left( a+5/8
 \right)  }},\\
\dfrac{{2}^{6\,n} \left( 9/8,n \right)  \left( 7/8,n \right) }{ {3}^{3\,n}
  \left( 2/3,n \right)   \left( 4/3,n \right) }
\quad {\text {if $a=-1/2-n$}}
\end{cases}
&
\end{flalign*}
(The second case is identical to Theorem 2 in [Ek]).
\begin{flalign*}
&{\rm (xiv)}\,{\mbox{$F$}(2\,a,1-a;\,3\,a+1/2;\,-1)}&\\
&=\begin{cases}
{\dfrac { {3}^{3\,a-3/8}\sqrt {\pi }\Gamma 
 \left( 3/4 \right) \Gamma  \left( a+1/6 \right)\Gamma  \left( a+5/6 \right)  }{ 2^{4\,a-1/4}\Gamma\left( 7/24 \right)\Gamma\left( 23/24 \right) \Gamma  \left( a+3/8 \right) \Gamma  \left( a+5/8
 \right)  }},\\
 \dfrac{-7\cdot {2}^{4\,n-6} \left( 17/8,n \right)  \left( 15/8,n \right)}{ {3}^{3\,n-1}
   \left( 5/3,n \right)   \left( 7/3,n \right)  }
\quad {\text {if $a=-3/2-n$}}.
\end{cases}
&
\end{flalign*}
\begin{flalign*}
&{\rm (xv)}\,{\mbox{$F$}(2\,a,4\,a-1/2;\,3\,a+1/2;\,1/2)}&\\
&=\begin{cases}
{\dfrac { {3}^{3\,a-3/8}\sqrt {\pi }\Gamma 
 \left( 3/4 \right) \Gamma  \left( a+1/6 \right)\Gamma  \left( a+5/6 \right)  }{ 2^{2\,a-1/4}\Gamma\left( 7/24 \right)\Gamma\left( 23/24 \right) \Gamma  \left( a+3/8 \right) \Gamma  \left( a+5/8
 \right)  }},\\
 \dfrac{-7\cdot {2}^{2\,n-9} \left( 17/8,n \right)  \left( 15/8,n \right)}{  {3}^{3\,n-1}
   \left( 5/3,n \right)   \left( 7/3,n \right) }
\quad {\text {if $a=-3/2-n$}}.
\end{cases}
&
\end{flalign*}
\begin{flalign*}
&{\rm (xvi)}\,{\mbox{$F$}(1-a,a+1/2;\,3\,a+1/2;\,1/2)}&\\
&=\begin{cases}
{\dfrac { {3}^{3\,a-3/8}\sqrt {\pi }\Gamma 
 \left( 3/4 \right) \Gamma  \left( a+1/6 \right)\Gamma  \left( a+5/6 \right)  }{ 2^{5\,a-5/4}\Gamma\left( 7/24 \right)\Gamma\left( 23/24 \right) \Gamma  \left( a+3/8 \right) \Gamma  \left( a+5/8
 \right)  }},\\
 \dfrac{{2}^{5\,n} \left( 9/8,n \right)  \left( 7/8,n \right)  }{{3}^{3\,n}
  \left( 2/3,n \right)   \left( 4/3,n \right)  }
 \quad {\text {if $a=-1/2-n$}}.
\end{cases}
&
\end{flalign*} 
The special values obtained from {\rm (xvii)-(xxiv)}
are contained in the above.

\paragraph{(1,4,3-3)}
The special values obtained from (1,4,3-3) coincide with 
those obtained from (1,4,3-2).

\subsection{$m=4$}
\subsubsection{$(k,l,m)=(0,4,4)$}
In this case, we have
\begin{align}
&(a,b,c,x)=(a,b,b+1-a,-1)\tag{0,4,4-1},\\
&(a,b,c,x)=(1,b,b,\lambda),\,
S^{(n)}=1 ,\tag{0,4,4-2}\\
&(a,b,c,x)=(0,b,b+1,\lambda),\,
S^{(n)}=1 \tag{0,4,4-3}
\end{align}
where, $\lambda$ is a solution of $x^2+1=0$.
\paragraph{(0,4,4-1)}
The special values obtained from (0,4,4-1) are 
evaluated in the case (0,2,2-1).
\paragraph{(0,4,4-2)}
The special values obtained from (0,4,4-2) except trivial values 
are the special cases of (\ref{F(1,b;2;x)}).

\paragraph{(0,4,4-3)}
The special values obtained from (0,4,4-3) coincide with those obtained from (0,4,4-2)
\subsubsection{$(k,l,m)=(1,3,4)$}
In this case, we have
\begin{align}
&(a,b,c,x)=(a,3\,a-1/2,4\,a,4),\,
S^{(n)}={\frac { \left( -3 \right) ^{3\,n} \left( a+2/3
,n \right)  \left( a+1/6,n \right) }{{2}^{4\,n} \left( a+1/4,n \right) \left( a+3/4,n
 \right) }},
\tag{1,3,4-1}\\
&(a,b,c,x)=(a,3\,a-3/2,4\,a-2,4)\,
S^{(n)}={\frac { \left( -3 \right) ^{3\,n} \left( a-1/3
,n \right)  \left( a+1/6,n \right) }{{2}^{4\,n} \left( a-1/4,n \right)  \left( a+1/4,n
 \right) }},
\tag{1,3,4-2}\\
&(a,b,c,x)=(a,3\,a-1/2,4\,a,-8),\,
S^{(n)}={\frac {{3}^{3\,n} \left( a+1/3,n \right)  \left( a+2/3,n \right) }{{2}^{2\,n}
 \left( a+1/4,n \right)  \left( a+3/4,n \right) }},
\tag{1,3,4-3}\\
&(a,b,c,x)=(a,3\,a-1/4,4\,a,-8),\,
S^{(n)}={\frac {{3}^{3\,n} \left( a+2/3,n \right)  \left( a+7/12,n \right) }{{2}^{2\,n}
 \left( a+1/2,n \right)  \left( a+3/4,n \right) }},
\tag{1,3,4-4}\\
&(a,b,c,x)=(a,3\,a-3/4,4\,a-1,-8),\,
S^{(n)}={\frac {{3}^{3\,n} \left( a+1/3,n \right)  \left( a+5/12,n \right) }{{2}^{2\,n}
 \left( a+1/2,n \right)  \left( a+1/4,n \right) }},
\tag{1,3,4-5}
\\
&(a,b,c,x)=(a,3\,a-5/4,4\,a-1,-8),\,
S^{(n)}={\frac {{3}^{3\,n} \left( a+1/3,n \right)  \left( a-1/12,n \right) }{{2}^{2\,n}
 \left( a+1/2,n \right)  \left( a-1/4,n \right) }},
\tag{1,3,4-6}
\\
&(a,b,c,x)=(a,3\,a-3/2,4\,a-2,-8),\,
S^{(n)}={\frac {{3}^{3\,n} \left( a-1/6,n \right)  \left( a+1/6,n \right) }{{2}^{2\,n}
 \left( a-1/4,n \right)  \left( a+1/4,n \right) }},
\tag{1,3,4-7}
\\
&(a,b,c,x)=(a,3\,a-7/4,4\,a-2,-8),\,
S^{(n)}={\frac {{3}^{3\,n} \left( a+1/3,n \right)  \left( a+5/12,n \right) }{{2}^{2\,n}
 \left( a+1/2,n \right)  \left( a+1/4,n \right) }}
\tag{1,3,4-8}.
\end{align}
\paragraph{(1,3,4-1)}
\begin{flalign*}
&{\rm (i)}\, {\mbox{$F$}(a,3\,a-1/2;\,4\,a;\,4)}
=\begin{cases}
\dfrac{-5\cdot ({-3})^{3\,n} \left( 4/3,n \right) 
 \left( 11/6,n \right)}{ {2}^{4\,n+1}  \left( 7/4,n \right)  \left( 5/4,n \right) }
 \quad {\text {if $a=-1-n$}},\\
\dfrac{\left( -3 \right) ^{3\,n} \left( 1/6,n \right)  \left( 2/3,n \right) }{ {2}^{4\,n}
 \left( 7/12,n \right)  \left( 1/12,n \right)  } 
 \quad {\text {if $a=1/6-n$}},\\
 0
 \quad {\text {if $a=-1/6-n$}},\\
\dfrac{-\left( -3 \right) ^{3\,n+1} \left( 11/6,n
 \right) \left( 7/3,n \right)}{ {2}^{4\,n-1}
  \left( 9/4,n \right)  \left( 7/4,n \right)} 
 \quad {\text {if $a=-3/2-n$}}
\end{cases}&
\end{flalign*}
(The first case is identical to Theorem 3 in [Ek]).
\begin{flalign*}
&{\rm (ii)}\,{\mbox{$F$}(3\,a,a+1/2;\,4\,a;\,4)}
=\begin{cases}
\dfrac{-5 \left( -3 \right) ^{3\,n} \left( 4/3,n \right)  \left( {11/6},n \right)}{ {2}^{4\,n+1}
  \left( 7/4,n \right)   \left( 5/4,n \right)  }
\quad {\text {if $a=-1-n$}},\\
\dfrac{-\left( -3 \right) ^{3\, n+1} \left( 2/3,n
 \right) \left( 7/6,n \right) }{{2}^{4\,n+1}
   \left( 13/12,n
 \right)   \left( 7/12
,n \right) }
\quad {\text {if $a=-1/3-n$}},\\
0
\quad {\text {if $a=-2/3-n$}},\\
{\dfrac { \left( -3 \right) ^{3\,n} \left( 5/6,n
 \right) \left( 4/3,n \right) }{{2}^{4\,n} \left( 5/4,n \right)  \left( 3/4,n
 \right) }}
\quad {\text {if $a=-1/2-n$}}\\
\end{cases}&
\end{flalign*}
(The fourth case is identical to Theorem 9 in [Ek]).
\begin{flalign*}
&{\rm (iii)}\,{\mbox{$F$}(a,a+1/2;\,4\,a;\,4/3)}
=\begin{cases}
\dfrac{5\cdot {3}^{2\,n-1} \left( 4/3,n \right) 
 \left( 11/6,n \right) }{ {2}^{4\,n+1}  \left( 7/4,n \right)    \left( 5/4,n \right)  }
\quad {\text {if $a=-1-n$}},\\
{\dfrac {{3}^{2\,n} \left( 5/6,n \right)  \left( 4/3,n \right) }{{2}^{4\,n} \left( 5
/4,n \right)  \left( 3/4,n \right) }}
\quad {\text {if $a=-1/2-n$}}.
\end{cases}&
\end{flalign*}
\begin{flalign*}
&{\rm (iv)}\,{\mbox{$F$}(3\,a,3\,a-1/2;\,4\,a;\,4/3)}&\\
&=\begin{cases}
{\dfrac {5}{27}}\,\dfrac{\left( 4/3,n \right)  \left(  11/6,n \right)}{ {2}^{4\,n+1} \left( 7/4,n \right)   \left( 5/4,n \right) }
\quad {\text {if $a=-1-n$}},\\
\dfrac{- \left( 2/3,n \right)  \left( 7/6,n
 \right)}{ {2}^{4\,n+1}
 \left( 13/12,n \right)  \left( 7/12,n \right) }
\quad {\text {if $a=-1/3-n$}},\\
0
\quad {\text {if $a=-2/3-n$}},\\
\dfrac{\left( 1/6,n \right)  \left( 2/3,n
 \right) }{ {2}^{4\,n}  
 \left( 7/12,n \right) \left( 1/12,n \right)  }
\quad {\text {if $a=1/6-n$}},\\
0
\quad {\text {if $a=-1/6-n$}},\\
-{\dfrac {1}{81}}\, \dfrac{\left( 11/6,n \right)  \left( 7/3,n \right)}{ {2}^{4\,n-1}  \left( 9/4,n \right)   \left( 7/4,n \right) }
\quad {\text {if $a=-3/2-n$}},\\
\end{cases}&
\end{flalign*}
\begin{flalign*}
&{\rm (v)}\,{\mbox{$F$}(a,3\,a-1/2;\,1/2;\,-3)}
=\begin{cases}
{\dfrac {\left( -16 \right) ^{n} \left( 5/6,n
 \right) }{ \left( 2/3,n \right) }}
\quad {\text {if $a=-n$}},\\
{\dfrac {\left( -16 \right) ^{n}  \left( 2/3,n
 \right) }{ \left( 1/2,n \right) }}
\quad {\text {if $a=1/6-n$}},\\
0
\quad {\text {if $a=-1/6-n$}},\\
{\dfrac { \left( -16 \right) ^{n+1} \left( 4
/3,n \right) }{ 2\left( 7/6,n \right) }}
\quad {\text {if $a=-1/2-n$}}.\\
\end{cases}&
\end{flalign*}
We find {\rm (vi)}$\leq ${\rm (v)}.
\begin{flalign*}
&{\rm (vii)}\,{\mbox{$F$}(a,1-3\,a;\,1/2;\,3/4)}
=
\dfrac{2^{2/3-2\,a}\sqrt{\pi}\Gamma (1/3)}{\Gamma (a+1/6)\Gamma (2/3-a)}
&
\end{flalign*}
(The above is identical to (3/4.4) in [Gos]).
We derive (vii) by using the algebraic transformation formula
(44) in 2.11  in [Erd] 
and connection formulae for the hypergeometric series(see (25)--(44) in 2.9 in [Erd]).
We find {\rm (viii)}$\leq ${\rm (vii)}.
\begin{flalign*}
&{\rm (ix)}\,{\mbox{$F$}(a,1-3\,a;\,3/2-2\,a;\,1/4)}
=
{\dfrac {{2}^{2-4\,a}\sqrt{\pi}\Gamma  \left(3/2 -2\,a \right) }{{3}^{3/2-3\,a}
\Gamma  \left( 2/3\right) \Gamma  \left( 4/3-2\,a \right) }}.
&
\end{flalign*}
\begin{flalign*}
&{\rm (x)}\,{\mbox{$F$}(a+1/2,3/2-3\,a;\,3/2-2\,a;\,1/4)}
=
{\dfrac {{2}^{3-4\,a}\sqrt{\pi}\Gamma  \left(3/2 -2\,a \right) }{{3}^{2-3\,a}
\Gamma  \left( 2/3\right) \Gamma  \left( 4/3-2\,a \right) }}.
&
\end{flalign*}
\begin{flalign*}
&{\rm (xi)}\,{\mbox{$F$}(a,a+1/2;\,3/2-2\,a;\,-1/3)}
=
{\dfrac {{2}^{2-6\,a}\sqrt{\pi}\Gamma  \left(3/2 -2\,a \right) }{{3}^{3/2-4\,a}
\Gamma  \left( 2/3\right) \Gamma  \left( 4/3-2\,a \right) }}.
&
\end{flalign*}
\begin{flalign*}
&{\rm (xii)}\,{\mbox{$F$}(1-3\,a,3/2-3\,a;\,3/2-2\,a;\,-1/3)}
=
{\dfrac {{2}^{2\,a}\sqrt{\pi}\Gamma  \left(3/2 -2\,a \right) }{\sqrt{3}
\Gamma  \left( 2/3\right) \Gamma  \left( 4/3-2\,a \right) }}.
&
\end{flalign*}
The special values obtained from {\rm (xiii)-(xx)} are contained in
the above those.
\begin{flalign*}
&{\rm (xxi)}\,{\mbox{$F$}(3\,a,a+1/2;\,3/2;\,-3)}
=\begin{cases}
{\dfrac { \left( -16 \right) ^{n} \left( 1/3,n
 \right) }{ \left( 7/6,n \right) }}
\quad {\text {if $a=-n$}},\\
{\dfrac { 4\left( -16 \right) ^{n}
 \left( 2/3,n \right) }{3 \left( 3/2,n \right) }}
\quad {\text {if $a=-1/3-n$}},\\
0
\quad {\text {if $a=-2/3-n$}},\\
{\dfrac { \left( -16 \right) ^{n} \left( 5/6,n
 \right) }{ \left( 5/3,n \right) }}
\quad {\text {if $a=-1/2-n$}}.
\end{cases}&
\end{flalign*}
We find {\rm (xxii)}$\leq ${\rm (xxi)}.
\begin{flalign*}
&{\rm (xxiii)}\,{\mbox{$F$}(3\,a,1-a;\,3/2;\,3/4)}
={\frac {{2}^{2\,a-1/3}\sqrt {\pi }\Gamma  \left( 4/3 \right) }{\Gamma 
 \left( 7/6-a \right) \Gamma  \left( a+2/3 \right) }}.&
\end{flalign*}
(The above is identical to (3/4.3) in [Gos]).
We derive (xxiii) by using the algebraic transformation formula
(46) in 2.11  in [Erd] 
and connection formulae for the hypergeometric series(see (25)--(44) in 2.9 in [Erd]).
We find {\rm (xxiv)}$\leq ${\rm (xxiii)}.

\paragraph{(1,3,4-2)}
The special values obtained from (1,3,4-2) coincide with those obtained
from (1,3,4-1).

\paragraph{(1,3,4-3)}
\begin{flalign*}
&{\rm (i)}\,{\mbox{$F$}(a,3\,a-1/2;\,4\,a;\,-8)}
=\begin{cases}
{\dfrac {{3}^{3\,n} \left( 4/3,n \right)  \left( 5/3,n \right) }{{2}^{2\,n-3}
 \left( 7/4,n \right)  \left( 5/4,n \right) }}
\quad {\text {if $a=-1-n$}},\\
\dfrac{{3}^{3\,n}\left( 1/6,n \right)
 \left( 1/2,n \right) }{ {2}^{2\,n}   \left( 7/12,n \right)   \left( 1/12,n \right)  }
\quad {\text {if $a=1/6-n$}},\\
\dfrac{{3}^{3\,n+1}\left( 1/2,n \right) 
 \left( 5/6,n \right)}{ {2}^{2\,n}   \left(  11/12,n \right)    \left( 5/12,n \right) }
\quad {\text {if $a=-1/6-n$}},\\
\dfrac{-7\cdot {3}^{3\,n+1} \left( 11/6,n \right)  \left( 13/6,n \right)}{{2}^{2\,n}
 \left( 9/4,n \right)  \left( 7/4,n \right)  }
\quad {\text {if $a=-3/2-n$}}
\end{cases}&
\end{flalign*}
(The first case is identical to Theorem 7 in [Ek]).
\begin{flalign*}
&{\rm (ii)}\,{\mbox{$F$}(3\,a,a+1/2;\,4\,a;\,-8)}
=\begin{cases}
{\dfrac {{3}^{3\,n} \left( 4/3,n \right)  \left( 5/3,n \right) }{{2}^{2\,n-3}
 \left( 7/4,n \right)  \left( 5/4,n \right) }}
\quad {\text {if $a=-1-n$}},\\
0
\quad {\text {if $a=-1/3-n$}},\\
0
\quad {\text {if $a=-2/3-n$}},\\
{\dfrac {{3}^{3\,n} \left( 5/6,n \right) \left( 7/6,n \right) }{{2}^{2\,n} \left( 5
/4,n \right)  \left( 3/4,n \right) }}
\quad {\text {if $a=-1/2-n$}}
\end{cases}&
\end{flalign*}
(The fourth case is identical to Theorem 5 in [Ek]).
\begin{flalign*}
&{\rm (iii)}\,{\mbox{$F$}(a,a+1/2;\,4\,a;\,8/9)}
=\begin{cases}
{\dfrac {{2}^{2\,a-1/2}\Gamma  \left( a+3/4 \right) \Gamma  \left( a+1/
4 \right) }{{3}^{a-1/2}\Gamma  \left( a+2/3 \right) \Gamma  \left( a+1
/3 \right) }},\\
{\dfrac {{3}^{n-2} \left( 4/3,n \right)  \left( 5/3,n \right) }{{2}^{2\,n-3} \left( 7/4,n
 \right)  \left( 5/4,n \right) }}
\quad {\text {if $a=-1-n$}},\\
{\dfrac {{3}^{n} \left( 5/6,n \right) 
 \left( 7/6,n \right) }{{2}^{2\,n} \left( 5/4,n
 \right)  \left( 3/4,n \right) }}
\quad {\text {if $a=-1/2-n$}}.
\end{cases}&
\end{flalign*}
\begin{flalign*}
&{\rm (iv)}\,{\mbox{$F$}(3\,a,3\,a-1/2;\,4\,a;\, 8/9)}
=\begin{cases}
{\dfrac {{2}^{2\,a-1/2}\Gamma  \left( a+3/4 \right) \Gamma  \left( a+1/
4 \right) }{{3}^{1/2-3\,a}\Gamma  \left( a+2/3 \right) \Gamma  \left( 
a+1/3 \right) }},\\
{\dfrac {\left( 4/3,n \right) \left( 
5/3,n \right) }{{2}^{2\,n-3}{3}^{3\,n+6} \left( 7/4,n
 \right)  \left( 5/4,n \right) }}
\quad {\text {if $a=-1-n$}},\\
{\dfrac {-7\left( 11/6,n \right) \left( 
13/6,n \right) }{{2}^{2\,n}{3}^{3\,n+9} \left( 9/4,n
 \right)  \left( 7/4,n \right) }}
\quad {\text {if $a=-3/2-n$}}
\end{cases}&
\end{flalign*}
(The first case is identical to (3.1) in [Ka]).
\begin{flalign*}
&{\rm (v)}\,{\mbox{$F$}(a,3\,a-1/2;\,1/2;\,9)}
=\begin{cases}
2^{6\,n}
\quad {\text {if $a=-n$}},\\
2^{6\,n}
\quad {\text {if $a=1/6-n$}},\\
2^{6\,n+2}
\quad {\text {if $a=-1/6-n$}},\\
-2^{6\, n+3}
\quad {\text {if $a=-1/2-n$}}.
\end{cases}&
\end{flalign*}
We find {\rm (vi)}$\leq ${\rm (v)}.
\begin{flalign*}
&{\rm (vii)}\,{\mbox{$F$}(a,1-3\,a;\,1/2;\,9/8)}
=\begin{cases}
\left( -2 \right) ^{3\,n}
\quad {\text {if $a=-n$}},\\
\left( -2 \right) ^{-3\,n}
\quad {\text {if $a=1/3+n$}},\\
\left( -2 \right) ^{-3\,n-1}
\quad {\text {if $a=2/3+n$}},\\
\left( -2 \right) ^{-3\,n-3}
\quad {\text {if $a=1+n$}}.
\end{cases}&
\end{flalign*}
We find {\rm (viii)}$\leq ${\rm (vii)}.
\begin{flalign*}
&{\rm (ix)}\,{\mbox{$F$}(a,1-3\,a;\,3/2-2\,a;\,-1/8)}
=\left(\frac{2}{3}\right)^{1-3\,a} \frac { \sqrt {\pi }\Gamma  \left(3/2 -2\,a \right) }{\Gamma  \left( 7/6-a \right) \Gamma  \left( 5/6-a
 \right) }.&
\end{flalign*}
\begin{flalign*}
&{\rm (x)}\,{\mbox{$F$}(a+1/2,3/2-3\,a;\,3/2-2\,a;\,-1/8)}
={\frac {{2}^{5/2-3\,a}\sqrt {\pi }\Gamma  \left( 3/2-2\,a \right) }{{
3}^{2-3\,a}\Gamma  \left( 7/6-a \right) \Gamma  \left( 5/6-a \right) }
}.&
\end{flalign*}
\begin{flalign*}
&{\rm (xi)}\,{\mbox{$F$}(a,a+1/2;\,3/2-2\,a;\,1/9)}
={\frac {{2}^{1-6\,a}\sqrt {\pi }\Gamma  \left( 3/2-2\,a \right) }{{3}
^{1-5\,a}\Gamma  \left( 7/6-a \right) \Gamma  \left( 5/6-a \right) }}&
\end{flalign*}
(The above is identical to (1/9.1) in [Gos]).
\begin{flalign*}
&{\rm (xii)}\,{\mbox{$F$}(1-3\,a,3/2-3\,a;\,3/2-2\,a;\,1/9)}
=\left(\frac{4}{3}\right)^{3\,a-1} {\frac { \sqrt {\pi }\Gamma  \left( 3/2-2\,a \right) }{\Gamma  \left( 7/6-a \right) \Gamma  \left( 5/6-a
 \right) }}
&
\end{flalign*}
(The above is identical to (1.1) in [Ka]).
The special values obtained from {\rm (xiii)-(xx)} are contained in
the above those.
\begin{flalign*}
&{\rm (xxi)}\,{\mbox{$F$}(3\,a,a+1/2;\,3/2;\,9)}
=\begin{cases}
{\dfrac {{2}^{6\,n} \left( 1/3,n \right)  \left( 2/3,n \right) }{\left( 5/6,n
 \right)  \left( 7/6,n \right) }}
\quad {\text {if $a=-n$}},\\
0
\quad {\text {if $a=-1/3-n$}},\\
0
\quad {\text {if $a=-2/3-n$}},\\
{\dfrac {{2}^{6\,n} \left( 5/6,n \right)  \left( 7/6,n \right) }{ \left( 4/3,n
 \right)  \left( 5/3,n \right) }}
\quad {\text {if $a=-1/2-n$}},\\
\end{cases}&
\end{flalign*}
We find {\rm (xxii)}$\leq ${\rm (xxi)}.
\begin{flalign*}
&{\rm (xxiii)}\,{\mbox{$F$}(3\,a,1-a;\,3/2;\, 9/8)}
=\begin{cases}
{\dfrac { \left( 1/3,n \right)  \left( 
2/3,n \right) }{ \left( -2 \right) ^{3\,n} \left( 5/6,
n \right)  \left( 7/6,n \right) }}
\quad {\text {if $a=-n$}},\\
0
\quad {\text {if $a=-1/3-n$}},\\
0
\quad {\text {if $a=-2/3-n$}},\\
{\dfrac { \left( -2 \right) ^{3\,n} \left( 5/6,n
 \right)  \left( 7/6,n \right) }{
 \left( 5/3,n \right)  \left( 4/3,n \right) }}
\quad {\text {if $a=1+n$}}.
\end{cases}&
\end{flalign*}
We find {\rm (xxiv)}$\leq ${\rm (xxiii)}.

\paragraph{(1,3,4-4)}
\begin{flalign*}
&{\rm (i)}\,{\mbox{$F$}(a,3\,a-1/4;\,4\,a;\,-8)}
=\begin{cases}
\dfrac{5\cdot {3}^{3\,n+1} \left( 4/3,n \right) 
 \left( 17/12,n \right)}{{2}^{2\,n+1} \left( 3/2,n \right)   \left( 5/4,n \right)  }
\quad {\text {if $a=-1-n$}},\\
\dfrac{{3}^{3\,n} \left( 1/4,n \right) 
 \left( 1/3,n \right)}{ {2}^{2\,n}   \left( 5/12,n \right)    \left( 1/6,n \right)  }
\quad {\text {if $a=1/12-n$}},\\
\dfrac{7\cdot {3}^{3\,n+1} \left( 19/12,n \right)  \left( 5/3,n \right) }{ {2}^{2\,n}   \left( 7/4,n \right) \left( 3/2,n \right)}
\quad {\text {if $a=-5/4-n$}},\\
0
\quad {\text {if $a=-7/12-n$}}\\
\end{cases}&
\end{flalign*}
(The first case is identical to Theorem 6 in [Ek]).
\begin{flalign*}
&{\rm (ii)}\,{\mbox{$F$}(3\,a,a+1/4;\,4\,a;\,-8)}
=\begin{cases}
\dfrac{5\cdot {3}^{3\,n+1} \left( 4/3,n \right)
 \left( 17/12,n \right) }{{2}^{2\,n+1}  \left( 3/2,n \right)   \left( 5/4,n \right) }
\quad {\text {if $a=-1-n$}},\\
{\dfrac {{3}^{3\,n+1} \left( 2/3,n \right)  \left( 3/4,n \right) }{{2}^{2\,n+1}
 \left( 5/6,n \right)  \left( 7/12,n \right) }}
\quad {\text {if $a=-1/3-n$}},\\
0
\quad {\text {if $a=-2/3-n$}},\\
{\dfrac {{3}^{3\,n}\left( 7/12,n \right)  \left( 2/3,n \right) }{{2}^{2\,n} \left( 3
/4,n \right)  \left( 1/2,n \right) }}
\quad {\text {if $a=-1/4-n$}}\\
\end{cases}&
\end{flalign*}
(The fourth case is identical to Theorem 4 in [Ek]).
\begin{flalign*}
&{\rm (iii)}\,{\mbox{$F$}(a,a+1/4;\,4\,a;\,8/9)}
=\begin{cases}
\dfrac{{2}^{2\,a-1/6}\Gamma  \left( 2/3 \right) \Gamma  \left( 3/4 \right) 
\Gamma  \left( a+3/4 \right) \Gamma  \left( a+1/2 \right) }{ {3}^
{a-1/4} \Gamma  \left( 5/6 \right)   \Gamma  \left( 7/12 \right) 
\Gamma  \left( a+2/3 \right)  \Gamma  \left( a+7/12 \right)},\\
\dfrac{5\cdot {3}^{n-1} \left( 4/3,n \right) 
 \left( 17/12,n \right)}{ {2}^{2\,n+1}   \left( 3/2,n \right)    \left( 5/4,n \right)}
\quad {\text {if $a=-1-n$}},\\
\dfrac{{3}^{n} \left( 7/12,n \right)  \left( 2/3,n \right)}{  {2}^{2\,n}   \left( 3/4,n \right)   \left( 1/2,n \right) }
 \quad {\text {if $a=-1/4-n$}}.\\
\end{cases}&
\end{flalign*}
\begin{flalign*}
&{\rm (iv)}\,{\mbox{$F$}(3\,a,3\,a-1/4;\,4\,a;\,8/9)}&\\
&=\begin{cases}
\dfrac{{2}^{2\,a-1/6}\Gamma  \left( 2/3 \right) \Gamma  \left( 3/4 \right) 
\Gamma  \left( a+3/4 \right) \Gamma  \left( a+1/2 \right) }{ {3}^
{1/4-3\,a} \Gamma  \left( 5/6 \right)   \Gamma  \left( 7/12 \right) 
\Gamma  \left( a+2/3 \right)  \Gamma  \left( a+7/12 \right)},\\
\dfrac{5 \left( 4/3,n \right)  \left( 
17/12,n \right)}{ {2}^{2\,n+1}  {3
}^{3\,n+5} \left( 3/2,n \right) 
 \left( 5/4,n \right) }
\quad {\text {if $a=-1-n$}},\\
\dfrac{7 \left( 19/12,n \right) 
 \left( 5/3,n \right)}{ {2}^{2\,n} {3}^{3\,
n+7}  \left( 7/4,n \right) 
  \left( 3/2,n \right)}
\quad {\text {if $a=-5/4-n$}}\\
\end{cases}&
\end{flalign*}
(The first case is identical to (3.4) in [Ka]).
\begin{flalign*}
&{\rm (v)}\,{\mbox{$F$}(a,3\,a-1/4;\,3/4;\,9)}
=\begin{cases}
\dfrac{{2}^{6\,n} \left( 5/12,n \right)}{  \left( 2/3,n \right)}
 \quad {\text {if $a=-n$}},\\
\dfrac{{2}^{6\,n} \left( 1/3,n \right) }{ \left( 7/12,n \right)  }
\quad {\text {if $a=1/12-n$}},\\
\dfrac{{2}^{6\,n+2} \left( 2/3,n \right) }{ \left( 11/12,n \right)  }
\quad {\text {if $a=-1/4-n$}},\\
0
\quad {\text {if $a=-7/12-n$}}.\\
\end{cases}&
\end{flalign*}
\begin{flalign*}
&{\rm (vi)}\,{\mbox{$F$}(1-3\,a,3/4-a;\,3/4;\,9)}
=\begin{cases}
\dfrac{{2}^{6\,n} \left( 2/3,n \right) }{ \left( 11/12,n \right)  }
\quad {\text {if $a=1/3+n$}},\\
0
\quad {\text {if $a=2/3+n$}},\\
\dfrac{-{2}^{6\,n+5} \left( 4/3,n \right)}{ 7  \left( 19/12,n \right) }
\quad {\text {if $a=1+n$}},\\
\dfrac{{2}^{6\,n} \left( 13/12,n \right) }{ \left( 4/3,n \right)  }
\quad {\text {if $a=3/4+n$}}.
\end{cases}&
\end{flalign*}
\begin{flalign*}
&{\rm (vii)}\,{\mbox{$F$}(a,1-3\,a;\,3/4;\,9/8)}
=\begin{cases}
\dfrac{({-2})^{3\,n} \left( 5/12,n \right) }{ \left( 2/3,n \right)  }
\quad {\text {if $a=-n$}},\\
\dfrac{ \left( 2/3,n \right) }{({-2})^{3\,n} \left( 11/12,n \right)  }
\quad {\text {if $a=1/3+n$}},\\
0
\quad {\text {if $a=2/3+n$}},\\
\dfrac{ \left( 4/3,n \right) }{7({-2})^{3\,n+1} \left( 19/12,n \right)  }
\quad {\text {if $a=1+n$}}.\\
\end{cases}&
\end{flalign*}
\begin{flalign*}
&{\rm (viii)}\,{\mbox{$F$}(3\,a-1/4,3/4-a;\,3/4;\,9/8)}&\\
&=\begin{cases}
\dfrac{ \left( 1/3,n \right) }{({-2})^{3\,n} \left( 7/12,n \right)  }
\quad {\text {if $a=1/12-n$}},\\
\dfrac{ \left( 2/3,n \right) }{({-2})^{3\,n+1} \left( 11/12,n \right)  }
\quad {\text {if $a=-1/4-n$}},\\
0
\quad {\text {if $a=-7/12-n$}},\\
\dfrac{ ({-2})^{3\,n}\left( 13/12,n \right) }{ \left( 4/3,n \right)  }
\quad {\text {if $a=3/4+n$}}.
\end{cases}&
\end{flalign*}
\begin{flalign*}
&{\rm (ix)}\,{\mbox{$F$}(a,1-3\,a;\,5/4-2\,a;\,-1/8)}&\\
&=
\dfrac{{3}^{3\,a+5/4}\Gamma  \left( 4/3 \right) \Gamma  \left( 11/12 \right) 
\Gamma  \left( 5/8-a \right) \Gamma  \left( 9/8-a \right) }{ {2}^
{5\,a+2} \Gamma  \left( 7/8 \right)   \Gamma  \left( 11/8 \right) 
\Gamma  \left( 13/12-a \right)  \Gamma  \left( 2/3-a \right)}.
&
\end{flalign*}
\begin{flalign*}
&{\rm (x)}\,{\mbox{$F$}(5/4-3\,a,a+1/4;\,5/4-2\,a;\,-1/8)}&\\
&=
\dfrac{{3}^{3\,a+3/4}\Gamma  \left( 4/3 \right) \Gamma  \left( 11/12 \right) 
\Gamma  \left( 5/8-a \right) \Gamma  \left( 9/8-a \right) }{ {2}^
{5\,a+5/4} \Gamma  \left( 7/8 \right)   \Gamma  \left( 11/8 \right) 
\Gamma  \left( 13/12-a \right)  \Gamma  \left( 2/3-a \right)}.
&
\end{flalign*}
\begin{flalign*}
&{\rm (xi)}\,{\mbox{$F$}(a,a+1/4;\,5/4-2\,a;\,1/9)}&\\
&=
\dfrac{{3}^{5\,a+5/4}\Gamma  \left( 4/3 \right) \Gamma  \left( 11/12 \right) 
\Gamma  \left( 5/8-a \right) \Gamma  \left( 9/8-a \right) }{ {2}^
{8\,a+2} \Gamma  \left( 7/8 \right)   \Gamma  \left( 11/8 \right) 
\Gamma  \left( 13/12-a \right)  \Gamma  \left( 2/3-a \right)}
&
\end{flalign*}
(The above is identical to (1/9.2) in [Gos], (6.5) in [GS] and (1.4) in [Ka]).
\begin{flalign*}
&{\rm (xii)}\,{\mbox{$F$}(1-3\,a,5/4-3\,a;\,5/4-2\,a;\,1/9)}&\\
&=
\dfrac{{3}^{13/4-3\,a}\Gamma  \left( 4/3 \right) \Gamma  \left( 11/12 \right) 
\Gamma  \left( 5/8-a \right) \Gamma  \left( 9/8-a \right) }{ {2}^
{5-4\,a} \Gamma  \left( 7/8 \right)   \Gamma  \left( 11/8 \right) 
\Gamma  \left( 13/12-a \right)  \Gamma  \left( 2/3-a \right)}.
&
\end{flalign*}
\begin{flalign*}
&{\rm (xiii)}\,{\mbox{$F$}(3/4-a,3\,a-1/4;\,2\,a+3/4;\,-1/8)}&\\
&=
\dfrac{{3}^{1/2-3\,a}\Gamma  \left( 2/3 \right) \Gamma  \left( 7/12 \right) 
\Gamma  \left( a+7/8 \right) \Gamma  \left( a+3/8 \right) }{ {2}^
{3/4-5\,a} \Gamma  \left( 7/8 \right)   \Gamma  \left( 3/8 \right) 
\Gamma  \left( a+2/3 \right)  \Gamma  \left( a+7/12 \right)}.
&
\end{flalign*}
\begin{flalign*}
&{\rm (xiv)}\, {\mbox{$F$}(3\,a,1-a;\,2\,a+3/4;\,-1/8)}&\\
&=
\dfrac{{3}^{-3\,a}\Gamma  \left( 2/3 \right) \Gamma  \left( 7/12 \right) 
\Gamma  \left( a+7/8 \right) \Gamma  \left( a+3/8 \right) }{ {2}^
{-5\,a} \Gamma  \left( 7/8 \right)   \Gamma  \left( 3/8 \right) 
\Gamma  \left( a+2/3 \right)  \Gamma  \left( a+7/12 \right)}.
&
\end{flalign*}
\begin{flalign*}
&{\rm (xv)}\, {\mbox{$F$}(3\,a,3\,a-1/4;\,2\,a+3/4;\,1/9)}&\\
&=
\dfrac{{3}^{3\,a}\Gamma  \left( 2/3 \right) \Gamma  \left( 7/12 \right) 
\Gamma  \left( a+7/8 \right) \Gamma  \left( a+3/8 \right) }{ {2}^
{4\,a} \Gamma  \left( 7/8 \right)   \Gamma  \left( 3/8 \right) 
\Gamma  \left( a+2/3 \right)  \Gamma  \left( a+7/12 \right)}.
&
\end{flalign*}
\begin{flalign*}
&{\rm (xvi)}\,{\mbox{$F$}(1-a,3/4-a;\,2\,a+3/4;\,1/9)}&\\
&=
\dfrac{{3}^{2-5\,a}\Gamma  \left( 2/3 \right) \Gamma  \left( 7/12 \right) 
\Gamma  \left( a+7/8 \right) \Gamma  \left( a+3/8 \right) }{ {2}^
{3-8\,a} \Gamma  \left( 7/8 \right)   \Gamma  \left( 3/8 \right) 
\Gamma  \left( a+2/3 \right)  \Gamma  \left( a+7/12 \right)}
&
\end{flalign*}
(The above is identical to (1/9.3) in [Gos], (6.6) in [GS] and (1.5) in [Ka]).
\begin{flalign*}
&{\rm (xvii)}\,{\mbox{$F$}(1-3\,a,3/4-a;\,2-4\,a;\,-8)}&\\
&=\begin{cases}
{\dfrac {{3}^{3\,n}\left( 1/4,n \right)  \left( 2/3,n \right) }{{2}^{2\,n}\left( 1
/12,n \right) \left( 5/6,n \right) }}
\quad {\text {if $a=1/3+n$}},\\
0
\quad {\text {if $a=2/3+n$}},\\
\dfrac{-11\cdot {3}^{3\,n} \left( 23/12,n \right)  \left( 7/3,n \right)}{  {2}^{2\,n-1}  \left( 7/4,n \right)  \left( 5/2,n \right)  }
 \quad {\text {if $a=2+n$}},\\
 \dfrac{ {3}^{3\,n} \left( 2/3,n \right) 
 \left( 13/12,n \right) }{ {2}^{2\,n}  \left( 1/2,n \right)    \left( 5/4,n \right) }
\quad {\text {if $a=3/4+n$}}\\
\end{cases}
&
\end{flalign*}
(The fourth case is identical to Theorem 8 in [Ek]).
\begin{flalign*}
&{\rm (xviii)}\,{\mbox{$F$}(1-a,5/4-3\,a;\,2-4\,a;\,-8)}&\\
&=\begin{cases}
\dfrac{{3}^{3\,n} \left( 11/12,n \right)  \left( 4/3,n \right)}{  {2}^{2\,n} 
 \left( 3/2,n \right) \left( 3/4,n \right) }
\quad {\text {if $a=1+n$}},\\
\dfrac{{3}^{3\,n} \left( 3/4,n \right) 
 \left( 1/3,n \right)}{{2}^{2\,n}  \left( 11/12,n \right)   \left( 1/6,n \right)  }
\quad {\text {if $a=5/12+n$}},\\
-{\dfrac {13}{5}}\,\dfrac{{3}^{3\,n+1} \left(  25/12
,n \right) \left( 5/3,n \right)}{{2}^{2\,n}
 \left( 9/4,n \right)  \left( 3/2,n \right)  }
  \quad {\text {if $a=7/4+n$}},\\
  0
  \quad {\text {if $a=13/12+n$}}
\end{cases}
\end{flalign*}
(The first case is identical to Theorem 10 in [Ek]).
\begin{flalign*}
&{\rm (xix)}\,{\mbox{$F$}(1-3\,a,5/4-3\,a;\,2-4\,a;\, 8/9)}&\\
&=\begin{cases}
\dfrac{{2}^{5/6-2\,a}\Gamma  \left( 1/4 \right) \Gamma  \left( 2/3 \right) 
\Gamma  \left( 1/2-a \right) \Gamma  \left( 5/4-a \right) }{ {3}^
{3\,a-5/4} \Gamma  \left( 1/12 \right)   \Gamma  \left( 5/6 \right) 
\Gamma  \left(2/3 -a \right)  \Gamma  \left(13/12- a \right)},\\
\dfrac{-11\, \left( 7/3,n \right)  \left( 
23/12,n \right)}{ {2}^{2\,n-1}  {3
}^{3\,n+10}   \left( 5/2,n
 \right)  \left( 7/4,n \right) }
 \quad {\text {if $a=2+n$}},\\
-{\dfrac {13}{5}}\,
\dfrac{\left( 5/3,n \right)  \left( 25/12,n \right)}{ {2}^{2\,n}
  {3}^{3\,n+7}  \left( 9/4,n \right)  \left( 3/2,n \right) }
  \quad {\text {if $a=7/4+n$}}\\
\end{cases}
&
\end{flalign*}
(The first case is identical to (3.5) in [Ka]).
\begin{flalign*}
&{\rm (xx)}\,
{\mbox{$F$}(1-a,3/4-a;\,2-4\,a;\,8/9)}&\\
&=\begin{cases}
\dfrac{{2}^{5/6-2\,a}\Gamma  \left( 1/4 \right) \Gamma  \left( 2/3 \right) 
\Gamma  \left( 1/2-a \right) \Gamma  \left( 5/4-a \right) }{ {3}^
{-a-3/4} \Gamma  \left( 1/12 \right)   \Gamma  \left( 5/6 \right) 
\Gamma  \left(2/3 -a \right)  \Gamma  \left(13/12- a \right)},\\
\dfrac{{3}^{n} \left( 4/3,n \right)  \left( 
11/12,n \right) }{ {2}^{2\,n}   \left( 3/2,n \right)   \left( 3/4,n \right) }
\quad {\text {if $a=1+n$}},\\
\dfrac{{3}^{n} \left( 13/12,n \right)  \left( 2/3,n \right)}{ {2}^{2\,n} \left( 5/4,n \right)  \left( 1/2,n \right) }
\quad {\text {if $a=3/4+n$}}.\\
\end{cases}
&
\end{flalign*}
\begin{flalign*}
&{\rm (xxi)}\,{\mbox{$F$}(3\,a,a+1/4;\,5/4;\,9)}
=\begin{cases}
\dfrac{{2}^{6\,n} \left( 1/3,n \right) }{\left( 13/12,n \right)  }
 \quad {\text {if $a=-n$}},\\
 \dfrac{{2}^{6\,n+3} \left( 2/3,n \right) }{5\left( 17/12,n \right)  }
 \quad {\text {if $a=-1/3-n$}},\\
0 
 \quad {\text {if $a=-2/3-n$}},\\
 \dfrac{{2}^{6\,n} \left( 7/12,n \right) }{\left( 4/3,n \right)  }
 \quad {\text {if $a=-1/4-n$}}.\\
\end{cases}
&
\end{flalign*}
\begin{flalign*}
&{\rm (xxii)}\, {\mbox{$F$}(1-a,5/4-3\,a;\,5/4;\,9)}
=\begin{cases}
\dfrac{{2}^{6\,n} \left( 11/12,n \right) }{\left( 5/3,n \right)  }
 \quad {\text {if $a=1+n$}},\\
\dfrac{{2}^{6\,n} \left( 1/3,n \right) }{\left( 13/12,n \right)  }
 \quad {\text {if $a=5/12+n$}},\\ 
\dfrac{-{2}^{6\,n+2} \left( 2/3,n \right) }{5\left( 17/12,n \right)  }
 \quad {\text {if $a=3/4+n$}},\\ 
 0
 \quad {\text {if $a=13/12+n$}}.\\ 
 \end{cases}
&
\end{flalign*}
\begin{flalign*}
&{\rm (xxiii)}\,{\mbox{$F$}(3\,a,1-a;\,5/4;\, 9/8)}
=\begin{cases}
\dfrac{ \left( 1/3,n \right) }{({-2})^{3\,n}\left( 13/12,n \right)  }
 \quad {\text {if $a=-n$}},\\
\dfrac{ -\left( 2/3,n \right) }{5\left({-2}\right)^{3\,n}\left( 17/12,n \right)  }
 \quad {\text {if $a=-1/3-n$}},\\
 0
 \quad {\text {if $a=-2/3-n$}},\\
 \dfrac{ ({-2})^{3\,n}\left( 11/12,n \right) }{\left( 5/3,n \right)  }
 \quad {\text {if $a=1+n$}}.\\
\end{cases}&
\end{flalign*}
\begin{flalign*}
&{\rm (xxiv)}\,{\mbox{$F$}(5/4-3\,a,a+1/4;\,5/4;\,9/8)}&\\
&=\begin{cases}
\dfrac{ \left( 1/3,n \right) }{({-2})^{3\,n}\left( 13/12,n \right)  }
 \quad {\text {if $a=5/12+n$}},\\
 \dfrac{ -\left( 2/3,n \right) }{5\left({-2}\right)^{3\,n+1}\left( 17/12,n \right)  }
 \quad {\text {if $a=3/4+n$}},\\
 0
 \quad {\text {if $a=13/12+n$}},\\
 \dfrac{ ({-2})^{3\,n}\left( 7/12,n \right) }{\left( 4/3,n \right)  }
 \quad {\text {if $a=-1/4-n$}}.\\
\end{cases}
&
\end{flalign*}
 
\paragraph{(1,3,4-5), (1,3,4-6), (1,3,4-8)}
The special values obtained from (1,3,4-5), (1,3,4-6) and (1,3,4-8) coincide with those obtained from (1,3,4-4).

\paragraph{(1,3,4-7)}
The special values obtained from (1,3,4-7) coincide with those obtained from (1,3,4-3).

\subsubsection{$(k,l,m)=(1,4,4)$}
In this case, we have
\begin{align}
&(a,b,c,x)=(a,4\,a+1,4\,a,4/3),\,
S^{(n)}=
{\frac { \left( -1 \right) ^{n}{2}^{8\,n} \left( a+1,n
 \right)  \left( 3\,a,3\,n \right) }{{3}^{4\,n}\left( 4\,a,4\,n \right) }},
\tag{1,4,4-1}\\
&(a,b,c,x)=(a,4\,a-4,4\,a-2,4/3),\,
S^{(n)}={\frac {4\,a-3}{ \left( -3 \right) ^{n} \left( 4\,a-3+4\,n \right) }}. \tag{1,4,4-2}
\end{align}
\paragraph{(1,4,4-1)}
\begin{flalign*}
&{\rm (i)}\,
{\mbox{$F$}(a,4\,a+1;\,4\,a;\,4/3)}&\\
&=\begin{cases}
0
\quad {\text {if $a=-1-n$}},\\
{\dfrac { \left( -1 \right) ^{n}{2}^{8\,n} \left( 1/4,n
 \right)  \left( 7/4,3\,n \right) }{{3}^{4\,n} \left( 2,4\,n \right) }}
\quad {\text {if $a=-1/4-n$}},\\
{\dfrac { \left( -1 \right) ^{n}{2}^{8\,n+1} \left( 1/2,n
 \right)  \left( 5/2,3\,n \right) }{{3}^{4\,n+1} \left( 3,4\,n \right) }}
 \quad {\text {if $a=-1/2-n$}},\\
{\dfrac { 5\left( -1 \right) ^{n}{2}^{8\,n-1} \left( 3/4,n
 \right)  \left( 13/4,3\,n \right) }{{3}^{4\,n+2} \left( 4,4\,n \right) }}
 \quad {\text {if $a=-3/4-n$}}.
\end{cases}&
\end{flalign*}
The special values obtained from {\rm (ii)} and {\rm (iii) } are trivial.
\begin{flalign*}
&{\rm (iv)}\,
{\mbox{$F$}(3\,a,4\,a+1;\,4\,a;\,4)}&\\
&=\begin{cases}
0
\quad {\text {if $a=-1/3-1/3\,n$}},\\
{\dfrac { \left( -1 \right) ^{n}{2}^{8\,n} \left( 1/4,n
 \right) \left( 7/4,3\,n \right) }{
 \left( 2,4\,n \right) }}
\quad {\text {if $a=-1/4-n$}},\\
{\dfrac { \left( -1 \right) ^{n+1}{2}^{8\,n+1} \left( 1/2,n
 \right) \left( 5/2,3\,n \right) }{
 \left( 3,4\,n \right) }}
\quad {\text {if $a=-1/2-n$}},\\
{\dfrac { 5\left( -1 \right) ^{n}{2}^{8\,n-1} \left( 3/4,n
 \right) \left( 13/4,3\,n \right) }{
 \left( 4,4\,n \right) }}
\quad {\text {if $a=-3/4-n$}}.
\end{cases}
\end{flalign*}
\begin{flalign*}
&{\rm (v)}\,
{\mbox{$F$}(a,4\,a+1;\,a+2;\,-1/3)}
=
{2}^{-8\,a}{3}^{4\,a} \left( a+1 \right)& 
\end{flalign*}
(The above is a special case of (\ref{1.6eb2})).
\begin{flalign*}
&{\rm (vi)}\,{\mbox{$F$}(2,1-3\,a;\,a+2;\,-1/3)}
=3/4\,a+3/4
&
\end{flalign*}
(The above is a special case of (\ref{1.5eb2})).
\begin{flalign*}
&{\rm (vii)}\,{\mbox{$F$}(a,1-3\,a;\,a+2;\,1/4)}
={2}^{-6\,a}{3}^{3\,a} \left( a+1 \right) 
&
\end{flalign*}
(The above is a special case of (\ref{1.6eb2})).
\begin{flalign*}
&{\rm (viii)}\,
{\mbox{$F$}(2,4\,a+1;\,a+2;\,1/4)}
=4/3\,a+4/3
&
\end{flalign*}
(The above is a special case of (\ref{1.5eb2})).
\begin{flalign*}
&{\rm (ix)}\,
{\mbox{$F$}(a,1-3\,a;\,-3\,a;\,3/4)}
=\begin{cases}
0,\\
{\dfrac {{3}^{3\,n} \left( 4/3,4\,n \right) }{{2}^{6\,n
} \left( 4/3,n \right)  \left( 1,3\,n \right) }}
\quad {\text {if $a=1/3+n$}},\\
{\dfrac {5\cdot {3}^{3\,n} \left( 8/3,4\,n \right) }{{2}^{6\,n+2
} \left( 5/3,n \right)  \left( 2,3\,n \right) }}
 \quad {\text {if $a=2/3+n$}},\\
{\dfrac { {3}^{3\,n+3} \left( 4,4\,n \right) }{{2}^{6\,n+4
} \left( 2,n \right)  \left( 3,3\,n \right) }}
 \quad {\text {if $a=1+n$}}.\\
\end{cases}
&
\end{flalign*}
The special values obtained from {\rm (x)} and {\rm (xi) } are trivial.
\begin{flalign*}
&{\rm (xii)}\,
{\mbox{$F$}(-4\,a,1-3\,a;\,-3\,a;\,-3)}
=\begin{cases}
0
 \quad {\text {if $a=1/4+n$}},\\
0
 \quad {\text {if $a=1/2+n$}},\\
0
 \quad {\text {if $a=3/4+n$}},\\
{\dfrac {{3}^{3\,n} \left( 4/3,4\,n \right) }{ \left( 4/3,n \right)  \left( 1,3\,n \right) }}
  \quad {\text {if $a=1/3+n$}},\\
{\dfrac {5\cdot {3}^{3\,n} \left( 8/3,4\,n \right) }{ \left( 5/3,n \right)  \left( 2,3\,n \right) }}
 \quad {\text {if $a=2/3+n$}},\\
{\dfrac {{3}^{3\,n+3} \left( 4,4\,n \right) }{ \left( 2,n \right)  \left( 3,3\,n \right) }}
  \quad {\text {if $a=1+n$}}.\\
 \end{cases}
&
\end{flalign*}
\begin{flalign*}
&{\rm (xiii)}\,
{\mbox{$F$}(2,4\,a+1;\,3\,a+2;\,3/4)} 
= 12\,a+4
 &
\end{flalign*}
(The above is a special case of (\ref{1.5eb2})).
\begin{flalign*}
&{\rm (xiv)}\,{\mbox{$F$}(3\,a,1-a;\,3\,a+2;\,3/4)}
={2}^{-2\,a} \left( 3\,a+1 \right) 
&
\end{flalign*}
(The above is a special case of (\ref{1.6eb2})).
\begin{flalign*}
&{\rm (xv)}\,
{\mbox{$F$}(3\,a,4\,a+1;\,3\,a+2;\,-3)}
=\begin{cases}
{2}^{8\,n} \left( 1-12\,n \right) 
 \quad {\text {if $a=-1/4-n$}},\\
 {2}^{8\,n+2} \left( -2-12\,n \right) 
 \quad {\text {if $a=-1/2-n$}},\\
{2}^{8\,n+4} \left( -5-12\,n \right) 
\quad {\text {if $a=-3/4-n$}}
\end{cases}
&
\end{flalign*}
(The above are special cases of (\ref{1.6eb2})).
\begin{flalign*}
&{\rm (xvi)}\,
{\mbox{$F$}(2,1-a;\,3\,a+2;\,-3)}
=3/4\,n+1
\quad {\text {if $a=1+n$}}&
\end{flalign*}
(The above is a special case of (\ref{1.5eb2})).
\begin{flalign*}
&{\rm (xvii)}\,
{\mbox{$F$}(2,1-3\,a;\,2-4\,a;\,4/3)}
=
4\,n+1
\quad {\text {if $a=1/3+1/3\,n$}}
&
\end{flalign*}
(The above is a special case of (\ref{1.5eb2})).
\begin{flalign*}
&{\rm (xviii)}\,
{\mbox{$F$}(-4\,a,1-a;\,2-4\,a;\,4/3)}
=
- \left( -3 \right) ^{-n-1} \left( 4\,n+3 \right) 
\quad {\text {if $a=1+n$}}&
\end{flalign*}
(The above is a special case of (\ref{1.6eb2})).
\begin{flalign*}
&{\rm (xix)}\,
{\mbox{$F$}(-4\,a,1-3\,a;\,2-4\,a;\,4)}
=
\left( -3 \right) ^{n} \left( 4\,n +1\right) 
\quad {\text {if $a=1/3+1/3\,n$}}
&
\end{flalign*}
(The above is a special case of (\ref{1.6eb2})).
\begin{flalign*}
&{\rm (xx)}\,
{\mbox{$F$}(2,1-a;\,2-4\,a;\,4)}
=
4/3\,n+1
\quad {\text {if $a=1+n$}}
&
\end{flalign*}
(The above is a special case of (\ref{1.5eb2})).
The special values obtained from {\rm (xxi)} are trivial.
\begin{flalign*}
&{\rm (xxii)}\,
{\mbox{$F$}(-4\,a,1-a;\,-a;\,-1/3)}
=\begin{cases}
0,\\
{\dfrac { \left( 4,4\,n \right) }{{3}^{n} \left( 2,n \right)  \left( 3,3\,n \right) }}
\quad {\text {if $a=1+n$}}.\\
\end{cases}
&
\end{flalign*}
\begin{flalign*}
&{\rm (xxiii)}\,
{\mbox{$F$}(3\,a,1-a;\,-a;\,1/4)}
=\begin{cases}
0,\\
{\dfrac { \left( 4,4\,n \right) }{{2}^{2\,n} \left( 2,n \right)  \left( 3,3\,n \right) }}
\quad {\text {if $a=1+n$}}
\end{cases}
&
\end{flalign*}
(The first case is identical to (29.3) in [Ge]).
The special values obtained from {\rm (xxiv)} are trivial.

\paragraph{(1,4,4-2)}
The special values obtained from (1,4,4-2) coincide with those obtained from (1,4,4-1).

\subsubsection{$(k,l,m)=(1,5,4)$}
In this case, there is no admissible quadruple.

\subsubsection{$(k,l,m)=(2,2,4)$}
In this case, we have
\begin{align}
(a,b,c,x)=(a,b,2\,a,2), \tag{2,2,4-1}\\
(a,b,c,x)=(a,b,2\,b,2). \tag{2,2,4-2}
\end{align}
\paragraph{(2,2,4-1), (2,2,4-2)}
The special values obtained from (2,2,4-1) and (2,2,4-2) coincide
with those obtained from (1,2,2-1).

\subsubsection{$(k,l,m)=(2,3,4)$}
In this case, we have
\begin{align}
&
\begin{cases}
(a,b,c,x)=(a,3/2\,a-1/4,2\,a,8+4\,\sqrt {3}),\,\\
S^{(n)}={\dfrac {{3}^{3/2\,n} \left( \sqrt {3}+2 \right) ^{3\,n} \left( 1/2\,a+7/12,n \right) }{{2}^{2\,n} \left( 1/2\,a+3/4,n
 \right) }},
\tag{2,3,4-1}
\end{cases}
\\
&
\begin{cases}
(a,b,c,x)=(a,3/2\,a-1/4,2\,a,8-4\,\sqrt {3}),\,\\
S^{(n)}={\dfrac {{3}^{3/2\,n} \left( \sqrt {3}-2 \right) ^{3\,n} \left( 1/2\,a+7/12,n \right) }{{2}^{2\,n} \left( 1/2\,a+3/4,n
 \right) }}.
\tag{2,3,4-2}
\end{cases}
\end{align}
\paragraph{(2,3,4-1)}
\begin{flalign*}
&{\rm (i)}\,{\mbox{$F$}(a,3/2\,a-1/4;\,2\,a;\,8+4\,\sqrt {3})}&\\
&=\begin{cases}
\dfrac{ 5\cdot {3}^{3/2\,n+1/2} \left( 2+\sqrt {3}
 \right) ^{3\,n+3} \left( 17/12,n \right) }{ {2}^{2\,n+2}  \left( 5/4,
n \right)  }
\quad {\text {if $a=-2-2\,n$}},\\
\dfrac{-  {3}^{3/2\,n+1/2} \left( 2+\sqrt {3}
 \right) ^{3\,n+2} \left( 11/12,n \right) }{ {2}^{2\,n+1}  \left( 3/4,
n \right)  }
\quad {\text {if $a=-1-2\,n$}},\\
\dfrac{ {3}^{3/2\,n} \left( 2+\sqrt {3}
 \right) ^{3\,n} \left( 1/3,n \right) }{ {2}^{2\,n}  \left( 1/6,
n \right)  }
\quad {\text {if $a=1/6-2\,n$}},\\
\dfrac{ -{3}^{3/2\,n+1} \left( 2+\sqrt {3}
 \right) ^{3\,n+4} \left( 5/3,n \right) }{ {2}^{2\,n}  \left( 3/2,
n \right)  }
\quad {\text {if $a=-5/2-2\,n$}},\\
0
\quad {\text {if $a=-7/6-2\,n$}}\\
\end{cases}&
\end{flalign*}
(The first case is identical to Theorem 28 in [Ek]).
\begin{flalign*}
&{\rm (ii)}\,{\mbox{$F$}(a,1/2\,a+1/4;\,2\,a;\,8+4\,\sqrt {3})}&\\
&=\begin{cases}
\dfrac{  5\cdot {3}^{3/2\,n+1/2} \left( -2-\sqrt {3}
 \right) ^{n+1} \left( 17/12,n \right) }{ {2}^{2\,n+2}  \left( 5/4,
n \right)  }
\quad {\text {if $a=-2-2\,n$}},\\
\dfrac{-  {3}^{3/2\,n+1/2} \left( -2-\sqrt {3}
 \right) ^{n} \left( 11/12,n \right) }{ {2}^{2\,n+1}  \left( 3/4,
n \right)  }
\quad {\text {if $a=-1-2\,n$}},\\
\dfrac{  {3}^{3/2\,n} \left( -2-\sqrt {3}
 \right) ^{n} \left( 2/3,n \right) }{ {2}^{2\,n}  \left( 1/2,
n \right)  }
\quad {\text {if $a=-1/2-2\,n$}}\\
\end{cases}&
\end{flalign*}
(The third case is identical to Theorem 16 in [Ek]).
\begin{flalign*}
&{\rm (iii)}\,{\mbox{$F$}(a,1/2\,a+1/4;\,2\,a;\,8-4\,\sqrt {3})}&\\
&=\begin{cases}
\dfrac{  5\cdot {3}^{3/2\,n+1/2} \left( 2-\sqrt {3}
 \right) ^{n+1} \left( 17/12,n \right) }{ {2}^{2\,n+2}  \left( 5/4,
n \right)  }
\quad {\text {if $a=-2-2\,n$}},\\
\dfrac{  {3}^{3/2\,n+1/2} \left( 2-\sqrt {3}
 \right) ^{n} \left( 11/12,n \right) }{ {2}^{2\,n+1}  \left( 3/4,
n \right)  }
\quad {\text {if $a=-1-2\,n$}},\\
\dfrac{  {3}^{3/2\,n} \left( 2-\sqrt {3}
 \right) ^{n} \left( 2/3,n \right) }{ {2}^{2\,n}  \left( 1/2,
n \right)  }
\quad {\text {if $a=-1/2-2\,n$}}
\end{cases}&
\end{flalign*}
(The third case is identical to Theorem 16 in [Ek]).
\begin{flalign*}
&{\rm (iv)}\,{\mbox{$F$}(a,3/2\,a-1/4;\,2\,a;\,8-4\,\sqrt {3})}&\\
&=\begin{cases}
\dfrac{  5\cdot {3}^{3/2\,n+1/2} \left( \sqrt {3}-2
 \right) ^{3\,n+3} \left( 17/12,n \right) }{ {2}^{2\,n+2}  \left( 5/4,
n \right)  }
\quad {\text {if $a=-2-2\,n$}},\\
\dfrac{  {3}^{3/2\,n+1/2} \left( \sqrt {3}-2
 \right) ^{3\,n+2} \left( 11/12,n \right) }{ {2}^{2\,n+1}  \left( 3/4,
n \right)  }
\quad {\text {if $a=-1-2\,n$}},\\
\dfrac{  {3}^{3/2\,n} \left( \sqrt {3}-2
 \right) ^{3\,n} \left( 1/3,n \right) }{ {2}^{2\,n}  \left( 1/6,
n \right)  }
\quad {\text {if $a=1/6-2\,n$}},\\
\dfrac{  -{3}^{3/2\,n+1} \left( \sqrt {3}-2
 \right) ^{3\,n+4} \left( 5/3,n \right) }{ {2}^{2\,n}  \left( 3/2,
n \right)  }
\quad {\text {if $a=-5/2-2\,n$}},\\
0
\quad {\text {if $a=-7/6-2\,n$}}\\
\end{cases}&
\end{flalign*}
(The first case is identical to Theorem 28 in [Ek]).
\begin{flalign*}
&{\rm (v)}\,{\mbox{$F$}(a,3/2\,a-1/4;\,1/2\,a+3/4;\,-7-4\,\sqrt {3})}&\\
&=\begin{cases}
\dfrac{  2^{4\,n}{3}^{3/2\,n} \left( -2-\sqrt {3}
 \right) ^{3\,n} \left( 5/12,n \right) }{ \left( 1/4,
n \right)  }
\quad {\text {if $a=-2\,n$}},\\
\dfrac{ - 2^{4\,n+2}{3}^{3/2\,n+1/2} \left( -2-\sqrt {3}
 \right) ^{3\,n+2} \left( 11/12,n \right) }{ \left( 3/4,
n \right)  }
\quad {\text {if $a=-1-2\,n$}},\\
\dfrac{  2^{4\,n}{3}^{3/2\,n} \left( -2-\sqrt {3}
 \right) ^{3\,n} \left( 1/3,n \right) }{ \left( 1/6,
n \right)  }
\quad {\text {if $a=1/6-2\,n$}},\\
\dfrac{  2^{4\,n+1}{3}^{3/2\,n+1/2} \left( -2-\sqrt {3}
 \right) ^{3\,n+1} \left( 2/3,n \right) }{ \left( 1/2,
n \right)  }
\quad {\text {if $a=-1/2-2\,n$}},\\
0
\quad {\text {if $a=-7/6-2\,n$}}\\
\end{cases}&
\end{flalign*}
(The first case is identical to Theorem 31 in [Ek]).
\begin{flalign*}
&{\rm (vi)}\,
{\mbox{$F$}(1-a,3/4-1/2\,a;\,1/2\,a+3/4;\,-7-4\,\sqrt {3})}&\\
&=\begin{cases}
\dfrac{ {2}^{4\,n} \left(-2 -\sqrt {3}
 \right) ^{n} \left( 5/4,n \right) }{ 3^{3/2\,n}\left( 13/12,
n \right)  }
\quad {\text {if $a=1+2\,n$}},\\
\dfrac{ -{2}^{4\,n+2} \left(-2 -\sqrt {3}
 \right) ^{n} \left( 7/4,n \right) }{ 7\cdot 3^{3/2\,n-1/2}\left( 19/12,
n \right)  }
\quad {\text {if $a=2+2\,n$}},\\
\dfrac{ {2}^{4\,n} \left(-2 -\sqrt {3}
 \right) ^{n} \left( 3/2,n \right) }{  3^{3/2\,n}\left( 4/3,
n \right)  }
\quad {\text {if $a=3/2+2\,n$}}.\\
\end{cases}&
\end{flalign*}
\begin{flalign*}
&{\rm (vii)}\,
{\mbox{$F$}(a,1-a;\,1/2\,a+3/4;\,1/2+1/4\,\sqrt {3})}&\\
&=\dfrac{4 \left( 2-\sqrt {3} \right) ^{1/2\,a-1/4}\sqrt {\pi }\sin \left( 
\pi \left( 1/2\,a+1/12 \right)  \right) \Gamma  \left( 1
/2\,a+3/4 \right)}{ {3}^{3/4\,a+3/8} 
\Gamma  \left( 2/3 \right)  \Gamma  \left( 1/2\,a
+7/12\right)  }
&
\end{flalign*}
(The above is a generalization of Theorem 37 in [Ek]).
We derive the above formula using the algebraic transformation
\begin{align*}
&F(a, 1-a, 1/2\,a+3/4, x)\\
&=(16x^2-16x+1)^{1/8-3/4\,a}(1-x)^{1/2\,a-1/4}\\
&\times F(1/4\,a-1/24, 1/4\,a+7/24, 1/2\,a+3/4, -108x(1-x)/(16x^2-16x+1)^3)
\label{algebraic_trans_2}
\end{align*}
and connection formulae for the hypergeometric series (see (25)--(44) in 2.9 in [Erd]).
\begin{flalign*}
&{\rm (viii)}\,
{\mbox{$F$}(3/4-1/2\,a,3/2\,a-1/4;\,1/2\,a+3/4;\,1/2+1/4\,\sqrt {3})}&\\
&=\dfrac{2^{a+3/2}\, \sqrt {\pi }\sin \left( 
\pi \left( 1/2\,a+1/12 \right)  \right) \Gamma  \left( 1
/2\,a+3/4 \right)}{ {3}^{3/4\,a+3/8} 
\Gamma  \left( 2/3 \right)  \Gamma  \left( 1/2\,a
+7/12\right)  }.
&
\end{flalign*}
\begin{flalign*}
&{\rm (ix)}\,
{\mbox{$F$}(a,1-a;\,5/4-1/2\,a;\,1/2-1/4\,\sqrt {3})}&\\
&=
\dfrac{\left( 2-\sqrt {3} \right) ^{1/2\,a-1/4}
\Gamma (1/3) \Gamma  \left( 5/4-1
/2\,a \right)}{ {3}^{5/8-3/4\,a} \sqrt{\pi}
\Gamma  \left( 13/12-1/2\,a \right)    }
&
\end{flalign*}
(The above is a generalization of Theorem 37 in [Ek]).
\begin{flalign*}
&{\rm (x)}\,
{\mbox{$F$}(5/4-3/2\,a,1/2\,a+1/4;\,5/4-1/2\,a;\,1/2-1/4\,\sqrt {3})}\\
&=
\dfrac{2^{1/2-a } \,
\Gamma (1/3) \Gamma  \left( 5/4-1
/2\,a \right)}{ {3}^{5/8-3/4\,a} \sqrt{\pi}
\Gamma  \left( 13/12-1/2\,a \right)  }.
&
\end{flalign*}
\begin{flalign*}
&{\rm (xi)}\,{\mbox{$F$}(a,1/2\,a+1/4;\,5/4-1/2\,a;\,-7+4\,\sqrt {3})}\\
&=\dfrac{2^{-2\,a } \,
 \left( 2+\sqrt {3} \right) ^{1/2\,a+1/4}\Gamma (1/3) \Gamma  \left( 5/4-1
/2\,a \right)}{ {3}^{5/8-3/4\,a} \sqrt{\pi}
\Gamma  \left( 13/12-1/2\,a \right)  }.
&
\end{flalign*}
\begin{flalign*}
&{\rm (xii)}\,
{\mbox{$F$}(1-a,5/4-3/2\,a;\,5/4-1/2\,a;\,-7+4\,\sqrt {3})}\\
&=\dfrac{2^{2\,a-2} \,
 \left( 2-\sqrt {3} \right) ^{3/2\,a-5/4}\Gamma (1/3) \Gamma  \left( 5/4-1
/2\,a \right)}{ {3}^{5/8-3/4\,a} \sqrt{\pi}
\Gamma  \left( 13/12-1/2\,a \right)  }
&
\end{flalign*}
(The above is a generalization of Theorem 31 in [Ek]).
The special values obtained from {\rm (xiii)-(xxiv)} are contained 
in the above those.
\paragraph{(2,3,4-2)}
The special values obtained from (2,3,4-2) coincide with those 
obtained from (2,3,4-1).
\subsubsection{$(k,l,m)=(2,4,4)$}
In this case, we have
\begin{align}
&(a,b,c,x)=(a,b,2\,a,2),\, 
 \tag{2,4,4-1}\\
 &
 \begin{cases}
 (a,b,c,x)=(a,2\,a-1/2,2\,a,-2+2\,\sqrt {2}), \,\\
 S^{(n)}=
 {\dfrac {{2}^{2\,n} \left( \sqrt {2}-1 \right) ^{4\,n} \left( 1/2\,a+3/8,n \right)  \left( 1/2\,a+5/8,n
 \right) }{ \left( a+1/2,2\,n \right) }},
\tag{2,4,4-2}
\end{cases}\\
&
\begin{cases}
(a,b,c,x)=(a,2\,a-1/2,2\,a,-2-2\,\sqrt {2}), \,\\
 S^{(n)}=
 {\dfrac {{2}^{2\,n} \left( \sqrt {2} +1\right) ^{4\,n} \left( 1/2\,a+3/8,n \right)  \left( 1/2\,a+5/8,n
 \right) }{ \left( a+1/2,2\,n \right) }}.
\tag{2,4,4-3}
\end{cases}
\end{align}
\paragraph{(2,4,4-1)}
The special values obtained from (2,4,4-1) are evaluated
in the case (1,2,2-1).
\paragraph{(2,4,4-2)}
\begin{flalign*}
&{\rm (i)}\,{\mbox{$F$}(a,2\,a-1/2;\,2\,a;\,-2+2\,\sqrt {2})}&\\
&=\begin{cases}
{\dfrac {{2}^{1/4-a}\left( \sqrt {2}-1 \right) ^{1/2-2\,a}\sqrt {\pi }\Gamma  \left( a+1/2 \right) 
 }{\Gamma  \left( 1/2\,a+3/8
 \right) \Gamma  \left( 1/2\,a+5/8 \right) }},\\
\dfrac{5 \left( \sqrt {2}-1 \right) ^{4\,n+4}
\left( 
13/8,n \right)  \left( 11/8,n \right) }{4\left( 7
/4,n \right)\left( 5/4,n \right) } 
\quad {\text {if $a=-2-2\,n$}},\\
\dfrac{
 - \left( \sqrt {2}-1 \right) ^{4\,n+3}
 \left( 9/8,n \right) \left( 7/8
,n \right)  }{2\left( 5/4,n \right)\left( 3/4,n \right) }
\quad {\text {if $a=-1-2\,n$}}.
\end{cases}
&
\end{flalign*}
\begin{flalign*}
&{\rm (ii)}\,{\mbox{$F$}(1/2,a;\,2\,a;\,-2+2\,\sqrt {2})}&\\
&=\begin{cases}
{\dfrac {{{2}^{1/4-a}\left(1+\,\sqrt {2}\right)}^{1/2}\sqrt {\pi }\Gamma  \left( a+1/2 \right) }{
\Gamma  \left( 1/2\,a+3/8 \right) \Gamma 
 \left( 1/2\,a+5/8 \right) }},\\
 \dfrac{5 \left( 13/8,n \right)  \left(  11/8,n \right)}{4 \left( 7/4,n \right)\left( 5/4,n \right)}    
 \quad {\text {if $a=-2-2\,n$}},\\
 \dfrac{\left( 1+\sqrt {2} \right) \left( 9/8,n \right) 
 \left( 7/8,n \right)  }{2\left( 5/4,n \right)\left( 3/4,n \right)}
 \quad {\text {if $a=-1-2\,n$}}
\end{cases}
&
\end{flalign*}
(The second case is identical to Theorem 33 in [Ek]).
\begin{flalign*}
&{\rm (iii)}\,{\mbox{$F$}(1/2,a;\,2\,a;\,-2-2\,\sqrt {2})}&\\
&=\begin{cases}
\dfrac{5 \left(  13/8,n \right)  \left( 11/8,n \right)}{ 4\left( 7/4,n \right)\left( 5/4,n \right)}    
\quad {\text {if $a=-2-2\,n$}},\\
\dfrac{  \left( 1-\sqrt {2} \right)
 \left( 9/8,n \right)
 \left(  7/8,n \right) }{2\left( 5/4,n \right)\left( 3/4,n \right)}
\quad {\text {if $a=-1-2\,n$}}
\end{cases}
&
\end{flalign*}
(The first case is identical to Theorem 33 in [Ek]).
\begin{flalign*}
&{\rm (iv)}\,{\mbox{$F$}(a,2\,a-1/2;\,2\,a;\,-2-2\,\sqrt {2})}&\\
&=\begin{cases}
\dfrac{5 \left( 1+\sqrt {2}
 \right) ^{4\,n+4}\left( 13/8,n \right)  \left( 11/8,n \right)}{4
 \left( 7/4,n \right)\left( 5/4,n  \right) }   
\quad {\text {if $a=-2-2\,n$}},\\
 \dfrac{  \left( 1+\sqrt {2} \right) ^{4\,n+3}
\left( 9/8,n \right) \left( 7/8,n \right)}{2 \left( 5/4,n \right)\left( 3/4,n \right)}     
\quad {\text {if $a=-1-2\,n$}},\\
\dfrac{{ \left( 1+\sqrt {2} \right) ^{4\,n}\left( 1/2,n \right) 
 \left( 1/4,n \right)  }}{{ \left( 5/8,n \right) \left( 1/8,n \right) }} 
\quad {\text {if $a=1/4-2\,n$}},\\
\dfrac{ \sqrt{2}\left( 1+\sqrt {2} \right) ^{4\,n+1} \left( 3/4,n \right) 
 \left( 1/2,n \right)    }{{ \left( 7/8,n \right)\left( 3/8,n \right) }}
\quad {\text {if $a=-1/4-2\,n$}},\\
0
\quad {\text {if $a=-3/4-2\,n$}},\\
0
\quad {\text {if $a=-5/4-2\,n$}}.
\end{cases}
&
\end{flalign*}
\begin{flalign*}
&{\rm (v)}\,{\mbox{$F$}(a,2\,a-1/2;\,a+1/2;\,3-2\,\sqrt {2})}
=
{\dfrac {{2}^{3/4-3\,a} \left( \sqrt {2}-1 \right) ^{1/2-2\,a}\sqrt {
\pi }\Gamma  \left( a+1/2 \right) }{\Gamma  \left( 1/2\,a+3/8 \right) 
\Gamma  \left( 1/2\,a+5/8 \right) }}
&
\end{flalign*}
(The above is a generalization of Theorem 25 in [Ek]).
\begin{flalign*}
&{\rm (vi)}\,{\mbox{$F$}(1/2,1-a;\,a+1/2;\,3-2\,\sqrt {2})}
=
{\dfrac {{2}^{-a-1/4} ({1+\sqrt {2}})^{1/2}\sqrt {\pi }\Gamma  \left( a+1/
2 \right) }{\Gamma  \left( 1/2\,a+3/8 \right) \Gamma  \left( 1/2\,a+5/
8 \right) }}.
&
\end{flalign*}
\begin{flalign*}
&{\rm (vii)}\,{\mbox{$F$}(a,1-a;\,a+1/2;\,1/2-1/2\,\sqrt {2})}
=
{\frac {{2}^{3/4-2\,a} \left( \sqrt {2}-1 \right) ^{1/2-a}\sqrt {\pi 
}\Gamma  \left( a+1/2 \right) }{\Gamma  \left( 1/2\,a+3/8 \right) 
\Gamma  \left( 1/2\,a+5/8 \right) }}
&
\end{flalign*}
(The above is a generalization of Theorem 36 in [Ek]).
\begin{flalign*}
&{\rm (viii)}\,{\mbox{$F$}(1/2,2\,a-1/2;\,a+1/2;\,1/2-1/2\,\sqrt {2})}
=
{\frac {{2}^{1/4-a}\sqrt {\pi }\Gamma  \left( a+1/2 \right) }{\Gamma 
 \left( 1/2\,a+3/8 \right) \Gamma  \left( 1/2\,a+5/8 \right) }}.
&
\end{flalign*}
\begin{flalign*}
&{\rm (ix)}\,
{\mbox{$F$}(a,1-a;\,3/2-a;\,1/2+1/2\,\sqrt {2})}&\\
&=
\begin{cases}
{\dfrac { \left( \sqrt {2}-1 \right) ^{2\,n}\left( 3/4,n \right)
\left( 5/4,n \right) 
 }{2^{2\,n}\,
 \left( 7/8,n \right) \left( 9/8,n \right)  }}
\quad {\text {if $a=-2\,n$}},\\
\dfrac{1}{5}\,\dfrac{\left( \sqrt {2}-1 \right) ^{2\,n+2} \left( 5/4,n \right) \left( 7/4,n \right) 
 }{ {2}^{2\,n} 
 \left( 11/8,n \right)  \left( 13/8,n \right)   }  
\quad {\text {if $a=-1-2\,n$}},\\
{\dfrac {{2}^{2\,n} \left( 1+\sqrt {2} \right) ^{2\,n}\left( 3/8,n \right)  \left( 5/8,n \right)  }{
 \left( 1/4,n \right)\left( 3/4,n \right) }}
\quad {\text {if $a=1+2\,n$}},\\
\dfrac{{2}^{2\,n} \left( 1+\sqrt {2}
 \right) ^{2\,n+2} \left( 7/8,n \right)  \left(9/8,n \right) }{\left( 3/4,n \right)\left( 5/4,n \right) }
\quad {\text {if $a=2+2\,n$}}
\end{cases}
&
\end{flalign*}
(The third case is identical to Theorem 36 in [Ek]).
\begin{flalign*}
&{\rm (x)}\, {\mbox{$F$}(1/2,3/2-2\,a;\,3/2-a;\,1/2+1/2\,\sqrt {2})}&\\
&=\begin{cases}
{\dfrac { \left( 1/4,n \right)\left( 1/2,n \right) }{ \left( 1/8,n
 \right)\left( 5/8,n\right) }}
\quad {\text {if $a=3/4+2\,n$}},\\
{\dfrac {-\sqrt{2}\left( 1/2,n \right) \left( 3/4,n \right) }{ \left( 3/8,
n \right)\left( 7/8,n \right) }}
\quad {\text {if $a=5/4+2\,n$}},\\
0
\quad {\text {if $a=7/4+2\,n$}},\\
0
\quad {\text {if $a=9/4+2\,n$}}.
\end{cases}
&
\end{flalign*}
\begin{flalign*}
&{\rm (xi)}\,{\mbox{$F$}(1/2,a;\,3/2-a;\,3+2\,\sqrt {2})}&\\
&=\begin{cases}
 \dfrac{\left( 3/4,n \right)\left( 5/4,n \right)}{    \left( 7/8,n \right) \left(  9/8,n \right) }
\quad {\text {if $a=-2\,n$}},\\
 \dfrac{2\left( 1-\sqrt {2}
 \right)\left( 5/4,n \right)\left( 7/4,n \right)   }{ 5 
 \left( 11/8,n \right)  \left( 13/8,n \right) }
\quad {\text {if $a=-1-2\,n$}}.
\end{cases}
&
\end{flalign*}
\begin{flalign*}
&{\rm (xii)}\, {\mbox{$F$}(1-a,3/2-2\,a;\,3/2-a;\,3+2\,\sqrt {2})}&\\
&=\begin{cases}
{\dfrac { \left( 4\,n +1\right) {2}^{4\,n} \left( 1+\sqrt {2}
 \right) ^{4\,n} \left( 3/8,n
 \right) \left( 5/8,n \right) }{ \left( 3/4,n \right)\left( 5/4,n \right) }}
\quad {\text {if $a=1+2\,n$}},\\
\dfrac{-{2}^{4\,n+1}\left( 1+\sqrt {2}
 \right) ^{4\,n+3} \left( 7/8,n \right)  \left( 9/8,n \right)   }{  \left( 3/4,n \right) \left( 5/4,n \right) }
\quad {\text {if $a=2+2\,n$}},\\
{\dfrac {{2}^{4\,n}\left( 1+\sqrt {2} \right) ^{4\,n} \left( 1/4,n \right) \left( 1/2,n \right)  }{ \left( 1/8,n \right)
\left( 5/8,n \right) }}
\quad {\text {if $a=3/4+2\,n$}},\\
{\dfrac {{2}^{4\,n+3/2}\left( 1+\sqrt {2} \right) ^{4\,n+1}\left( 1/2,n \right)  \left( 3/4,n \right)  
 }{ \left( 3/8,n
 \right)\left( 7/8,n
 \right) }}
\quad {\text {if $a=5/4+2\,n$}},\\
0
\quad {\text {if $a=7/4+2\,n$}},\\
0
\quad {\text {if $a=9/4+2\,n$}}\\
\end{cases}
&
\end{flalign*}
(The first case is identical to Theorem 25 in [Ek]).
The special values obtained from {\rm (xiii)-(xxiv)} coincide with
the above.
\paragraph{(2,4,4-3)}
The special values obtained from (2,4,4-3) coincide with
those obtained from (2,4,4-2).

\subsubsection{$(k,l,m)=(2,5,4)$}
In this case, we have
\begin{align}
&\begin{cases}
(a,b,c,x)=(a,5/2\,a-1,2\,a,-1/2+1/2\,\sqrt {5}),\\
S^{(n)}={\dfrac {{5}^{5/2\,n}  \left( \sqrt {5}
-1 \right) ^{5\,n} \left( 1/2\,a+2/5,n \right)  \left( 1/2\,a+3/5,n \right) }{{2}^{11\,n}
 \left( 1/2\,a+1/4,n \right)  \left( 1/2\,a+3/4,n
 \right) }},
 \tag{2,5,4-1}
\end{cases}\\
&\begin{cases}
(a,b,c,x)=(a,5/2\,a-1,2\,a,-1/2-1/2\,\sqrt {5}),\\
S^{(n)}={\dfrac {{5}^{5/2\,n} \left( \sqrt {5}
+1 \right) ^{5\,n} \left( 1/2\,a+2/5,n \right)  \left( 1/2\,a+3/5,n \right) }{{2}^{11\,n}
 \left( 1/2\,a+1/4,n \right)  \left( 1/2\,a+3/4,n
 \right) }},
 \tag{2,5,4-2}
\end{cases}\\
&\begin{cases}
(a,b,c,x)=(a,5/2\,a-1/2,2\,a,-1/2+1/2\,\sqrt {5}),\\
S^{(n)}={\dfrac {{5}^{5/2\,n}  \left( \sqrt {5}
-1 \right) ^{5\,n} \left( 1/2\,a+3/10,n \right)  \left( 1/2\,a+7/10,n \right) }{{2}^{11\,n}
 \left( 1/2\,a+1/4,n \right)  \left( 1/2\,a+3/4,n
 \right) }},
 \tag{2,5,4-3}
\end{cases}\\
&\begin{cases}
(a,b,c,x)=(a,5/2\,a-1/2,2\,a,-1/2-1/2\,\sqrt {5}),\\
S^{(n)}={\dfrac {{5}^{5/2\,n}  \left( \sqrt {5}
+1 \right) ^{5\,n} \left( 1/2\,a+3/10,n \right)  \left( 1/2\,a+7/10,n \right) }{{2}^{11\,n}
 \left( 1/2\,a+1/4,n \right)  \left( 1/2\,a+3/4,n
 \right) }}
 \tag{2,5,4-4}.
\end{cases}
\end{align}
\paragraph{(2,5,4-1)}
\begin{flalign*}
&{\rm (i)}\, {\mbox{$F$}(a,5/2\,a-1;\,2\,a;\,-1/2+1/2\,\sqrt {5})}&\\
&=\begin{cases}
{\dfrac {{5}^{1/2-5/4\,a} 
 \left( \sqrt {5}-1 \right) ^{1-5/2\,a}\Gamma\left(3/5\right)
 \Gamma\left(4/5\right)\Gamma\left(a+1/2\right)  }{{2}^{9/5-9/2
\,a}\Gamma\left(9/10\right)\Gamma\left(1/2\,a+2/5\right)
\Gamma\left(1/2\,a+3/5\right)}},\\
\dfrac{{5}^{5/2\,n+1/2} \left( \sqrt {5}-
1 \right) ^{5\,n+5} \left( 8/5,n \right) 
 \left( 7/5,n \right) }{ {2}^{11\,n+6}
   \left( 7/4,n \right)  \left( 5/4,n \right) }
\quad {\text {if $a=-2-2\,n$}},\\
\dfrac{-{5}^{5/2\,n+1/2}  \left( \sqrt {5}-
1 \right) ^{5\,n+4} \left( 11/10,n \right) 
 \left( 9/10,n \right) }{ {2}^{11\,n+6}
   \left( 5/4,n \right)  \left( 3/4,n \right) }
\quad {\text {if $a=-1-2\,n$}}.
\end{cases}
&
\end{flalign*}
\begin{flalign*}
&{\rm (ii)}\, 
{\mbox{$F$}(a,1-1/2\,a;\,2\,a;\,-1/2+1/2\,\sqrt {5})}&\\
&=\begin{cases}
{\dfrac {{5}^{1/2-5/4\,a} 
 \left( \sqrt {5}-1 \right) ^{1/2\,a-1}\Gamma\left(3/5\right)
 \Gamma\left(4/5\right)\Gamma\left(a+1/2\right)  }{{2}^{-3/2
\,a-1/5}\Gamma\left(9/10\right)\Gamma\left(1/2\,a+2/5\right)
\Gamma\left(1/2\,a+3/5\right)}},\\
\dfrac{{5}^{5/2\,n+1/2}  \left( \sqrt {5}+
1 \right) ^{n+1} \left( 8/5,n \right) 
 \left( 7/5,n \right) }{ {2}^{7\,n+2}
   \left( 7/4,n \right)  \left( 5/4,n \right) }
\quad {\text {if $a=-2-2\,n$}},\\
\dfrac{{5}^{5/2\,n+1/2}  \left( \sqrt {5}+
1 \right) ^{n+2} \left( 11/10,n \right) 
 \left( 9/10,n \right) }{ {2}^{7\,n+4}
   \left( 5/4,n \right)  \left( 3/4,n \right) }
\quad {\text {if $a=-1-2\,n$}}
\end{cases}
&
\end{flalign*}
(The first case is a generalization of Theorem 24 in [Ek]).
\begin{flalign*}
&{\rm (iii)}\,{\mbox{$F$}(a,1-1/2\,a;\,2\,a;\,-1/2-1/2\,\sqrt {5}})&\\
&=\begin{cases}
\dfrac{{5}^{5/2\,n+1/2}  \left( \sqrt {5}-
1 \right) ^{n+1} \left( 8/5,n \right) 
 \left( 7/5,n \right) }{ {2}^{7\,n+2}
   \left( 7/4,n \right)  \left( 5/4,n \right) }
\quad {\text {if $a=-2-2\,n$}},\\
\dfrac{-{5}^{5/2\,n+1/2}  \left( \sqrt {5}-
1 \right) ^{n+2} \left( 11/10,n \right) 
 \left( 9/10,n \right) }{ {2}^{7\,n+4}
   \left( 5/4,n \right)  \left( 3/4,n \right) }
\quad {\text {if $a=-1-2\,n$}},\\
\dfrac{ {2}^{5\,n}  \left( \sqrt {5}+
1 \right) ^{n} \left( 5/4,n \right) 
 \left( 7/4,n \right) }{{5}^{5/2\,n}
   \left( 7/5,n \right)  \left( 8/5,n \right) }
\quad {\text {if $a=2+2\,n$}}.
\end{cases}
&
\end{flalign*}
\begin{flalign*}
&{\rm (iv)}\,
{\mbox{$F$}(a,5/2\,a-1;\,2\,a;\,-1/2-1/2\,\sqrt {5})}&\\
&=\begin{cases}
\dfrac{{5}^{5/2\,n+1/2}  \left( \sqrt {5}+
1 \right) ^{5\,n+5} \left( 8/5,n \right) 
 \left( 7/5,n \right) }{ {2}^{11\,n+6}
   \left( 7/4,n \right)  \left( 5/4,n \right) }
\quad {\text {if $a=-2-2\,n$}},\\
\dfrac{{5}^{5/2\,n+1/2}  \left( \sqrt {5}+
1 \right) ^{5\,n+4} \left( 11/10,n \right) 
 \left( 9/10,n \right) }{ {2}^{11\,n+6}
   \left( 5/4,n \right)  \left( 3/4,n \right) }
\quad {\text {if $a=-1-2\,n$}},\\
\dfrac{{5}^{5/2\,n}  \left( \sqrt {5}+
1 \right) ^{5\,n} \left( 1/5,n \right) 
 \left( 2/5,n \right) }{ {2}^{11\,n}
   \left( 11/20,n \right)  \left( 1/20,n \right) }
\quad {\text {if $a=2/5-2\,n$}},\\
\dfrac{{5}^{5/2\,n+1}
 \left( \sqrt {5}+
1 \right) ^{5\,n+2} \left( 4/5,n \right) 
 \left( 3/5,n \right) }{ {2}^{11\,n+3}
   \left( 19/20,n \right)  \left( 9/20,n \right) }
\quad {\text {if $a=-2/5-2\,n$}},\\
0
\quad {\text {if $a=-4/5-2\,n$}},\\
0
\quad {\text {if $a=-6/5-2\,n$}}.
\end{cases}
&
\end{flalign*}
\begin{flalign*}
&{\rm (v)}\,
{\mbox{$F$}(a,5/2\,a-1;\,3/2\,a;\,3/2-1/2\,\sqrt {5})}&\\
&=\begin{cases}
{\dfrac {{3}^{3/2\,a-3/5} 
 \left( \sqrt {5}+1 \right) ^{5/2\,a-1}\Gamma  \left( 3/5 \right) 
\Gamma  \left( 4/5 \right) \Gamma  \left( 1/2\,a+1/3 \right) \Gamma 
 \left( 1/2\,a+2/3 \right) }{{2}^{5/2\,a-1}{5}^{5/4\,a-1/2}\Gamma 
 \left( 8/15 \right) \Gamma  \left( 13/15
 \right) \Gamma  \left( 1/2\,a+2/5 \right) \Gamma  \left( 1/2\,a+3/5
 \right)}},\\
\dfrac{{5}^{5/2\,n+1/2} \left( \sqrt {5}-
1 \right) ^{5\,n+5} \left( 8/5,n \right) 
 \left( 7/5,n \right) }{ {2}^{5\,n+5}3^{3\,n}
   \left( 5/3,n \right)  \left( 4/3,n \right) }
   \quad {\text {if $a=-2-2\,n$}}.
\end{cases}
 &
\end{flalign*}
\begin{flalign*}
&{\rm (vi)}\,
{\mbox{$F$}(1-a,1/2\,a;\,3/2\,a;\,3/2-1/2\,\sqrt {5})}&\\
&=\begin{cases}
{\dfrac {{3}^{3/2\,a-3/5}
 \left( \sqrt {5}+1 \right) ^{1/2\,a}\Gamma  \left( 3/5 \right) 
\Gamma  \left( 4/5 \right) \Gamma  \left( 1/2\,a+1/3 \right) \Gamma 
 \left( 1/2\,a+2/3 \right) }{{2}^{1/2\,a}{5}^{5/4\,a-1/2}\Gamma 
 \left( 8/15 \right) \Gamma  \left( 13/15
 \right) \Gamma  \left( 1/2\,a+2/5 \right) \Gamma  \left( 1/2\,a+3/5
 \right)}},\\
 \dfrac{{5}^{5/2\,n+1/2} \left( \sqrt {5}-1 \right) ^{n+1}  \left( 8/5,n \right) 
 \left( 7/5,n \right) }{ {2}^{n+1}3^{3\,n}
   \left( 5/3,n \right)  \left( 4/3,n \right) }
   \quad {\text {if $a=-2-2\,n$}}
\end{cases}
 &
\end{flalign*}
(The second case is identical to Theorem 21 in [Ek]).
\begin{flalign*}
&{\rm (vii)}\,
{\mbox{$F$}(a,1-a;\,3/2\,a;\,1/2-1/2\,\sqrt {5})}&\\
&=\begin{cases}
{\dfrac {{3}^{3/2\,a-3/5}
 \left( \sqrt {5}+1 \right) ^{3/2\,a-1}\Gamma  \left( 3/5 \right) 
\Gamma  \left( 4/5 \right) \Gamma  \left( 1/2\,a+1/3 \right) \Gamma 
 \left( 1/2\,a+2/3 \right) }{{2}^{3/2\,a-1}{5}^{5/4\,a-1/2}\Gamma 
 \left( 8/15 \right) \Gamma  \left( 13/15
 \right) \Gamma  \left( 1/2\,a+2/5 \right) \Gamma  \left( 1/2\,a+3/5
 \right)}},\\
 \dfrac{{5}^{5/2\,n+1/2} \left( \sqrt {5}-1 \right) ^{3\,n+3}  \left( 8/5,n \right) 
 \left( 7/5,n \right) }{ 2^{3\,n+3}3^{3\,n}
   \left( 5/3,n \right)  \left( 4/3,n \right) }
   \quad {\text {if $a=-2-2\,n$}}.
\end{cases}
 &
\end{flalign*}
\begin{flalign*}
&{\rm (viii)}\,
{\mbox{$F$}(1/2\,a,5/2\,a-1;\,3/2\,a;\,1/2-1/2\,\sqrt {5})}&\\
&=\begin{cases}
{\dfrac {{3}^{3/2\,a-3/5}
\Gamma  \left( 3/5 \right) 
\Gamma  \left( 4/5 \right) \Gamma  \left( 1/2\,a+1/3 \right) \Gamma 
 \left( 1/2\,a+2/3 \right) }{{5}^{5/4\,a-1/2}\Gamma 
 \left( 8/15 \right) \Gamma  \left( 13/15
 \right) \Gamma  \left( 1/2\,a+2/5 \right) \Gamma  \left( 1/2\,a+3/5
 \right)}},\\
 \dfrac{{5}^{5/2\,n+1/2} \left( 8/5,n \right) 
 \left( 7/5,n \right) }{3^{3\,n}
   \left( 5/3,n \right)  \left( 4/3,n \right) }
   \quad {\text {if $a=-2-2\,n$}}.
\end{cases}
 &
\end{flalign*}
\begin{flalign*}
&{\rm (ix)}\,
{\mbox{$F$}(a,1-a;\,2-3/2\,a;\,1/2+1/2\,\sqrt {5})}&\\
&=\begin{cases}
{\dfrac {{3}^{3\,n} \left( \sqrt {5}-1 \right) ^{3\,n}
 \left( 4/3,n \right)  \left( 2/3,n \right) }{2^{3\,n}{5}^{5/2
\,n} \left( 6/5,n \right)  \left( 4/5,
n \right) }}
\quad {\text {if $a=-2\,n$}},\\
{\dfrac {{3}^{3\,n} \left( \sqrt {5}-1 \right) ^{3\,n+3}
 \left( 11/6,n \right)  \left( 7/6,n \right) }{7\cdot 2^{3\,n+3}{5}^{5/2
\,n-1/2} \left( 17/10,n \right)  \left( 13/10,
n \right) }}
\quad {\text {if $a=-1-2\,n$}},\\
\dfrac{ {5}^{5/2\,n} \left( \sqrt {5}+1 \right) ^{3\,n} \left( 3/
10,n \right) \left( 7/10,n \right) }{2^{3\,n} {3}^{3\,n} \left( 1/6
,n \right)  \left( 5/6,n
 \right)  }
\quad {\text {if $a=1+2\,n$}},\\
\dfrac{{5}^{5/2\,n+1/2} \left( \sqrt {5}+1 \right) ^{3\,n+6} \left( 9/
5,n \right) \left( 11/5,n \right) }{ 2^{3\,n+6}{3}^{3\,n}  \left( 5/3
,n \right)  \left( 7/3,n
 \right)  }
\quad {\text {if $a=4+2\,n$}}
\end{cases}
 &
\end{flalign*}
(The third case is identical to Theorem 34 in [Ek]).
\begin{flalign*}
&{\rm (x)}\,{\mbox{$F$}(2-5/2\,a,1-1/2\,a;\,2-3/2\,a;\,1/2+1/2\,\sqrt {5})}&\\
&=\begin{cases}
\dfrac{{5}^{5/2\,n}  \left( 1/5,n \right)
 \left( 3/5,n \right) }{ {\left(-3\right)}^{3\,n}    \left( 1/15,n
 \right)   \left( 11/15
,n \right)  }
 \quad {\text {if $a=4/5+2\,n$}},\\
 \dfrac{-{5}^{5/2\,n+1/2}  \left( 2/5,n \right)
 \left( 4/5,n \right) }{ {\left(-3\right)}^{3\,n}   \left( 4/15,n
 \right)   \left( 14/15
,n \right)  }
 \quad {\text {if $a=6/5+2\,n$}},\\
 0
  \quad {\text {if $a=8/5+2\,n$}},\\
{\dfrac {{5}^{5/2\,n}  \left( 4/5,n \right)  \left( 6/5,n \right) }{ \left( -3 \right) ^{3\,n
} \left( 2/3,n \right) \left( 4/3,n
 \right) }}
\quad {\text {if $a=2+2\,n$}},\\
0
\quad {\text {if $a=12/5+2\,n$}}.
\end{cases}
 &
\end{flalign*}
\begin{flalign*}
&{\rm (xi)}\,
{\mbox{$F$}(a,1-1/2\,a;\,2-3/2\,a;\,3/2+1/2\,\sqrt {5})}&\\
&=\begin{cases}
{\dfrac {{3}^{3\,n}  \left( \sqrt {5}-1
 \right) ^{n} \left( 4/3,n \right) 
 \left( 2/3,n \right) }{{2}^{n}{5}^{5/2\,n} \left( 
6/5,n \right)  \left( 4/5,n \right) }}
\quad {\text {if $a=-2\,n$}},\\
\dfrac{-{3}^{3\,n}  \left( \sqrt {5}-1
 \right) ^{n+2} \left( 11/6,n \right)  \left( 7/6,n \right)}{7\cdot  {2}^{n+2}  {5}^{5/2\,n-1/2}
 \left( 17/10,n \right)   \left( 13/10,n \right)  }
\quad {\text {if $a=-1-2\,n$}},\\
{\dfrac {{5}^{5/2\,n} \left( \sqrt {
5}+1 \right) ^{n} \left( 4/5,n \right)  \left( 6/5,n \right) }{{2}^{n}{3}^{3\,n} \left( 2/3,n \right)  \left( 4/3,n \right) }}
\quad {\text {if $a=2+2\,n$}}\\
\end{cases}
 &
\end{flalign*}
(The third case is identical to Theorem 22 in [Ek]).
\begin{flalign*}
&{\rm (xii)}\,
{\mbox{$F$}(1-a,2-5/2\,a;\,2-3/2\,a;\,3/2+1/2\,\sqrt {5})}&\\
&=\begin{cases}
\dfrac{{5}^{5/2\,n} \left( \sqrt {5}+1 \right) ^{5\,n}  \left( 3/10,n \right)  \left( 7/10,n \right)}{ {2}^{5\,n} {3}^{3\,n} \left( 1/6,n
 \right) \left( 5/6,n \right) 
 }
\quad {\text {if $a=1+2\,n$}},\\
\dfrac{-{5}^{5/2\,n+1/2}  \left( \sqrt {5
}+1 \right) ^{5\,n+9} \left( 9/5,n \right)  \left( 11/5,n \right)}{{2}^{5\,n+9}
  {3}^{3\,n} 
 \left( 5/3,n \right)   \left( 7/
3,n \right)  }
\quad {\text {if $a=4+2\,n$}},\\
\dfrac{{5}^{5/2\,n} \left( \sqrt {5}+1 \right) ^{5\,n}  \left( 1/5,n \right) 
 \left( 3/5,n \right) }{  {2}^{5\,n} {3}^{3\,
n}   \left( 1/15,n \right) 
  \left(  11/15,n
 \right)  }
\quad {\text {if $a=4/5+2\,n$}},\\
\dfrac{{5}^{5/2\,n+1/2} \left( \sqrt {5}+1 \right) ^{5\,n+1} \left( 2/5,n \right)  \left( 4/5,n \right) }{ {2}^{5\,n+1}  {3}^{3\,n}  \left( 4/15,n \right)  
 \left( 14/15,n \right) }
\quad {\text {if $a=6/5+2\,n$}},\\
0
\quad {\text {if $a=8/5+2\,n$}},\\
0
\quad {\text {if $a=12/5+2\,n$}}.\\
\end{cases}
 &
\end{flalign*}
\begin{flalign*}
&{\rm (xiii)}\,
{\mbox{$F$}(1/2\,a,5/2\,a-1;\,3/2\,a;\,1/2+1/2\,\sqrt {5})}&\\
&=\begin{cases}
{\dfrac {-{5}^{5/2\,n+1/2} \left( 8/5,n \right)  \left( 7/5,n \right) }{ \left( -3 \right) ^{3\,n} \left( 5/3,n \right) \left( 4/3,n
 \right) }}
\quad {\text {if $a=-2-2\,n$}},\\
\dfrac{{5}^{5/2\,n} \left( 2/5,n \right) 
 \left( 1/5,n \right) }{ \left( -3 \right) ^{3\,n}    \left( 7/15,n \right)  
 \left( 2/15,n \right) }
\quad {\text {if $a=2/5-2\,n$}},\\
\dfrac{-{5}^{5/2\,n+1} \left( 4/5,n \right) 
 \left( 3/5,n \right) }{ \left( -3 \right) ^{3\,n+1}  \left( 13/15,n \right)  
 \left( 8/15,n \right) }
\quad {\text {if $a=-2/5-2\,n$}},\\
0
\quad {\text {if $a=-4/5-2\,n$}},\\
0
\quad {\text {if $a=-6/5-2\,n$}}.
\end{cases}
 &
\end{flalign*}
\begin{flalign*}
&{\rm (xiv)}\,
{\mbox{$F$}(a,1-a;\,3/2\,a;\,1/2+1/2\,\sqrt {5})}&\\
&=\begin{cases}
{\dfrac {{5}^{5/2\,n+1/2} \left( \sqrt {5}+1 \right) ^{3\,n+3}  \left( 8/5,n \right)\left( 7/5,n \right) }{{2}^{3\,n+3}{3}^{3\,n} \left( 5/3,n \right)  \left( 4/3,n \right) }}
\quad {\text {if $a=-2-2\,n$}},\\
\dfrac{{5}^{5/2\,n+1/2} \left( \sqrt {5}+1 \right) ^{3\,n+3}  \left( 11/10,n \right)  \left( 9/10,n \right) }{{2}^{3\,n+3}
  {3}^{3\,n+1}  \left( 7/6,n \right)   \left( 5/6,n \right) }
\quad {\text {if $a=-1-2\,n$}},\\
\dfrac{{3}^{3\,n} \left( \sqrt {5}-1 \right) ^{3\,n} \left( 5/6,n \right) 
 \left( 7/6,n \right)}{{2}^{3\,n}  {5}^{5/2
\,n} \left( 9/10,n
 \right)   \left( 11/10
,n \right)  }
\quad {\text {if $a=1+2\,n$}},\\
{\dfrac {-{3}^{3\,n-1} \left( \sqrt {5}-1 \right) ^{3\,n+3}  \left( 4/3,n \right)  \left( 5/3,n \right) }{{2}^{3\,n+3}{5}^{5/2\,n} \left( 7/5,n \right)  \left( 8/5,n
 \right) }}
\quad {\text {if $a=2+2\,n$}}\\
\end{cases}
 &
\end{flalign*}
(The first case is identical to Theorem 35 in [Ek]).
\begin{flalign*}
&{\rm (xv)}\,
{\mbox{$F$}(a,5/2\,a-1;\,3/2\,a;\,3/2+1/2\,\sqrt {5})}&\\
&=\begin{cases}
{\dfrac {{5}^{5/2\,n+1/2} \left( \sqrt {5}+1 \right) ^{5\,n+5}  \left( 8/5,n \right) \left( 7/5,n \right) }{{2}^{5\,n+5}{3}^{3\,n} \left( 5/3,n \right)  \left( 4/3,n
 \right) }}
\quad {\text {if $a=-2-2\,n$}},\\
\dfrac{-{5}^{5/2\,n+1/2} \left( \sqrt {5}+1 \right) ^{5\,n+4} \left( 11/10,n
 \right)  \left( 9/10,n \right) }{ {2}
^{5\,n+4} {3}^{3\,n+1}   \left( 7/6,n \right)  \left( 5/6,n \right) } 
\quad {\text {if $a=-1-2\,n$}},\\
\dfrac{{5}^{5/2\,n} \left( \sqrt {5}+1 \right) ^{5\,n}  \left( 2/5,n \right)
 \left( 1/5,n \right) }{ {2}^{5\,n} {3}^{3\,
n}  \left( 7/15,n
 \right)   \left( 2/15,n \right) }
\quad {\text {if $a=2/5-2\,n$}},\\
\dfrac{{5}^{5/2\,n+1} \left( \sqrt {5}+1 \right) ^{5\,n+2}  \left( 4/5,n \right)  \left( 3/5,n \right)}{ {2}^{5\,n+2}  {3}^{3\,n+1}  \left( 13/15,n \right) 
 \left( 8/15,n \right) }
\quad {\text {if $a=-2/5-2\,n$}},\\
0
\quad {\text {if $a=-4/5-2\,n$}},\\
0
\quad {\text {if $a=-6/5-2\,n$}}.\\
\end{cases}
&
\end{flalign*}
\begin{flalign*}
&{\rm (xvi)}\,
{\mbox{$F$}(1/2\,a,1-a;\,3/2\,a;\,3/2+1/2\,\sqrt {5})}&\\
&=\begin{cases}
{\dfrac {{5}^{5/2\,n+1/2} \left( \sqrt {5}+1 \right) ^{n+1} \left( 8/5,n \right)  \left( 7/5,n
 \right) }{{2}^{n+1}{3}^{3\,n} \left( 5/3,n \right)  \left( 4/3,n \right) }}
\quad {\text {if $a=-2-2\,n$}},\\
\dfrac{{3}^{3\,n} \left( \sqrt {5}-1 \right) ^{n} \left( 5/6,
n \right)  \left( 7/6,n \right) }{ {2}^{n}
  {5}^{5/2\,n}   \left( 9/10,n \right)   \left( 11/10,n \right)  }
\quad {\text {if $a=1+2\,n$}},\\
{\dfrac {{3}^{3\,n-1}  \left( \sqrt {5}-1
 \right) ^{n+2} \left( 4/3,n \right) 
 \left( 5/3,n \right) }{{2}^{n+2}{5}^{5/2\,n} \left( 7
/5,n \right)  \left( 8/5,n \right) }}
\quad {\text {if $a=2+2\,n$}}\\
\end{cases}
&
\end{flalign*}
(The first case is identical to Theorem 21 in [Ek]).
\begin{flalign*}
&{\rm (xvii)}\,
{\mbox{$F$}(1/2\,a,1-a;\,2-2\,a;\,-1/2+1/2\,\sqrt {5})}&\\
&=\begin{cases}
{\dfrac {{2}^{-3/2\,a} \left( \sqrt {5}-1 \right) ^{-1/2\,a} \Gamma\left(6/5\right)\Gamma\left(4/5\right) \Gamma  \left( 
3/2-a \right)  }{{5}^{-
5/4\,a}\Gamma\left(3/2\right)
\Gamma  \left( 6/5-1/2\,a \right) \Gamma  \left( 4/5-1/2\,a
 \right) }},\\
\dfrac{7\cdot {5}^{5/2\,n+1/2} \left( \sqrt {5}+1 \right) ^{n+1}  \left( 13/10,n \right)  \left( 17/10,n \right) }{{2}^{7\,n+5}
 \left( 5/4,n \right)  \left( 7/4,n \right) }
\quad {\text {if $a=3+2\,n$}},\\
{\dfrac {{5}^{5/2\,n} \left( \sqrt {5}+1 \right) ^{n+2}  \left( 4/5,n \right)  \left( 6/5,n \right) }{{2}^{7\,n+3}
 \left( 3/4,n \right)  \left( 5/4,n \right) }}
\quad {\text {if $a=2+2\,n$}}\\
\end{cases}
&
\end{flalign*}
(The first case is a generalization of Theorem 23 in [Ek]).
\begin{flalign*}
&{\rm (xviii)}\,
{\mbox{$F$}(1-a,2-5/2\,a;\,2-2\,a;\,-1/2+1/2\,\sqrt {5})}&\\
&=\begin{cases}
{\dfrac {{2}^{2-9/2\,a} \left( \sqrt {5}-1 \right) ^{5/2\,a-2} \Gamma\left(6/5\right)\Gamma\left(4/5\right) \Gamma  \left( 
3/2-a \right)  }{{5}^{-
5/4\,a}\Gamma(3/2)
\Gamma  \left( 6/5-1/2\,a \right) \Gamma  \left( 4/5-1/2\,a
 \right) }},\\
\dfrac{7\cdot {5}^{5/2\,n+1/2} \left( \sqrt {5}-1 \right) ^{5\,n+5}  \left( 13/10,n
 \right)  \left( 17/10,n \right)}{ {2
}^{11\,n+9}   \left( 5/4,n \right) 
  \left( 7/4,n \right) }
\quad {\text {if $a=3+2\,n$}},\\
{\dfrac {{5}^{5/2\,n} \left( \sqrt {5}-1 \right) ^{5\,n+4} \left( 4/5,n \right)  \left( 6/5,n \right) }{{2}^{11\,n+5}\left( 3
/4,n \right) \left( 5/4,n \right) }}
\quad {\text {if $a=2+2\,n$}},\\
\end{cases}
&
\end{flalign*}
\begin{flalign*}
&{\rm (xix)}\,
{\mbox{$F$}(1-a,2-5/2\,a;\,2-2\,a;\,-1/2-1/2\,\sqrt {5})}&\\
&=\begin{cases}
\dfrac{7\cdot {5}^{5/2\,n+1/2} \left( \sqrt {5}+1 \right) ^{5\,n+5}  \left( 13/10
,n \right)  \left( 17/10,n \right)}{ 
{2}^{11\,n+9} \left( 5/4,n
 \right)   \left( 7/4,n \right) }
\quad {\text {if $a=3+2\,n$}},\\
{\dfrac {{5}^{5/2\,n} \left( \sqrt {5}+1 \right) ^{5\,n+4} \left( 4/5,n \right)  \left( 6/5,n \right) }{{2}^{11\,n+5}
 \left( 3/4,n \right)  \left( 5/4,n \right) }}
\quad {\text {if $a=2+2\,n$}},\\
\dfrac{{5}^{5/2\,n} \left( \sqrt {5}+1 \right) ^{5\,n}  \left( 1/5,n \right) 
 \left( 3/5,n \right) }{  {2}^{11\,n}  \left( 3/20,n \right)   \left(13/20,n \right)}
\quad {\text {if $a=4/5+2\,n$}},\\
\dfrac{{5}^{5/2\,n+1/2} \left( \sqrt {5}+1 \right) ^{5\,n+1} \left( 2/5,n \right)  \left( 4/5,n \right) }{ {2}^{11\,n+2}  \left( 7/20,n \right)    \left( 17/20,n \right)  }
 \quad {\text {if $a=6/5+2\,n$}},\\
 0
 \quad {\text {if $a=8/5+2\,n$}},\\
 0
 \quad {\text {if $a=12/5+2\,n$}}.
\end{cases}
&
\end{flalign*}
\begin{flalign*}
&{\rm (xx)}\,{\mbox{$F$}(1/2\,a,1-a;\,2-2\,a;\,-1/2-1/2\,\sqrt {5})}&\\
&=\begin{cases}
{\dfrac {{2}^{5\,n} \left( \sqrt {5}+1 \right) ^{n}
 \left( 5/4,n \right) \left( 3/4,n \right) }{{5}^{5/2
\,n} \left( 6/5,n \right)  \left( 4/5,
n \right) }}
 \quad {\text {if $a=-2\,n$}},\\
\dfrac{7\cdot {5}^{5/2\,n+1/2} \left( \sqrt {5}-1 \right) ^{n+1}
 \left( 13/10,n \right) \left( 17
/10,n \right) }{{2}^{7\,n+5}   \left( 5/4,n \right)   \left( 7/4,n \right)  }
\quad {\text {if $a=3+2\,n$}},\\
{\dfrac {{5}^{5/2\,n}   \left( \sqrt {5}-1
 \right) ^{n+2} \left( 4/5,n \right) 
 \left( 6/5,n \right) }{{2}^{7\,n+3} \left( 3/4,n
 \right)  \left( 5/4,n \right) }}
\quad {\text {if $a=2+2\,n$}}.
\end{cases}
&
\end{flalign*}
\begin{flalign*}
&{\rm (xxi)}\,
{\mbox{$F$}(a,1-1/2\,a;\,2-3/2\,a;\,3/2-1/2\,\sqrt {5})}&\\
&=\begin{cases}
{\dfrac {{3}^{-3/2\,a} \left( \sqrt {5}+1 \right) ^{-1/2\,a} \Gamma\left(6/5\right)\Gamma\left(4/5\right) \Gamma 
 \left( 4/3-1/2\,a \right) \Gamma  \left( 2/3-1/2\,a \right)
}{{2}^{-1/2\,a}{5}^{-5/4\,a}\Gamma \left(4/3\right) \Gamma\left(2/3\right)
\Gamma  \left( 6/5-1/2\,a \right) \Gamma 
 \left( 4/5-1/2\,a \right)}},\\
 {\dfrac {{5}^{5/2\,n} \left( \sqrt {5}-1 \right) ^{n}\left( 4/5,n \right)  \left( 6/5,n \right) }{{2}^{n}{3}^{3\,n} \left( 2/3,n \right)  \left( 4/3,n \right) }}
\quad {\text {if $a=2+2\,n$}}
\end{cases}
&
\end{flalign*}
(The second case is identical to Theorem 22 in [Ek]).
\begin{flalign*}
&{\rm (xxii)}\,
{\mbox{$F$}(1-a,2-5/2\,a;\,2-3/2\,a;\,3/2-1/2\,\sqrt {5})}&\\
&=\begin{cases}
{\dfrac {{3}^{-3/2\,a} \left( \sqrt {5}+1 \right) ^{1-5/2\,a} \Gamma\left(6/5\right)\Gamma\left(4/5\right) \Gamma 
 \left( 4/3-1/2\,a \right) \Gamma  \left( 2/3-1/2\,a \right)
}{{2}^{1-5/2\,a}{5}^{-5/4\,a}\Gamma \left(4/3\right) \Gamma\left(2/3\right)
\Gamma  \left( 6/5-1/2\,a \right) \Gamma 
 \left( 4/5-1/2\,a \right)}},\\
\dfrac{-{5}^{5/2\,n+1/2} \left( \sqrt {5}-1 \right) ^{5\,n+9} \left( 9/5,n \right) 
 \left( 11/5,n \right)  }{ {2}^{5\,n+9} {3}^{3\,n} \left( 5/3,n
 \right)   \left( 7/3,n \right) 
  }
 \quad {\text {if $a=4+2\,n$}}.
\end{cases}
&
\end{flalign*}
\begin{flalign*}
&{\rm (xxiii)}\,
{\mbox{$F$}(a,1-a;\,2-3/2\,a;\,1/2-1/2\,\sqrt {5})}&\\
&=\begin{cases}
{\dfrac {{3}^{-3/2\,a} \left( \sqrt {5}+1 \right) ^{-3/2\,a} \Gamma\left(6/5\right)\Gamma\left(4/5\right) \Gamma 
 \left( 4/3-1/2\,a \right) \Gamma  \left( 2/3-1/2\,a \right)
}{{2}^{-3/2\,a}{5}^{-5/4\,a}\Gamma \left(4/3\right) \Gamma\left(2/3\right)
\Gamma  \left( 6/5-1/2\,a \right) \Gamma 
 \left( 4/5-1/2\,a \right)}},\\
 \dfrac{-{5}^{5/2\,n+1/2} \left( \sqrt {5}-1 \right) ^{3\,n+6} \left( 9/5,n \right)  \left( 11/5,n \right)}{{2}^{3\,n+6}
  {3}^{3\,n} 
 \left( 5/3,n \right)  \left( 7/
3,n \right)  }
 \quad {\text {if $a=4+2\,n$}}.
 \end{cases}
&
\end{flalign*}
\begin{flalign*}
&{\rm (xxiv)}\,
{\mbox{$F$}(2-5/2\,a,1-1/2\,a;\,2-3/2\,a;\,1/2-1/2\,\sqrt {5})}&\\
&=\begin{cases}
{\dfrac {{3}^{-3/2\,a} \left( \sqrt {5}-1 \right)  \Gamma\left(6/5\right)\Gamma\left(4/5\right) \Gamma 
 \left( 4/3-1/2\,a \right) \Gamma  \left( 2/3-1/2\,a \right)
}{{2}\cdot {5}^{-5/4\,a}\Gamma \left(4/3\right) \Gamma\left(2/3\right)
\Gamma  \left( 6/5-1/2\,a \right) \Gamma 
 \left( 4/5-1/2\,a \right)}},\\
{\dfrac {{5}^{5/2\,n} \left( 4/5,n \right)  \left( 6/5,n \right) }{{3}^{3\,n} \left( 2/3,n \right)  \left( 4/3,n \right) }}
\quad {\text {if $a=2+2\,n$}}.
  \end{cases}
&
\end{flalign*}
\paragraph{(2,5,4-2), (2,5,4-3), (2,5,4-4)}
The special values obtained from (2,5,4-2), (2,5,4-3) and (2,5,4-4)
coincide with those obtained from (2,5,4-1).

\subsubsection{$(k,l,m)=(2,6,4)$}
In this case, we have
\begin{align}
&(a,b,c,x)=(a,b,2\,a,2),\, 
\tag{2,6,4-1}\\
&(a,b,c,x)=(a,b,b+1-a,-1),\,
\tag{2,6,4-2}\\
&(a,b,c,x)=(a,b,1/2\,a+1/2\,b+1/2,1/2),\,
\tag{2,6,4-3}\\
&(a,b,c,x)=(a,3\,a-1,2\,a,1/2+1/2\,i\sqrt {3}),
\tag{2,6,4-4}\\
&(a,b,c,x)=(a,3\,a-1,2\,a,1/2-1/2\,i\sqrt {3})
\tag{2,6,4-5}.
\end{align}
\paragraph{(2,6,4-1), (2,6,4-2), (2,6,4-3)}
The special values obtained from (2,6,4-1), (2,6,4-2) and (2,6,4-3)
are evaluated in paragraphs (1,2,2-1) and (0,2,2-1).
\paragraph{(2,6,4-4), (2,6,4-5)}
The special values obtained from (2,6,4-4) and (2,6,4-5) coincide with those obtained from (1,3,2-1).

\subsection{$m=5$}
\subsubsection{$(k,l,m)=(0,5,5)$}
In this case, we have
\begin{align}
&(a,b,c,x)=(1,b,b,\lambda),\, S^{(n)}=1, \tag{0,5,5-1}\\
&(a,b,c,x)=(0,b,b+1,\lambda),\, S^{(n)}=1, \tag{0,5,5-2}
\end{align}
where $\lambda$ is a solution of $x^4+x^3+x^2+x+1=0$.
\paragraph{(0,5,5-1)}
The special values obtained from (0,5,5-1) except trivial values 
are the special cases of (\ref{F(1,b;2;x)}).

\paragraph{(0,5,5-2)}
The special values obtained from (0,5,5-2) coincide with
those obtained from (0,5,5-1).

\subsubsection{$(k,l,m)=(1,4,5)$}
In this case, there is no admissible quadruple.

\subsubsection{$(k,l,m)=(1,5,5)$}
In this case, we have
\begin{align}
&
(a,b,c,x)=(a,5\,a+1,5\,a,5/4),\, 
S^{(n)}={\frac { \left( -5 \right) ^{5\,n} \left( a+1,n
 \right)  \left( 4\,a,4\,n \right) }{{2}^{10\,n} \left( 5\,a,5\,n \right) }},
\tag{1,5,5-1}
\\
&
(a,b,c,x)=(a,5\,a-5,5\,a-3,5/4),\, 
S^{(n)}={\frac {5\,a-4}{ \left( -4 \right) ^{n} \left( 5\,a-4+5\,n \right) }}
\tag{1,5,5-2}.
\end{align}
\paragraph{(1,5,5-1)}
\begin{flalign*}
&{\rm (i)}\,{\mbox{$F$}(a,5\,a+1;\,5\,a;\,5/4)}
=\begin{cases}
0
\quad {\text {if $a=-1-n$}},\\
{\dfrac { \left( -5 \right) ^{5\,n} \left( 1/5,n
 \right) \left( 9/5,4\,n \right) }{{2}^{10\,n} \left( 2,5\,n \right) }}
\quad {\text {if $a=-1/5-n$}},\\
{\dfrac { 3\left( -5 \right) ^{5\,n} \left( 2/5,n
 \right) \left( 13/5,4\,n \right) }{{2}^{10\,n+2} \left( 3,5\,n \right) }}
\quad {\text {if $a=-2/5-n$}},\\
{\dfrac { 7\left( -5 \right) ^{5\,n} \left( 3/5,n
 \right) \left( 17/5,4\,n \right) }{{2}^{10\,n+4} \left( 4,5\,n \right) }}
\quad {\text {if $a=-3/5-n$}},\\
{\dfrac { 11\left( -5 \right) ^{5\,n} \left( 4/5,n
 \right) \left( 21/5,4\,n \right) }{{2}^{10\,n+6} \left( 5,5\,n \right) }}
\quad {\text {if $a=-4/5-n$}}.\\
\end{cases}
&
\end{flalign*}
The special values obtained from {\rm {(ii) and (iii)}} are trivial.
\begin{flalign*}
&{\rm (iv)}\,{\mbox{$F$}(4\,a,5\,a+1;\,5\,a;\,5)}
=\begin{cases}
0
\quad {\text {if $a=-1/4-1/4\,n$}},\\
{\dfrac {{5}^{5\,n}\left( 1/5,n \right)  \left( 9/5,4\,n \right) }{ \left( 2,5\,n\right) }}
\quad {\text {if $a=-1/5-n$}},\\
{\dfrac {-3\cdot {5}^{5\,n}\left( 2/5,n \right)  \left( 13/5,4\,n \right) }{ \left( 3,5\,n\right) }}
\quad {\text {if $a=-2/5-n$}},\\
{\dfrac {7\cdot {5}^{5\,n}\left( 3/5,n \right)  \left( 17/5,4\,n \right) }{ \left( 4,5\,n\right) }}
\quad {\text {if $a=-3/5-n$}},\\
{\dfrac {-11\cdot {5}^{5\,n}\left( 4/5,n \right)  \left( 21/5,4\,n \right) }{ \left( 5,5\,n\right) }}
\quad {\text {if $a=-4/5-n$}}.
\end{cases}
&
\end{flalign*}
\begin{flalign*}
&{\rm (v)}\,{\mbox{$F$}(a,5\,a+1;\,a+2;\,-1/4)}
={2}^{10\,a}{5}^{-5\,a} \left( a+1 \right) &
\end{flalign*}
(The above is a special case of (\ref{1.6eb2})).
\begin{flalign*}
&{\rm (vi)}\,
{\mbox{$F$}(2,1-4\,a;\,a+2;\,-1/4)}
=4/5\,a+4/5
&
\end{flalign*}
(The above is a special case of (\ref{1.5eb2})).
\begin{flalign*}
&{\rm (vii)}\,
{\mbox{$F$}(a,1-4\,a;\,a+2;\,1/5)}
={2}^{8\,a}{5}^{-4\,a} \left( a+1 \right) 
&
\end{flalign*}
(The above is a special case of (\ref{1.6eb2})).
\begin{flalign*}
&{\rm (viii)}\,
{\mbox{$F$}(2,5\,a+1;\,a+2;\,1/5)}
=5/4\,a+5/4
&
\end{flalign*}
(The above is a special case of (\ref{1.5eb2})).
\begin{flalign*}
&{\rm (ix)}\,
{\mbox{$F$}(a,1-4\,a;\,-4\,a;\,4/5)}
=\begin{cases}
0,\\
{\dfrac {{2}^{8\,n}\left( 9/4,5\,n \right) }{{5}^{4\,n
}\left( 5/4,n \right)\left( 2,4\,n\right) }}
\quad {\text {if $a=1/4+n$}},\\
{\dfrac {3\cdot {2}^{8\,n+1}\left( 7/2,5\,n \right) }{{5}^{4\,n+1
}\left( 3/2,n \right)\left( 3,4\,n\right) }}
\quad {\text {if $a=1/2+n$}},\\
{\dfrac {77\cdot {2}^{8\,n-1}\left( 19/4,5\,n \right) }{{5}^{4\,n+2
}\left( 7/4,n \right)\left( 4,4\,n\right) }}
\quad {\text {if $a=3/4+n$}},\\
{\dfrac {{2}^{8\,n+8}\left( 6,5\,n \right) }{{5}^{4\,n+3
}\left( 2,n \right)\left( 5,4\,n\right) }}
\quad {\text {if $a=1+n$}}.\\
\end{cases}&
\end{flalign*}
The special values obtained from {\rm {(x) and (xi)}} are trivial.
\begin{flalign*}
&{\rm (xii)}\,{\mbox{$F$}(-5\,a,1-4\,a;\,-4\,a;\,-4)}
=\begin{cases}
0
\quad {\text {if $a=1/5+n$}},\\
0
\quad {\text {if $a=2/5+n$}},\\
0
\quad {\text {if $a=3/5+n$}},\\
0
\quad {\text {if $a=4/5+n$}},\\
{\dfrac {{2}^{8\,n} \left( 9/4,5\,n \right) }{ \left( 5/4,n \right)  \left( 2,4\,n \right) }}
\quad {\text {if $a=1/4+n$}},\\
{\dfrac {3\cdot {2}^{8\,n+1} \left( 7/2,5\,n \right) }{ \left( 3/2,n \right)  \left( 3,4\,n \right) }}
\quad {\text {if $a=1/2+n$}},\\
{\dfrac {77\cdot {2}^{8\,n-1} \left( 19/4,5\,n \right) }{ \left( 7/4,n \right)  \left( 4,4\,n \right) }}
\quad {\text {if $a=3/4+n$}},\\
{\dfrac { {2}^{8\,n+8} \left( 6,5\,n \right) }{ \left( 2,n \right)  \left( 5,4\,n \right) }}
\quad {\text {if $a=1+n$}}.\\
\end{cases}&
\end{flalign*}
\begin{flalign*}
&{\rm (xiii)}\,
{\mbox{$F$}(2,5\,a+1;\,4\,a+2;\,4/5)}=5\left(4\,a+1\right)
&
\end{flalign*}
(The above is a special case of (\ref{1.5eb2})).
\begin{flalign*}
&{\rm (xiv)}\,
{\mbox{$F$}(4\,a,1-a;\,4\,a+2;\,4/5)}={5}^{-a} \left( 4\,a+1 \right)
&
\end{flalign*}
(The above is a special case of (\ref{1.6eb2})).
\begin{flalign*}
&{\rm (xv)}\,
{\mbox{$F$}(4\,a,5\,a+1;\,4\,a+2;\,-4)}
=\begin{cases}
-{5}^{5\,n} \left( 20\,n -1\right)
\quad {\text {if $a=-1/5-n$}},\\
-{5}^{5\,n+1} \left( 20\,n+3 \right) 
\quad {\text {if $a=-2/5-n$}},\\
-{5}^{5\,n+2} \left( 20\,n +7\right) 
\quad {\text {if $a=-3/5-n$}},\\
-{5}^{5\,n+3} \left( 20\,n +11\right) 
\quad {\text {if $a=-4/5-n$}}\\
\end{cases}&
\end{flalign*}
(The above are special cases of (\ref{1.6eb2})).
\begin{flalign*}
&{\rm (xvi)}\,
{\mbox{$F$}(2,1-a;\,4\,a+2;\,-4)}
=4/5\,n+1
\quad {\text {if $a=1+n$}}
&
\end{flalign*}
(The above is a special case of (\ref{1.5eb2})).
\begin{flalign*}
&{\rm (xvii)}\,{\mbox{$F$}(2,1-4\,a;\,2-5\,a;\,5/4)}
=
5\,n+1
\quad {\text {if $a=1/4+1/4\,n$}}
&
\end{flalign*}
(The above is a special case of (\ref{1.5eb2})).
\begin{flalign*}
&{\rm (xviii)}\,
{\mbox{$F$}(-5\,a,1-a;\,2-5\,a;\,5/4)}
=- \left( -4 \right) ^{-n-1} \left( 5\,n +4\right) 
\quad {\text {if $a=1+n$}}
&
\end{flalign*}
(The above is a special case of (\ref{1.6eb2})).
\begin{flalign*}
&{\rm (xix)}\,
{\mbox{$F$}(-5\,a,1-4\,a;\,2-5\,a;\,5)}
=\begin{cases}
{2}^{8\,n} \left( 20\,n+1 \right)
\quad {\text {if $a=1/4+n$}},\\
-{2}^{8\,n+2} \left( 20\,n+6 \right) 
\quad {\text {if $a=1/2+n$}},\\
{2}^{8\,n+4} \left( 20\,n+11 \right) 
\quad {\text {if $a=3/4+n$}},\\
-{2}^{8\,n+14} \left( 20\,n+36 \right) 
\quad {\text {if $a=2+n$}}
\end{cases}&
\end{flalign*}
(The above are special cases of (\ref{1.6eb2})).
\begin{flalign*}
&{\rm (xx)}\,
{\mbox{$F$}(2,1-a;\,2-5\,a;\,5)}
=5/4\,n+1
\quad {\text {if $a=1+n$}}
&
\end{flalign*}
(The above is a special case of (\ref{1.5eb2})).
The special values obtained from {\rm (xxi)} are trivial.
\begin{flalign*}
&{\rm (xxii)}\,
{\mbox{$F$}(-5\,a,1-a;\,-a;\,-1/4)}
=\begin{cases}
0,\\
{\dfrac { \left( 5,5\,n \right) }{{2}^{2\,n} \left( 2,n \right)  \left( 4,4\,n \right) }}
\quad {\text {if $a=1+n$}}\\
\end{cases}&
\end{flalign*}
(The first case is identical to (29.6) in [Ge]).
\begin{flalign*}
&{\rm (xxiii)}\,
{\mbox{$F$}(4\,a,1-a;\,-a;\,1/5)}
=\begin{cases}
0,\\
{\dfrac {\left( 5,5\,n \right) }{{5}^{n} \left( 2,n \right)  \left( 4,4\,n \right) }}
\quad {\text {if $a=1+n$}}.\\
\end{cases}&
\end{flalign*}
The special values obtained from {\rm (xxiv)} are trivial.
\paragraph{(1,5,5-2)}
The special values obtained from (1,5,5-2) coincide with
those obtained from (1,5,5-1).
\subsubsection{$(k,l,m)=(1,6,5)$}
In this case, we have
\begin{gather*}
(a,b,c,x)=(a,b,b+1-a,-1). \tag{1,6,5-1}
\end{gather*}
\paragraph{(1,6,5-1)}
The special values obtained from (1,6,5-1)
are evaluated in paragraphs (1,2,2-1) and (0,2,2-1). 

\subsubsection{$(k,l,m)=(2,3,5)$}
In this case, there is no admissible quadruple.

\subsubsection{$(k,l,m)=(2,4,5)$}
In this case, there is no admissible quadruple.

\subsubsection{$(k,l,m)=(2,5,5)$}
In this case, we have
\begin{align}
&
(a,b,c,x)=(a,5/2\,a+1,5/2\,a,5/3),\,
S^{(n)}=
{\dfrac {{5}^{5\,n} \left( a+1,2\,n \right)  \left( 3/2\,a,3\,n \right) }{{3}^{5\,n}
 \left( 5/2\,a,5\,n \right) }},
\tag{2,5,5-1}
\\
&
(a,b,c,x)=(a,5/2\,a-5/2,5/2\,a-1/2,5/3),\,
S^{(n)}=
\dfrac{2^{2\,n} (5\,a-3)}{5\,a-3+10\,n}.
\tag{2,5,5-2}
\end{align}
\paragraph{(2,5,5-1)}
\begin{flalign*}
&{\rm (i)}\,
{\mbox{$F$}(a,5/2\,a+1;\,5/2\,a;\,5/3)}&\\
&=\begin{cases}
0 
\quad {\text {if $a=-1-n$}},\\
{\dfrac {{5}^{5\,n} \left( 2/5,2\,n \right)  \left( 8/5,3\,n \right) }{{3}^{5\,n}
 \left( 2,5\,n \right) }}
\quad {\text {if $a=-2/5-2\,n$}},\\
{\dfrac {{5}^{5\,n} \left( 4/5,2\,n \right) \left( 11/5,3\,n \right) }{{3}^{5\,n+1}
 \left( 3,5\,n \right) }}
\quad {\text {if $a=-4/5-2\,n$}},\\
{\dfrac {-2\cdot {5}^{5\,n} \left( 6/5,2\,n \right) \left( 14/5,3\,n \right) }{{3}^{5\,n+2}
 \left( 4,5\,n \right) }}
\quad {\text {if $a=-6/5-2\,n$}},\\
{\dfrac {-7\cdot {5}^{5\,n} \left( 8/5,2\,n \right) \left( 17/5,3\,n \right) }{{3}^{5\,n+3}
 \left( 5,5\,n \right) }}
\quad {\text {if $a=-8/5-2\,n$}}.\\
\end{cases}&
\end{flalign*}
The special values obtained from {\rm (ii)} and {\rm (iii)} are trivial.
\begin{flalign*}
&{\rm (iv)}\,
{\mbox{$F$}(3/2\,a,5/2\,a+1;\,5/2\,a;\,5/2)}&\\
&=\begin{cases}
0
\quad {\text {if $a=-2/3-2/3\,n$}},\\
{\dfrac { \left( -5 \right) ^{5\,n} \left( 2/5,2\,n
 \right)\left( 8/5,3\,n \right) }{{2}^{5\,n} \left( 2,5\,n \right) }}
\quad {\text {if $a=-2/5-2\,n$}},\\
{\dfrac { -\left( -5 \right) ^{5\,n} \left( 4/5,2\,n
 \right)\left( 11/5,3\,n \right) }{{2}^{5\,n+1} \left( 3,5\,n \right) }}
\quad {\text {if $a=-4/5-2\,n$}},\\
{\dfrac { -\left( -5 \right) ^{5\,n} \left( 6/5,2\,n
 \right)\left( 14/5,3\,n \right) }{{2}^{5\,n+1} \left( 4,5\,n \right) }}
\quad {\text {if $a=-6/5-2\,n$}},\\
{\dfrac { 7 \left( -5 \right) ^{5\,n} \left( 8/5,2\,n
 \right)\left( 17/5,3\,n \right) }{{2}^{5\,n+3} \left( 5,5\,n \right) }}
\quad {\text {if $a=-8/5-2\,n$}}.\\
\end{cases}&
\end{flalign*}
\begin{flalign*}
&{\rm (v)}\,
{\mbox{$F$}(a,5/2\,a+1;\,a+2;\,-2/3)}
=
{3}^{5/2\,a}{5}^{-5/2\,a} \left( a+1 \right)
&
\end{flalign*}
(The above is a special case of (\ref{1.6eb2})).
\begin{flalign*}
&{\rm (vi)}\,
{\mbox{$F$}(2,1-3/2\,a;\,a+2;\,-2/3)}
=3/5\,a+3/5
&
\end{flalign*}
(The above is a special case of (\ref{1.5eb2})).
\begin{flalign*}
&{\rm (vii)}\,
{\mbox{$F$}(a,1-3/2\,a;\,a+2;\,2/5)}
=
{3}^{3/2\,a}{5}^{-3/2\,a} \left( a+1 \right) 
&
\end{flalign*}
(The above is a special case of (\ref{1.6eb2}) and is a generalization of Theorem 38 in [Ek]).
\begin{flalign*}
&{\rm (viii)}\,
{\mbox{$F$}(2,5/2\,a+1;\,a+2;\,2/5)}
=
5/3\,a+5/3
&
\end{flalign*}
(The above is a special case of (\ref{1.5eb2})).
\begin{flalign*}
&{\rm (ix)}\,
{\mbox{$F$}(a,1-3/2\,a;\,-3/2\,a;\,3/5)}&\\
&=\begin{cases}
0,\\
{\dfrac {{3}^{3\,n} \left( 8/3,5\,n \right) }{{5}^{3\,n
}\left( 5/3,2\,n \right)  \left( 2,3
\,n \right) }}
\quad {\text {if $a=2/3+2\,n$}},\\
{\dfrac {7\cdot {3}^{3\,n} \left( 13/3,5\,n \right) }{{5}^{3\,n+1
}\left( 7/3,2\,n \right)  \left( 3,3
\,n \right) }}
\quad {\text {if $a=4/3+2\,n$}},\\
{\dfrac {18\cdot {3}^{3\,n+1} \left( 6,5\,n \right) }{{5}^{3\,n+2
}\left( 3,2\,n \right)  \left( 4,3
\,n \right) }}
\quad {\text {if $a=2+2\,n$}}.
\end{cases}&
\end{flalign*}
The special values obtained from {\rm (x)} and {\rm (xi)} are trivial.
\begin{flalign*}
&{\rm (xii)}\,
{\mbox{$F$}(-5/2\,a,1-3/2\,a;\,-3/2\,a;\,-3/2)}&\\
&=\begin{cases}
0
\quad {\text {if $a=2/5+2\,n$}},&\\
0
\quad {\text {if $a=4/5+2\,n$}},&\\
0
\quad {\text {if $a=6/5+2\,n$}},&\\
0
\quad {\text {if $a=8/5+2\,n$}},&\\
{\dfrac {{3}^{3\,n}\left( 8/3,5\,n \right) }{{2}^{3\,n
} \left( 5/3,2\,n \right)  \left( 2,3
\,n \right) }}
\quad {\text {if $a=2/3+2\,n$}},&\\
{\dfrac {7\cdot {3}^{3\,n}\left( 13/3,5\,n \right) }{{2}^{3\,n+1
} \left( 7/3,2\,n \right)  \left( 3,3
\,n \right) }}
\quad {\text {if $a=4/3+2\,n$}},&\\
{\dfrac { {3}^{3\,n+3}\left( 6,5\,n \right) }{{2}^{3\,n+1
} \left( 3,2\,n \right)  \left( 4,3
\,n \right) }}
\quad {\text {if $a=2+2\,n$}}.&
\end{cases}&
\end{flalign*}
\begin{flalign*}
&{\rm (xiii)}\,{\mbox{$F$}(2,5/2\,a+1;\,3/2\,a+2;\,3/5)}
=15/4\,a+5/2
&
\end{flalign*}
(The above is a special case of (\ref{1.5eb2})).
\begin{flalign*}
&{\rm (xiv)}\,
{\mbox{$F$}(3/2\,a,1-a;\,3/2\,a+2;\,3/5)}
={2}^{a-1}{5}^{-a} \left( 3\,a+2 \right)
&
\end{flalign*}
(The above is a special case of (\ref{1.6eb2})).
\begin{flalign*}
&{\rm (xv)}\,
{\mbox{$F$}(3/2\,a,5/2\,a+1;\,3/2\,a+2;\,-3/2)}&\\
&=\begin{cases}
{2}^{-5\,n-1}{5}^{5\,n} \left( 2-15\,n \right) 
\quad {\text {if $a=-2/5-2\,n$}},&\\
{2}^{-5\,n-2}{5}^{5\,n+1} \left( -1-15\,n \right) 
\quad {\text {if $a=-4/5-2\,n$}},&\\
{2}^{-5\,n-3}{5}^{5\,n+2} \left( -4-15\,n \right) 
\quad {\text {if $a=-6/5-2\,n$}},&\\
{2}^{-5\,n-4}{5}^{5\,n+3} \left( -7-15\,n \right) 
\quad {\text {if $a=-8/5-2\,n$}}&\\
\end{cases}&
\end{flalign*}
(The above are special cases of (\ref{1.6eb2})).
\begin{flalign*}
&{\rm (xvi)}\,
{\mbox{$F$}(2,1-a;\,3/2\,a+2;\,-3/2)}
=
3/5\,n+1
\quad {\text {if $a=1+n$}}
&
\end{flalign*}
(The above is a special case of (\ref{1.5eb2})).
\begin{flalign*}
&{\rm (xvii)}\,
{\mbox{$F$}(2,1-3/2\,a;\,2-5/2\,a;\,5/3)}
=
5/2\,n+1
\quad {\text {if $a=2/3+2/3\,n$}}
&
\end{flalign*}
(The above is a special case of (\ref{1.5eb2})).
\begin{flalign*}
&{\rm (xviii)}\,
{\mbox{$F$}(-5/2\,a,1-a;\,2-5/2\,a;\,5/3)}
=
\left(-2\right)^{n}3^{-n-1}\left(5\,n+3\right)
\quad {\text {if $a=1+n$}}
&
\end{flalign*}
(The above is a special case of (\ref{1.6eb2})).
\begin{flalign*}
&{\rm (xix)}\,
{\mbox{$F$}(-5/2\,a,1-3/2\,a;\,2-5/2\,a;\,5/2)}&\\
&=
{2}^{-n-1} \left( -3 \right) ^{n} \left( 5\,n +2\right)
 \quad {\text {if $a=2/3+2/3\,n$}}
&
\end{flalign*}
(The above is a special case of (\ref{1.6eb2})).
\begin{flalign*}
&{\rm (xx)}\,
{\mbox{$F$}(2,1-a;\,2-5/2\,a;\,5/2)}
=
5/3\,n+1
 \quad {\text {if $a=1+n$}}
&
\end{flalign*}
(The above is a special case of (\ref{1.5eb2})).
The special values obtained from {\rm (xxi)} are trivial.
\begin{flalign*}
&{\rm (xxii)}\,
{\mbox{$F$}(-5/2\,a,1-a;\,-a;\,-2/3)}
=\begin{cases}
0,\\
{\dfrac {{2}^{2\,n}\left( 5/2,5\,n \right) }{{3}^{2\,n
}\left( 2,2\,n \right)  \left( 3/2,3
\,n \right) }}
\quad {\text {if $a=1+2\,n$}},&\\
{\dfrac {{2}^{2\,n+3}\left( 5,5\,n \right) }{{3}^{2\,n+1
}\left( 3,2\,n \right)  \left( 3,3
\,n \right) }}
 \quad {\text {if $a=2+2\,n$}}.&\\
\end{cases}
&
\end{flalign*}
\begin{flalign*}
&{\rm (xxiii)}\,
{\mbox{$F$}(3/2\,a,1-a;\,-a;\,2/5)}
=\begin{cases}
0,\\
{\dfrac {{2}^{2\,n} \left( 5/2,5\,n \right) }{{5}^{2\,n
} \left( 2,2\,n \right)  \left( 3/2,3
\,n \right) }}
\quad {\text {if $a=1+2\,n$}},&\\
{\dfrac {{2}^{2\,n+3} \left( 5,5\,n \right) }{{5}^{2\,n+1
} \left( 3,2\,n \right)  \left( 3,3
\,n \right) }}
 \quad {\text {if $a=2+2\,n$}}.&\\
\end{cases}
&
\end{flalign*}
The special values obtained from {\rm (xxiv)} are trivial.
\paragraph{(2,5,5-2)}
The special values obtained from (2,5,5-2)
coincide with those obtained from (2,5,5-1).

\subsubsection{$(k,l,m)=(2,6,5)$}
In this case, there is no admissible quadruple.

\subsubsection{$(k,l,m)=(2,7,5)$}
In this case, there is no admissible quadruple.

\subsection{$m=6$}
\subsubsection{$(k,l,m)=(0,6,6)$}
In this case, we have
\begin{align}
&(a,b,c,x)=(a,b,b+1-a,-1), \tag{0,6,6-1}\\
&(a,b,c,x)=(1,b,b,\lambda), \tag{0,6,6-2}\\
&(a,b,c,x)=(0,b,b+1,\lambda), \tag{0,6,6-3}\\
&(a,b,c,x)=(1,b,b,\mu), \tag{0,6,6-4}\\
&(a,b,c,x)=(0,b,b+1,\mu). \tag{0,6,6-5},
\end{align}
where $\lambda$ and $\mu$ are solutions of 
$x^2+x+1=0$ and $x^2-x+1=0$, respectively. 
\paragraph{(0,6,6-1)}
The special values obtained from (0,6,6-1)
are contained in those from (1,2,2-1) and (0,2,2-1).
\paragraph{(0,6,6-2)}
The special values obtained from (0,6,6-2) except trivial values
are special cases of (\ref{F(1,b;2;x)}).
\paragraph{(0,6,6-3)}
The special values obtained from (0,6,6-3) coincide with 
those obtained from (0,6,6-2).
\paragraph{(0,6,6-4)}
The special values obtained from (0,6,6-4) except trivial values
are special cases of (\ref{F(1,b;2;x)}).
\paragraph{(0,6,6-5)}
The special values obtained from (0,6,6-5) coincide with 
those obtained from (0,6,6-4).
\subsubsection{$(k,l,m)=(1,5,6)$}
In this case, we have
\begin{align}
&\begin{cases}
(a,b,c,x)=(a,5\,a-1/2,6\,a,-4),\\
S^{(n)}=\dfrac{{5}^{5\,n} \left( a+1/5,n \right)
 \left( a+4/5,n \right)  \left( a+3/10,n \right)  \left( a+7/10,n \right) }{ {3}^{6\,n}
  \left( a+1/3,n \right) 
  \left( a+2/3,n \right) 
  \left( a+1/6,n \right) 
   \left( a+5/6,n \right) 
 },\tag{1,5,6-1}
\end{cases}\\
&\begin{cases}
(a,b,c,x)=(a,5\,a-3/2,6\,a-1,-4),\\
S^{(n)}=\dfrac{{5}^{5\,n} \left( a+2/5,n \right)
 \left( a+3/5,n \right)  \left( a-1/10,n \right)  \left( a+1/10,n \right) }{ {3}^{6\,n}
  \left( a+1/3,n \right) 
  \left( a+2/3,n \right) 
  \left( a-1/6,n \right) 
   \left( a+1/6,n \right) 
 },\tag{1,5,6-2}
\end{cases}\\
&\begin{cases}
(a,b,c,x)=(a,5\,a-5/2,6\,a-3,-4),\\
S^{(n)}=\dfrac{{5}^{5\,n} \left( a-1/5,n \right)
 \left( a+1/5,n \right)  \left( a-3/10,n \right)  \left( a+3/10,n \right) }{ {3}^{6\,n}
  \left( a-1/3,n \right) 
  \left( a+1/3,n \right) 
  \left( a-1/6,n \right) 
   \left( a+1/6,n \right) 
 },\tag{1,5,6-3}
\end{cases}\\
&\begin{cases}
(a,b,c,x)=(a,5\,a-7/2,6\,a-4,-4),\\
S^{(n)}=\dfrac{{5}^{5\,n} \left( a-3/5,n \right)
 \left( a-2/5,n \right)  \left( a-1/10,n \right)  \left( a+1/10,n \right) }{ {3}^{6\,n}
  \left( a-2/3,n \right) 
  \left( a-1/3,n \right) 
  \left( a-1/6,n \right) 
   \left( a+1/6,n \right) 
 }.\tag{1,5,6-4}
\end{cases}
\end{align}
\paragraph{(1,5,6-1)}
\begin{flalign*}
&{\rm (i)}\,
{\mbox{$F$}(a,5\,a-1/2;\,6\,a;\,-4)}&\\
&=\begin{cases}
\dfrac{14\cdot {5}^{5\,n} \left( 9/5,n \right) 
 \left( 6/5,n \right)  \left( 17/10,n
 \right)  \left( 13/10,n \right)}{ {3
}^{6\,n+1} \left( 5/3,n \right) 
\left( 4/3,n \right) \left( 11/6,n \right) 
  \left( 7/6,n \right) }
\quad {\text {if $a=-1-n$}},&\\
\dfrac{{5}^{5\,n} \left( 7/10,n \right)  \left( 1/10,n \right)  \left( 3/5,n
 \right)  \left( 1/5,n \right) }{{3}^{6\,n}
   \left( 17/30,n
 \right)   \left( 7/30
,n \right)  \left( 11/
15,n \right) \left( 1/15,n
 \right)  }
\quad {\text {if $a=1/10-n$}},&\\
\dfrac{{5}^{5\,n+1} \left( 9/10,n \right)  \left( 3/10,n \right)  \left( 4/5,n
 \right)  \left( 2/5,n \right) }{{3}^{6\,n+1}
   \left( 23/30,n
 \right)   \left( 13/30
,n \right)  \left( 14/
15,n \right) \left( 4/15,n
 \right)  }
\quad {\text {if $a=-1/10-n$}},&\\
0
\quad {\text {if $a=-3/10-n$}},&\\
\dfrac{13\cdot {5}^{5\,n+1} \left( 23/10,n \right)  \left( 17/10,n \right)
 \left( 11/5,n \right) \left( 9/5,n
 \right)}{ {3}^{6\,n+1}
 \left( 13/6,n \right)  \left( 11/6,n \right)   \left( 7/3,n \right)  \left( 5/3,n \right)  }
\quad {\text {if $a=-3/2-n$}},&\\
0
\quad {\text {if $a=-7/10-n$}}.&\\
\end{cases}
&
\end{flalign*}
\begin{flalign*}
&{\rm (ii)}\,
{\mbox{$F$}(5\,a,a+1/2;\,6\,a;\,-4)}&\\
&=\begin{cases}
\dfrac{14\cdot {5}^{5\,n} \left( 9/5,n \right)
 \left( 6/5,n \right) \left( 17/10,n
 \right)  \left( 13/10,n \right)}{{3
}^{6\,n+1}  \left( 5/3,n \right) 
   \left( 4/3,n \right)   \left( 11/6,n \right) 
\left( 7/6,n \right) }
\quad {\text {if $a=-1-n$}},&\\
0
\quad {\text {if $a=-1/5-n$}},&\\
\dfrac{{5}^{5\,n+2} \left( 6/5,n \right)  \left( 3/5,n \right)  \left( 11/10,n
 \right)\left( 7/10,n \right) }{7\cdot  {3}
^{6\,n+1}  \left( 16/15
,n \right)  \left( 11/
15,n \right)  \left( 
37/30,n \right)   \left(17/30,n \right)  }
\quad {\text {if $a=-2/5-n$}},&\\
\dfrac{{5}^{5\,n+3} \left( 7/5,n \right)  \left( 4/5,n \right)  \left( 13/10,n
 \right)\left( 9/10,n \right) }{26\cdot  {3}
^{6\,n+1}  \left( 19/15
,n \right)  \left( 14/
15,n \right)  \left( 
43/30,n \right)   \left(23/30,n \right)  }
\quad {\text {if $a=-3/5-n$}},&\\
0
\quad {\text {if $a=-4/5-n$}},&\\
\dfrac{{5}^{5\,n} \left( 13/10,n \right)  \left(7/10,n \right)  \left( 6
/5,n \right)  \left( 4/5,n \right)}{ {3}^{6\,n}
\left( 7/6,n \right)   \left( 5/6,n \right)   \left( 4/3,n \right)   \left( 2/3,n \right) }  
\quad {\text {if $a=-1/2-n$}}.&\\
\end{cases}
&
\end{flalign*}
\begin{flalign*}
&{\rm (iii)}\,
{\mbox{$F$}(a,a+1/2;\,6\,a;\,4/5)}&\\
&\begin{cases}
\dfrac{{3}^{6\,a-3/5} \Gamma  \left( 3/5 \right)\Gamma \left(4/5\right) 
\Gamma \left(2\,a+1/3\right)
\Gamma  \left( 2\,a+2/3
 \right)  }{ 5^{4\,a-1/2}
\Gamma  \left( 8/15 \right)\Gamma \left(13/15\right) 
\Gamma \left(2\,a+2/5\right)
\Gamma  \left( 2\,a+3/5
 \right)}, \\
\dfrac{14\cdot 5^{4\,n-1}\left( 9/5,n \right)  \left( 6/5,n \right) \left( 17/
10,n \right)  \left( 13/10,n \right)}{{3}^{6\,n+1}  \left( 5/3,
n \right)   \left( 4/3,n \right) 
 \left( 11/6,n
 \right)   \left( 7/6,n \right) } 
\quad {\text {if $a=-1-n$}},&\\
\dfrac{{5}^{4\,n} \left(13/10,n \right)  \left( 7/10,n \right)  \left( 6
/5,n \right)\left( 4/5,n \right) }{ {3}^{6\,n}
 \left( 7/6,n \right)   \left( 5/6,n \right)   \left( 4/3,n \right)   \left( 2/3,n \right) }
\quad {\text {if $a=-1/2-n$}}.&
\end{cases}&
\end{flalign*}

\begin{flalign*}
&{\rm (iv)}\,
{\mbox{$F$}(5\,a,5\,a-1/2;\,6\,a;\,4/5)}&\\
&\begin{cases}
\dfrac{{3}^{6\,a-3/5} \Gamma  \left( 3/5 \right)\Gamma \left(4/5\right) 
\Gamma \left(2\,a+1/3\right)
\Gamma  \left( 2\,a+2/3
 \right)  }{ \Gamma  \left( 8/15 \right)\Gamma \left(13/15\right) 
\Gamma \left(2\,a+2/5\right)
\Gamma  \left( 2\,a+3/5
 \right)},\\
{\dfrac {14}{5^5}}\,\dfrac{\left( 9/5,n \right)  \left( 6/5,n \right) \left( 17/
10,n \right)  \left( 13/10,n \right)}{{3}^{6\,n+1}  \left( 5/3,
n \right)   \left( 4/3,n \right) 
 \left( 11/6,n
 \right)   \left( 7/6,n \right) } 
\quad {\text {if $a=-1-n$}},&\\
{\dfrac {13}{5^7}}\,\dfrac{\left( 11/5,n \right)  \left( 9/5,n \right) \left( 23/
10,n \right)  \left( 17/10,n \right)}{{3}^{6\,n+1}  \left( 7/3,
n \right)   \left( 5/3,n \right) 
 \left( 13/6,n
 \right)   \left( 11/6,n \right) } 
\quad {\text {if $a=-3/2-n$}}.&\\
\end{cases}
&
\end{flalign*}
\begin{flalign*}
&{\rm (v)}\,
{\mbox{$F$}(a,5\,a-1/2;\,1/2;\,5)}
=\begin{cases}
{\dfrac {{2}^{6\,n} \left( 7/10,n \right)  \left( 3/10,n \right) }{ \left( 3/5,n
 \right)  \left( 2/5,n \right) }}
\quad {\text {if $a=-n$}},\\
{\dfrac {{2}^{6\,n} \left( 3/5,n \right) \left( 1/5,n \right) }{ \left( 1/2,n
 \right)  \left( 3/10,n \right) }}
\quad {\text {if $a=1/10-n$}},\\
{\dfrac {{2}^{6\,n+1} \left( 4/5,n \right) \left( 2/5,n \right) }{ \left( 7/10,n
 \right)  \left( 1/2,n \right) }}
\quad {\text {if $a=-1/10-n$}},\\
0
\quad {\text {if $a=-3/10-n$}},\\
{\dfrac {{2}^{6\,n+4} \left( 6/5,n \right) \left( 4/5,n \right) }{ \left( 11/10,n
 \right)  \left( 9/10,n \right) }}
\quad {\text {if $a=-1/2-n$}},\\
0
\quad {\text {if $a=-7/10-n$}}.\\
\end{cases}&
\end{flalign*}
\begin{flalign*}
&{\rm (vi)}\,
{\mbox{$F$}(1-5\,a,1/2-a;\,1/2;\,5)}
=\begin{cases}
{\dfrac {{2}^{6\,n}\left( 3/5,n \right)  \left( 4/5,n \right) }{ \left( 1/2,n
 \right)  \left( 9/10,n \right) }}
\quad {\text {if $a=1/5+n$}},\\
0
\quad {\text {if $a=2/5+n$}},\\
0
\quad {\text {if $a=3/5+n$}},\\
{\dfrac {{2}^{6\,n+6}\left( 6/5,n \right)  \left( 7/5,n \right) }{5 \left( 11/10,n
 \right)  \left( 3/2,n \right) }}
\quad {\text {if $a=4/5+n$}},\\
{\dfrac {-{2}^{6\,n+7}\left( 7/5,n \right)  \left( 8/5,n \right) }{7 \left( 13/10,n
 \right)  \left( 17/10,n \right) }}
\quad {\text {if $a=1+n$}},\\
{\dfrac {{2}^{6\,n}\left( 9/10,n \right)  \left( 11/10,n \right) }{ \left( 4/5,n
 \right)  \left( 6/5,n \right) }}
\quad {\text {if $a=1/2+n$}}.\\
\end{cases}&
\end{flalign*}
\begin{flalign*}
&{\rm (vii)}\,
{\mbox{$F$}(a,1-5\,a;\,1/2;\,5/4)}&\\
&=\begin{cases}
{\dfrac { \left( -16 \right) ^{n}\left( 7/10,
n \right) \left( 3/10,n \right) }{
 \left( 3/5,n \right)  \left( 2/5,n \right) }}
\quad {\text {if $a=-n$}},\\
{\dfrac { \left( 3/5,
n \right) \left( 4/5,n \right) }{\left( -16 \right) ^{n}
 \left( 1/2,n \right)  \left( 9/10,n \right) }}
\quad {\text {if $a=1/5+n$}},\\
0
\quad {\text {if $a=2/5+n$}},\\
0
\quad {\text {if $a=3/5+n$}},\\
{\dfrac { -\left( 6/5,
n \right) \left( 7/5,n \right) }{5\left( -16 \right) ^{n}
 \left( 11/10,n \right)  \left( 3/2,n \right) }}
\quad {\text {if $a=4/5+n$}},\\
{\dfrac { -\left( 7/5,
n \right) \left( 8/5,n \right) }{14\left( -16 \right) ^{n}
 \left( 13/10,n \right)  \left( 17/10,n \right) }}
\quad {\text {if $a=1+n$}}.\\
\end{cases}&
\end{flalign*}
\begin{flalign*}
&{\rm (viii)}\,
{\mbox{$F$}(5\,a-1/2,1/2-a;\,1/2;\,5/4)}&\\
&=\begin{cases}
{\dfrac { \left( 3/5,n \right)\left( 
1/5,n \right) }{ \left( -16 \right) ^{n}
 \left( 1/2,n \right) \left( 3/10,n \right) }}
\quad {\text {if $a=1/10-n$}},\\
{\dfrac { -\left( 4/5,n \right)\left( 
2/5,n \right) }{2 \left( -16 \right) ^{n}
 \left( 7/10,n \right) \left( 1/2,n \right) }}
\quad {\text {if $a=-1/10-n$}},\\
0
\quad {\text {if $a=-3/10-n$}},\\
{\dfrac { -\left( 6/5,n \right)\left( 
4/5,n \right) }{ 4\left( -16 \right) ^{n}
 \left(11/10,n \right) \left( 9/10,n \right) }}
\quad {\text {if $a=-1/2-n$}},\\
0
\quad {\text {if $a=-7/10-n$}},\\
{\dfrac {\left( -16 \right) ^{n} \left( 9/10,n \right)\left( 
11/10,n \right) }{ 
 \left(4/5,n \right) \left( 6/5,n \right) }}
\quad {\text {if $a=1/2+n$}}.\\
\end{cases}&
\end{flalign*}
\begin{flalign*}
&{\rm (ix)}\,
{\mbox{$F$}(a,1-5\,a;\,3/2-4\,a;\,-1/4)}&\\
&=
\dfrac{5^{5\,a} \Gamma \left(4/5\right)\Gamma \left(6/5\right)
\Gamma \left(3/2-4\,a\right)}
{2^{8\,a} \Gamma \left(3/2\right)\Gamma \left(4/5-2\,a\right)
\Gamma \left(6/5-2\,a\right)}.&
\end{flalign*}
\begin{flalign*}
&{\rm (x)}\,
{\mbox{$F$}(3/2-5\,a,a+1/2;\,3/2-4\,a;\,-1/4)}&\\
&=
\dfrac{5^{5\,a-1/2} \Gamma \left(4/5\right)\Gamma \left(6/5\right)
\Gamma \left(3/2-4\,a\right)}
{2^{8\,a-1} \Gamma \left(3/2\right)\Gamma \left(4/5-2\,a\right)
\Gamma \left(6/5-2\,a\right)}.&
\end{flalign*}
\begin{flalign*}
&{\rm (xi)}\,
{\mbox{$F$}(a,a+1/2;\,3/2-4\,a;\,1/5)}=
\dfrac{5^{6\,a} \Gamma \left(4/5\right)\Gamma \left(6/5\right)
\Gamma \left(3/2-4\,a\right)}
{2^{10\,a} \Gamma \left(3/2\right)\Gamma \left(4/5-2\,a\right)
\Gamma \left(6/5-2\,a\right)}&
\end{flalign*}
(The above is a generalization of Theorem 20 in [Ek]).
\begin{flalign*}
&{\rm (xii)}\,
{\mbox{$F$}(1-5\,a,3/2-5\,a;\,3/2-4\,a;\,1/5)}&\\
&=
\dfrac{5\, \Gamma \left(4/5\right)\Gamma \left(6/5\right)
\Gamma \left(3/2-4\,a\right)}
{2^{2-2\,a} \Gamma \left(3/2\right)\Gamma \left(4/5-2\,a\right)
\Gamma \left(6/5-2\,a\right)}.&
\end{flalign*}
\begin{flalign*}
&{\rm (xiii)}\,
{\mbox{$F$}(5\,a-1/2,1/2-a;\,4\,a+1/2;\,-1/4)}&\\
&=
\dfrac{2^{8\,a-4}\Gamma \left(7/5\right)\Gamma \left(8/5\right)
\Gamma \left(4\,a+1/2\right)}
{5^{5\,a-5/2} \Gamma \left(5/2\right)\Gamma \left(2\,a+2/5\right)
\Gamma \left(2\,a+3/5\right)}.&
\end{flalign*}
\begin{flalign*}
&{\rm (xiv)}\,
{\mbox{$F$}(5\,a,1-a;\,4\,a+1/2;\,-1/4)}=
\dfrac{2^{8\,a-3}\Gamma \left(7/5\right)\Gamma \left(8/5\right)
\Gamma \left(4\,a+1/2\right)}
{5^{5\,a-2} \Gamma \left(5/2\right)\Gamma \left(2\,a+2/5\right)
\Gamma \left(2\,a+3/5\right)}.&
\end{flalign*}
\begin{flalign*}
&{\rm (xv)}\,
{\mbox{$F$}(5\,a,5\,a-1/2;\,4\,a+1/2;\,1/5)}=
\dfrac{25\cdot 2^{-2\,a-3}\Gamma \left(7/5\right)\Gamma \left(8/5\right)
\Gamma \left(4\,a+1/2\right)}
{ \Gamma \left(5/2\right)\Gamma \left(2\,a+2/5\right)
\Gamma \left(2\,a+3/5\right)}.&
\end{flalign*}
\begin{flalign*}
&{\rm (xvi)}\,
{\mbox{$F$}(1-a,1/2-a;\,4\,a+1/2;\,1/5)}=
\dfrac{ 2^{10\,a-5}\Gamma \left(7/5\right)\Gamma \left(8/5\right)
\Gamma \left(4\,a+1/2\right)}
{ 5^{6\,a-3}\Gamma \left(5/2\right)\Gamma \left(2\,a+2/5\right)
\Gamma \left(2\,a+3/5\right)}&
\end{flalign*}
(The above is a generalization of Theorem 19 in [Ek]).
\begin{flalign*}
&{\rm (xvii)}\,
{\mbox{$F$}(1-5\,a,1/2-a;\,2-6\,a;\,-4)}&\\
&=\begin{cases}
\dfrac{{5}^{5\,n} \left( 1/10,n \right) 
 \left( 3/10,n \right)  \left( 3/5,n \right)  \left( 4/5,n \right)}{  {3}^{6\,n}   \left( 1/30,n \right)   \left( 11/30,n \right)   \left( 8/15,n \right)  \left( 13/15,n \right)  }
 \quad {\text {if $a=1/5+n$}},\\
 0
 \quad {\text {if $a=2/5+n$}},\\
  0
 \quad {\text {if $a=3/5+n$}},\\
 \dfrac{{5}^{5\,n+3} \left( 7/10,n \right) 
 \left( 9/10,n \right)  \left( 6/5,n \right)  \left( 7/5,n \right)}{7\cdot  {3}^{6\,n+1}   \left( 19/30,n \right)   \left( 29/30,n \right)   \left( 17/15,n \right)  \left( 22/15,n \right)  }
 \quad {\text {if $a=4/5+n$}},\\
 \dfrac{-22\cdot {5}^{5\,n} \left( 19/10,n \right) 
 \left( 21/10,n \right)  \left( 12/5,n \right)  \left( 13/5,n \right)}{  {3}^{6\,n}   \left( 11/6,n \right)   \left( 13/6,n \right)   \left( 7/3,n \right)  \left( 8/3,n \right)  }
 \quad {\text {if $a=2+n$}},\\
{\dfrac {{5}^{5\,n} \left( 2/5,n \right)  \left( 3/5,n \right)   \left( 9/10,
n \right)  \left( 11/10\right) }{{3}^{6\,n}\left( 1/3,n
 \right) \left( 2/3,n \right) 
 \left( 5/6,n \right) \left( 7/6,n \right) }}
\quad {\text {if $a=1/2+n$}}.\\
\end{cases}&
\end{flalign*}
\begin{flalign*}
&{\rm (xviii)}\,
{\mbox{$F$}(1-a,3/2-5\,a;\,2-6\,a;\,-4)}&\\
&=\begin{cases}
{\dfrac {{5}^{5\,n} \left( 9/10,n \right)  \left( 11/10,n \right)  \left( 7/5,n
 \right) \left( 8/5,n \right) }{{3}^{6\,n} \left( 5/6,n \right)  \left( 7/6,n
 \right)  \left( 4/3,n \right) 
 \left( 5/3,n \right) }}
\quad {\text {if $a=1+n$}},\\
{\dfrac {{5}^{5\,n} \left( 1/5,n \right)  \left( 2/5,n \right)  \left( 7/10,n
 \right) \left( 9/10,n \right) }{{3}^{6\,n} \left( 2/15,n \right)  \left( 7/15,n
 \right)  \left( 19/30,n \right) 
 \left( 29/30,n \right) }}
\quad {\text {if $a=3/10+n$}},\\
{\dfrac {-11\cdot {5}^{5\,n} \left( 7/5,n \right)  \left( 8/5,n \right)  \left( 19/10,n
 \right) \left( 21/10,n \right) }{7\cdot {3}^{6\,n-1} \left( 4/3,n \right)  \left( 5/3,n
 \right)  \left( 11/6,n \right) 
 \left( 13/6,n \right) }}
\quad {\text {if $a=3/2+n$}},\\
{\dfrac {{5}^{5\,n+2} \left( 3/5,n \right)  \left( 4/5,n \right)  \left( 11/10,n
 \right) \left( 13/10,n \right) }{11\cdot {3}^{6\,n} \left( 8/15,n \right)  \left( 13/15,n
 \right)  \left( 31/30,n \right) 
 \left( 41/30,n \right) }}
\quad {\text {if $a=7/10+n$}},\\
0
\quad {\text {if $a=9/10+n$}},\\
0
\quad {\text {if $a=11/10+n$}}.\\
\end{cases}&
\end{flalign*}
\begin{flalign*}
&{\rm (xix)}\,
{\mbox{$F$}(1-5\,a,3/2-5\,a;\,2-6\,a;\,4/5)}&\\
&=\begin{cases}
\dfrac{ 3^{6/5-6\,a}\Gamma \left(2/5\right)\Gamma \left(4/5\right)
\Gamma \left(2/3-2\,a\right)\Gamma \left( 4/3-2\,a \right)}
{ \Gamma \left(4/15\right)\Gamma \left(14/15\right)
\Gamma \left(4/5-2\,a\right)\Gamma \left( 6/5-2\,a \right) },\\
\dfrac{-22\left( 19/10,n
 \right) \left( 21/10,n \right)  \left( 12/5,n \right)  \left( 13/5,n \right)}{ 5^9{3}^{6\,n}   \left( 11/6,n \right)   \left( 13/6,n \right)  \left( 7/3,n \right)  \left( 8/3,n \right) }
\quad {\text {if $a=2+n$}},\\
\dfrac{-11\left( 7/5,n
 \right) \left( 8/5,n \right)  \left( 19/10,n \right)  \left( 21/10,n \right)}{ 5^6\cdot 7\cdot {3}^{6\,n-1}   \left( 4/3,n \right)   \left( 5/3,n \right)  \left( 11/6,n \right)  \left( 13/6,n \right) }
\quad {\text {if $a=3/2+n$}}.\\
\end{cases}&
\end{flalign*}
\begin{flalign*}
&{\rm (xx)}\,
{\mbox{$F$}(1-a,1/2-a;\,2-6\,a;\,4/5)}&\\
&=\begin{cases}
\dfrac{ 3^{6/5-6\,a}\Gamma \left(2/5\right)\Gamma \left(4/5\right)
\Gamma \left(2/3-2\,a\right)\Gamma \left( 4/3-2\,a \right)}
{5^{1/2-4\,a} \Gamma \left(4/15\right)\Gamma \left(14/15\right)
\Gamma \left(4/5-2\,a\right)\Gamma \left( 6/5-2\,a \right) },\\
{\dfrac {{5}^{4\,n} \left( 9/10,n \right)  \left( 11/10,n \right)  \left( 7/5,n
 \right)  \left( 8/5,n \right) }{{3}^{6\,n} \left( 5/6,n \right)  \left( 7/6,n
 \right) \left( 4/3,n \right) 
 \left( 5/3,n \right) }}
\quad {\text {if $a=1+n$}},\\
{\dfrac {{5}^{4\,n}\left( 2/5,n \right) \left( 3/5,n \right)  \left( 9/10,n
 \right) \left( 11/10,n \right) }{{3}^{6\,n} \left( 1/3,n \right)  \left( 2/3,n
 \right)  \left( 5/6,n \right) 
 \left( 7/6,n \right) }}
\quad {\text {if $a=1/2+n$}}.\\
\end{cases}&
\end{flalign*}
\begin{flalign*}
&{\rm (xxi)}\,
{\mbox{$F$}(5\,a,a+1/2;\,3/2;\,5)}
=\begin{cases}
{\dfrac {{2}^{6\,n} \left( 4/5,n \right)  \left( 1/5,n \right) }{ \left( 11/10,n
 \right)  \left( 9/10,n \right) }}
\quad {\text {if $a=-n$}},\\
0
\quad {\text {if $a=-1/5-n$}},\\
{\dfrac {{2}^{6\,n+4} \left( 6/5,n \right)  \left( 3/5,n \right) }{15 \left( 3/2,n
 \right)  \left( 13/10,n \right) }}
\quad {\text {if $a=-2/5-n$}},\\
{\dfrac {{2}^{6\,n+6} \left( 7/5,n \right)  \left( 4/5,n \right) }{35 \left( 17/10,n
 \right)  \left( 3/2,n \right) }}
\quad {\text {if $a=-3/5-n$}},\\
0
\quad {\text {if $a=-4/5-n$}},\\
{\dfrac {{2}^{6\,n} \left( 13/10,n \right)  \left( 7/10,n \right) }{ \left( 8/5,n
 \right)  \left( 7/5,n \right) }}
\quad {\text {if $a=-1/2-n$}}.
\end{cases}&
\end{flalign*}
\begin{flalign*}
&{\rm (xxii)}\,
{\mbox{$F$}(1-a,3/2-5\,a;\,3/2;\,5)}
=\begin{cases}
{\dfrac {{2}^{6\,n} \left( 9/10,n \right)  \left( 11/10,n \right) }{ \left( 6/5,n
 \right)  \left( 9/5,n \right) }}
 \quad {\text {if $a=1+n$}},\\
{\dfrac {{2}^{6\,n} \left( 1/5,n \right)  \left( 2/5,n \right) }{ \left( 1/2,n
 \right)  \left( 11/10,n \right) }}
 \quad {\text {if $a=3/10+n$}},\\
 {\dfrac {-{2}^{6\,n+1} \left( 2/5,n \right)  \left( 3/5,n \right) }{ 3\left( 7/10,n
 \right)  \left( 13/10,n \right) }}
 \quad {\text {if $a=1/2+n$}},\\
 {\dfrac {{2}^{6\,n+3} \left( 3/5,n \right)  \left( 4/5,n \right) }{ 5\left( 9/10,n
 \right)  \left( 3/2,n \right) }}
\quad {\text {if $a=7/10+n$}},\\
0
\quad {\text {if $a=9/10+n$}},\\
0
\quad {\text {if $a=11/10+n$}}.\\
\end{cases}&
\end{flalign*}
\begin{flalign*}
&{\rm (xxiii)}\,
{\mbox{$F$}(5\,a,1-a;\,3/2;\,5/4)}&\\
&=\begin{cases}
{\dfrac { \left( 4/5,n \right)  \left( 1/5,n \right) }
{\left(-16\right)^{n} \left( 11/10,n
 \right)  \left( 9/10,n \right) }}
\quad {\text {if $a=-n$}},\\
0
\quad {\text {if $a=-1/5-n$}},\\ 
{\dfrac { \left( 6/5,n \right)  \left( 3/5,n \right) }
{15\left( -16\right)^{n} \left( 3/2,n
 \right)  \left( 13/10,n \right) }}
\quad {\text {if $a=-2/5-n$}},\\
{\dfrac {- \left( 7/5,n \right)  \left( 4/5,n \right) }
{35\left( -16\right)^{n} \left( 17/10,n
 \right)  \left( 3/2,n \right) }}
\quad {\text {if $a=-3/5-n$}},\\
0
\quad {\text {if $a=-4/5-n$}},\\
{\dfrac {{\left(-16\right)}^{n} \left( 9/10,n \right)  \left(11/10,n \right) }{  \left( 6/5,n
 \right)  \left( 9/5,n \right) }}
\quad {\text {if $a=1+n$}}.\\
\end{cases}&
\end{flalign*}
\begin{flalign*}
&{\rm (xxiv)}\,
{\mbox{$F$}(a+1/2,3/2-5\,a;\,3/2;\,5/4)}&\\
&=\begin{cases}
{\dfrac {{\left(-16\right)}^{n} \left( 13/10,n \right)  \left(7/10,n \right) }{  \left( 8/5,n
 \right)  \left( 7/5,n \right) }}
\quad {\text {if $a=-1/2-n$}},\\
{\dfrac { \left( 1/5,n \right)  \left(2/5,n \right) }
{ \left(-16\right)^{n} \left( 1/2,n
 \right)  \left( 11/10,n \right) }}
\quad {\text {if $a=3/10+n$}},\\
{\dfrac { \left( 2/5,n \right)  \left(3/5,n \right) }
{ 6\left(-16\right)^{n} \left( 7/10,n
 \right)  \left( 13/10,n \right) }}
\quad {\text {if $a=1/2+n$}},\\
{\dfrac { \left( 3/5,n \right)  \left(4/5,n \right) }
{ 10 \left(-16\right)^{n} \left( 9/10,n
 \right)  \left( 3/2,n \right) }}
\quad {\text {if $a=7/10+n$}},\\
0
\quad {\text {if $a=9/10+n$}},\\
0
\quad {\text {if $a=11/10+n$}}.
\end{cases}&
\end{flalign*}

\paragraph{(1,5,6-2), (1,5,6-3), (1,5,6-4)}
The special values obtained from (1,5,6-2), (1,5,6-3) and (1,5,6-4)
coincide with those obtained from (1,5,6-1).

\subsubsection{$(k,l,m)=(1,6,6)$}
In this case, we have
\begin{align*}
&
(a,b,c,x)=(a,6\,a+1,6\,a,6/5),\,
S^{(n)}=
\dfrac {\left(-1\right)^n 2^{6\,n}3^{6\,n}
\left( a+1,n \right)\left( 5\,a,5\,n \right)  
  }{5^{6\,n}\left(6\,a,6\,n\right) },
\tag{1,6,6-1}\\
&(a,b,c,x)=(a,6\,a-6,6\,a-4,6/5),\,
S^{(n)}=
{\frac {6\,a-5}{ \left( -5 \right) ^{n} \left( 6\,a-5+6\,n \right) }}.
\tag{1,6,6-2}
\end{align*}
\paragraph{(1,6,6-1)}
\begin{flalign*}
&{\rm (i)}\,{\mbox{$F$}(a,6\,a+1;\,6\,a;\,6/5)}&\\
&=\begin{cases}
0
\quad {\text {if $a=-1-n$}},\\
\dfrac{\left( -1 \right) ^{n}{2}^{6\,n}{3}^{6\,n} \left( 1/6
,n \right)  \left( 11/6,5\,n \right) }{ {5}^{6\,n} \left( 2,6\,n \right) }
\quad {\text {if $a=-1/6-n$}},\\
\dfrac{\left( -1 \right) ^{n}{2}^{6\,n+2}{3}^{6\,n} \left( 1/3
,n \right)  \left( 8/3,5\,n \right) }{ {5}^{6\,n+1} \left( 3,6\,n \right) }
\quad {\text {if $a=-1/3-n$}},\\
\dfrac{\left( -1 \right) ^{n}{2}^{6\,n-1}{3}^{6\,n+3} \left( 1/2
,n \right)  \left( 7/2,5\,n \right) }{ {5}^{6\,n+2} \left( 4,6\,n \right) }
\quad {\text {if $a=-1/2-n$}},\\
\dfrac{7\left( -1 \right) ^{n}{2}^{6\,n+4}{3}^{6\,n-1} \left( 2/3
,n \right)  \left( 13/3,5\,n \right) }{ {5}^{6\,n+3} \left( 5,6\,n \right) }
\quad {\text {if $a=-2/3-n$}},\\
\dfrac{1729\left( -1 \right) ^{n}{2}^{6\,n-3}{3}^{6\,n-1} \left( 5/6
,n \right)  \left( 31/6,5\,n \right) }{ {5}^{6\,n+4} \left( 6,6\,n \right) }
\quad {\text {if $a=-5/6-n$}}.\\
\end{cases}&
\end{flalign*}
The special values obtained from {\rm (ii)} and {\rm (iii)}
are trivial.
\begin{flalign*}
&{\rm (iv)}\,{\mbox{$F$}(5\,a,6\,a+1;\,6\,a;\,6)}&\\
&=\begin{cases}
0
\quad {\text {if $a=-1/5-1/5\,n$}},\\
\dfrac{\left( -1 \right) ^{n}{2}^{6\,n}{3}^{6\,n} \left( 1/6
,n \right)  \left( 11/6,5\,n \right) }{ \left( 2,6\,n \right) }
\quad {\text {if $a=-1/6-n$}},\\
\dfrac{\left( -1 \right) ^{n+1}{2}^{6\,n+2}{3}^{6\,n} \left( 1/3
,n \right)  \left( 8/3,5\,n \right) }{ \left( 3,6\,n \right) }
\quad {\text {if $a=-1/3-n$}},\\
\dfrac{\left( -1 \right) ^{n}{2}^{6\,n-1}{3}^{6\,n+3} \left( 1/2
,n \right)  \left( 7/2,5\,n \right) }{ \left( 4,6\,n \right) }
\quad {\text {if $a=-1/2-n$}},\\
\dfrac{7\left( -1 \right) ^{n+1}{2}^{6\,n+4}{3}^{6\,n-1} \left( 2/3
,n \right)  \left( 13/3,5\,n \right) }{ \left( 5,6\,n \right) }
\quad {\text {if $a=-2/3-n$}},\\
\dfrac{1729\left( -1 \right) ^{n}{2}^{6\,n-3}{3}^{6\,n-1} \left( 5/6
,n \right)  \left( 31/6,5\,n \right) }{ \left( 6,6\,n \right) }
\quad {\text {if $a=-5/6-n$}}.
\end{cases}&
\end{flalign*}
\begin{flalign*}
&{\rm (v)}\,
{\mbox{$F$}(a,6\,a+1;\,a+2;\,-1/5)}
={5}^{6\,a}{6}^{-6\,a} \left( a+1 \right) 
&
\end{flalign*}
(The above is a special case of (\ref{1.6eb2})).
\begin{flalign*}
&{\rm (vi)}\,
{\mbox{$F$}(2,1-5\,a;\,a+2;\,-1/5)}
=5/6\,a+5/6
&
\end{flalign*}
(The above is a special case of (\ref{1.5eb2})).
\begin{flalign*}
&{\rm (vii)}\,
{\mbox{$F$}(a,1-5\,a;\,a+2;\,1/6)}
={5}^{5\,a}{6}^{-5\,a} \left( a+1 \right) 
&
\end{flalign*}
(The above is a special case of (\ref{1.6eb2})).
\begin{flalign*}
&{\rm (viii)}\,{\mbox{$F$}(2,6\,a+1;\,a+2;\,1/6)}
=6/5\,a+6/5
&
\end{flalign*}
(The above is a special case of (\ref{1.5eb2})).
\begin{flalign*}
&{\rm (ix)}\,
{\mbox{$F$}(a,1-5\,a;\,-5\,a;\,5/6)}&\\
&=\begin{cases}
0,\\
{\dfrac {{5}^{5\,n} \left( 6/5,6\,n \right) }{{2}^{5\,n
}{3}^{5\,n} \left( 6/5,n \right) 
 \left( 1,5\,n \right) }}
\quad {\text {if $a=1/5+n$}},\\
{\dfrac {7\cdot {5}^{5\,n} \left( 12/5,6\,n \right) }{{2}^{5\,n+1
}{3}^{5\,n+1} \left( 7/5,n \right) 
 \left( 2,5\,n \right) }}
\quad {\text {if $a=2/5+n$}},\\
{\dfrac {13\cdot {5}^{5\,n} \left( 18/5,6\,n \right) }{{2}^{5\,n
}{3}^{5\,n+2} \left( 8/5,n \right) 
 \left( 3,5\,n \right) }}
\quad {\text {if $a=3/5+n$}},\\
{\dfrac {133\cdot {5}^{5\,n} \left( 24/5,6\,n \right) }{{2}^{5\,n+3
}{3}^{5\,n+2} \left( 9/5,n \right) 
 \left( 4,5\,n \right) }}
\quad {\text {if $a=4/5+n$}},\\
{\dfrac {{5}^{5\,n+5} \left( 6,6\,n \right) }{{2}^{5\,n+4
}{3}^{5\,n+4} \left( 2,n \right) 
 \left( 5,5\,n \right) }}
\quad {\text {if $a=1+n$}}.
\end{cases}&
\end{flalign*}
The special values obtained from {\rm (x)} and {\rm (xi)} are trivial.
\begin{flalign*}
&{\rm (xii)}\,
{\mbox{$F$}(-6\,a,1-5\,a;\,-5\,a;\,-5)}
=\begin{cases}
0
\quad {\text {if $a=1/6+n$}},\\
0
\quad {\text {if $a=1/3+n$}},\\
0
\quad {\text {if $a=1/2+n$}},\\
0
\quad {\text {if $a=2/3+n$}},\\
0
\quad {\text {if $a=5/6+n$}},\\
{\dfrac {{5}^{5\,n} \left( 6/5,6\,n \right) }{ \left( 6/5,n \right)  \left( 1,5\,n \right) }}
\quad {\text {if $a=1/5+n$}},\\
{\dfrac {7\cdot {5}^{5\,n} \left( 12/5,6\,n \right) }{ \left( 7/5,n \right)  \left( 2,5\,n \right) }}
\quad {\text {if $a=2/5+n$}},\\
{\dfrac {52\cdot {5}^{5\,n} \left( 18/5,6\,n \right) }{ \left( 8/5,n \right)  \left( 3,5\,n \right) }}
\quad {\text {if $a=3/5+n$}},\\
{\dfrac {399\cdot {5}^{5\,n} \left( 24/5,6\,n \right) }{ \left( 9/5,n \right)  \left( 4,5\,n \right) }}
\quad {\text {if $a=4/5+n$}},\\
{\dfrac { {5}^{5\,n+5} \left( 6,6\,n \right) }{ \left( 2,n \right)  \left( 5,5\,n \right) }}
\quad {\text {if $a=1+n$}}.
\end{cases}&
\end{flalign*}
\begin{flalign*}
&{\rm (xiii)}\,{\mbox{$F$}(2,6\,a+1;\,5\,a+2;\,5/6)}=
6\left(5\,a+1\right)
&
\end{flalign*}
(The above is a special case of (\ref{1.5eb2})).
\begin{flalign*}
&{\rm (xiv)}\,
{\mbox{$F$}(5\,a,1-a;\,5\,a+2;\,5/6)}
={6}^{-a} \left( 5\,a+1 \right) 
&
\end{flalign*}
(The above is a special case of (\ref{1.6eb2})).
\begin{flalign*}
&{\rm (xv)}\,
{\mbox{$F$}(5\,a,6\,a+1;\,5\,a+2;\,-5)}
=\begin{cases}
6^{6\,n}\left(1-30\,n\right)
\quad {\text {if $a=-1/6-n$}},\\
6^{6\,n+1}\left(-4-30\,n\right)
\quad {\text {if $a=-1/3-n$}},\\
6^{6\,n+2}\left(-9-30\,n\right)
\quad {\text {if $a=-1/2-n$}},\\
6^{6\,n+3}\left(-14-30\,n\right)
\quad {\text {if $a=-2/3-n$}},\\
6^{6\,n+4}\left(-19-30\,n\right)
\quad {\text {if $a=-5/6-n$}}
\end{cases}&
\end{flalign*}
(The above are special cases of (\ref{1.6eb2})).
\begin{flalign*}
&{\rm (xvi)}\,
{\mbox{$F$}(2,1-a;\,5\,a+2;\,-5)}
=5/6\,n+1
\quad {\text {if $a=1+n$}}
&
\end{flalign*}
(The above is a special case of (\ref{1.5eb2})).
\begin{flalign*}
&{\rm (xvii)}\,
{\mbox{$F$}(2,1-5\,a;\,2-6\,a;\,6/5)}
=
6\,n+1
\quad {\text {if $a=1/5+1/5\,n$}}&
\end{flalign*}
(The above is a special case of (\ref{1.5eb2})).
\begin{flalign*}
&{\rm (xviii)}\,
{\mbox{$F$}(-6\,a,1-a;\,2-6\,a;\,6/5)}
=-5^{-n-1}\left(6\,n+5\right)
\quad {\text {if $a=1+n$}}
&
\end{flalign*}
(The above is a special case of (\ref{1.6eb2})).
\begin{flalign*}
&{\rm (xix)}\,
{\mbox{$F$}(-6\,a,1-5\,a;\,2-6\,a;\,6)}
=
\left(-5\right)^{n}\left(6\,n+1\right)
\quad {\text {if $a=1/5+1/5\,n$}}
&
\end{flalign*}
(The above is a special case of (\ref{1.6eb2})).
\begin{flalign*}
&{\rm (xx)}\,
{\mbox{$F$}(2,1-a;\,2-6\,a;\,6)}
=6/5\,n+1
\quad {\text {if $a=1+n$}}
&
\end{flalign*}
(The above is a special case of (\ref{1.5eb2})).
The special values obtained from {\rm (xxi)} are trivial.
\begin{flalign*}
&{\rm (xxii)}\,
{\mbox{$F$}(-6\,a,1-a;\,-a;\,-1/5)}
=\begin{cases}
0,\\
{\dfrac { \left( 6,6\,n \right) }{{5}^{n} \left( 2,n \right) \left( 5,5\,n \right) }}
\quad {\text {if $a=1+n$}}.
\end{cases}&
\end{flalign*}
\begin{flalign*}
&{\rm (xxiii)}\,
{\mbox{$F$}(5\,a,1-a;\,-a;\,1/6)}
=\begin{cases}
0,\\
{\dfrac { \left( 6,6\,n \right) }{{2}^{n}3^n \left( 2,n \right) \left( 5,5\,n \right) }}
\quad {\text {if $a=1+n$}}.
\end{cases}&
\end{flalign*}
The special values obtained from {\rm (xxiv)} are trivial.

\paragraph{(1,6,6-2)}
The special values obtained from (1,6,6-2) coincide with 
those obtained from (1,6,6-1).

\subsubsection{$(k,l,m)=(1,7,6)$}
In this case, there is no admissible quadruple.

\subsubsection{$(k,l,m)=(2,4,6)$}
In this case, we have
\begin{align}
&(a,b,c,x)=(a,2\,a-1/3,3\,a,9),\tag{2,4,6-1}\\
&(a,b,c,x)=(a,2\,a-2/3,3\,a-1,9),\tag{2,4,6-2}\\
&\begin{cases}
(a,b,c,x)=(a,2\,a-1/3,3\,a,-9+6\,\sqrt {3}),\\
S^{(n)}=
\dfrac{{\left(-2\right)}^{n} \left( \sqrt {3}-1 \right) ^{6\,n} \left( 1
/2\,a+3/4,n \right)  \left( 1/2\,a+5/12,n
 \right) }{ {3}^{3/2\,n}  
 \left( 1/2\,a+1/3,n \right) 
 \left( 1/2\,a+5/6,n \right)} ,
\tag{2,4,6-3}
\end{cases}
\\
&
\begin{cases}
(a,b,c,x)=(a,2\,a-1/3,3\,a,-9-6\,\sqrt {3}),\\
S^{(n)}=
\dfrac{{2}^{n} \left( \sqrt {3}+1 \right) ^{6\,n} \left( 1
/2\,a+3/4,n \right)  \left( 1/2\,a+5/12,n
 \right) }{ {3}^{3/2\,n}  
 \left( 1/2\,a+1/3,n \right) 
 \left( 1/2\,a+5/6,n \right)} ,
\tag{2,4,6-4}
\end{cases}\\
&
\begin{cases}
(a,b,c,x)=(a,2\,a-2/3,3\,a-1,-9+6\,\sqrt {3}),\\
S^{(n)}=
\dfrac{{\left(-2\right)}^{n} \left( \sqrt {3}-1 \right) ^{6\,n} \left( 1
/2\,a+1/4,n \right)  \left( 1/2\,a+7/12,n
 \right) }{ {3}^{3/2\,n}  
 \left( 1/2\,a+2/3,n \right) 
 \left( 1/2\,a+1/6,n \right)} ,
\tag{2,4,6-5}
\end{cases}
\\
&\begin{cases}
(a,b,c,x)=(a,2\,a-2/3,3\,a-1,-9-6\,\sqrt {3}),\\
S^{(n)}=
\dfrac{{2}^{n} \left( \sqrt {3}+1 \right) ^{6\,n} \left( 1
/2\,a+1/4,n \right)  \left( 1/2\,a+7/12,n
 \right) }{ {3}^{3/2\,n}  
 \left( 1/2\,a+2/3,n \right) 
 \left( 1/2\,a+1/6,n \right)}.
\tag{2,4,6-6}
\end{cases}
\end{align}
\paragraph{(2,4,6-1), (2,4,6-2)}
The special values obtained from (2,4,6-1) and (2,4,6-2)
coincide with those obtained from (1,2,3-1).
\paragraph{(2,4,6-3)}
\begin{flalign*}
&{\rm (i)}\,
{\mbox{$F$}(a,2\,a-1/3;\,3\,a;\,-9+6\,\sqrt {3})}&\\
&=\begin{cases}
{\dfrac {7\left(-2\right)^{n-3} \left( \sqrt {3}-1 \right) ^{6\,n+6} \left( 5/4,n \right) \left( 19/12,n
 \right) }{{3}^{3/2\,n+1/2} \left( 5/3,n \right)  \left( 7/6,n \right) }}
\quad {\text {if $a=-2-2\,n$}},\\
{\dfrac {{\left(-2\right)}^{n-1} \left( \sqrt {3}-1 \right) ^{6\,n+4} \left( 3/4,n \right) \left( 13/12,n
 \right) }{{3}^{3/2\,n+1/2} \left( 7/6,n \right)  \left( 2/3,n \right) }}
\quad {\text {if $a=-1-2\,n$}},\\
{\dfrac {{\left(-2\right)}^{n} \left( \sqrt {3}-1 \right) ^{6\,n}
 \left( 1/6,n \right) \left( 1/2,n \right) }{{3}^{3/2
\,n}\left( 7/12,n \right)  \left( 1/12
,n \right) }}
\quad {\text {if $a=1/6-2\,n$}},\\
{\dfrac {-{\left(-2\right)}^{n-2} \left( \sqrt {3}-1 \right) ^{6\,n+8}
 \left( 17/12,n \right) \left( 7/4,n \right) }{{3}^{3/2
\,n-1/2}\left(11/6,n \right)  \left( 4/3
,n \right) }}
\quad {\text {if $a=-7/3-2\,n$}}\\
0
\quad {\text {if $a=-5/6-2\,n$}},\\
{\dfrac {{\left(-2\right)}^{n-1} \left( \sqrt {3}-1 \right) ^{6\,n+4}
 \left( 11/12,n \right) \left( 5/4,n \right) }{{3}^{3/2
\,n}\left( 4/3,n \right)  \left( 5/6
,n \right) }}
\quad {\text {if $a=-4/3-2\,n$}}.
\end{cases}&
\end{flalign*}
\begin{flalign*}
&{\rm (ii)}\,
{\mbox{$F$}(2\,a,a+1/3;\,3\,a;\,-9+6\,\sqrt {3})}&\\
&=\begin{cases}
{\dfrac {7\left(-2\right)^{n-3} \left( \sqrt {3}-1 \right) ^{6\,n+6} \left( 5/4,n \right) \left( 19/12,n
 \right) }{{3}^{3/2\,n+1/2} \left( 5/3,n \right)  \left( 7/6,n \right) }}
\quad {\text {if $a=-2-2\,n$}},\\
{\dfrac { - {\left(-2\right)}^{n+1} \left( \sqrt {3}-1 \right) ^{6\,n+1} \left( 1/2,n \right) \left( 5/6,n
 \right) }{{3}^{3/2\,n+1/2} \left( 11/12,n \right)  \left( 5/12,n \right) }}
\quad {\text {if $a=-1/2-2\,n$}},\\
{\dfrac {{\left(-2\right)}^{n} \left( \sqrt {3}-1 \right) ^{6\,n+2} \left( 3/4,n \right) \left( 13/12,n
 \right) }{{3}^{3/2\,n+1/2} \left( 7/6,n \right)  \left( 2/3,n \right) }}
\quad {\text {if $a=-1-2\,n$}},\\
0
\quad {\text {if $a=-3/2-2\,n$}},\\
{\dfrac {{\left(-2\right)}^{n} \left( \sqrt {3}-1 \right) ^{6\,n} \left( 5/12,n \right) \left( 3/4,n
 \right) }{{3}^{3/2\,n} \left( 5/6,n \right)  \left( 1/3,n \right) }}
\quad {\text {if $a=-1/3-2\,n$}},\\
{\dfrac {{\left(-2\right)}^{n-2} \left( \sqrt {3}-1 \right) ^{6\,n+4} \left( 11/12,n \right) \left( 5/4,n
 \right) }{{3}^{3/2\,n} \left( 4/3,n \right)  \left( 5/6,n \right) }}
\quad {\text {if $a=-4/3-2\,n$}}.
\end{cases}&
\end{flalign*}
\begin{flalign*}
&{\rm (iii)}\,
{\mbox{$F$}(a,a+1/3;\,3\,a;\,9/4+3/4\,\sqrt {3})}&\\
&=\begin{cases}
{\dfrac { 7 \left(-2\right)^{n-3} \left( 5/4,n \right) \left( 19/12,n
 \right) }{{3}^{3/2\,n+1/2} \left( 5/3,n \right)  \left( 7/6,n \right) }}
\quad {\text {if $a=-2-2\,n$}},\\
{\dfrac {  -\left(-2\right)^{n-1}\left(\sqrt{3}-1\right) \left( 3/4,n \right) \left( 13/12,n
 \right) }{{3}^{3/2\,n+1/2} \left( 7/6,n \right)  \left( 2/3,n \right) }}
\quad {\text {if $a=-1-2\,n$}},\\
{\dfrac { \left(-2\right)^{n} \left( 5/12,n \right) \left( 3/4,n
 \right) }{{3}^{3/2\,n} \left( 5/6,n \right)  \left( 1/3,n \right) }}
\quad {\text {if $a=-1/3-2\,n$}},\\
{\dfrac {  -\left(-2\right)^{n-2}\left(\sqrt{3}-1\right) \left( 11/12,n \right) \left( 5/4,n
 \right) }{{3}^{3/2\,n} \left( 4/3,n \right)  \left( 5/6,n \right) }}
\quad {\text {if $a=-4/3-2\,n$}}.\\
\end{cases}&
\end{flalign*}
\begin{flalign*}
&{\rm (iv)}\,
{\mbox{$F$}(2\,a,2\,a-1/3;\,3\,a;\,9/4+3/4\,\sqrt {3})}&\\
&=\begin{cases}
{\dfrac { 7 \left( \sqrt {3}+1 \right) ^{6\,n+6}  \left( 5/4,n \right) \left( 19/12,n
 \right) }{\left(-2\right)^{5\,n+9}{3}^{3/2\,n+1/2} \left( 5/3,n \right)  \left( 7/6,n \right) }}
\quad {\text {if $a=-2-2\,n$}},\\
{\dfrac { \left( \sqrt {3}+1 \right) ^{6\,n+2}  \left( 1/2,n \right) \left( 5/6,n
 \right) }{\left(-2\right)^{5\,n+1}{3}^{3/2\,n+1/2} \left( 11/12,n \right)  \left( 5/12,n \right) }}
\quad {\text {if $a=-1/2-2\,n$}},\\
{\dfrac { \left( \sqrt {3}+1 \right) ^{6\,n+4}  \left( 3/4,n \right) \left( 13/12,n
 \right) }{\left(-2\right)^{5\,n+4}{3}^{3/2\,n+1/2} \left( 7/6,n \right)  \left( 2/3,n \right) }}
\quad {\text {if $a=-1-2\,n$}},\\
0
\quad {\text {if $a=-3/2-2\,n$}},\\
{\dfrac { \left( \sqrt {3}+1 \right) ^{6\,n}  \left( 1/6,n \right) \left( 1/2,n
 \right) }{\left(-2\right)^{5\,n}{3}^{3/2\,n} \left( 7/12,n \right)  \left( 1/12,n \right) }}
\quad {\text {if $a=1/6-2\,n$}},\\
{\dfrac { -\left( \sqrt {3}+1 \right) ^{6\,n+7}  \left( 17/12,n \right) \left( 7/4,n
 \right) }{\left(-2\right)^{5\,n+9}{3}^{3/2\,n-1/2} \left( 11/6,n \right)  \left( 4/3,n \right) }}
\quad {\text {if $a=-7/3-2\,n$}},\\
0
\quad {\text {if $a=-5/6-2\,n$}},\\
{\dfrac { \left( \sqrt {3}+1 \right) ^{6\,n+5}  \left( 11/12,n \right) \left( 5/4,n
 \right) }{\left(-2\right)^{5\,n+6}{3}^{3/2\,n} \left( 4/3,n \right)  \left( 5/6,n \right) }}
\quad {\text {if $a=-4/3-2\,n$}}.\\
\end{cases}&
\end{flalign*}
\begin{flalign*}
&{\rm (v)}\,
{\mbox{$F$}(a,2\,a-1/3;\,2/3;\,10-6\,\sqrt {3})}
=\dfrac{{3}^{3/8-9/4\,a} \left( \sqrt {3}-1 \right) ^{1/2-3\,a}\sqrt {\pi }
\Gamma  \left( 2/3 \right) }{2^{1/4-3/2\,a} \Gamma  \left( 3/4-1/2\,a \right) 
 \Gamma  \left( 1/2\,a+5/12 \right) }
&
\end{flalign*}
(The above is a generalization of Theorem 26 in [Ek]).
We derive (v) by using (131) in [Gour]
and connection formulae for the hypergeometric series(see (25)--(44) in 2.9 in [Erd]).
We find {\rm (vi)}$\leq${\rm (v)}.
\begin{flalign*}
&{\rm (vii)}\,
{\mbox{$F$}(a,1-2\,a;\,2/3;\,2/3-2/9\,\sqrt {3})}
=\dfrac{{3}^{3/8-3/4\,a} \left( \sqrt {3}-1 \right) ^{1/2-a}\sqrt {\pi }
\Gamma  \left( 2/3 \right) }{ 2^{1/4-1/2\,a}\Gamma  \left( 3/4-1/2\,a \right) 
 \Gamma  \left( 1/2\,a+5/12 \right) }
&
\end{flalign*}
(The above is a generalization of Theorem 40 in [Ek]).
We derive (vii) by using (130) in [Gour]
and connection formulae for the hypergeometric series(see (25)--(44) in 2.9 in [Erd]).
We find {\rm (viii)}$\leq${\rm (vii)}.
\begin{flalign*}
&{\rm (ix)}\,
{\mbox{$F$}(a,1-2\,a;\,4/3-a;\,1/3+2/9\,\sqrt {3})}&\\
&=
\dfrac{
{3}^{-3/4\,a} \left( \sqrt {3}-1 \right) ^{-a}\sin \left( \pi 
 \left( 5/12-1/2\,a \right)  \right) \Gamma  \left( 3/4
 \right) \Gamma \left(5/12\right) \Gamma 
 \left( 4/3-a \right) }{
  {2}^{1/2-7/2\,a}   \Gamma  \left( 2/3 \right)  \Gamma 
 \left( 13/12-1/2\,a \right)  \Gamma 
 \left( 3/4-1/2\,a \right)  }.
&
\end{flalign*}
\begin{flalign*}
&{\rm (x)}\,
{\mbox{$F$}(4/3-2\,a,a+1/3;\,4/3-a;\,1/3+2/9\,\sqrt {3})}&\\
&=
\dfrac{
{3}^{1/2-3/4\,a} \left( \sqrt {3}-1 \right) ^{-a-1/3}\sin \left( \pi 
 \left( 5/12-1/2\,a \right)  \right) \Gamma  \left( 3/4
 \right) \Gamma \left(5/12\right) \Gamma 
 \left( 4/3-a \right) }{
  {2}^{5/6-7/2\,a}  \Gamma  \left( 2/3 \right)  \Gamma 
 \left( 13/12-1/2\,a \right)  \Gamma 
 \left( 3/4-1/2\,a \right)  }.
&
\end{flalign*}
We derive {\rm (ix) and (x)} using formula (129) in [Gour]
and connection formulae for the hypergeometric series (see (25)--(44) in 2.9 in [Erd]).
\begin{flalign*}
&{\rm (xi)}\,
{\mbox{$F$}(a,a+1/3;\,4/3-a;\,-5/4-3/4\,\sqrt {3})}&\\
&=\begin{cases}
\dfrac{{3}^{9/2\,n}  \left( 7/6,n \right)
 \left( 2/3,n \right) }{ {\left(-2\right)}^{7\,n} \left( 13/12,n \right)    \left( 3/4,n \right) } 
\quad {\text {if $a=-2\,n$}},\\
\dfrac{-{3}^{9/2\,n+3/2} \left( \sqrt {3}-1 \right)  \left( 5/3,n \right)
 \left( 7/6,n \right) }{ 7 \left(-2\right)^{7\,n+1} \left( 19/12,n \right)    \left( 5/4,n \right) } 
\quad {\text {if $a=-1-2\,n$}},\\
\dfrac{{3}^{9/2\,n}  \left( 4/3,n \right)
 \left( 5/6,n \right) }{ \left(-2\right)^{7\,n} \left( 5/4,n \right)    \left( 11/12,n \right) } 
\quad {\text {if $a=-1/3-2\,n$}},\\
\dfrac{{3}^{9/2\,n+1} \left( \sqrt {3}-1 \right) \left(11/6,n \right)
 \left( 4/3,n \right) }{ \left(-2\right)^{7\,n+3} \left( 7/4,n \right)    \left( 17/12,n \right) } 
\quad {\text {if $a=-4/3-2\,n$}}.\\
\end{cases}&
\end{flalign*}
\begin{flalign*}
&{\rm (xii)}\,
{\mbox{$F$}(1-2\,a,4/3-2\,a;\,4/3-a;\,-5/4-3/4\,\sqrt {3})}&\\
&=\begin{cases}
{\dfrac {{3}^{9/2\,n} \left( \sqrt {3}+1 \right) ^{6\,n} \left( 1/6,n \right)  \left( 1/2,n \right) }{\left(-2\right)^{5
\,n}\left( 1/12,n \right) \left( 7/12
,n \right) }}
\quad {\text {if $a=1/2+2\,n$}},\\
{\dfrac {{3}^{9/2\,n+1}\left(\sqrt{3}+1\right)^{6\,n+1}  \left( 5/12,n \right)  \left( 3/4,n \right) }{\left(-2\right)^{5
\,n+1}\left( 1/3,n \right) \left( 5/6
,n \right) }}
\quad {\text {if $a=1+2\,n$}},\\
0
\quad {\text {if $a=3/2+2\,n$}},\\
{\dfrac {-{3}^{9/2\,n+4}\left(\sqrt{3}+1\right)^{6\,n+5} \left( 11/12,n \right)  \left( 5/4,n \right) }{\left(-2\right)^{5
\,n+6}\left( 5/6,n \right) \left( 4/3
,n \right) }}
\quad {\text {if $a=2+2\,n$}},\\
{\dfrac {{3}^{9/2\,n} \left( \sqrt {3}+1 \right) ^{6\,n} \left( 1/4,n \right)  \left( 7/12,n \right) }{\left(-2\right)^{5
\,n}\left( 1/6,n \right) \left( 2/3
,n \right) }}
\quad {\text {if $a=2/3+2\,n$}},\\
{\dfrac {{3}^{9/2\,n+3/2} \left( \sqrt {3}+1 \right) ^{6\,n+2} \left( 1/2,n \right)  \left( 5/6,n \right) }{\left(-2\right)^{5
\,n+1}\left( 5/12,n \right) \left(11/12
,n \right) }}
\quad {\text {if $a=7/6+2\,n$}},\\
{\dfrac {-{3}^{9/2\,n+5/2} \left( \sqrt {3}+1 \right) ^{6\,n+4} \left( 3/4,n \right)  \left( 13/12,n \right) }{\left(-2\right)^{5
\,n+4}\left( 2/3,n \right) \left(7/6
,n \right) }}
\quad {\text {if $a=5/3+2\,n$}},\\
0
\quad {\text {if $a=13/6+2\,n$}}.
\end{cases}&
\end{flalign*}
The special values obtained from {\rm (xiii)-(xx)}
are contained in the above.
\begin{flalign*}
&{\rm (xxi)}\,
{\mbox{$F$}(2\,a,a+1/3;\,4/3;\,10-6\,\sqrt {3})} &\\
&=\dfrac{{3}^{-9/4\,a-5/8} \left( \sqrt {3}-1 \right)^{-3\,a-1/2} {\sqrt{\pi} }\Gamma\left(1/3\right)}{ 2^{5/4-3/2\,a}\, \Gamma  \left( 13/12-1/2\,a
 \right) \Gamma  \left( 1/2\,a+3/4 \right) 
}.
\end{flalign*}
(The above is a generalization of Theorem 27 in [Ek]). 
We derive {\rm (xxi)} using formula (133) in [Gour]
and connection formulae for the hypergeometric series (see (25)--(44) in 2.9 in [Erd]).
We find {\rm (xxii)}$\leq$ {\rm (xxi)}.
\begin{flalign*}
&{\rm (xxiii)}\,
{\mbox{$F$}(2\,a,1-a;\,4/3;\,2/3-2/9\,\sqrt {3})}&\\
&=\dfrac{{3}^{3/4\,a-5/8} \left( \sqrt {3}-1 \right)^{a-1/2} {\sqrt{\pi} }\Gamma\left(1/3\right)}{ 2^{1/2\,a+5/4}\, \Gamma  \left( 13/12-1/2\,a
 \right) \Gamma  \left( 1/2\,a+3/4 \right) 
}
\end{flalign*}
(The above is a generalization of Theorem 39 in [Ek]).
We derive {\rm (xxiii)} using formula (132) in [Gour]
and connection formulae for the hypergeometric series (see (25)--(44) in 2.9 in [Erd]).
We find {\rm (xxiv)}$\leq$ {\rm (xxiii)}.
\paragraph{(2,4,6-4)}
\begin{flalign*}
&{\rm (i)}\,
{\mbox{$F$}(a,2\,a-1/3;\,3\,a;\,-9-6\,\sqrt {3})}&\\
&=\begin{cases}
{\dfrac {7\cdot {2}^{n-3} \left( \sqrt {3}+1 \right) ^{6\,n+6} \left( 5/4,n \right) \left( 19/12,n
 \right) }{{3}^{3/2\,n+1/2} \left( 5/3,n \right)  \left( 7/6,n \right) }}
\quad {\text {if $a=-2-2\,n$}},\\
{\dfrac {{2}^{n-1} \left( \sqrt {3}+1 \right) ^{6\,n+4} \left( 3/4,n \right) \left( 13/12,n
 \right) }{{3}^{3/2\,n+1/2} \left( 7/6,n \right)  \left( 2/3,n \right) }}
\quad {\text {if $a=-1-2\,n$}},\\
{\dfrac {{2}^{n} \left( \sqrt {3}+1 \right) ^{6\,n}
 \left( 1/6,n \right) \left( 1/2,n \right) }{{3}^{3/2
\,n}\left( 7/12,n \right)  \left( 1/12
,n \right) }}
\quad {\text {if $a=1/6-2\,n$}},\\
{\dfrac {{2}^{n-2} \left( \sqrt {3}+1 \right) ^{6\,n+8}
 \left( 17/12,n \right) \left( 7/4,n \right) }{{3}^{3/2
\,n-1/2}\left(11/6,n \right)  \left( 4/3
,n \right) }}
\quad {\text {if $a=-7/3-2\,n$}}\\
0
\quad {\text {if $a=-5/6-2\,n$}},\\
{\dfrac {-{2}^{n-1} \left( \sqrt {3}+1 \right) ^{6\,n+4}
 \left( 11/12,n \right) \left( 5/4,n \right) }{{3}^{3/2
\,n}\left( 4/3,n \right)  \left( 5/6
,n \right) }}
\quad {\text {if $a=-4/3-2\,n$}}.
\end{cases}&
\end{flalign*}
\begin{flalign*}
&{\rm (ii)}\,
{\mbox{$F$}(2\,a,a+1/3;\,3\,a;\,-9-6\,\sqrt {3})}&\\
&=\begin{cases}
{\dfrac {7\cdot {2}^{n-3} \left( \sqrt {3}+1 \right) ^{6\,n+6} \left( 5/4,n \right) \left( 19/12,n
 \right) }{{3}^{3/2\,n+1/2} \left( 5/3,n \right)  \left( 7/6,n \right) }}
\quad {\text {if $a=-2-2\,n$}},\\
{\dfrac { {2}^{n+1} \left( \sqrt {3}+1 \right) ^{6\,n+1} \left( 1/2,n \right) \left( 5/6,n
 \right) }{{3}^{3/2\,n+1/2} \left( 11/12,n \right)  \left( 5/12,n \right) }}
\quad {\text {if $a=-1/2-2\,n$}},\\
{\dfrac {-{2}^{n} \left( \sqrt {3}+1 \right) ^{6\,n+2} \left( 3/4,n \right) \left( 13/12,n
 \right) }{{3}^{3/2\,n+1/2} \left( 7/6,n \right)  \left( 2/3,n \right) }}
\quad {\text {if $a=-1-2\,n$}},\\
0
\quad {\text {if $a=-3/2-2\,n$}},\\
{\dfrac {{2}^{n} \left( \sqrt {3}+1 \right) ^{6\,n} \left( 5/12,n \right) \left( 3/4,n
 \right) }{{3}^{3/2\,n} \left( 5/6,n \right)  \left( 1/3,n \right) }}
\quad {\text {if $a=-1/3-2\,n$}},\\
{\dfrac {{2}^{n-2} \left( \sqrt {3}+1 \right) ^{6\,n+4} \left( 11/12,n \right) \left( 5/4,n
 \right) }{{3}^{3/2\,n} \left( 4/3,n \right)  \left( 5/6,n \right) }}
\quad {\text {if $a=-4/3-2\,n$}}.
\end{cases}&
\end{flalign*}
\begin{flalign*}
&{\rm (iii)}\,
{\mbox{$F$}(a,a+1/3;\,3\,a;\,9/4-3/4\,\sqrt {3})}&\\
&=\begin{cases}
{\dfrac {{3}^{3/4\,a-1/8} \left( \sqrt {3}+1 \right) ^{1/2}\sqrt {\pi }\Gamma  \left( a+2/3
 \right) }{{2}^{3/2\,a-1/4}\Gamma  \left( 1/2\,a+5/12 \right) \Gamma 
 \left( 1/2\,a+3/4 \right) }},\\
{\dfrac { 7 \cdot 2^{n-3}\left( 5/4,n \right) \left( 19/12,n
 \right) }{{3}^{3/2\,n+1/2} \left( 5/3,n \right)  \left( 7/6,n \right) }}
\quad {\text {if $a=-2-2\,n$}},\\
{\dfrac { 2^{n-1}\left(\sqrt{3}+1\right)\left( 3/4,n \right) \left( 13/12,n
 \right) }{{3}^{3/2\,n+1/2} \left( 7/6,n \right)  \left( 2/3,n \right) }}
\quad {\text {if $a=-1-2\,n$}},\\
{\dfrac { 2^{n} \left( 5/12,n \right) \left( 3/4,n
 \right) }{{3}^{3/2\,n} \left( 5/6,n \right)  \left( 1/3,n \right) }}
\quad {\text {if $a=-1/3-2\,n$}},\\
{\dfrac { 2^{n-2} \left( \sqrt {3}+1 \right) \left( 11/12,n \right) \left( 5/4,n
 \right) }{{3}^{3/2\,n} \left( 4/3,n \right)  \left( 5/6,n \right) }}
\quad {\text {if $a=-4/3-2\,n$}}.\\
\end{cases}
\end{flalign*}
\begin{flalign*}
&{\rm (iv)}\,
{\mbox{$F$}(2\,a,2\,a-1/3;\,3\,a;\,9/4-3/4\,\sqrt {3})}&\\
&=\begin{cases}
{\dfrac {{3}^{3/4\,a-1/8} \left( \sqrt {3}+1 \right) ^{3\,a-1/2}\sqrt {\pi }\Gamma  \left( a+2/3
 \right) }{{2}^{3/2\,a-1/4}\Gamma  \left( 1/2\,a+5/12 \right) \Gamma 
 \left( 1/2\,a+3/4 \right) }},\\
{\dfrac { 7 \left( \sqrt {3}-1 \right) ^{6\,n+6}  \left( 5/4,n \right) \left( 19/12,n
 \right) }{2^{5\,n+9}{3}^{3/2\,n+1/2} \left( 5/3,n \right)  \left( 7/6,n \right) }}
\quad {\text {if $a=-2-2\,n$}},\\
{\dfrac {- \left( \sqrt {3}-1 \right) ^{6\,n+4}  \left( 3/4,n \right) \left( 13/12,n
 \right) }{2^{5\,n+4}{3}^{3/2\,n+1/2} \left( 7/6,n \right)  \left( 2/3,n \right) }}
\quad {\text {if $a=-1-2\,n$}},\\
{\dfrac { \left( \sqrt {3}-1 \right) ^{6\,n+7}  \left( 17/12,n \right) \left( 7/4,n
 \right) }{2^{5\,n+9}{3}^{3/2\,n-1/2} \left( 11/6,n \right)  \left( 4/3,n \right) }}
\quad {\text {if $a=-7/3-2\,n$}},\\
{\dfrac { -\left( \sqrt {3}-1 \right) ^{6\,n+5}  \left( 11/12,n \right) \left( 5/4,n
 \right) }{2^{5\,n+6}{3}^{3/2\,n} \left( 4/3,n \right)  \left( 5/6,n \right) }}
\quad {\text {if $a=-4/3-2\,n$}}.\\
\end{cases}&
\end{flalign*}
\begin{flalign*}
&{\rm (v)}\,
{\mbox{$F$}(a,2\,a-1/3;\,2/3;\,10+6\,\sqrt {3})}&\\
&=\begin{cases}
{\dfrac {{3}^{9/2\,n} \left( \sqrt {3}+1 \right) ^{6\,n} \left( 7/12,n \right) }{{2}^{3\,n} \left( 3/4,n
 \right) }}
\quad {\text {if $a=-2\,n$}},\\
{\dfrac {{3}^{9/2\,n+3/2} \left( \sqrt {3}+1 \right) ^{6\,n+4} \left( 13/12,n \right) }{{2}^{3\,n+2} \left( 5/4,n
 \right) }}
\quad {\text {if $a=-1-2\,n$}},\\
{\dfrac {{3}^{9/2\,n} \left( \sqrt {3}+1 \right) ^{6\,n} \left( 1/2,n \right) }{{2}^{3\,n} \left( 2/3,n
 \right) }}
\quad {\text {if $a=1/6-2\,n$}},\\
{\dfrac {{3}^{9/2\,n+1} \left( \sqrt {3}+1 \right) ^{6\,n+2} \left( 3/4,n \right) }{{2}^{3\,n+1} \left( 11/12,n
 \right) }}
\quad {\text {if $a=-1/3-2\,n$}},\\
0
\quad {\text {if $a=-5/6-2\,n$}},\\
{\dfrac {-{3}^{9/2\,n+4} \left( \sqrt {3}+1 \right) ^{6\,n+4} \left( 5/4,n \right) }{5\cdot {2}^{3\,n+2} \left( 17/12,n
 \right) }}
\quad {\text {if $a=-4/3-2\,n$}}.
\end{cases}&
\end{flalign*}
(The first case is identical to Theorem 26 in [Ek]). We find {\rm (vi)}$\leq${\rm (v)}. 
\begin{flalign*}
&{\rm (vii)}\,
{\mbox{$F$}(a,1-2\,a;\,2/3;\,2/3+2/9\,\sqrt {3})}&\\
&=\begin{cases}
{\dfrac {{3}^{3/2\,n} \left( \sqrt {3}+1 \right) ^{2\,n} \left( 7/12,n \right) }{{2}^{n} \left( 3/4,n
 \right) }}
\quad {\text {if $a=-2\,n$}},\\
{\dfrac {-{3}^{3/2\,n} \left( \sqrt {3}+1 \right) ^{2\,n+2} \left( 13/12,n \right) }{{2}^{n+1} \left( 5/4,n
 \right) }}
\quad {\text {if $a=-1-2\,n$}},\\
{\dfrac { \left( \sqrt {3}-1 \right) ^{2\,n} \left( 1/2
,n \right) }{{2}^{n}{3}^{3/2\,n} \left( 2/3,n \right) }}
\quad {\text {if $a=1/2+2\,n$}},\\
{\dfrac { -\left( \sqrt {3}-1 \right) ^{2\,n} \left( 3/4
,n \right) }{{2}^{n}{3}^{3/2\,n+1/2} \left( 11/12,n \right) }}
\quad {\text {if $a=1+2\,n$}},\\
0
\quad {\text {if $a=3/2+2\,n$}},\\
{\dfrac { \left( \sqrt {3}-1 \right) ^{2\,n+2} \left( 5/4
,n \right) }{5\cdot {2}^{n+1}{3}^{3/2\,n+1/2} \left( 17/12,n \right) }}
\quad {\text {if $a=2+2\,n$}}
\end{cases}&
\end{flalign*}
(The first case is identical to Theorem 40 in [Ek]).
We find {\rm (viii)}$\leq${\rm (vii)}.
\begin{flalign*}
&{\rm (ix)}\,
{\mbox{$F$}(a,1-2\,a;\,4/3-a;\,1/3-2/9\,\sqrt {3})}&\\
&={\frac {{3}^{-3/4\,a} \left( \sqrt {3}+1 \right) ^{-a}\Gamma 
 \left( 3/4 \right) \Gamma \left( 13/12  \right) \Gamma  \left( 2/3-1/2
\,a \right) \Gamma  \left( 7/6-1/2\,a \right) }
{{2}^{-5/2\,a}\Gamma \left(2/3\right)
\Gamma  \left( 7/6 \right) \Gamma  \left( 3/4-1/2\,a \right) \Gamma 
 \left(13/12 -1/2\,a\right) }}.
\end{flalign*}
\begin{flalign*}
&{\rm (x)}\,
{\mbox{$F$}(a+1/3,4/3-2\,a;\,4/3-a;\,1/3-2/9\,\sqrt {3})}&\\
&={\frac {{3}^{1/2-3/4\,a} \left( \sqrt {3}+1 \right) ^{-a-1/3}\Gamma 
 \left( 3/4 \right) \Gamma \left( 13/12  \right) \Gamma  \left( 2/3-1/2
\,a \right) \Gamma  \left( 7/6-1/2\,a \right) }
{{2}^{1/3-5/2\,a}\Gamma \left(2/3\right)
\Gamma  \left( 7/6 \right) \Gamma  \left( 3/4-1/2\,a \right) \Gamma 
 \left(13/12 -1/2\,a\right) }}.
\end{flalign*}
\begin{flalign*}
&{\rm (xi)}\,
{\mbox{$F$}(a,a+1/3;\,4/3-a;\,-5/4+3/4\,\sqrt {3})}&\\
&={\frac {{3}^{-9/4\,a} \Gamma 
 \left( 3/4 \right) \Gamma \left( 13/12  \right) \Gamma  \left( 2/3-1/2
\,a \right) \Gamma  \left( 7/6-1/2\,a \right) }
{{2}^{-7/2\,a}\Gamma \left(2/3\right)
\Gamma  \left( 7/6 \right) \Gamma  \left( 3/4-1/2\,a \right) \Gamma 
 \left(13/12 -1/2\,a\right) }}.
\end{flalign*}
\begin{flalign*}
&{\rm (xii)}\,
{\mbox{$F$}(1-2\,a,4/3-2\,a;\,4/3-a;\,-5/4+3/4\,\sqrt {3})}&\\
&={\frac {{3}^{9/4\,a-3/2} \left( \sqrt {3}+1 \right) ^{1-3\,a}\Gamma 
 \left( 3/4 \right) \Gamma \left( 13/12  \right) \Gamma  \left( 2/3-1/2
\,a \right) \Gamma  \left( 7/6-1/2\,a \right) }
{{2}^{-1/2\,a-1}\Gamma \left(2/3\right)
\Gamma  \left( 7/6 \right) \Gamma  \left( 3/4-1/2\,a \right) \Gamma 
 \left(13/12 -1/2\,a\right) }}.
\end{flalign*}
The special values obtained from {\rm (xiii)-(xx)}
are contained in the above.
\begin{flalign*}
&{\rm (xxi)}\,
{\mbox{$F$}(2\,a,a+1/3;\,4/3;\,10+6\,\sqrt {3})}&\\
&=\begin{cases}
{\dfrac {{3}^{9/2\,n} \left( \sqrt {3}+1 \right) ^{6\,n} \left( 1/4,n \right) }{{2}^{3\,n}\left( 13/12,n
 \right) }}
\quad {\text {if $a=-2\,n$}},\\
{\dfrac {{3}^{9/2\,n+3/2} \left( \sqrt {3}+1 \right) ^{6\,n+1} \left( 1/2,n \right) }{{2}^{3\,n+2}\left( 4/3,n
 \right) }}
\quad {\text {if $a=-1/2-2\,n$}},\\
{\dfrac {-{3}^{9/2\,n+5/2} \left( \sqrt {3}+1 \right) ^{6\,n+2} \left( 3/4,n \right) }{7\cdot {2}^{3\,n+1}\left( 19/12,n
 \right) }}
\quad {\text {if $a=-1-2\,n$}},\\
0
\quad {\text {if $a=-3/2-2\,n$}},\\
{\dfrac {{3}^{9/2\,n} \left( \sqrt {3}+1 \right) ^{6\,n} \left( 5/12,n \right) }{{2}^{3\,n}\left( 5/4,n
 \right) }}
\quad {\text {if $a=-1/3-2\,n$}},\\
{\dfrac {{3}^{9/2\,n+1} \left( \sqrt {3}+1 \right) ^{6\,n+4} \left( 11/12,n \right) }{{2}^{3\,n+2}\left( 7/4,n
 \right) }}
\quad {\text {if $a=-4/3-2\,n$}}
\end{cases}&
\end{flalign*}
(The fifth case is identical to Theorem 27 in [Ek]).
We find {\rm (xxii)}$\leq${\rm (xxi)}.
\begin{flalign*}
&{\rm (xxiii)}\,
{\mbox{$F$}(2\,a,1-a;\,4/3;\,2/3+2/9\,\sqrt {3})}&\\
&=\begin{cases}
{\dfrac { \left( \sqrt {3}-1 \right) ^{2\,n} \left( 1/4
,n \right) }{{2}^{n}{3}^{3/2\,n} \left( 13/12,n
 \right) }}
\quad {\text {if $a=-2\,n$}},\\
{\dfrac { -\left( \sqrt {3}-1 \right) ^{2\,n+1} \left( 1/2
,n \right) }{{2}^{n+2}{3}^{3/2\,n} \left( 4/3,n
 \right) }}
\quad {\text {if $a=-1/2-2\,n$}},\\
{\dfrac { -\left( \sqrt {3}-1 \right) ^{2\,n+2} \left( 3/4
,n \right) }{7\cdot {2}^{n+1}{3}^{3/2\,n+1/2} \left( 19/12,n
 \right) }}
\quad {\text {if $a=-1-2\,n$}},\\
0
\quad {\text {if $a=-3/2-2\,n$}},\\
{\dfrac {{3}^{3/2\,n} \left( \sqrt {3}+1 \right) ^{2\,n} \left( 5/12,n \right) }{{2}^{n} \left( 5/4,n
 \right) }}
\quad {\text {if $a=1+2\,n$}},\\
{\dfrac {-{3}^{3/2\,n-1/2} \left( \sqrt {3}+1 \right) ^{2\,n+2} \left( 11/12,n \right) }{{2}^{n+1} \left( 7/4,n
 \right) }}
\quad {\text {if $a=2+2\,n$}}
\end{cases}&
\end{flalign*}
(The fifth case is identical to Theorem 39 in [Ek]).
We find {\rm (xxiv)}$\leq${\rm (xxiii)}.
\paragraph{(2,4,6-5)}
The special values obtained from (2,4,6-5) coincide with
those (2,4,6-3).
\paragraph{(2,4,6-6)}
The special values obtained from (2,4,6-6) coincide with
those (2,4,6-4).

\subsubsection{$(k,l,m)=(2,5,6)$}
In this case, there is no admissible quadruple.

\subsubsection{$(k,l,m)=(2,6,6)$}
In this case, we get
\begin{align}
&(a,b,c,x)=(a,3\,a-1/2,3\,a,-3),
\tag{2,6,6-1}\\
&(a,b,c,x)=(a,3\,a-3/2,3\,a-1,-3),
\tag{2,6,6-2}\\
&(a,b,c,x)=(a,3\,a+1,3\,a,3/2),
\tag{2,6,6-3}\\
&(a,b,c,x)=(a,3\,a-3,3\,a-1,3/2).
\tag{2,6,6-4}
\end{align}
\paragraph{(2,6,6-1), (2,6,6-2)}
The special values obtained from (2,6,6-1) and (2,6,6-2)
coincide those obtained from (1,3,3-1).
\paragraph{(2,6,6-3), (2,6,6-4)}
The special values obtained from (2,6,6-3) and (2,6,6-4) 
coincide those obtained from (1,3,3-3).

\subsubsection{$(k,l,m)=(2,7,6)$}
In this case, there is no admissible quadruple.

\subsubsection{$(k,l,m)=(2,8,6)$}
In this case, we have
\begin{align}
&(a,b,c,x)=(a,b,b+1-a,-1),
\tag{2,8,6-1}\\
&(a,b,c,x)=(a,4\,a-1/2,3\,a,-1),
\tag{2,8,6-2}&\\
&(a,b,c,x)=(a,4\,a-5/2,3\,a-1,-1).
\tag{2,8,6-3}&
\end{align}
\paragraph{(2,8,6-1)}
The special values obtained from (2,8,6-1) are evaluated in 
paragraphs (1,2,2-1) and (0,2,2-1).
\paragraph{(2,8,6-2), (2,8,6-3)}
The special values obtained from (2,8,6-2) and (2,8,6-3) coincide with those obtained from (1,4,3-2).

\subsubsection{$(k,l,m)=(3,3,6)$}
In this case, we have
\begin{align}
&(a,b,c,x)=(a,a-1/2,2\,a,4),\,
S^{(n)}=1,
\tag{3,3,6-1}\\
&(a,b,c,x)=(a,a+1/2,2\,a,4),\,
S^{(n)}=1,
\tag{3,3,6-2}\\
&(a,b,c,x)=(a,a-1/2,2\,a-1,4),\,
S^{(n)}=1,
\tag{3,3,6-3}\\
&(a,b,c,x)=(a,a+1/2,2\,a+1,4),\,
S^{(n)}=1,
\tag{3,3,6-4}\\
&(a,b,c,x)=(a,a+3/2,2\,a-1,4),\,
S^{(n)}=\frac{\left(2\,a-3\right)\left(2\,a+1+6\,n\right)}{\left(2\,a+1\right)\left(2\,a-3+6\,n\right)},
\tag{3,3,6-5}\\
&(a,b,c,x)=(a,a-3/2,2\,a+1,4),\,
S^{(n)}=\frac{\left(2\,a-1\right)\left(2\,a+3+6\,n\right)}{\left(2\,a+3\right)\left(2\,a-1+6\,n\right)},
\tag{3,3,6-6}\\
&(a,b,c,x)=(a,a-3/2,2\,a-4,4),\,
S^{(n)}=\frac{\left(a-3\right)\left(a-1+3\,n\right)}{\left(a-1\right)\left(a-3+3\,n\right)},
\tag{3,3,6-7}\\
&(a,b,c,x)=(a,a+3/2,2\,a+4,4),\,
S^{(n)}=\frac{\left(a+1\right)\left(a+3+3\,n\right)}{\left(a+3\right)\left(a+1+3\,n\right)},
\tag{3,3,6-8}\\
&(a,b,c,x)=(a,a-1/2,2\,a,4/3),\,
S^{(n)}=\left(-3\right)^{-3\,n},
\tag{3,3,6-9}\\
&(a,b,c,x)=(a,a+1/2,2\,a,4/3),\,
S^{(n)}=\left(-3\right)^{-3\,n},
\tag{3,3,6-10}\\
&(a,b,c,x)=(a,a-1/2,2\,a-1,4/3),\,
S^{(n)}=\left(-3\right)^{-3\,n},
\tag{3,3,6-11}\\
&(a,b,c,x)=(a,a+1/2,2\,a+1,4/3),\,
S^{(n)}=\left(-3\right)^{-3\,n},
\tag{3,3,6-12}\\
&(a,b,c,x)=(a,a-5/2,2\,a-1,4/3),\,
S^{(n)}=\dfrac{\left(2\,a-3\right)\left(2\,a+1+6\,n\right)}{\left(-3\right)^{3\,n}\left(2\,a+1\right)\left(2\,a-3+6\,n\right)},
\tag{3,3,6-13}\\
&(a,b,c,x)=(a,a+5/2,2\,a+1,4/3),\,
S^{(n)}=\dfrac{\left(2\,a-1\right)\left(2\,a+3+6\,n\right)}{\left(-3\right)^{3\,n}\left(2\,a+3\right)\left(2\,a-1+6\,n\right)},
\tag{3,3,6-14}\\
&(a,b,c,x)=(a,a-5/2,2\,a-4,4/3),\,
S^{(n)}=\frac{\left(a-3\right)\left(a-1+3\,n\right)}{\left(-3\right)^{3\,n}\left(a-1\right)\left(a-3+3\,n\right)},
\tag{3,3,6-15}\\
&(a,b,c,x)=(a,a+5/2,2\,a+4,4/3),\,
S^{(n)}=\frac{\left(a+1\right)\left(a+3+3\,n\right)}{\left(-3\right)^{3\,n}\left(a+3\right)\left(a+1+3\,n\right)}.
\tag{3,3,6-16}
\end{align}
\paragraph{(3,3,6-1)}
\begin{flalign*}
&{\rm (i)}\,{\mbox{$F$}(a,a-1/2;\,2\,a;\,4)}
=\begin{cases}
1
\quad {\text {if $a=-3-3\,n$}},\\
-2
\quad {\text {if $a=-1-3\,n$}},\\
1
\quad {\text {if $a=-2-3\,n$}},\\
2
\quad {\text {if $a=-5/2-3\,n$}},\\
-1
\quad {\text {if $a=-7/2-3\,n$}},\\
-1
\quad {\text {if $a=-3/2-3\,n$}}.
\end{cases}&
\end{flalign*}
\begin{flalign*}
&{\rm (ii)}\,{\mbox{$F$}(a,a+1/2;\,2\,a;\,4)}
=\begin{cases}
1
\quad {\text {if $a=-3-3\,n$}},\\
0
\quad {\text {if $a=-1-3\,n$}},\\
-1
\quad {\text {if $a=-2-3\,n$}},\\
1
\quad {\text {if $a=-1/2-3\,n$}},\\
-1
\quad {\text {if $a=-3/2-3\,n$}},\\
0
\quad {\text {if $a=-5/2-3\,n$}}.
\end{cases}&
\end{flalign*}
\begin{flalign*}
&{\rm (iii)}\,
{\mbox{$F$}(a,a+1/2;\,2\,a;\,4/3)}
=\begin{cases}
\left( -3 \right) ^{-3\,n-3}
\quad {\text {if $a=-3-3\,n$}},\\
-2\left( -3 \right) ^{-3\,n-1}
\quad {\text {if $a=-1-3\,n$}},\\
\left( -3 \right) ^{-3\,n-2}
\quad {\text {if $a=-2-3\,n$}},\\
\left( -3 \right) ^{-3\,n}
\quad {\text {if $a=-1/2-3\,n$}},\\
-\left( -3 \right) ^{-3\,n-1}
\quad {\text {if $a=-3/2-3\,n$}},\\
0
\quad {\text {if $a=-5/2-3\,n$}}.
\end{cases}&
\end{flalign*}
\begin{flalign*}
&{\rm (iv)}\,
{\mbox{$F$}(a,a-1/2;\,2\,a;\,4/3)}
=\begin{cases}
\left( -3 \right) ^{-3\,n-3}
\quad {\text {if $a=-3-3\,n$}},\\
0
\quad {\text {if $a=-1-3\,n$}},\\
-\left( -3 \right) ^{-3\,n-2}
\quad {\text {if $a=-2-3\,n$}},\\
2\left( -3 \right) ^{-3\,n-3}
\quad {\text {if $a=-5/2-3\,n$}},\\
-\left( -3 \right) ^{-3\,n-4}
\quad {\text {if $a=-7/2-3\,n$}},\\
-\left( -3 \right) ^{-3\,n-2}
\quad {\text {if $a=-3/2-3\,n$}}.
\end{cases}&
\end{flalign*}
\begin{flalign*}
&{\rm (v)}\,
{\mbox{$F$}(a,a-1/2;\,1/2;\,-3)}
=\begin{cases}
2^{6\,n}
\quad {\text {if $a=-3\,n$}},\\
-2^{6\,n+3}
\quad {\text {if $a=-1-3\,n$}},\\
2^{6\,n+4}
\quad {\text {if $a=-2-3\,n$}},\\
2^{6\,n}
\quad {\text {if $a=1/2-3\,n$}},\\
-2^{6\,n+1}
\quad {\text {if $a=-1/2-3\,n$}},\\
-2^{6\,n+3}
\quad {\text {if $a=-3/2-3\,n$}}.
\end{cases}&
\end{flalign*}
We find {\rm (vi)}$\leq${\rm (v)}.
\begin{flalign*}
&{\rm (vii)}\,
{\mbox{$F$}(a,1-a;\,1/2;\,3/4)}=
2\,\sin \left( 1/6 \left( 4\,a+1 \right) \pi  \right).&
\end{flalign*}
\begin{flalign*}
&{\rm (viii)}\,
{\mbox{$F$}(1/2-a,a-1/2;\,1/2;\,3/4)}
=\sin \left( 1/6 \left( 4\,a+1 \right) \pi  \right).& 
\end{flalign*}
\begin{flalign*}
&{\rm (ix)}\,
{\mbox{$F$}(a,1-a;\,3/2;\,1/4)}
={\frac {2\,\sin \left( 1/6 \left( 2\,a-1 \right) \pi  \right) }{2\,a
-1}}.&
\end{flalign*}
\begin{flalign*}
&{\rm (x)}\,
{\mbox{$F$}(3/2-a,a+1/2;\,3/2;\,1/4)}
={\frac {4\,\sin \left( 1/6 \left( 2\,a-1 \right) \pi  \right) }{\sqrt{3}\left(2\,a-1\right)}}.&
\end{flalign*}
\begin{flalign*}
&{\rm (xi)}\,
{\mbox{$F$}(a,a+1/2;\,3/2;\,-1/3)}
={\frac {{3}^{a}\sin \left( 1/6 \left( 2\,a-1 \right) \pi  \right) }{
{2}^{2\,a-1} \left( 2\,a-1 \right) }}.&
\end{flalign*}
We find {\rm (xii)}$\leq${\rm (xi)}.
\begin{flalign*}
&{\rm (xiii)}\,
{\mbox{$F$}(1/2-a,a-1/2;\,1/2;\,1/4)}
=\sin \left( 1/3 \left( 2-a \right) \pi  \right).&
\end{flalign*}
\begin{flalign*}
&{\rm (xiv)}\,
{\mbox{$F$}(a,1-a;\,1/2;\,1/4)}
=\frac{2}{\sqrt{3}}\sin \left( 1/3 \left( 2-a \right) \pi  \right).&
\end{flalign*}
\begin{flalign*}
&{\rm (xv)}\,
{\mbox{$F$}(a,a-1/2;\,1/2;\,-1/3)}
={\frac {{3}^{a-1/2}\sin \left( 1/3 \left( 2-a \right) \pi  \right) 
}{{2}^{2\,a-1}}}.&
\end{flalign*}
We find {\rm (xvi)}$\leq${\rm (xv)}, 
{\rm (xvii)}$\leq${\rm (i)},
{\rm (xviii)}$\leq${\rm (ii)},
{\rm (xix)}$\leq${\rm (iii)} and
{\rm (xx)}$\leq${\rm (iv)}.
\begin{flalign*}
&{\rm (xxi)}\,
{\mbox{$F$}(a,a+1/2;\,3/2;\,-3)}
=\begin{cases}
{\dfrac {{2}^{6\,n}}{6\,n+1}}
\quad {\text {if $a=-3\,n$}},\\
0
\quad {\text {if $a=-1-3\,n$}},\\
{\dfrac {-{2}^{6\,n+4}}{6\,n+5}}
\quad {\text {if $a=-2-3\,n$}},\\
{\dfrac {{2}^{6\,n}}{3\,n+1}}
\quad {\text {if $a=-1/2-3\,n$}},\\
{\dfrac {-{2}^{6\,n+2}}{3\,n+2}}
\quad {\text {if $a=-3/2-3\,n$}},\\
0
\quad {\text {if $a=-5/2-3\,n$}}.
\end{cases}&
\end{flalign*}
We find {\rm (xxii)}$\leq${\rm (xxi)}.
\begin{flalign*}
&{\rm (xxiii)}\,
{\mbox{$F$}(a,1-a;\,3/2;\,3/4)} 
={\frac {2\,\sin \left( 1/3 \left( 2\,a-1 \right) \pi 
 \right) }{\sqrt{3}\left(2\,a-1\right)}}.&
\end{flalign*}
\begin{flalign*}
&{\rm (xxiv)}\,
{\mbox{$F$}(3/2-a,a+1/2;\,3/2;\,3/4)}
={\frac {4\,\sin \left( 1/3 \left( 2\,a-1 \right) \pi 
 \right) }{\sqrt{3}\left(2\,a-1\right)}}.&
\end{flalign*}

\paragraph{(3,3,6-2), (3,3,6-3), (3,3,6-4), (3,3,6-9), (3,3,6-10), (3,3,6-11), (3,3,6-12)}
The special values obtained from (3,3,6-2), (3,3,6-3),  (3,3,6-4),
(3,3,6-9), (3,3,6-10), (3,3,6-11) and (3,3,6-12)
coincide with those obtained from (3,3,6-1).

\paragraph{(3,3,6-5)}
\begin{flalign*}
&{\rm (i)}\,
{\mbox{$F$}(a,a+3/2;\,2\,a-1;\,4)}
=\begin{cases}
-{\dfrac {2\,n+1}{6\,n-1}}\vspace{3pt}
\quad {\text {if $a=-3\,n$}},\\
{\dfrac {6\,n+5}{18\,n+3}}\vspace{3pt}
\quad {\text {if $a=-1-3\,n$}},\\
0
\quad {\text {if $a=-2-3\,n$}},\\
{\dfrac {n+1}{3\,n+1}}\vspace{3pt}
\quad {\text {if $a=-3/2-3\,n$}},\\
-\,{\dfrac {3\,n+4}{9\,n+6}}\vspace{3pt}
\quad {\text {if $a=-5/2-3\,n$}},\\
0
\quad {\text {if $a=-7/2-3\,n$}}.
\end{cases}&
\end{flalign*}
\begin{flalign*}
&{\rm (ii)}\,
{\mbox{$F$}(a-1,a-5/2;\,2\,a-1;\,4)}
=\begin{cases}
-{\dfrac {12\,n+14}{2\,n+1}}\vspace{3pt}
\quad {\text {if $a=-2-3\,n$}},\\
{\dfrac {18\,n+27}{6\,n+5}}\vspace{3pt}
\quad {\text {if $a=-3-3\,n$}},\\
{\dfrac {18\,n+15}{6\,n+1}}\vspace{3pt}
\quad {\text {if $a=-1-3\,n$}},\\
{\dfrac {6\,n+10}{n+1}}\vspace{3pt}
\quad {\text {if $a=-7/2-3\,n$}},\\
-{\dfrac {9\,n+18}{3\,n+4}}\vspace{3pt}
\quad {\text {if $a=-9/2-3\,n$}},\\
-{\dfrac {9\,n+12}{3\,n+2}}
\quad {\text {if $a=-5/2-3\,n$}}.
\end{cases}&
\end{flalign*}
\begin{flalign*}
&{\rm (iii)}\,
{\mbox{$F$}(a,a-5/2;\,2\,a-1;\,4/3)}
=\begin{cases}
-{\dfrac {2\,n+1}{ \left( -3 \right) ^{3\,n} \left( 6\,n-1 \right) }}
\quad {\text {if $a=-3\,n$}},\\ 
-{\dfrac {6\,n+5}{ \left( -3 \right) ^{3\,n+2} \left( 6\,n+1 \right) }}
\quad {\text {if $a=-1-3\,n$}},\\
0
\quad {\text {if $a=-2-3\,n$}},\\
{\dfrac {6\,n+10}{ \left( -3 \right) ^{3\,n+6} \left( n+1 \right) }}
\quad {\text {if $a=-7/2-3\,n$}},\\
{\dfrac {3\,n+6}{ \left( -3 \right) ^{3\,n+6} \left( 3\,n+4 \right) }}
\quad {\text {if $a=-9/2-3\,n$}},\\
{\dfrac {3\,n+4}{ \left( -3 \right) ^{3\,n+4} \left( 3\,n+2 \right) }}
\quad {\text {if $a=-5/2-3\,n$}}.
\end{cases}&
\end{flalign*}
\begin{flalign*}
&{\rm (iv)}\,
{\mbox{$F$}(a-1,a+3/2;\,2\,a-1;\,4/3)}
=\begin{cases}
-{\dfrac {12\,n+14}{ \left( -3 \right) ^{3\,n+3} \left( 2\,n+1 \right) 
}}
\quad {\text {if $a=-2-3\,n$}},\\ 
{\dfrac {2\,n+3}{ \left( -3 \right) ^{3\,n+2} \left( 6\,n+5 \right) }}
\quad {\text {if $a=-3-3\,n$}},\\ 
-{\dfrac {6\,n+5}{ \left( -3 \right) ^{3\,n+1} \left( 6\,n+1 \right) }}
\quad {\text {if $a=-1-3\,n$}},\\ 
{\dfrac {n+1}{ \left( -3 \right) ^{3\,n} \left( 3\,n+1 \right) }}
\quad {\text {if $a=-3/2-3\,n$}},\\ 
{\dfrac {3\,n+4}{ \left( -3 \right) ^{3\,n+2} \left( 3\,n+2 \right) }}
\quad {\text {if $a=-5/2-3\,n$}},\\ 
0
\quad {\text {if $a=-7/2-3\,n$}}.
\end{cases}&
\end{flalign*}
\begin{flalign*}
&{\rm (v)}\,{\mbox{$F$}(a,a+3/2;\,7/2;\,-3)}
=\begin{cases}
{\dfrac {-5\cdot {2}^{6\,n}}{ \left( 6\,n+5 \right)  \left( 3\,n+1 \right) 
 \left( 6\,n-1 \right) }}
\quad {\text {if $a=-3\,n$}},\\ 
{\dfrac {5\cdot {2}^{6\,n+2}}{ \left( 6\,n+7 \right)  \left( 3\,n+2
 \right)  \left( 6\,n+1 \right) }}
\quad {\text {if $a=-1-3\,n$}},\\
0
\quad {\text {if $a=-2-3\,n$}},\\
{\dfrac {5\cdot {2}^{6\,n+2}}{ \left( 3\,n+4 \right)  \left( 6\,n+5
 \right)  \left( 3\,n+1 \right) }}
\quad {\text {if $a=-3/2-3\,n$}},\\
{\dfrac {-5\cdot {2}^{6\,n+4}}{ \left( 3\,n+5 \right)  \left( 6\,n+7
 \right)  \left( 3\,n+2 \right) }}
\quad {\text {if $a=-5/2-3\,n$}},\\
0
\quad {\text {if $a=-7/2-3\,n$}}.\\
\end{cases}&
\end{flalign*}
We find {\rm (vi)}$\leq${\rm (v)}.
\begin{flalign*}
&{\rm (vii)}\,
{\mbox{$F$}(a,2-a;\,7/2;\,3/4)}
=
{\dfrac {10\sin \left( 1/3 \left( 2\,a+1 \right) \pi 
 \right) }{\sqrt{3} \left( a-1 \right)  \left( 2\,a+1 \right)  \left( 2\,a-5
 \right) }}.
&
\end{flalign*}
\begin{flalign*}
&{\rm (viii)}\,
{\mbox{$F$}(a+3/2,7/2-a;\,7/2;\,3/4)}
=
{\dfrac {80}{\sqrt{3}}}\,{\dfrac {\sin \left( 1/3\left( 2\,a+1
 \right) \pi  \right) }{ \left( 2\,a+1 \right)  \left( 2\,a-5 \right) 
 \left( a-1 \right) }}.
&
\end{flalign*}
\begin{flalign*}
&{\rm (ix)}\,
{\mbox{$F$}(a,2-a;\,-1/2;\,1/4)}
=\dfrac{-2}{3^{3/2}}\,{\dfrac { \left( 2\,a-1 \right)  \left( 2\,a-3 \right) \sin
 \left( 1/3 \left( a+2 \right) \pi  \right) }{a-1}}.
&
\end{flalign*}
\begin{flalign*}
&{\rm (x)}\,
{\mbox{$F$}(a-5/2,-a-1/2;\,-1/2;\,1/4)}
=\dfrac{-3}{16}\,{\dfrac { \left( 2\,a-1 \right)  \left( 2\,a-3 \right) \sin
 \left( 1/3 \left( a+2 \right) \pi  \right) }{a-1}}.
&
\end{flalign*}
\begin{flalign*}
&{\rm (xi)}\,
{\mbox{$F$}(a,a-5/2;\,-1/2;\,-1/3)}
=\dfrac{-3^{a-3/2}}{2^{2\,a-1}}\,{\dfrac { \left( 2\,a-1 \right)  \left( 2\,a-3 \right) \sin
 \left( 1/3 \left( a+2 \right) \pi  \right) }{a-1}}.
&
\end{flalign*}
We find {\rm (xii)}$\leq${\rm (xi)}. 
\begin{flalign*}
&{\rm (xiii)}\,
{\mbox{$F$}(a+3/2,7/2-a;\,5/2;\,1/4)}
={\dfrac {-32}{\sqrt{3}}}\,{\dfrac {\sin \left( 1/3 \left( a+1/2
 \right) \pi  \right) }{ \left( 2\,a-5 \right)  \left( 2\,a+1 \right) 
}}.
&
\end{flalign*}
\begin{flalign*}
&{\rm (xiv)}\,
{\mbox{$F$}(a-1,1-a;\,5/2;\,1/4)}
={\dfrac {-9\sin \left( 1/3 \left( a+1/2
 \right) \pi  \right) }{ \left( 2\,a-5 \right)  \left( 2\,a+1 \right) 
}}.
&
\end{flalign*}
\begin{flalign*}
&{\rm (xv)}\,
{\mbox{$F$}(a-1,a+3/2;\,5/2;\,-1/3)}
=\dfrac {-3^{a+1}}{2^{2\,a-2}}\,{\dfrac {\sin \left( 1/3 \left( a+1/2
 \right) \pi  \right) }{ \left( 2\,a-5 \right)  \left( 2\,a+1 \right) 
}}.
&
\end{flalign*}
We find {\rm (xvi)}$\leq${\rm (xv)}, 
{\rm (xvii)}$\leq${\rm (i)},
{\rm (xviii)}$\leq${\rm (ii)},  
{\rm (xix)}$\leq${\rm (iii)}, 
{\rm (xx)}$\leq${\rm (iv)}.
\begin{flalign*}
&{\rm (xxi)}\,
{\mbox{$F$}(a-1,a-5/2;\,-3/2;\,-3)} &\\
&=\begin{cases}
-{2}^{6\,n+6} \left( 6\,n+5 \right)  \left( 6\,n+7 \right) 
\quad {\text {if $a=-2-3\,n$}},\\
{2}^{6\,n+1} \left( 6\,n+1 \right)  \left( 6\,n+3 \right)
\quad {\text {if $a=-3\,n$}},\\
{2}^{6\,n+3} \left( 6\,n+3 \right)  \left( 6\,n+5 \right)
\quad {\text {if $a=-1-3\,n$}},\\
{2}^{6\,n-3} \left( 6\,n-2 \right)  \left( 6\,n-4 \right)
\quad {\text {if $a=5/2-3\,n$}},\\
-{2}^{6\,n+4} \left( 6\,n+4 \right)  \left( 6\,n+6 \right)
\quad {\text {if $a=-3/2-3\,n$}},\\
-{2}^{6\,n+6} \left( 6\,n+6 \right)  \left( 6\,n+8 \right)
\quad {\text {if $a=-5/2-3\,n$}}.
\end{cases}&
\end{flalign*}
We find {\rm (xxii)}$\leq${\rm (xxi)}.
\begin{flalign*}
&{\rm (xxiii)}\,
{\mbox{$F$}(a-1,1-a;\,-3/2;\,3/4)}
=- \left( 2\,a-1 \right)  \left( 2\,a-3 \right) \sin \left( 1/6
 \left( 4\,a-1 \right) \pi  \right).&
\end{flalign*}
\begin{flalign*}
&{\rm (xxiv)}\,
{\mbox{$F$}(a-5/2,-a-1/2;\,-3/2;\,3/4)}&\\
&=-1/8\,\left( 2\,a-1 \right)  \left( 2\,a-3 \right) \sin \left( 1/6
 \left( 4\,a-1 \right) \pi  \right).&
\end{flalign*}

\paragraph{(3,3,6-6), (3,3,6-7), (3,3,6-8), (3,3,6-13), (3,3,6-14), (3,3,6-15), (3,3,6-16)}
The special values obtained from (3,3,6-6), (3,3,6,-7), (3,3,6-8),
(3,3,6-13), (3,3,6-14), (3,3,6-15) and (3,3,6-16)
coincide with those obtained from (3,3,6,5).

\subsubsection{$(k,l,m)=(3,4,6)$}
In this case, we have
\begin{align}
(a,b,c,x)=(a,b,2\,a,2). \tag{3,4,6-1}
\end{align}
\paragraph{(3,4,6-1)}
The special values obtained from (3,4,6-1) are evaluated 
in paragraphs (1,2,2-1) and (0,2,2-1).

\subsubsection{$(k,l,m)=(3,5,6)$}
In this case, we have
\begin{align}
&\begin{cases}
(a,b,c,x)=(a, 5/3\,a-1/2, 2\,a,-8+4\,\sqrt{5}),\\
S^{(n)}=\dfrac{{5}^{5/2\,n} \left( \sqrt {5}-1 \right) ^{15\,n}
 \left( 1/3\,a+3/10,n \right) \left( 1/3\,a+7/10,n \right)}{{2}^{15\,n} {3}^{3\,n}
 \left( 1/3\,a+1/6,n \right) 
 \left( 1/3\,a+5/6,n \right)},
\tag{3,5,6-1}
\end{cases}\\
&\begin{cases}
(a,b,c,x)=(a, 5/3\,a-1/2, 2\,a,-8-4\,\sqrt{5}),\\
S^{(n)}=\dfrac{{5}^{5/2\,n} \left( \sqrt {5}+1 \right) ^{15\,n}
 \left( 1/3\,a+3/10,n \right) \left( 1/3\,a+7/10,n \right)}{{2}^{15\,n} {3}^{3\,n}
 \left( 1/3\,a+1/6,n \right) 
 \left( 1/3\,a+5/6,n \right)},
\tag{3,5,6-2}
\end{cases}\\
&\begin{cases}
(a,b,c,x)=(a, 5/3\,a-1/6, 2\,a,-8+4\,\sqrt{5}),\\
S^{(n)}=\dfrac{{5}^{5/2\,n} \left( \sqrt {5}-1 \right) ^{15\,n}
 \left( 1/3\,a+17/30,n \right) \left( 1/3\,a+23/30,n \right)}{{2}^{15\,n} {3}^{3\,n}
 \left( 1/3\,a+1/2,n \right) 
 \left( 1/3\,a+5/6,n \right)},
\tag{3,5,6-3}
\end{cases}\\
&\begin{cases}
(a,b,c,x)=(a, 5/3\,a-1/6, 2\,a,-8-4\,\sqrt{5}),\\
S^{(n)}=\dfrac{{5}^{5/2\,n} \left( \sqrt {5}+1 \right) ^{15\,n}
 \left( 1/3\,a+17/30,n \right) \left( 1/3\,a+23/30,n \right)}{{2}^{15\,n} {3}^{3\,n}
 \left( 1/3\,a+1/2,n \right) 
 \left( 1/3\,a+5/6,n \right)}.
\tag{3,5,6-4}
\end{cases}
\end{align}
\paragraph{(3,5,6-1)}
\begin{flalign*}
&{\rm (i)}\,
{\mbox{$F$}(a,5/3\,a-1/2;\,2\,a;\,-8+4\,\sqrt {5})}&\\
&=\begin{cases}
\dfrac{5^{1/4-5/6\,a}\left(\sqrt{5}-1\right)^{3/2-5\,a}
\Gamma\left(2/5 \right)\Gamma\left(4/5 \right)
\Gamma\left(1/3\,a+1/6 \right)\Gamma\left(1/3\,a+5/6 \right)
}{
2^{3/2-5\,a}3^{3/10-a}
\Gamma\left(4/15 \right)\Gamma\left(14/15 \right)
\Gamma\left(1/3\,a+3/10 \right)\Gamma\left(1/3\,a+7/10 \right)
},\\
\dfrac{7\cdot {5}^{5/2\,n-1/2} \left( \sqrt {5}-1 \right) ^{15\,n+15} \left( 17/10,n \right)  \left( 
13/10,n \right)}{  {2}^{15\,n+15} {3}^{3\,n}  \left( 11/6,n \right) 
 \left( 7/6,n \right) }
\quad {\text {if $a=-3-3\,n$}},\\
\dfrac{- {5}^{5/2\,n} \left( \sqrt {5}-1 \right) ^{15\,n+7} \left( 31/30,n \right)  \left( 
19/30,n \right)}{  {2}^{15\,n+6} {3}^{3\,n+1}  \left( 7/6,n \right) 
 \left( 1/2,n \right) }
\quad {\text {if $a=-1-3\,n$}},\\
\dfrac{-11\cdot {5}^{5/2\,n+1/2} \left( \sqrt {5}-1 \right) ^{15\,n+12} \left( 41/30,n \right)  \left( 
29/30,n \right)}{  {2}^{15\,n+12} {3}^{3\,n+3}  \left( 3/2,n \right) 
 \left( 5/6,n \right) }
\quad {\text {if $a=-2-3\,n$}},\\
\dfrac{ -{5}^{5/2\,n+1} \left( \sqrt {5}-1 \right) ^{15\,n+24} \left( 11/5,n \right)  \left( 
9/5,n \right)}{  {2}^{15\,n+24} {3}^{3\,n}  \left( 7/3,n \right) 
 \left(5/3,n \right) }
\quad {\text {if $a=-9/2-3\,n$}},\\
\end{cases}&
\end{flalign*}
\begin{flalign*}
&{\rm (ii)}\,
{\mbox{$F$}(a,1/3\,a+1/2;\,2\,a;\,-8+4\,\sqrt {5})}&\\
&=\begin{cases}
\dfrac{5^{1/4-5/6\,a}\left(\sqrt{5}-1\right)^{-a-3/2}
\Gamma\left(2/5 \right)\Gamma\left(4/5 \right)
\Gamma\left(1/3\,a+1/6 \right)\Gamma\left(1/3\,a+5/6 \right)
}{
2^{-a-3/2}3^{3/10-a}
\Gamma\left(4/15 \right)\Gamma\left(14/15 \right)
\Gamma\left(1/3\,a+3/10 \right)\Gamma\left(1/3\,a+7/10 \right)
},\\
\dfrac{7\cdot {5}^{5/2\,n-1/2} \left( \sqrt {5}-1 \right) ^{3\,n+3} \left( 17/10,n \right)  \left( 
13/10,n \right)}{  {2}^{3\,n+3} {3}^{3\,n}  \left( 11/6,n \right) 
 \left( 7/6,n \right) }
\quad {\text {if $a=-3-3\,n$}},\\
\dfrac{ {5}^{5/2\,n} \left( \sqrt {5}-1 \right) ^{3\,n-1} \left( 31/30,n \right)  \left( 
19/30,n \right)}{  {2}^{3\,n-2} {3}^{3\,n+1}  \left( 7/6,n \right) 
 \left( 1/2,n \right) }
\quad {\text {if $a=-1-3\,n$}},\\
\dfrac{11\cdot {5}^{5/2\,n+1/2} \left( \sqrt {5}-1 \right) ^{3\,n} \left( 41/30,n \right)  \left( 
29/30,n \right)}{  {2}^{3\,n} {3}^{3\,n+3}  \left( 3/2,n \right) 
 \left( 5/6,n \right) }
\quad {\text {if $a=-2-3\,n$}},\\
\dfrac{ {5}^{5/2\,n} \left( \sqrt {5}-1 \right) ^{3\,n} \left( 6/5,n \right)  \left( 
4/5,n \right)}{  {2}^{3\,n} {3}^{3\,n}  \left( 4/3,n \right) 
 \left( 2/3,n \right) }
\quad {\text {if $a=-3/2-3\,n$}}.
\end{cases}&
\end{flalign*}
\begin{flalign*}
&{\rm (iii)}\,
{\mbox{$F$}(a,1/3\,a+1/2;\,2\,a;\,-8-4\,\sqrt {5})}&\\
&=\begin{cases}
\dfrac{7\cdot {5}^{5/2\,n-1/2} \left( \sqrt {5}+1 \right) ^{3\,n+3} \left( 17/10,n \right)  \left( 
13/10,n \right)}{  {2}^{3\,n+3} {3}^{3\,n}  \left( 11/6,n \right) 
 \left( 7/6,n \right) }
\quad {\text {if $a=-3-3\,n$}},\\
\dfrac{ -{5}^{5/2\,n} \left( \sqrt {5}+1 \right) ^{3\,n-1} \left( 31/30,n \right)  \left( 
19/30,n \right)}{  {2}^{3\,n-2} {3}^{3\,n+1}  \left( 7/6,n \right) 
 \left( 1/2,n \right) }
\quad {\text {if $a=-1-3\,n$}},\\
\dfrac{-11\cdot {5}^{5/2\,n+1/2} \left( \sqrt {5}+1 \right) ^{3\,n} \left( 41/30,n \right)  \left( 
29/30,n \right)}{  {2}^{3\,n} {3}^{3\,n+3}  \left( 3/2,n \right) 
 \left( 5/6,n \right) }
\quad {\text {if $a=-2-3\,n$}},\\
\dfrac{ {5}^{5/2\,n} \left( \sqrt {5}+1 \right) ^{3\,n} \left( 6/5,n \right)  \left( 
4/5,n \right)}{  {2}^{3\,n} {3}^{3\,n}  \left( 4/3,n \right) 
 \left( 2/3,n \right) }
\quad {\text {if $a=-3/2-3\,n$}}.
\end{cases}&
\end{flalign*}
\begin{flalign*}
&{\rm (iv)}\,
{\mbox{$F$}(a,5/3\,a-1/2;\,2\,a;\,-8-4\,\sqrt {5})}&\\
&=\begin{cases}
\dfrac{7\cdot {5}^{5/2\,n-1/2} \left( \sqrt {5}+1 \right) ^{15\,n+15} \left( 17/10,n \right)  \left( 
13/10,n \right)}{  {2}^{15\,n+15} {3}^{3\,n}  \left( 11/6,n \right) 
 \left( 7/6,n \right) }
\quad {\text {if $a=-3-3\,n$}},\\
\dfrac{{5}^{5/2\,n} \left( \sqrt {5}+1 \right) ^{15\,n+7} \left( 31/30,n \right)  \left( 
19/30,n \right)}{  {2}^{15\,n+6} {3}^{3\,n+1}  \left( 7/6,n \right) 
 \left( 1/2,n \right) }
\quad {\text {if $a=-1-3\,n$}},\\
\dfrac{11\cdot {5}^{5/2\,n+1/2} \left( \sqrt {5}+1 \right) ^{15\,n+12} \left( 41/30,n \right)  \left( 
29/30,n \right)}{  {2}^{15\,n+12} {3}^{3\,n+3}  \left( 3/2,n \right) 
 \left( 5/6,n \right) }
\quad {\text {if $a=-2-3\,n$}},\\
\dfrac{ {5}^{5/2\,n} \left( \sqrt {5}+1 \right) ^{15\,n} \left( 3/5,n \right)  \left( 
1/5,n \right)}{  {2}^{15\,n} {3}^{3\,n}  \left( 11/15,n \right) 
 \left( 1/15,n \right) }
\quad {\text {if $a=3/10-3\,n$}},\\
\dfrac{ {5}^{5/2\,n+1/2} \left( \sqrt {5}+1 \right) ^{15\,n+3} \left( 4/5,n \right)  \left( 
2/5,n \right)}{  {2}^{15\,n+3} {3}^{3\,n}  \left( 14/15,n \right) 
 \left( 4/15,n \right) }
\quad {\text {if $a=-3/10-3\,n$}},\\
0
\quad {\text {if $a=-9/10-3\,n$}},\\
\dfrac{ -{5}^{5/2\,n+1} \left( \sqrt {5}+1 \right) ^{15\,n+24} \left( 11/5,n \right)  \left( 
9/5,n \right)}{  {2}^{15\,n+24} {3}^{3\,n}  \left( 7/3,n \right) 
 \left(5/3,n \right) }
\quad {\text {if $a=-9/2-3\,n$}},\\
0
\quad {\text {if $a=-21/10-3\,n$}}.
\end{cases}&
\end{flalign*}
\begin{flalign*}
&{\rm (v)}\,
{\mbox{$F$}(a,5/3\,a-1/2;\,2/3\,a+1/2;\,9-4\,\sqrt {5})}&\\
&=
\dfrac{{5}^{-5/6\,a} \left( \sqrt {5}-1 \right) ^{1-5\,a}\Gamma  \left( 1/3\,a
+1/4 \right) \Gamma  \left( 1/3\,a+3/4 \right)}{{2}^{1/2-11/3\,
a}    \Gamma  \left( 1/3\,a+3/10 \right) 
\Gamma  \left( 1/3\,a+7/10\right) 
}.
&
\end{flalign*}
\begin{flalign*}
&{\rm (vi)}\,
{\mbox{$F$}(1-a,1/2-1/3\,a;\,2/3\,a+1/2;\,9-4\,\sqrt {5})}&\\
&=
\dfrac{{5}^{-5/6\,a} \left( \sqrt {5}-1 \right) ^{a-2}\Gamma  \left( 1/3\,a
+1/4 \right) \Gamma  \left( 1/3\,a+3/4 \right)}{{2}^{-5/3\,
a-1/2}   \Gamma  \left( 1/3\,a+3/10 \right) 
\Gamma  \left( 1/3\,a+7/10\right) 
}
&
\end{flalign*}
(The above is a generalization of Theorem 14 in [Ek]).
\begin{flalign*}
&{\rm (vii)}\,
{\mbox{$F$}(a,1-a;\,2/3\,a+1/2;\,1/2-1/4\,\sqrt {5})}&\\
&=
\dfrac{{5}^{-5/6\,a} \left( \sqrt {5}-1 \right) ^{1-2\,a}\Gamma  \left( 1/3\,a
+1/4 \right) \Gamma  \left( 1/3\,a+3/4 \right)}{{2}^{1/2-8/3\,
a}    \Gamma  \left( 1/3\,a+3/10 \right) 
\Gamma  \left( 1/3\,a+7/10\right) 
}.
&
\end{flalign*}
\begin{flalign*}
&{\rm (viii)}\,
{\mbox{$F$}(1/2-1/3\,a,5/3\,a-1/2;\,2/3\,a+1/2;\,1/2-1/4\,\sqrt {5})}&\\
&=
\dfrac{{5}^{-5/6\,a} \left( \sqrt {5}-1 \right) ^{-1/2}\Gamma  \left( 1/3\,a
+1/4 \right) \Gamma  \left( 1/3\,a+3/4 \right)}{{2}^{-2\,
a}    \Gamma  \left( 1/3\,a+3/10 \right) 
\Gamma  \left( 1/3\,a+7/10\right) 
}.
&
\end{flalign*}
\begin{flalign*}
&{\rm (ix)}\,
{\mbox{$F$}(a,1-a;\,3/2-2/3\,a;\,1/2+1/4\,\sqrt {5})}&\\
&=\begin{cases}
{\dfrac { \left( \sqrt {5}-1 \right) ^{6\,n} \left( 5/4
,n \right)  \left( 3/4,n \right) }{ \left( -16
 \right) ^{n}{5}^{5/2\,n} \left( 11/10,n \right)\left( 9/10,n \right) }}
\quad {\text {if $a=-3\,n$}},\\
{\dfrac { \left( \sqrt {5}-1 \right) ^{6\,n+4} \left( 19/12
,n \right)  \left( 13/12,n \right) }{ 104\left( -16
 \right) ^{n}{5}^{5/2\,n} \left( 43/30,n \right)\left( 37/30,n \right) }}
\quad {\text {if $a=-1-3\,n$}},\\
{\dfrac { 11\left( \sqrt {5}-1 \right) ^{6\,n+6} \left( 23/12
,n \right)  \left( 17/12,n \right) }{ 1564\left( -16
 \right) ^{n+1}{5}^{5/2\,n-1/2} \left( 53/30,n \right)\left( 47/30,n \right) }}
\quad {\text {if $a=-2-3\,n$}},\\
{\dfrac { 5^{5/2\,n}\left( \sqrt {5}+1 \right) ^{6\,n} \left( 7/30
,n \right)  \left( 13/30,n \right) }{ \left( -256
 \right) ^{n} \left(1/12,n \right)\left( 7/12,n \right) }}
\quad {\text {if $a=1+3\,n$}},\\
{\dfrac { -5^{5/2\,n+1/2}\left( \sqrt {5}+1 \right) ^{6\,n+2} \left( 17/30
,n \right)  \left( 23/30,n \right) }{ 2\left( -256
 \right) ^{n} \left(5/12,n \right)\left( 11/12,n \right) }}
\quad {\text {if $a=2+3\,n$}},\\
{\dfrac { -5^{5/2\,n+1/2}\left( \sqrt {5}+1 \right) ^{6\,n+6} \left( 9/10
,n \right)  \left( 11/10,n \right) }{ 64\left( -256
 \right) ^{n} \left(3/4,n \right)\left( 5/4,n \right) }}
\quad {\text {if $a=3+3\,n$}}.
\end{cases}&
\end{flalign*}
\begin{flalign*}
&{\rm (x)}\,
{\mbox{$F$}(1/3\,a+1/2,3/2-5/3\,a;\,3/2-2/3\,a;\,1/2+1/4\,\sqrt {5})}&\\
&=\begin{cases}
{\dfrac { \left( -64 \right) ^{n} \left( 7/4,n \right) 
 \left( 5/4,n \right) }{{5}^{5/2\,n}
 \left( 8/5,n \right) \left( 7/5,n \right) }}
\quad {\text {if $a=-3/2-3\,n$}},\\
{\dfrac { {5}^{5/2\,n} \left( 1/5,n \right) 
 \left( 2/5,n \right) }{\left( -64 \right) ^{n}
 \left( 1/20,n \right) \left( 11/20,n \right) }}
\quad {\text {if $a=9/10+3\,n$}},\\
{\dfrac { -{5}^{5/2\,n+1/2} \left( 2/5,n \right) 
 \left( 3/5,n \right) }{2\left( -64 \right) ^{n}
 \left( 1/4,n \right) \left( 3/4,n \right) }}
\quad {\text {if $a=3/2+3\,n$}},\\
{\dfrac { {5}^{5/2\,n+1} \left( 3/5,n \right) 
 \left( 4/5,n \right) }{2\left( -64 \right) ^{n}
 \left( 9/20,n \right) \left( 19/20,n \right) }}
\quad {\text {if $a=21/10+3\,n$}},\\
0
\quad {\text {if $a=27/10+3\,n$}},\\
0
\quad {\text {if $a=33/10+3\,n$}}.
\end{cases}&
\end{flalign*}
\begin{flalign*}
&{\rm (xi)}\,
{\mbox{$F$}(a,1/3\,a+1/2;\,3/2-2/3\,a;\,9+4\,\sqrt {5})}&\\
&=\begin{cases}
{\dfrac {{2}^{5\,n} \left( \sqrt {5}+1 \right) ^{3\,n}
 \left( 5/4,n \right) \left( 3/4,n \right) }{{5}^{5/2
\,n}\left( 11/10,n \right) \left( 9/
10,n \right) }}
\quad {\text {if $a=-3\,n$}},\\
{\dfrac {-{2}^{5\,n+4} \left( \sqrt {5}+1 \right) ^{3\,n-1}
 \left( 19/12,n \right) \left( 13/12,n \right) }{13\cdot {5}^{5/2
\,n}\left( 43/30,n \right) \left( 37/
30,n \right) }}
\quad {\text {if $a=-1-3\,n$}},\\
{\dfrac {-11\cdot {2}^{5\,n+4} \left( \sqrt {5}+1 \right) ^{3\,n}
 \left( 23/12,n \right) \left( 17/12,n \right) }{391\cdot {5}^{5/2
\,n-1/2}\left( 53/30,n \right) \left( 47/
30,n \right) }}
\quad {\text {if $a=-2-3\,n$}},\\
{\dfrac { {2}^{5\,n} \left( \sqrt {5}+1 \right) ^{3\,n}
 \left( 7/4,n \right) \left( 5/4,n \right) }{ {5}^{5/2
\,n}\left( 8/5,n \right) \left( 7/
5,n \right) }}
\quad {\text {if $a=-3/2-3\,n$}}
\end{cases}&
\end{flalign*}
(The fourth case is identical to Theorem 15 in [Ek]).
\begin{flalign*}
&{\rm (xii)}\,
{\mbox{$F$}(1-a,3/2-5/3\,a;\,3/2-2/3\,a;\,9+4\,\sqrt {5})}&\\
&=\begin{cases}
\dfrac{{5}^{5/2\,n} \left( \sqrt {5}+1 \right) ^{15\,n}
 \left( 7/30,n \right)  \left( 13/
30,n \right)}{ {2}^{11\,n} \left( 1/12,n \right)   \left(7/12,n \right) }  
\quad {\text {if $a=1+3\,n$}},\\
\dfrac{{5}^{5/2\,n+1/2} \left( \sqrt {5}+1 \right) ^{15\,n+5}
 \left( 17/30,n \right)  \left( 23/
30,n \right)}{ {2}^{11\,n+2} \left( 5/12,n \right)   \left(11/12,n \right) }  
\quad {\text {if $a=2+3\,n$}},\\
\dfrac{-{5}^{5/2\,n+1/2} \left( \sqrt {5}+1 \right) ^{15\,n+12}
 \left( 9/10,n \right)  \left( 11/
10,n \right)}{ {2}^{11\,n+8} \left( 3/4,n \right)   \left(5/4,n \right) }  
\quad {\text {if $a=3+3\,n$}},\\
\dfrac{{5}^{5/2\,n} \left( \sqrt {5}+1 \right) ^{15\,n}
 \left( 1/5,n \right)  \left( 2/
5,n \right)}{ {2}^{11\,n} \left( 1/20,n \right)   \left(11/20,n \right) }  
\quad {\text {if $a=9/10+3\,n$}},\\
\dfrac{{5}^{5/2\,n+1/2} \left( \sqrt {5}+1 \right) ^{15\,n+3}
 \left( 2/5,n \right)  \left( 3/
5,n \right)}{ {2}^{11\,n+2} \left( 1/4,n \right)   \left(3/4,n \right) }  
\quad {\text {if $a=3/2+3\,n$}},\\
\dfrac{{5}^{5/2\,n+1} \left( \sqrt {5}+1 \right) ^{15\,n+6}
 \left( 3/5,n \right)  \left( 4/
5,n \right)}{ {2}^{11\,n+3} \left( 9/20,n \right)   \left(19/20,n \right) }  
\quad {\text {if $a=21/10+3\,n$}},\\
0
\quad {\text {if $a=27/10+3\,n$}},\\
0
\quad {\text {if $a=33/10+3\,n$}}.
\end{cases}&
\end{flalign*}
\begin{flalign*}
&{\rm (xiii)}\,
{\mbox{$F$}(1/2-1/3\,a,5/3\,a-1/2;\,2/3\,a+1/2;\,1/2+1/4\,\sqrt {5})}&\\
&=\begin{cases}
{\dfrac { \left( -64 \right) ^{n} \left( 3/4,n \right) 
 \left( 5/4,n \right) }{{5}^{5/2\,n}
 \left( 4/5,n \right)  \left( 6/5,n \right) }}
\quad {\text {if $a=3/2+3\,n$}},\\
{\dfrac { {5}^{5/2\,n} \left( 3/5,n \right) 
 \left( 1/5,n \right) }{\left( -64 \right) ^{n}
 \left( 13/20,n \right)  \left( 3/20,n \right) }}
\quad {\text {if $a=3/10-3\,n$}},\\
{\dfrac { -{5}^{5/2\,n+1/2} \left( 4/5,n \right) 
 \left( 2/5,n \right) }{2\left( -64 \right) ^{n}
 \left( 17/20,n \right)  \left( 7/20,n \right) }}
\quad {\text {if $a=-3/10-3\,n$}},\\
0
\quad {\text {if $a=-9/10-3\,n$}},\\
{\dfrac { -{5}^{5/2\,n+1/2} \left( 6/5,n \right) 
 \left( 4/5,n \right) }{4\left( -64 \right) ^{n}
 \left( 5/4,n \right)  \left( 3/4,n \right) }}
\quad {\text {if $a=-3/2-3\,n$}},\\
0
\quad {\text {if $a=-21/10-3\,n$}}.
\end{cases}&
\end{flalign*}
\begin{flalign*}
&{\rm (xiv)}\,
{\mbox{$F$}(a,1-a;\,2/3\,a+1/2;\,1/2+1/4\,\sqrt {5})}&\\
&=\begin{cases}
\dfrac{{5}^{5/2\,n} \left( \sqrt {5}+1 \right) ^{6\,n}
 \left( 7/10,n \right) \left( 3/10,n
 \right)}{ \left( -256 \right) ^{n} \left( 3/4,n \right)   \left( 1/4,n \right) }
\quad {\text {if $a=-3\,n$}},\\
\dfrac{{5}^{5/2\,n} \left( \sqrt {5}+1 \right) ^{6\,n+4}
 \left( 31/30,n \right) \left( 19/30,n
 \right)}{ 8\left( -256 \right) ^{n} \left( 13/12,n \right)   \left( 7/12,n \right) }
\quad {\text {if $a=-1-3\,n$}},\\
\dfrac{-11\cdot {5}^{5/2\,n+1/2} \left( \sqrt {5}+1 \right) ^{6\,n+6}
 \left( 41/30,n \right) \left( 29/30,n
 \right)}{ 320\left( -256 \right) ^{n} \left( 17/12,n \right)   \left( 11/12,n \right) }
\quad {\text {if $a=-2-3\,n$}},\\
\dfrac{\left( \sqrt {5}-1 \right) ^{6\,n} \left( 7/
12,n \right) \left( 13/12,n \right) }{\left( -16 \right) ^{n}{5}^{5/2\,n}
 \left( 19/30,n
 \right)  \left( 31/30,n \right) }
\quad {\text {if $a=1+3\,n$}},\\
\dfrac{-\left( \sqrt {5}-1 \right) ^{6\,n+2} \left( 11/
12,n \right) \left( 17/12,n \right) }{22\left( -16 \right) ^{n}{5}^{5/2\,n-1/2}
 \left( 29/30,n
 \right)  \left( 41/30,n \right) }
\quad {\text {if $a=2+3\,n$}},\\
\dfrac{-\left( \sqrt {5}-1 \right) ^{6\,n+6} \left( 5/
4,n \right) \left( 7/4,n \right) }{2240\left( -16 \right) ^{n}{5}^{5/2\,n-1/2}
 \left( 13/10,n
 \right)  \left( 17/10,n \right) }
\quad {\text {if $a=3+3\,n$}}.
\end{cases}&
\end{flalign*}
\begin{flalign*}
&{\rm (xv)}\,
{\mbox{$F$}(a,5/3\,a-1/2;\,2/3\,a+1/2;\,9+4\,\sqrt {5})}&\\
&=\begin{cases}
\dfrac{{5}^{5/2\,n} \left( \sqrt {5}+1 \right) ^{15\,n}
 \left( 7/10,n \right)\left( 3/10,n
 \right) }{ {2}^{11\,n} 
 \left( 3/4,n \right)  \left( 1/4,n \right)}  
\quad {\text {if $a=-3\,n$}},\\
\dfrac{-{5}^{5/2\,n} \left( \sqrt {5}+1 \right) ^{15\,n+7}
 \left( 31/30,n \right)\left( 19/30,n
 \right) }{ {2}^{11\,n+4} 
 \left( 13/12,n \right)  \left( 7/12,n \right)}  
\quad {\text {if $a=-1-3\,n$}},\\
\dfrac{-11\cdot {5}^{5/2\,n-1/2} \left( \sqrt {5}+1 \right) ^{15\,n+12}
 \left( 41/30,n \right)\left( 29/30,n
 \right) }{ {2}^{11\,n+8} 
 \left( 17/12,n \right)  \left( 11/12,n \right)}  
\quad {\text {if $a=-2-3\,n$}},\\
\dfrac{ {5}^{5/2\,n} \left( \sqrt {5}+1 \right) ^{15\,n}
 \left( 3/5,n \right)\left( 1/5,n
 \right) }{ {2}^{11\,n} 
 \left( 13/20,n \right)  \left( 3/20,n \right)}  
\quad {\text {if $a=3/10-3\,n$}},\\
\dfrac{ {5}^{5/2\,n+1/2} \left( \sqrt {5}+1 \right) ^{15\,n+3}
 \left( 4/5,n \right)\left( 2/5,n
 \right) }{ {2}^{11\,n+2} 
 \left( 17/20,n \right)  \left( 7/20,n \right)}  
\quad {\text {if $a=-3/10-3\,n$}},\\
0
\quad {\text {if $a=-9/10-3\,n$}},\\
\dfrac{ {5}^{5/2\,n+1/2} \left( \sqrt {5}+1 \right) ^{15\,n+9}
 \left( 6/5,n \right)\left( 4/5,n
 \right) }{ {2}^{11\,n+5} 
 \left( 5/4,n \right)  \left( 3/4,n \right)}  
\quad {\text {if $a=-3/2-3\,n$}},\\
0
\quad {\text {if $a=-21/10-3\,n$}}.
\end{cases}&
\end{flalign*}
\begin{flalign*}
&{\rm (xvi)}\,
{\mbox{$F$}(1-a,1/2-1/3\,a;\,2/3\,a+1/2;\,9+4\,\sqrt {5})}&\\
&=\begin{cases}
\dfrac{{2}^{5\,n} \left( \sqrt {5}+1 \right) ^{3\,n} \left( 7/12,n \right)  \left( 13/12,n
 \right)}{{5}^{5/2\,n} 
 \left( 19/30,n \right)  \left( 31/30,n \right)  }
\quad {\text {if $a=1+3\,n$}},\\
\dfrac{{2}^{5\,n+2} \left( \sqrt {5}+1 \right) ^{3\,n+1} \left( 11/12,n \right)  \left( 17/12,n
 \right)}{11\cdot {5}^{5/2\,n-1/2} 
 \left( 29/30,n \right)  \left( 41/30,n \right)  }
\quad {\text {if $a=2+3\,n$}},\\
\dfrac{-{2}^{5\,n+4} \left( \sqrt {5}+1 \right) ^{3\,n} \left( 5/4,n \right)  \left( 7/4,n
 \right)}{7\cdot {5}^{5/2\,n+1/2} 
 \left( 13/10,n \right)  \left( 17/10,n \right)  }
\quad {\text {if $a=3+3\,n$}},\\
\dfrac{{2}^{5\,n} \left( \sqrt {5}+1 \right) ^{3\,n} \left( 3/4,n \right)  \left( 5/4,n
 \right)}{ {5}^{5/2\,n} 
 \left( 4/5,n \right)  \left( 6/5,n \right)  }
\quad {\text {if $a=3/2+3\,n$}}
\end{cases}&
\end{flalign*}
(The fourth case is identical to Theorem 14 in [Ek]).
\begin{flalign*}
&{\rm (xvii)}\,
{\mbox{$F$}(1-a,1/2-1/3\,a;\,2-2\,a;\,-8+4\,\sqrt {5})}&\\
&=\begin{cases}
\dfrac{{3}^{9/10-a} \left( \sqrt {5}+1 \right) ^{3/2-a}
 \Gamma \left( 4/5 \right) \Gamma  \left( 3/5 \right)
\Gamma  \left( 7/6-1/3\,a \right) \Gamma  \left( 5/6-1/3\,a \right)}{ 
{2}^{3/2-a}{5}^{3/4-5/6\,a} 
\Gamma  \left( 13/15 \right)  \Gamma  \left( 8/15 \right)
\Gamma  \left( 11/10-1/3\,a \right) \Gamma  \left( 9/10-1/3\,a \right) },\\
\dfrac{91\cdot{5}^{5/2\,n+1/2} \left( \sqrt {5}-1 \right) ^{3\,n+3} \left( 37/30,n \right) \left( 
43/30,n \right)}{{2}^{3\,n+3}  {
3}^{3\,n+4} \left( 7/6,n
 \right)   \left( 3/2,n \right)
}  
\quad {\text {if $a=4+3\,n$}},\\
\dfrac{{5}^{5/2\,n+1/2} \left( \sqrt {5}-1 \right) ^{3\,n+1} \left( 17/30,n \right) \left( 
23/30,n \right)}{{2}^{3\,n}  {
3}^{3\,n+1} \left( 1/2,n
 \right)   \left( 5/6,n \right)
}  
\quad {\text {if $a=2+3\,n$}},\\
\dfrac{{5}^{5/2\,n+1/2} \left( \sqrt {5}-1 \right) ^{3\,n} \left( 9/10,n \right) \left( 
11/10,n \right)}{{2}^{3\,n}  {
3}^{3\,n+1} \left( 5/6,n
 \right)   \left( 7/6,n \right)
}  
\quad {\text {if $a=3+3\,n$}},\\
\dfrac{{5}^{5/2\,n} \left( \sqrt {5}-1 \right) ^{3\,n} \left( 2/5,n \right) \left( 
3/5,n \right)}{{2}^{3\,n}  {
3}^{3\,n} \left( 1/3,n
 \right)   \left( 2/3,n \right)
}  
\quad {\text {if $a=3/2+3\,n$}}.
\end{cases}&
\end{flalign*}
\begin{flalign*}
&{\rm (xviii)}\,
{\mbox{$F$}(1-a,3/2-5/3\,a;\,2-2\,a;\,-8+4\,\sqrt {5})}&\\
&=\begin{cases}
\dfrac{{3}^{9/10-a} \left( \sqrt {5}+1 \right) ^{9/2-5\,a}
 \Gamma \left( 4/5 \right) \Gamma  \left( 3/5 \right)
\Gamma  \left( 7/6-1/3\,a \right) \Gamma  \left( 5/6-1/3\,a \right)}{ 
{2}^{9/2-5\,a}{5}^{3/4-5/6\,a} 
\Gamma  \left( 13/15 \right)  \Gamma  \left( 8/15 \right)
\Gamma  \left( 11/10-1/3\,a \right) \Gamma  \left( 9/10-1/3\,a \right) },\\
\dfrac{91\cdot {5}^{5/2\,n+1/2} \left( \sqrt {5}-1 \right) ^{15\,n+15} \left(37/30,n \right) \left( 
43/30,n \right)}{ {2}^{15\,n+15} {3}^{3\,n+4} \left( 7/6,
n \right)  \left( 3/2,n \right) }
\quad {\text {if $a=4+3\,n$}},\\
\dfrac{ {5}^{5/2\,n+1/2} \left( \sqrt {5}-1 \right) ^{15\,n+5} \left(17/30,n \right) \left( 
23/30,n \right)}{ {2}^{15\,n+4} {3}^{3\,n+1} \left( 1/2,
n \right)  \left( 5/6,n \right) }
\quad {\text {if $a=2+3\,n$}},\\
\dfrac{ -{5}^{5/2\,n+1/2} \left( \sqrt {5}-1 \right) ^{15\,n+12} \left(9/10,n \right) \left( 
11/10,n \right)}{ {2}^{15\,n+12} {3}^{3\,n+1} \left( 5/6,
n \right)  \left( 7/6,n \right) }
\quad {\text {if $a=3+3\,n$}},\\
\dfrac{ {5}^{5/2\,n+1} \left( \sqrt {5}-1 \right) ^{15\,n+18} \left(7/5,n \right) \left( 
8/5,n \right)}{ {2}^{15\,n+18} {3}^{3\,n} \left( 4/3,
n \right)  \left( 5/3,n \right) }
\quad {\text {if $a=9/2+3\,n$}}.
\end{cases}&
\end{flalign*}
\begin{flalign*}
&{\rm (xix)}\,
{\mbox{$F$}(1-a,3/2-5/3\,a;\,2-2\,a;\,-8-4\,\sqrt {5})}&\\
&=\begin{cases}
\dfrac{91\cdot {5}^{5/2\,n+1/2} \left( \sqrt {5}+1 \right) ^{15\,n+15} \left(37/30,n \right) \left( 
43/30,n \right)}{ {2}^{15\,n+15} {3}^{3\,n+4} \left( 7/6,
n \right)  \left( 3/2,n \right) }
\quad {\text {if $a=4+3\,n$}},\\
\dfrac{ {5}^{5/2\,n+1/2} \left( \sqrt {5}+1 \right) ^{15\,n+5} \left(17/30,n \right) \left( 
23/30,n \right)}{ {2}^{15\,n+4} {3}^{3\,n+1} \left( 1/2,
n \right)  \left( 5/6,n \right) }
\quad {\text {if $a=2+3\,n$}},\\
\dfrac{ {5}^{5/2\,n+1/2} \left( \sqrt {5}+1 \right) ^{15\,n+12} \left(9/10,n \right) \left( 
11/10,n \right)}{ {2}^{15\,n+12} {3}^{3\,n+1} \left( 5/6,
n \right)  \left( 7/6,n \right) }
\quad {\text {if $a=3+3\,n$}},\\
\dfrac{ {5}^{5/2\,n} \left( \sqrt {5}+1 \right) ^{15\,n} \left(1/5,n \right) \left( 
2/5,n \right)}{ {2}^{15\,n} {3}^{3\,n} \left( 2/15,
n \right)  \left( 7/15,n \right) }
\quad {\text {if $a=9/10+3\,n$}},\\
\dfrac{ {5}^{5/2\,n+1} \left( \sqrt {5}+1 \right) ^{15\,n+18} \left(7/5,n \right) \left( 
8/5,n \right)}{ {2}^{15\,n+18} {3}^{3\,n} \left( 4/3,
n \right)  \left( 5/3,n \right) }
\quad {\text {if $a=9/2+3\,n$}},\\
\dfrac{ {5}^{5/2\,n+1} \left( \sqrt {5}+1 \right) ^{15\,n+6} \left(3/5,n \right) \left( 
4/5,n \right)}{ {2}^{15\,n+6} {3}^{3\,n+1} \left( 8/15,
n \right)  \left( 13/15,n \right) }
\quad {\text {if $a=21/10+3\,n$}},\\
0
\quad {\text {if $a=27/10+3\,n$}},\\
0
\quad {\text {if $a=33/10+3\,n$}}.
\end{cases}&
\end{flalign*}
\begin{flalign*}
&{\rm (xx)}\,
{\mbox{$F$}(1-a,1/2-1/3\,a;\,2-2\,a;\,-8-4\,\sqrt {5})}&\\
&=\begin{cases}
\dfrac{91\cdot{5}^{5/2\,n+1/2} \left( \sqrt {5}+1 \right) ^{3\,n+3} \left( 37/30,n \right) \left( 
43/30,n \right)}{{2}^{3\,n+3}  {
3}^{3\,n+4} \left( 7/6,n
 \right)   \left( 3/2,n \right)
}  
\quad {\text {if $a=4+3\,n$}},\\
\dfrac{{5}^{5/2\,n+1/2} \left( \sqrt {5}+1 \right) ^{3\,n+1} \left( 17/30,n \right) \left( 
23/30,n \right)}{{2}^{3\,n}  {
3}^{3\,n+1} \left( 1/2,n
 \right)   \left( 5/6,n \right)
}  
\quad {\text {if $a=2+3\,n$}},\\
\dfrac{-{5}^{5/2\,n+1/2} \left( \sqrt {5}+1 \right) ^{3\,n} \left( 9/10,n \right) \left( 
11/10,n \right)}{{2}^{3\,n}  {
3}^{3\,n+1} \left( 5/6,n
 \right)   \left( 7/6,n \right)
}  
\quad {\text {if $a=3+3\,n$}},\\
\dfrac{{5}^{5/2\,n} \left( \sqrt {5}+1 \right) ^{3\,n} \left( 2/5,n \right) \left( 
3/5,n \right)}{{2}^{3\,n}  {
3}^{3\,n} \left( 1/3,n
 \right)   \left( 2/3,n \right)
}  
\quad {\text {if $a=3/2+3\,n$}}.
\end{cases}&
\end{flalign*}
\begin{flalign*}
&{\rm (xxi)}\,
{\mbox{$F$}(a,1/3\,a+1/2;\,3/2-2/3\,a;\,9-4\,\sqrt {5})}&\\
&=
\dfrac{{2}^{5/2-5/3\,a} \left( \sqrt {5}-1 \right) ^{-a-1}\Gamma  \left(5/4 -1/3\,
a \right) \Gamma  \left( 3/4-1/3\,a \right)}{{5}^{1-5/6\,a
}   \Gamma  \left( 11/10-1/3\,a \right)   \Gamma  \left(9/10-1/3\,a \right)  }&
\end{flalign*} 
(The above is a generalization of Theorem 15 in [Ek]).
\begin{flalign*}
&{\rm (xxii)}\,
{\mbox{$F$}(1-a,3/2-5/3\,a;\,3/2-2/3\,a;\,9-4\,\sqrt {5})}&\\
&=
\dfrac{{2}^{7/2-11/3\,a} \left( \sqrt {5}-1 \right) ^{5\,a-4}\Gamma  \left(5/4 -1/3\,
a \right) \Gamma  \left( 3/4-1/3\,a \right)}{{5}^{1-5/6\,a
}  \Gamma  \left( 11/10-1/3\,a \right)   \Gamma  \left(9/10-1/3\,a \right)  }.&
\end{flalign*} 
\begin{flalign*}
&{\rm (xxiii)}\,
{\mbox{$F$}(a,1-a;\,3/2-2/3\,a;\,1/2-1/4\,\sqrt {5})}&\\
&=
\dfrac{{2}^{5/2-8/3\,a} \left( \sqrt {5}-1 \right) ^{2\,a-1}\Gamma  \left(5/4 -1/3\,
a \right) \Gamma  \left( 3/4-1/3\,a \right)}{{5}^{1-5/6\,a
}    \Gamma  \left( 11/10-1/3\,a \right)   \Gamma  \left(9/10-1/3\,a \right)  }.&
\end{flalign*} 
\begin{flalign*}
&{\rm (xxiv)}\,
{\mbox{$F$}(1/3\,a+1/2,3/2-5/3\,a;\,3/2-2/3\,a;\,1/2-1/4\,\sqrt {5})}
&\\
&=
\dfrac{{2}^{2-2\,a} \left( \sqrt {5}-1 \right) ^{1/2}\Gamma  \left(5/4 -1/3\,
a \right) \Gamma  \left( 3/4-1/3\,a \right)}{{5}^{1-5/6\,a
}  \Gamma  \left( 11/10-1/3\,a \right)   \Gamma  \left(9/10-1/3\,a \right)  }.&
\end{flalign*} 
\paragraph{(3,5,6-2), (3,5,6-3), (3,5,6-4)}
The special values obtained from (3,5,6-2), (3,5,6-3), (3,5,6-4)
coincide with those obtained from (3,5,6-1).
\subsubsection{$(k,l,m)=(3,6,6)$}
In this case, we have
\begin{gather*}
(a,b,c,x)=(a,b,2\,a,2)
\tag{3,6,6-1}
\end{gather*}
\paragraph{(3,6,6-1)}
The special values obtained from (3,6,6-1) are evaluated in 
paragraphs (1,2,2-1) and (0,2,2-1).

\subsubsection{$(k,l,m)=(3,7,6)$}
In this case, there is no admissible quadruple.

\subsubsection{$(k,l,m)=(3,8,6)$}
In this case, we have
\begin{gather*}
(a,b,c,x)=(a,b,2\,a,2)
\tag{3,8,6-1}
\end{gather*}
\paragraph{(3,8,6-1)}
The special values obtained from (3,8,6-1) are evaluated in 
paragraphs (1,2,2-1) and (0,2,2-1).
\subsubsection{$(k,l,m)=(3,9,6)$}
In this case, we have
\begin{align}
(a,b,c,x)&=(a,3\,a-1,2\,a,1/2+1/2\,i\sqrt {3}), \tag {3,9,6-1}\\
(a,b,c,x)&=(a,3\,a-1,2\,a,1/2-1/2\,i\sqrt {3}).\tag {3,9,6-2}
\end{align}
\paragraph{(3,9,6-1), (3,9,6-2)}
The special values obtained from (3,9,6-1) and (3,9,6-2) coincide with those obtained from (1,3,2-1).

\appendix
\section{Some hypergeometric identities 
for generalized hypergeometric series and Appell-Lauricella hypergeometric series}
In this appendix, 
we give some identities for generalized hypergeometric series 
and Appell-Lauricella hypergeometric series using the method introduced in section 1. 

\subsection{Some examples for generalized hypergeometric series}
Many identities for generalized hypergeometric series 
$$
\pFq{p}{q}{a_1, a_2, \cdots , a_p}{b_1, b_2, \cdots , b_q}{x}
:=\sum _{n=0} ^{\infty}\frac{(a_1,n)(a_2,n)\cdots  (a_p,n)}{(b_1,n)(b_2,n)\cdots  (b_q,n)(1,n)}x^n
$$
have been discovered 
by making use of various methods.
For example, see ``Hypergeometric Database'' in 36 page,
Table 6.1 in 108 page and Table 8.1 in 158 page of [Ko].
The author feels that
applying the method introduced in section 1 of this paper 
to generalized hypergeometric series,
we can obtain all identities appearing in those tables of [Ko] and,
of course, also new identities.
So, the author is planning to 
tabulate hypergeometric identities for generalized hypergeometric series 
obtained from our method.
Here, we only rediscover some identities appearing in tables of [Ko],
since the purpose of this article is not to tabulate 
hypergeomeric identities for generalized hypergeometric series 
but to tabulate those for the hypergeometric series .

Now, we apply our method to generalized hypergeometric series $_3F_2$.
It is known that for a given quintuple of integers $(k_1, k_2, k_3, l_1, l_2) \in\Z^5$,
there exists a unique triple of rational functions 
$(Q_2(x), Q_1(x), Q_0(x)) \in (\Q (a_1,a_2,a_3,b_1,b_2,x))^3$ 
satisfying 
\begin{align}
\begin{split}
\pFq{3}{2}{a_1+k_1,a_2+k_2,a_3+k_3}{b_1+l_1,b_2+l_2}{x}
&=Q_2(x)\pFq{3}{2}{a_1+2,a_2+2,a_3+2}{b_1+2,b_2+2}{x}\\
&+Q_1(x)\pFq{3}{2}{a_1+1,a_2+1,a_3+1}{b_1+1,b_2+1}{x}\\
&+Q_0(x)\pFq{3}{2}{a_1,a_2,a_3}{b_1,b_2}{x}.
\end{split}
\label{four_term_g}
\end{align}
This relation is called the contiguity relation of $_3F_2$ (or the four term relation of $_3F_2$). 
We remark that 
coefficients $Q_2(x)$, $Q_1(x)$, $Q_0(x)$ can be exactly computed 
for a given $(k_1, k_2, k_3, l_1, l_2)$.
As to a calculation method of these coefficients, for example, see [Ta1].
From (\ref{four_term_g}), we have
\begin{align}
\begin{split}
&\pFq{3}{2}{a_1+nk_1,a_2+nk_2,a_3+nk_3}{b_1+nl_1,b_2+nl_2}{x}\\
&=Q_2 ^{(n)}(x)\pFq{3}{2}{a_1+(n-1)k_1+2,a_2+(n-1)k_2+2,a_3+(n-1)k_3+2}
{b_1+(n-1)l_1+2,b_2+(n-1)l_2+2}{x}\\
&+Q_1 ^{(n)}(x)\pFq{3}{2}{a_1+(n-1)k_1+1,a_2+(n-1)k_2+1,a_3+(n-1)k_3+1}
{b_1+(n-1)l_1+1,b_2+(n-1)l_2+1}{x}\\
&+Q_0 ^{(n)}(x)\pFq{3}{2}{a_1+(n-1)k_1,a_2+(n-1)k_2,a_3+(n-1)k_3}{b_1+(n-1)l_1,b_2+(n-1)l_2}{x},
\end{split}
\label{four_term_g1}
\end{align}
where 
$$
Q_i ^{(n)}(x):=Q_i (x)|_{(a_1,a_2,a_3,b_1,b_2)\rightarrow
(a_1+(n-1)k_1,a_2+(n-1)k_2,a_3+(n-1)k_3, b_1+(n-1)l_1,b_2+(n-1)l_2)}
$$
for $i=0,1,2$.
Let sextuple $(a_1,a_2,a_3,b_1,b_2,x)$ satisfy 
\begin{gather}
\begin{cases}
Q_1^{(n)}(x)=0,\\
Q_2^{(n)}(x)=0\quad
\label{case_g}
\end{cases}
 {\text {for every integer $n$}}.
\end{gather}
Then, for this sextuple, we have the first order difference equation
\begin{gather}
\pFq{3}{2}{a_1,a_2,a_3}{b_1,b_2}{x}=\frac{1}{Q_0 ^{(1)}(x)}\times
\pFq{3}{2}{a_1+k_1,a_2+k_2,a_3+k_3}{b_1+l_1,b_2+l_2}{x}
\label{diff_eq_g}
\end{gather}
or
\begin{align}
\begin{split}
&\pFq{3}{2}{a_1,a_2,a_3}{b_1,b_2}{x}\\
&=\frac{1}{Q_0 ^{(1)}(x)Q_0 ^{(2)}(x)\cdots Q_0 ^{(n)}(x)}\times
\pFq{3}{2}{a_1+nk_1,a_2+nk_2,a_3+nk_3}{b_1+nl_1,b_2+nl_2}{x}
\end{split}
\label{diff_eq_g1}
\end{align}
Finally, solving (\ref{diff_eq_g}) or (\ref{diff_eq_g1}), 
we will yield hypergeometric identity for $_3F_2$.
We summarize the whole process:
\begin{itemize}
\item[] Step 1': Give $(k_1,k_2,k_3,l_1,l_2)\in \Z ^5$. 
\item[] Step 2': Obtain sextuple  $(a_1,a_2,a_3,b_1,b_2,x)$ satisfying the system (\ref{case_g}).
\item[] Step 3': Solve (\ref{diff_eq_g}) or (\ref{diff_eq_g1}). 
\end{itemize} 
Thus, the method introduced in the main part of section 1 
has broad utility. 
Although stated in section 1, 
we remark once again that 
implementing Step 2' is equivalent to solving polynomial systems in $a_1,a_2,a_3,b_1,b_2,x$.
We now do the above process for two quintuple of integers $(k_1,k_2,k_3,l_1,l_2)$
to give examples for the generalized hypergeometric series 
obtained from our method.
\subsubsection*{Example 1': $(k_1,k_2,k_3,l_1,l_2)=(-1,-1,0,1,1)$}
When $(k_1,k_2,k_3,l_1,l_2)=(-1,-1,0,1,1)$, 
$$
(a_1,a_2,a_3,b_1,b_2,x)=(a,b,c,c+1-a,c+1-b,1)
$$
is given as one of solutions of the system (\ref{case_g}).
In this case, we have
\begin{align}
\begin{split}
&\pFq{3}{2}{a,b,c}{c+1-a,c+1-b}{1}\\
&=\frac{(1/2\,c+1-a,n)(1/2\,c+1-b,n)(c+1-a-b,2n)}{(c+1-a,n)(c+1-b,n)(1/2\,c+1-a-b,2n)}\\
&\times\pFq{3}{2}{a-n,b-n,c}{c+1-a+n,c+1-b+n}{1}.
\end{split}
\label{ex_eq_g1}
\end{align}
By solving this equation (\ref{ex_eq_g1}), we obtain
\begin{align}
\begin{split}
&\pFq{3}{2}{a,b,c}{c+1-a,c+1-b}{1}\\
&=\frac{\G(1/2\, c+1)\G(c+1-a)\G(c+1-b)\G(1/2\, c+1-a-b)}
{\G(c+1)\G(1/2\,c+1-a)\G(1/2\,c+1-b)\G(c+1-a-b)}.
\end{split}
\label{ex_g1}
\end{align}
This was discovered by Dixon 
(see no.4 of ``Hypergeometric Database'' in page 36 of [Ko]).
Moreover, remarking that
$$
\lim _{b\rightarrow -\infty}\pFq{3}{2}{a,b,c}{c+1-a,c+1-b}{1}
=F(a,c;c+1-a;-1),
$$
equation (\ref{ex_eq_g1}) degenerates into
\begin{gather}
F(a,c;c+1-a;-1)=\frac{(1/2\,c+1-a,n)}{(c+1-a,n)}F(a-n,c;c+1-a+n;-1)
\label{red_ex_eq_g1}
\end{gather}
(the equation with equivalent to (\ref{red_ex_eq_g1}) also appears 
in the process of getting formulae (1,2,2-1)(v)).
By solving this equation (\ref{red_ex_eq_g1}), 
we can find formulae which are equivalent to (1,2,2-1)(v).
We also note that taking $b\rightarrow -\infty$ in formula (\ref{ex_g1}), 
the first case of (1,2,2-1)(v) is obtained. 
\subsubsection*{Example 2': $(k_1,k_2,k_3,l_1,l_2)=(0,1,2,1,2)$}
When $(k_1,k_2,k_3,l_1,l_2)=(0,1,2,1,2)$, we yield
$$
(a_1,a_2,a_3,b_1,b_2,x)=(a,b,c,1/2\, (a+c+1),2b,1)
$$
as one of solutions of the system (\ref{case_g}).
In this case, we have
\begin{align}
\begin{split}
&\pFq{3}{2}{a,b,c}{1/2\, (a+c+1),2b}{1}\\
&=\frac{(1/2\,(c+1),n)(b+1/2\,(1-a),n)}{(1/2\,(a+c+1),n)(b+1/2,n)}
\times\pFq{3}{2}{a,b+n,c+2n}{1/2\, (a+c+1)+n,2b+2n}{1}.
\end{split}
\label{ex_eq_g2}
\end{align}
Solving this equation (\ref{ex_eq_g1}), we yield
\begin{align}
\begin{split}
&\pFq{3}{2}{a,b,c}{1/2\, (a+c+1),2b}{1}\\
&=\frac{\G(1/2)\G(b+1/2\, (1-a-c))\G(1/2\, (a+c+1))\G(b+1/2)}
{\G(1/2\, (a+1))\G(b+1/2\,(1-c))\G(1/2\,(c+1))\G(b+1/2\, (1-a))}.
\end{split}
\label{ex_g2}
\end{align}
This formula is due to Watson and Whipple
(see no.5 of ``Hypergeometric Database'' in page 36 of [Ko]).
Moreover, paying attention to
$$
\lim _{b\rightarrow +\infty}\pFq{3}{2}{a,b,c}{1/2\, (a+c+1),2b}{1}
=F(a,c;1/2\, (a+c+1);1/2),
$$
equation (\ref{ex_eq_g2}) degenerates into
\begin{align}
\begin{split}
&F(a,c;1/2\, (a+c+1);1/2)\\
&=\frac{(1/2\,(c+1),n)}{(1/2\,(a+c+1),n)}F(a,c+2n;1/2\,(a+c+1)+n;1/2)
\label{red_ex_eq_g2}
\end{split}
\end{align}
(the equation with equivalent to (\ref{red_ex_eq_g2}) also appears 
in the process of getting formulae (1,2,2-1)(viii)).
Solving this equation (\ref{red_ex_eq_g2}), 
we can find formulae which are equivalent to (1,2,2-1)(viii).
We also remark that taking $b\rightarrow +\infty$ in formula (\ref{ex_g2}), 
the first case of (1,2,2-1)(viii) is obtained.

In a similar manner, we can make use of our method 
for generalized hypergeometric series $_pF_q$.
However, as stated before, we save this study for another day.

\subsection{Some examples for Appell-Lauricella hypereometric series}
Hypergeometric identities for Appell-Lauricella hypergeometric series
have not been studied yet,
while many transformation formulae for these series have been given
(for example, see formulae (8.3.1)--(8.4.1.2) in [Sl]).
The reason  comes from the fact 
that it is difficult to use previously known methods for these series,
because these are multiple series.
However, our method can be applied to these series,
because these have contiguity relations like generalized hypergeometric series.
Therefore, our method will become a powerful tool 
for investigating hypergeometric identities for Appell-Lauricella hypereometric series. 

We now apply our method to Appell hypergeometric series
$$
\pFq{}{1}{a;b_1,b_2}{c}{x,y}:=\sum _{m=0}^{\infty} \sum _{n=0} ^{\infty}
\frac{(a,m+n)(b_1,m)(b_2,n)}{(c,m+n)(1,m)(1,n)}{x^my^n}.
$$
We discuss the detail of applying our method to $F_1$, 
although we perform in much the same way as $_3F_2$ case.
It is known that for a given quadruple of integers $(k, l_1, l_2, m) \in\Z^4$,
there exists a unique triple of rational functions 
$(Q_{10}(x,y), Q_{01}(x,y), Q_{00}(x,y)) \in (\Q (a,b_1,b_2,c,x,y))^3$ 
satisfying 
\begin{align}
\begin{split}
\pFq{}{1}{a+k;b_1+l_1,b_2+l_2}{c+m}{x,y}
&=Q_{10}(x,y)\pFq{}{1}{a+1;b_1+1,b_2}{c+1}{x,y}\\
&+Q_{01}(x,y)\pFq{}{1}{a+1;b_1,b_2+1}{c+1}{x,y}\\
&+Q_{00}(x,y)\pFq{}{1}{a;b_1,b_2}{c}{x,y}.
\end{split}
\label{four_term_a}
\end{align}
This relation is called the contiguity relation of $F_1$ (or the four term relation of $F_1$). 
As to a calculation method of coefficients for the above relation, for example, 
see [Ta1] and [Ta2].
We define $Q_{ij} ^{(n)}(x,y)$ by
$$
Q_{ij} ^{(n)}(x,y):=Q_{ij} (x,y)|_{(a,b_1,b_2,c)\rightarrow
(a+(n-1)k,b_1+(n-1)l_1,b_2+(n-1)l_2, c+(n-1)m)}.
$$
Let sextuple $(a,b_1,b_2,c,x,y)$ satisfy 
\begin{gather}
\begin{cases}
Q_{10}^{(n)}(x,y)=0,\\
Q_{01}^{(n)}(x,y)=0\quad
\label{case_a}
\end{cases}
 {\text {for every integer $n$}}.
\end{gather}
Then, for this sextuple, we have the first order difference equation
\begin{gather}
\pFq{}{1}{a;b_1,b_2}{c}{x,y}=\frac{1}{Q_{00} ^{(1)}(x)}\times
\pFq{}{1}{a+k;b_1+l_1,b_2+l_2}{c+m}{x,y}
\label{diff_eq_a}
\end{gather}
or
\begin{align}
\begin{split}
&\pFq{}{1}{a;b_1,b_2}{c}{x,y}\\
&=\frac{1}{Q_{00} ^{(1)}(x,y)Q_{00} ^{(2)}(x,y)\cdots Q_{00} ^{(n)}(x,y)}\times
\pFq{}{1}{a+nk;b_1+nl_1,b_2+nl_2}{c+nm}{x,y}
\end{split}
\label{diff_eq_a1}
\end{align}
Solving (\ref{diff_eq_a}) or (\ref{diff_eq_a1}) as the final step, 
we will find hypergeometric identity for $F_1$.
We summarize the whole process:
\begin{itemize}
\item[] Step 1'': Give $(k,l_1,l_2,m)\in \Z ^4$. 
\item[] Step 2'': Obtain sextuple  $(a,b_1,b_2,c,x,y)$ satisfying the system (\ref{case_a}).
\item[] Step 3'': Solve (\ref{diff_eq_a}) or (\ref{diff_eq_a1}). 
\end{itemize} 
Now, we perform the above process for two quadruple of integers $(k,l_1,l_2,m)$
to obtain some examples of hypergeometric identities for $F_1$.
\subsubsection*{Example 1'': $(k,l_1,l_2,m)=(2,-1,6,3)$}
When $(k,l_1,l_2,m)=(2,-1,6,3)$, we have
$$
(a,b_1,b_2,c,x,y)=(2\alpha,1-\alpha,6\beta,\alpha+2\beta+2/3,1/9,1/3)
$$
as one of solutions of the system (\ref{case_a}).
In this case, the equation (\ref{diff_eq_a1}) becomes 
\begin{align}
\begin{split}
&\pFq{}{1}{2\alpha;1-\alpha,6\beta}{\alpha+2\beta+2/3}{\frac{1}{9},\frac{1}{3}}\\
&=\frac{2^{4n}(\alpha+1/2,n)(\beta+1/3,n)(\beta+1/6,n)}{3^{n}(\alpha+2\beta+2/3,3n)}
\pFq{}{1}{2\alpha+2n;1-\alpha-n,6\beta+6n}{\alpha+2\beta+2/3+3n}{\frac{1}{9},\frac{1}{3}}.
\end{split}
\label{ex_eq_a1}
\end{align}
By solving this equation (\ref{ex_eq_a1}), we find
\begin{gather}
\pFq{}{1}{2\alpha;1-\alpha,6\beta}{\alpha+2\beta+2/3}{\frac{1}{9},\frac{1}{3}}
=\frac{3^{\alpha}\sqrt{\pi}\G(\alpha+2\beta+2/3)
}{2^{2\alpha}\G(\alpha+1/2)\G(2\beta+2/3).}
\label{ex_a1}
\end{gather}
Pay attention to the fact that 
\begin{gather*}
\pFq{}{1}{a;b_1,0}{c}{x,y}=F(a,b_1;c;x),\ \pFq{}{1}{a;0,b_2}{c}{x,y}=F(a,b_2;c;y).
\end{gather*}
From the fact, 
we can degenerate formula (\ref{ex_a1}) 
into the following two hypergeometric identities for $F$:
\begin{align}
F(2\alpha,1-\alpha;\alpha+2/3;1/9)&=
\frac{3^{\alpha}\sqrt{\pi}\G(\alpha+2/3)
}{2^{2\alpha}\G(2/3)\G(\alpha+1/2)}, \label{ex_a1_d1}\\
F(2,6\beta;2\beta+5/3;1/3)&=3\beta+1. \label{ex_a1_d2}
\end{align}
Formulae (\ref{ex_a1_d1}) and (\ref{ex_a1_d2}) are equivalent to
formulae (1,2,3-1)(x) and (1,3,3-3)(viii), respectively.
Hence, we can regard (\ref{ex_a1}) as a fusion of 
hypergeometric identities for the hypergeometric series  (1,2,3-1)(x) and (1,3,3-3)(viii).
\subsubsection*{Example 2'': $(k,l_1,l_2,m)=(3,0,2,2)$}
When $(k,l_1,l_2,m)=(3,0,2,2)$, we have
$$
(a,b_1,b_2,c,x,y)=(3\alpha,2\beta-1, 2\alpha+1-2\beta,2\alpha+\beta,3/4,1/4)
$$
as one of solutions of the system (\ref{case_a}).
In this case, the equation (\ref{diff_eq_a1}) becomes 
\begin{align}
\begin{split}
&\pFq{}{1}{3\alpha; 2\beta-1, 2\alpha+1-2\beta}{2\alpha+\beta}{\frac{3}{4},\frac{1}{4}}\\
&=\frac{3^{3n}(\alpha+1/3,n)(\alpha+2/3,n)}{2^{4n}(2\alpha+\beta,2n)}
\pFq{}{1}{3\alpha+3n; 2\beta-1, 2\alpha+1-2\beta+2n}{2\alpha+\beta+2n}{\frac{3}{4},\frac{1}{4}}.
\end{split}
\label{ex_eq_a2}
\end{align}
By solving this equation (\ref{ex_eq_a2}), we find
\begin{gather}
\pFq{}{1}{3\alpha; 2\beta-1, 2\alpha+1-2\beta}{2\alpha+\beta}{\frac{3}{4},\frac{1}{4}}
=\frac{2^{4\alpha+1}{\pi}\G(2\alpha+\beta)
}{3^{3\alpha+1/2}\G(\alpha+1/3)\G(\alpha+2/3)\G(\beta).}
\label{ex_a2}
\end{gather}
In a similar manner to example 1'', 
we can degenerate formula (\ref{ex_a2}) into 
the following two hypergeometric identities for $F$:
\begin{align}
F(3\alpha,2\alpha;3\alpha+1/2;3/4)&=
\frac{2^{4\alpha} \Gamma(\alpha+1/6)\Gamma (\alpha+5/6)}
{\sqrt{3}\Gamma(\alpha+1/3)\Gamma (\alpha+2/3)},
\label{ex_a2_d1}\\
F(3\alpha,2\alpha;2\alpha+1/2;1/4)&=
\frac{2^{4\alpha+1}\sqrt{\pi}\Gamma(2\alpha+1/2)}
{3^{3\alpha+1/2}\Gamma(\alpha+1/3)\Gamma(\alpha+2/3)}.
\label{ex_a2_d2}
\end{align}
Here, we assume that $\alpha \notin \{ -3/2,-5/2,-7/2,\cdots \}$ in formula (\ref{ex_a2_d1}), 
because we have extended $F(a,b;c;x)$ as (\ref{extend_hgf}).
Formulae (\ref{ex_a2_d1}) and (\ref{ex_a2_d2}) are equivalent to
the above case of formula (1,3,3-1)(xix) and formula (1,3,3-1)(xii), respectively.
Hence, we can regard (\ref{ex_a2}) as a fusion of 
hypergeometric identities for the hypergeometric series  (1,3,3-1)(xix) and (1,3,3-1)(xii),
while we are not able to obtain the case below of  formula (1,3,3-1)(xix) 
from equation (\ref{ex_eq_a2}) or formula (\ref{ex_a2}).

In a similar manner, we can make use of our method for
Appell series $F_2$, $F_3$, $F_4$ and Lauricella series $F_A$, $F_B$, $F_C$, $F_D$.
We have never investigated hypergeometric identities 
for Appell-Lauricella series,
whereas we have advanced the study of those 
for generalized hypergeometric series so far.
This reason is due to the difficulty 
in applying previously known methods to the former series directly.
However, our method can be applied to
Appell-Lauricella series as well as generalized hypergeometric series.
The author feels worthy enough to 
launch the study of hypergeometric identities for Appell-Lauricella series 
from the point of view of fusions
(or generalizations) of those identities for the hypergeometric series.

\medskip
\begin{flushleft}
Akihito Ebisu\\
Department of Mathematics\\
Hokkaido University\\
Kita 10, Nishi 8, Kita-ku, Sapporo, 060-0810\\
Japan\\
a-ebisu@math.sci.hokudai.ac.jp
\end{flushleft}

\end{document}